\documentclass[opre,nonblindrev]{oo} % current default for manuscript submission

\usepackage{amsmath,amssymb,ifthen,url,graphicx,color,array,theorem}
\usepackage[table]{xcolor}
\usepackage{enumitem}
\usepackage[toc,page]{appendix}
\usepackage{pgfplots}
 \usepgfplotslibrary{groupplots,dateplot}
 \usetikzlibrary{patterns,shapes.arrows,fit}
 \pgfplotsset{compat=newest}

\usepackage{amsfonts}  
\usepackage{extarrows}
\usepackage{mathdots}
\usepackage{yhmath}
\usepackage{cancel}
\usepackage{tabularx}
\usepackage{comment}
\usepackage{eurosym}
\usepackage{booktabs}
\usepackage{multicol}
\usepackage{multirow}
\usepackage{longtable}
\usepackage{float}
\usepackage{caption}
\usepackage{subfigure}
\usepackage{rotating}
\usepackage{algorithm}
\usepackage[compatible]{algpseudocode}
 \usepackage{natbib}
    \bibpunct[, ]{(}{)}{,}{a}{}{,}%
    \def\bibfont{\small}%
\usepackage[colorlinks=true,linkcolor=blue,citecolor=blue, hypertexnames=false]{hyperref}

\allowdisplaybreaks

\OneAndAHalfSpacedXI % current default line spacing

\usepackage{natbib}
 \bibpunct[, ]{(}{)}{,}{a}{}{,}%
 \def\bibfont{\small}%

\TheoremsNumberedThrough     % Preferred (Theorem 1, Lemma 1, Theorem 
\EquationsNumberedThrough    % Default: (1), (2), ...
\ECRepeatTheorems

\newcommand{\convx}{\textup{conv}}
\newcommand{\conv}[1]{\convx (#1)}

\newcommand{\Exp}{\mathbb{E}}
\newcommand{\R}{\mathbb{R}}
\newcommand{\Z}{\mathbb{Z}}

\newcommand{\dist}{\mathbb{P}}

\newcommand{\ntminusone}{\npar_{t-1}}
\newcommand{\mpar}{m}
\newcommand{\mstatet}{\mpar^\textsc{s}_t}
\newcommand{\mrecourset}{\mpar^\textsc{r}_t}
\newcommand{\npar}{n}
\newcommand{\nt}{\npar_t}

\newcommand{\xipar}{\xi}
\newcommand{\xiparrand}{\xipar} %{\boldsymbol{\xipar}}
\newcommand{\xiparval}{\xipar}

\newcommand{\xicolor}{black} % LimeGreen
\newcommand{\xiprimecolor}{black} % Magenta
\newcommand{\xiT}{{\color{\xicolor}{\xiparrand^T}}}

\newcommand{\xiprimeT}{{\color{\xiprimecolor}{{\xiparrand'}^T}}}
\newcommand{\xisupt}{{\color{\xicolor}{\xiparrand^t}}}

\newcommand{\xiprimesupt}{{\color{\xiprimecolor}{{\xiparrand'}^t}}}
\newcommand{\xisuptminusone}{{\color{\xicolor}{\xiparrand^{t-1}}}}

\newcommand{\xiTval}{{\color{\xicolor}{\xiparval^T}}}
\newcommand{\xiprimeTval}{{\color{\xiprimecolor}{{\xiparval'}^T}}}
\newcommand{\xisuptval}{{\color{\xicolor}{\xiparval^t}}}

\newcommand{\xiprimesuptval}{{\color{\xiprimecolor}{{\xiparval'}^t}}}
\newcommand{\xisuptminusoneval}{{\color{\xicolor}{\xiparval^{t-1}}}}
\newcommand{\xisuptplusoneval}{{\color{\xicolor}{\xiparval^{t+1}}}}

\newcommand{\xisupp}{\Xi^T}
\newcommand{\xisuppt}{\Xi^t}
\newcommand{\naldr}{\alpha}

\newcommand{\set}[1] {[#1]}
\newcommand{\bsub}{\begin{subequations}}
\newcommand{\esub}{\end{subequations}}

\newcommand{\bsubeq}{\begin{subequations}}
\newcommand{\esubeq}{\end{subequations}}
\newcommand{\BI}{\begin{itemize}}
\newcommand{\EI}{\end{itemize}}
\newcommand{\BE}{\begin{enumerate}}
\newcommand{\EE}{\end{enumerate}}
\newcommand{\I}{\item}
\newcommand{\BSE}{\begin{subequations}}
\newcommand{\ESE}{\end{subequations}}

\newcommand{\node}{v}

\def \R {\mathbb{R}}
\def \Z {\mathbb{Z}}

\newcommand{\xiTminusone}{{\color{black}{\xiparrand^{T-1}}}}

\newcommand{\XiOpt}{\hat{\Xi}^\textsc{Feas}}
\newcommand{\XiExt}{\Xi^{\textsc{Ext}}}
\newcommand{\XiMP}{\hat{\Xi}}

\newcommand{\Ct}{c_t(\xisupt)}
\newcommand{\Ctval}{c_t(\xisuptval)}
\newcommand{\At}{A_t(\xisupt)}
\newcommand{\Atval}{A_t(\xisuptval)}
\newcommand{\Bt}{B_t(\xisupt)}

\newcommand{\Btplusoneval}{B_{t+1}(\xisuptplusoneval)}
\newcommand{\Dt}{D_t(\xisupt)}
\newcommand{\Dtval}{D_t(\xisuptval)}

\newcommand{\statevar}{x^\texttt{s}}
\newcommand{\statevart}{x^\texttt{s}_t(\xisupt)}
\newcommand{\statevartminusone}{x^\texttt{s}_{t-1}(\xisuptminusone)}
\newcommand{\statevarT}{x^\texttt{s}_T(\xiT)}
\newcommand{\recvar}{x^\texttt{r}}
\newcommand{\recvart}{x^\texttt{r}_t(\xisupt)}
\newcommand{\recvartminusone}{x^\texttt{r}_{t-1}(\xisuptminusone)}
\newcommand{\recvarT}{x^\texttt{r}_T(\xiT)}
\newcommand{\recvartval}{x^\texttt{r}_{t,\xiparval^{\crev T}}}

\newcommand{\stateCt}{c^\texttt{s}_t(\xisupt)}
\newcommand{\recCt}{c^\texttt{r}_t(\xisupt)}

\newcommand{\stateAt}{A^\texttt{s}_t(\xisupt)}
\newcommand{\recAt}{A^\texttt{r}_t(\xisupt)}

\newcommand{\stateBt}{B^\texttt{s}_t(\xisupt)}

\newcommand{\stateCtFixed}{c^\texttt{s}_t}
\newcommand{\recCtFixed}{c^\texttt{r}_t}
\newcommand{\stateAtFixed}{A^\texttt{s}_t}
\newcommand{\recAtFixed}{A^\texttt{r}_t}
\newcommand{\stateBtFixed}{B^\texttt{s}_t}

\newcommand{\stateDtFixed}{D^\texttt{s}_t}
\newcommand{\recDtFixed}{D^\texttt{r}_t}

\newcommand{\twostageLDR}{\beta}
\newcommand{\twostageLDRVart}{\twostageLDR_t}
\newcommand{\twostageLDRVartminusone}{\twostageLDR_{t-1}}
\newcommand{\twostageLDRBF}{\Phi}
\newcommand{\twostageLDRBFt}{\twostageLDRBF_t(\xisupt)}
\newcommand{\twostageLDRBFtminusone}{\twostageLDRBF_{t-1}(\xisuptminusone)}
\newcommand{\twostageLDRBFtVal}{\twostageLDRBF_t(\xisuptval)}
\newcommand{\twostageLDRBFtminusoneVal}{\twostageLDRBF_{t-1}(\xisuptminusoneval)}
\newcommand{\twostageBF}{\Theta}
\newcommand{\twostageBFt}{\twostageBF_t(\xisupt,\twostageLDRVart)}

\newcommand{\twostageBFtminusone}{\twostageBF_{t-1}(\xisuptminusone,\twostageLDRVartminusone)}
\newcommand{\twostageKBF}{\Upsilon}

\newcommand{\twostageKBFti}{\twostageKBF_{ti}(\xisupt,\twostageLDR_{t}^i)}

\newcommand{\K}{\mathcal{K}}
\newcommand{\twostageKBFAff}{{\crev \hat{\Upsilon}}}

\newcommand{\twostageKBFAffti}{\twostageKBFAff_{ti}(\xisupt)}

\newcommand{\MP}{\mathcal{MP}(\XiMP)}
\newcommand{\SPname}{\mathcal{SP}}
\newcommand{\SP}{\SPname(\hat{\twostageLDR},{\crev \hat{x}^\texttt{s}_1})}
\newcommand{\SPPLDR}{\SPname^\text{PCDR}(\twostageLDR,{\crev \statevar_1})}
\newcommand{\twostagePLDR}{x^\texttt{s}}
\newcommand{\twostagePLDRVart}{\twostagePLDR_{ti}}

\newcommand{\naBF}{\Psi}
\newcommand{\naBFt}{\naBF_{t}(\xiT)}

\newcommand{\nadualobj}{u_{\xiTval}}
\newcommand{\nadualstate}{v_{t\xiTval}}
\newcommand{\nadualstateplusone}{v_{t+1,\xiTval}}
\newcommand{\nadualrec}{w_{t\xiTval}}

\newcommand{\dnavar}{z}

\newcommand{\dualSPstate}{\pi^\texttt{S}}
\newcommand{\dualSPrec}{\pi^\texttt{R}}
\newcommand{\dualSPstatet}{\dualSPstate_t}
\newcommand{\dualSPrect}{\dualSPrec_t}

\newcommand{\MSARO}{MSARO}

\newcommand{\bound}{\nu}
\newcommand{\exactB}{\bound^\star}

\newcommand{\PILB}{\bound^\text{PI}}

\newcommand{\twostageUB}{\bound^\text{2S}}
\newcommand{\twostageLDRUB}{\bound^\text{2S-LDR}}
\newcommand{\twostagePLDRUB}{\bound^\text{2S-PCDR}}
\newcommand{\LDRUB}{\bound^\text{LDR}}
\newcommand{\pDistVar}{\rho}
\newcommand{\pDist}{\pDistVar^{\dist}_{\xiTval}}
\newcommand{\pDistprime}{\pDistVar^{\dist}_{\xiprimeTval}}
\newcommand{\pNormalization}{\rho_{|\xisuptval}^{\dist}}

\newcommand{\optT}{\textsc{opt}^\texttt{T}}
\newcommand{\optN}{\textsc{opt}^\texttt{N}}

\newcommand{\crev}{\color{black}} %blue

\begin{document}
%%%%%%%%%%%%%%%%
\RUNAUTHOR{Daryalal, Arslan, and Bodur}
\RUNTITLE{Primal and Dual Decision Rules For MSARO}
\TITLE{
% Primal and Dual Bounding Techniques
Two-stage and Lagrangian Dual Decision Rules for Multistage Adaptive Robust Optimization}
\ARTICLEAUTHORS{}
\ABSTRACT{In this work, we design primal and dual bounding methods for multistage adaptive robust optimization (\MSARO) problems {\crev motivated by} two decision rules rooted in the stochastic programming literature. From the primal perspective, this is achieved by applying decision rules that restrict the functional forms of {\crev only} a certain subset of decision variables {\crev resulting in an approximation of {\MSARO} as a two-stage adjustable robust optimization problem. We leverage the two-stage robust optimization literature in the solution of this approximation. From the dual perspective, decision rules are applied to the Lagrangian multipliers of a Lagrangian dual of {\MSARO}, resulting in a two-stage stochastic optimization problem. As the quality of the resulting dual bound depends on the distribution chosen when developing the dual formulation, we define a distribution optimization problem with the aim of optimizing the obtained bound and develop solution methods tailored to the nature of the recourse variables.} Our framework is general-purpose and does not require strong assumptions such as a stage-wise independent uncertainty set, and can consider integer recourse variables. Computational experiments on newsvendor, location-transportation, and capital budgeting problems show that our bounds yield considerably smaller optimality gaps compared to the existing methods.

}

\ARTICLEAUTHORS{
	\AUTHOR{Maryam Daryalal}
	\AFF{Department of Decision Sciences, HEC Montr\'{e}al, Montr\'{e}al, Qu\'{e}bec H3T 2A7, Canada, \\\EMAIL{maryam.daryalal@hec.ca}}
	\AUTHOR{Ay\c{s}e N. Arslan}
	\AFF{Univ. Bordeaux, CNRS, INRIA, Bordeaux INP, IMB, UMR 5251, F-33400 Talence, France,\\
 Centre Inria de l'Universite de Bordeaux, F-33405 Talence, France,\\
\EMAIL{ayse-nur.arslan@inria.fr}}
	\AUTHOR{Merve Bodur}
	\AFF{School of Mathematics, University of Edinburgh, Edinburgh  EH9 3FD, UK, \\ \EMAIL{merve.bodur@ed.ac.uk}}
}

\maketitle
%%%%%%%%%%%%%%%%%%%%%%%%%%%%%%%%
% ------------------------------
\section{Introduction}
Many practical planning, design and operational problems involve making decisions under uncertainty at consecutive stages, where the decisions in one stage affect the decisions of the future stages. In such sequential decision-making problems, first-stage (\emph{here-and-now}) decisions are the ones that are immediately implementable. Subsequent recourse (\emph{wait-and-see}) decisions depend on the state of the system, which is a result of previous decisions and observations of the uncertain parameters. A solution is then an adaptable \emph{policy} or \emph{decision rule} that takes the previous decisions and history of uncertainty realizations as an input, and returns a new implementable decision. The dynamics of a sequential decision-making problem is depicted in Figure \ref{fig:seq-decision}. 
\begin{figure}[htbp]
    \centering
    \small
    \scalebox{0.9}{
    \tikzset{every picture/.style={line width=0.75pt}} %set default line width to 0.75pt        

\begin{tikzpicture}[x=0.75pt,y=0.75pt,yscale=-1,xscale=1]
%uncomment if require: \path (0,300); %set diagram left start at 0, and has height of 300

%Rounded Rect [id:dp024760691964394232] 
\draw  [color={rgb, 255:red, 0; green, 0; blue, 0 }  ,draw opacity=1 ][line width=0.75]  (60.48,111.95) .. controls (60.48,105.87) and (65.41,100.93) .. (71.49,100.93) -- (137.77,100.93) .. controls (143.85,100.93) and (148.78,105.87) .. (148.78,111.95) -- (148.78,144.99) .. controls (148.78,151.07) and (143.85,156) .. (137.77,156) -- (71.49,156) .. controls (65.41,156) and (60.48,151.07) .. (60.48,144.99) -- cycle ;
%Straight Lines [id:da5351498880138911] 
\draw [color={rgb, 255:red, 0; green, 0; blue, 0 }  ,draw opacity=1 ]   (149.47,127.41) -- (194.87,127.52) ;
\draw [shift={(196.87,127.52)}, rotate = 180.13] [color={rgb, 255:red, 0; green, 0; blue, 0 }  ,draw opacity=1 ][line width=0.75]    (10.93,-3.29) .. controls (6.95,-1.4) and (3.31,-0.3) .. (0,0) .. controls (3.31,0.3) and (6.95,1.4) .. (10.93,3.29)   ;
%Straight Lines [id:da31548847674679037] 
\draw [color={rgb, 255:red, 155; green, 155; blue, 155 }  ,draw opacity=1 ] [dash pattern={on 0.84pt off 2.51pt}]  (297,127) -- (326,127) ;
%Straight Lines [id:da2346366602602119] 
\draw [color={rgb, 255:red, 155; green, 155; blue, 155 }  ,draw opacity=1 ]   (423,128) -- (469.75,127.54) ;
\draw [shift={(471.75,127.52)}, rotate = 179.44] [color={rgb, 255:red, 155; green, 155; blue, 155 }  ,draw opacity=1 ][line width=0.75]    (10.93,-3.29) .. controls (6.95,-1.4) and (3.31,-0.3) .. (0,0) .. controls (3.31,0.3) and (6.95,1.4) .. (10.93,3.29)   ;
%Straight Lines [id:da9997302461308888] 
\draw [color={rgb, 255:red, 155; green, 155; blue, 155 }  ,draw opacity=1 ]   (561.44,127.47) -- (602.5,127.47) ;
\draw [shift={(604.5,127.47)}, rotate = 180] [color={rgb, 255:red, 155; green, 155; blue, 155 }  ,draw opacity=1 ][line width=0.75]    (10.93,-3.29) .. controls (6.95,-1.4) and (3.31,-0.3) .. (0,0) .. controls (3.31,0.3) and (6.95,1.4) .. (10.93,3.29)   ;
%Rounded Rect [id:dp6156219639621894] 
\draw  [color={rgb, 255:red, 155; green, 155; blue, 155 }  ,draw opacity=1 ][line width=0.75]  (198.28,111.95) .. controls (198.28,105.87) and (203.21,100.93) .. (209.29,100.93) -- (275.57,100.93) .. controls (281.65,100.93) and (286.58,105.87) .. (286.58,111.95) -- (286.58,144.99) .. controls (286.58,151.07) and (281.65,156) .. (275.57,156) -- (209.29,156) .. controls (203.21,156) and (198.28,151.07) .. (198.28,144.99) -- cycle ;
%Rounded Rect [id:dp31351413199947864] 
\draw  [color={rgb, 255:red, 155; green, 155; blue, 155 }  ,draw opacity=1 ][line width=0.75]  (333.92,110.95) .. controls (333.92,104.87) and (338.85,99.93) .. (344.94,99.93) -- (411.22,99.93) .. controls (417.3,99.93) and (422.23,104.87) .. (422.23,110.95) -- (422.23,143.99) .. controls (422.23,150.07) and (417.3,155) .. (411.22,155) -- (344.94,155) .. controls (338.85,155) and (333.92,150.07) .. (333.92,143.99) -- cycle ;
%Rounded Rect [id:dp9927491202268093] 
\draw  [color={rgb, 255:red, 155; green, 155; blue, 155 }  ,draw opacity=1 ][line width=0.75]  (472.44,110.95) .. controls (472.44,104.87) and (477.37,99.93) .. (483.45,99.93) -- (549.73,99.93) .. controls (555.82,99.93) and (560.75,104.87) .. (560.75,110.95) -- (560.75,143.99) .. controls (560.75,150.07) and (555.82,155) .. (549.73,155) -- (483.45,155) .. controls (477.37,155) and (472.44,150.07) .. (472.44,143.99) -- cycle ;
%Curve Lines [id:da0532240535388806] 
\draw [color={rgb, 255:red, 155; green, 155; blue, 155 }  ,draw opacity=1 ]   (156.5,69.03) .. controls (162.5,85.67) and (168.5,103.9) .. (172.5,88.84) .. controls (176.4,74.16) and (180.3,84.34) .. (198.1,100.28) ;
\draw [shift={(199.5,101.52)}, rotate = 221.21] [color={rgb, 255:red, 155; green, 155; blue, 155 }  ,draw opacity=1 ][line width=0.75]    (10.93,-3.29) .. controls (6.95,-1.4) and (3.31,-0.3) .. (0,0) .. controls (3.31,0.3) and (6.95,1.4) .. (10.93,3.29)   ;
%Curve Lines [id:da34084432675297904] 
\draw [color={rgb, 255:red, 155; green, 155; blue, 155 }  ,draw opacity=1 ]   (292.8,69.03) .. controls (298.8,85.67) and (304.8,103.9) .. (308.8,88.84) .. controls (312.7,74.16) and (316.6,84.34) .. (334.4,100.28) ;
\draw [shift={(335.8,101.52)}, rotate = 221.21] [color={rgb, 255:red, 155; green, 155; blue, 155 }  ,draw opacity=1 ][line width=0.75]    (10.93,-3.29) .. controls (6.95,-1.4) and (3.31,-0.3) .. (0,0) .. controls (3.31,0.3) and (6.95,1.4) .. (10.93,3.29)   ;
%Curve Lines [id:da7560458339022407] 
\draw [color={rgb, 255:red, 155; green, 155; blue, 155 }  ,draw opacity=1 ]   (431.1,69.03) .. controls (437.1,85.67) and (443.1,103.9) .. (447.1,88.84) .. controls (451,74.16) and (454.9,84.34) .. (472.71,100.28) ;
\draw [shift={(474.1,101.52)}, rotate = 221.21] [color={rgb, 255:red, 155; green, 155; blue, 155 }  ,draw opacity=1 ][line width=0.75]    (10.93,-3.29) .. controls (6.95,-1.4) and (3.31,-0.3) .. (0,0) .. controls (3.31,0.3) and (6.95,1.4) .. (10.93,3.29)   ;

% Text Node
\draw (104.63,128.26) node [align=center] {
Decisions at\\stage 1};
% Text Node
\draw (242.43,128.26) node [align=center] {
Decisions at\\stage 2};
% Text Node
\draw (378.08,128.26) node [align=center] {
Decisions at\\stage $t-1$
};
% Text Node
\draw (516.59,128.26) node  [align=center] {Decisions at\\stage $t$};
% Text Node
\draw (100,52.03) node [anchor=north west][inner sep=0.75pt]   {Uncertainty};
% Text Node
\draw (241,52.82) node [anchor=north west][inner sep=0.75pt]  {Uncertainty};
% Text Node
\draw (382,52.82) node [anchor=north west][inner sep=0.75pt]   {Uncertainty};

\end{tikzpicture}
    }
    \caption{Sequential decision-making under uncertainty}
    \label{fig:seq-decision}
\end{figure}

There are several modeling frameworks for sequential decision-making problems under uncertainty. {\crev When the probability distribution governing the uncertain parameters is known}, these problems may be addressed by the \emph{multistage stochastic programming} (MSP) paradigm, with the goal of optimizing some statistical performance measure over the planning horizon. There is an extensive body of research on MSP problems of various structures, with a rich literature on problems with continuous decision variables. {\crev However, \cite{bertsimas2006robust}  point out that implemented solutions may perform poorly if the probability distribution used in the MSP model is different than the \emph{true} distribution,  even if both distributions share the same first and second moments. }To mitigate this effect, \emph{distributionally robust optimization} (DRO) models are proposed for making decisions that are based on a family of probability distributions, often defined by using historical data  \citep{goh2010distributionally}. These models aim to hedge against tuning decisions to a perceived distribution. {\crev While the DRO framework has received significant attention (\cite{mohajerinesfahani2018datadriven,cheramin2022computationally}) from the research community and some recent studies have proposed tractable solution methods for linear DRO problems under certain conditions \citep{philpott2018distributionally,bertsimas2019adaptive}, they remain largely challenging to solve especially in the multi-stage setting}.

{\crev\emph{Multistage adaptive/adjustable robust optimization} (MSARO) is another framework for modeling sequential decision-making problems under uncertainty that does not require any knowledge about the probability distribution governing the uncertain parameters. This framework is also adapted to contexts where the underlying uncertainty is not stochastic in nature, for instance, in the case of adversarial participants. In MSARO, the uncertainty is represented as belonging to a pre-structured (often compact) set, called the \emph{uncertainty set}, and the decisions are optimized with respect to the worst-case outcome in this set. 

In this paper, we focus on the MSARO framework.
Throughout, we use $ [a] := \{1,2, \dots, a\} $ and $ [a,b] := \{a, a+1, \dots, b\} $ for positive integers $a$ and $b$ (with $a\leq b $), and $ (\cdot)^{\top}$ for the transpose operator. 
For a problem with $T$ decision-making stages, we denote with $\xisupp\subseteq \R^{\ell^T}$ the uncertainty set which governs the set of uncertain parameters $(\xi_1,\xi_2, \hdots, \xi_T)$, where $\xi_t$ denotes the vector of parameters associated with stage $t$, and  $\xi_1=1$, by convention. For ease of presentation, we also define the sequence of uncertain parameter vectors up to stage $t$ along with their (projected) support as $\xi^t := (\xi_1,\hdots,\xi_t)\in \Xi^t := \text{proj}_{\xi^t}(\xisupp)\subseteq \R^{\ell^t}$. We assume $\Xi^t$ is compact for all $t \in [T]$. Then, we study the following general MSARO problem: 
}
{\crev
\begin{align}
\label{eq:MSAROnested}
    & \hspace*{-0.22cm} \min_{x_1 \in X_1(\xipar_1)}  c_1(\xipar_1)^\top x_1  + \sup_{\xi_2: (\xi^1,\xi_2)\in \Xi^2} 
 \min_{\substack{x_2 \in X_2(\xipar^2) : \\ A_2(\xipar^2) x_2 + B_2(\xipar^2) x_{1} \leq b_2(\xipar^2)}} \hspace*{-0.8cm} c_2(\xipar^2)^\top x_2 + \cdots 
\\*[0.15cm]
& \hspace*{0.1cm} \cdots + \hspace*{-0.17cm} \sup_{\xipar_t : (\xipar^{t-1},\xipar_t)\in \Xi^t} \min_{\substack{x_t \in X_t(\xipar^t): \\ A_t(\xipar^t) x_t + B_t (\xipar^t) x_{t-1} \leq b_t(\xipar^t)}} \hspace*{-0.9cm} c_t(\xipar^t)^\top x_t + \ \cdots \ + 
\hspace*{-0.17cm} \sup_{\xipar_T:(\xipar^{T-1},\xipar_T) \in \xisupp}
\min_{\substack{x_T \in X_T(\xipar^T) : \\ A_T(\xipar^T) x_T + B_T(\xipar^T) x_{T-1} \leq b_T(\xipar^T)}} \hspace*{-0.9cm} c_T(\xipar^T)^\top x_T  \nonumber
\end{align}
where $X_t(\xisupt) := \big\{x_t\in\R^{n_t-n^\texttt{i}_t}\times\Z^{n^\texttt{i}_t}:\  \Dt x_t\leq d_t(\xisupt)\big\}$ and $ c_t : \R^{\ell^t} \rightarrow \R^{\nt}, 
b_t : \R^{\ell^t} \rightarrow \R^{\mstatet}, 
d_t : \R^{\ell^t} \rightarrow \R^{\mrecourset}, 
A_t : \R^{\ell^t} \rightarrow \R^{\mstatet\times \nt}, B_t : \R^{\ell^t} \rightarrow \R^{\mstatet\times \ntminusone}, 
D_t : \R^{\ell^t} \rightarrow \R^{\mrecourset\times \nt}$ for $t\in [T]$. 
The main output of model \eqref{eq:MSAROnested} is the first-stage deterministic (here-and-now) decisions $x_1$ which minimize the worst-case objective value over $T$ stages taking into account the sequential uncertainty realizations and optimal wait-and-see decisions. In the sequential framework, at each stage $t\in T$, the worst-case realization $\xi_t$ that is consistent with the history of the realizations up to stage $t-1$, $\xi^{t-1}$, is revealed. The vector $\xi_t$ combined with the history $\xi^{t-1}$ yields the history of realizations up to stage $t$, that is, $\xi^t=(\xi^{t-1},\xi_t)$ from the support $\Xi^t$. This determines the parameters of the stage-$t$ minimization problem, from which the optimal wait-and-see decision vector $x_t$ is obtained. As such, the wait-and-see decisions, $x_2,\ldots,x_T$, also known as recourse decisions, are adapted to the history of the uncertain parameter realizations up to their decision-making stage, $\xi^2,\ldots, \xi^T$, respectively. We remark that, by definitions of $X_t$ and $\xisuppt$, we allow for the possibility of mixed-integer wait-and-see decisions and dependence between the uncertain parameters of different stages. In the following, we assume, for the data, that all uncertain vectors and matrices are affine functions of the associated uncertain parameters $\xipar^t$, as the majority of the literature mentioned makes this assumption; further assumptions will be specified explicitly when necessary. 
When the set $X_{t}(\xisupt), t\in [2,T]$ does not have integrality restrictions ($n_t^{\texttt{i}}=0$), MSARO problem has continuous recourse, otherwise, it has (mixed-)integer recourse.

MSARO problems are highly challenging to solve. Indeed, as has been recently proven by \cite{goerigk2024complexity}, they are, in general, harder than \texttt{NP}-hard problems, lying at the higher levels of the polynomial hierarchy, and their complexity increases with the number of decision stages. Specifically, $T$-stage MSARO problems with certain uncertainty set structures are  $\Sigma^\texttt{P}_{2T-1}$-hard. However, some sub-classes and special cases of MSARO problems are theoretically and/or computationally more tractable, depending on the number of stages, the structure of the uncertainty set, and the nature of recourse decisions. Among these, the most well-known is \emph{static} robust optimization, which considers that all decisions are here-and-now.  
For a considerable number of problem structures, e.g., when the uncertainty set is a polyhedron or an ellipsoid, it is possible to derive a monolithic reformulation of the static robust optimization problem through a compact reformulation of the adversarial problem, usually relying on duality techniques \citep{ben1999robust,ben2009robust,bertsimas2011theory,bertsimas2015tight}.
However, this paradigm cannot capture the flexibility offered by the possibility of adapting some decisions to the realization of uncertainty, thus often producing overly conservative decisions. As such, there has been significant research effort on developing exact and approximate solution methods for two-stage adjustable robust optimization (2ARO) problems, i.e., MSARO with $T=2$ \citep{ben2004adjustable,zeng2013solving,postek2016multistage,subramanyam2020k}. In developing these methods, the presence of discrete recourse variables, known as (mixed-)integer recourse, poses additional challenges in ensuring exact or high-quality solutions within a reasonable computational effort compared to the continuous recourse case. 

On the other hand, scientific progress on general MSARO has been much more limited. Given the aforementioned theoretical complexity of these problems, the focus of most existing studies is the approximate solution of these problems. Approximations proposed for MSARO mostly rely on reducing the multi-stage problem to a static problem, with a view to leverage the tractability of these problems. 
While these approximations can produce feasible solutions for MSARO problems, they can lead to a significant degradation in solution quality. Furthermore, some of these methods are quite restrictive being only applicable to special classes such as MSARO with continuous recourse.

To address these limitations, this paper aims to develop approximations for general MSARO problems of form  \eqref{eq:MSAROnested}. More specifically, 
inspired by the recent developments in the MSP literature, namely the two-stage linear decision rules \citep{bodur2018two}, 
we propose applying decision rule approximations to only a certain subset of decision variables, resulting in an approximation of MSARO problems as 2ARO problems. In so doing we have three motivations: (i) the strength of the added flexibility in adapting recourse decisions to uncertainty in 2ARO compared to static robust optimization, (ii) the significantly reduced theoretical complexity of 2ARO compared to MSARO, and (iii) the progress made in computationally viable exact and approximate solution methodologies for 2ARO. As a result of point (iii), our proposed framework is capable of considering a large variety of MSARO classes, most notably the mixed-integer recourse case, and will directly benefit from the developments in the highly active 2ARO literature in the future, e.g., the incorporation of machine learning for computational enhancements \citep{julien2022machine, dumouchelle2023neur2ro}.

While the aforementioned ideas are aimed at providing feasible solutions for MSARO problems, an important question arises as to the quality of the obtained solution. In order to evaluate the quality of a feasible policy, one could use a dual bound on the optimal value of the MSARO problem. Unfortunately, obtaining dual bounds for MSARO problems is a largely unexplored topic in the literature, especially in the case of mixed-integer recourse. To fill this gap, we propose to develop dual approximations for MSARO. In particular, we propose a Lagrangian dual for the MSARO problem and apply decision rules to the Lagrangian multipliers, leveraging ideas rooted in the MSP literature, namely, Lagrangian dual decision rules \citep{daryalal2020lagrangian}. In deriving a Lagrangian dual of the MSARO problem, we assign a probability distribution with the support as the uncertainty set and use the assigned distribution in dualizing a subset of constraints along with their Lagrangian multipliers. As a result, we obtain a dual approximation of MSARO in the form of a two-stage stochastic optimization problem, which can be solved with the help of state-of-the-art methods for two-stage stochastic problems. Since the quality of the resulting dual bound depends on the assigned distribution used while developing the dual formulation, we define a distribution optimization problem with the aim of identifying the strongest such dual bound. We develop appropriate solution methods tailored to the nature of the recourse variables for the resulting distribution optimization problem.

The contributions of our work are summarized as follows:
\BI
\I We develop a solution framework for MSARO problems that returns adaptable policies as well as a dual bound measuring the quality of these policies. We do this by employing two-stage and Lagrangian dual  decision rules, leading to novel techniques that reveal new theoretical and practical avenues. 

\I We propose to approximate MSARO problems via 2ARO problems to leverage existing solution methodologies and future developments for the latter in designing adaptable policies for the former. To this end, we present two-stage decision rules for MSARO, the first adaptation of a generalization of two-stage linear decision rules from the MSP literature in robust optimization, which can be applied to a broad range of problems. We employ, for an illustration of our approach, a tailored constraint-and-column generation algorithm to solve the resulting 2ARO approximation. The optimal solution of this approximation not only provides a primal policy but can also contribute to the calculation of a dual bound through identification of critical realizations in the uncertainty set. 

\I With a similar motivation, we derive a dual approximation of MSARO in the form of a two-stage stochastic optimization problem, which we show to be a strong dual in certain cases. Moreover, in order to obtain the strongest possible such dual bound, we study the numerical solution of a distribution optimization problem.
More specifically, we apply decision rules to dual variables, and design a cutting-plane algorithm to solve the obtained restricted dual problem. We also show that in the special case of continuous recourse, the restricted dual problem can be reformulated as a monolithic bilinear program. Additionally, we present an alternative decomposable dual problem which offers an improved numerical performance. These novel techniques contribute to the scarce literature for obtaining dual bounds for \MSARO{} problems with mixed-integer recourse.

\I We evaluate the performance of our solution framework over multistage versions of 
three 
classical problems from the \MSARO{} literature: $(i)$ the newsvendor problem, 
$(ii)$ the location-transportation problem, and 
$(iii)$ the capital budgeting problem. 
Each of these problem classes is suitable for a different solution method developed in this work and  our analysis over various instances attests to the quality of the returned primal and dual bounds.
\EI
}

{\crev

The remainder of the paper is organized as follows. In Section \ref{sec:intro-lit} we review the literature relevant to our work. In Section \ref{sec:primal} we introduce two-stage decision rules for obtaining primal adaptable policies. In Section \ref{sec:NA-dual} we present our approach to deriving dual bounds. 
This is followed by numerical experiments in Section \ref{sec:experiments} and concluding remarks. We remark that all proofs are deferred to Appendix~\ref{sec:proofs}. 

}

%%%%%%%%%%%%%%%%%%%%%%%%%%%%%%%%%
%%%%%%%%%%%%%%%%%%%%%%%%%%%%%%%%%
%%%%%%%%%%%%%%%%%%%%%%%%%%%%%%%%%
%%%%%%%%%%%%%%%%%%%%%%%%%%%%%%%%
\section{Literature Review}\label{sec:intro-lit}
Figure \ref{fig:literature} presents a summary of existing solution methods for obtaining exact/approximate solutions  and dual bounds for \MSARO{}, with methods developed specifically for 2ARO separately categorized. In the following, we briefly discuss each method and the specific problem structure it can address.  
\begin{figure}[htbp]
    \centering
    \scalebox{0.85}{
    \tikzset{every picture/.style={line width=0.75pt}} %set default line width to 0.75pt        

\begin{tikzpicture}[x=0.75pt,y=0.75pt,yscale=-1,xscale=1]
%uncomment if require: \path (0,300); %set diagram left start at 0, and has height of 300

%Rounded Rect [id:dp2675717185191404] 
\draw  [draw opacity=0][fill={rgb, 255:red, 230; green, 226; blue, 226 }  ,fill opacity=1 ] (73.42,-54.28) .. controls (73.42,-62.4) and (80,-68.98) .. (88.11,-68.98) -- (610.72,-68.98) .. controls (618.84,-68.98) and (625.42,-62.4) .. (625.42,-54.28) -- (625.42,-10.2) .. controls (625.42,-2.08) and (618.84,4.5) .. (610.72,4.5) -- (88.11,4.5) .. controls (80,4.5) and (73.42,-2.08) .. (73.42,-10.2) -- cycle ;
%Rounded Rect [id:dp8861250869570559] 
\draw  [draw opacity=0][fill={rgb, 255:red, 230; green, 226; blue, 226 }  ,fill opacity=1 ] (59.42,54.4) .. controls (59.42,40.1) and (71.01,28.5) .. (85.32,28.5) -- (613.52,28.5) .. controls (627.82,28.5) and (639.42,40.1) .. (639.42,54.4) -- (639.42,132.1) .. controls (639.42,146.4) and (627.82,158) .. (613.52,158) -- (85.32,158) .. controls (71.01,158) and (59.42,146.4) .. (59.42,132.1) -- cycle ;
%Rounded Rect [id:dp9447202998965468] 
\draw  [draw opacity=0][fill={rgb, 255:red, 230; green, 226; blue, 226 }  ,fill opacity=1 ] (130.42,186.63) .. controls (130.42,174.64) and (140.14,164.92) .. (152.13,164.92) -- (288.7,164.92) .. controls (300.69,164.92) and (310.42,174.64) .. (310.42,186.63) -- (310.42,251.78) .. controls (310.42,263.78) and (300.69,273.5) .. (288.7,273.5) -- (152.13,273.5) .. controls (140.14,273.5) and (130.42,263.78) .. (130.42,251.78) -- cycle ;
%Rounded Rect [id:dp5035240359182185] 
\draw  [draw opacity=0][fill={rgb, 255:red, 230; green, 226; blue, 226 }  ,fill opacity=1 ] (388.42,186.63) .. controls (388.42,174.64) and (398.14,164.92) .. (410.13,164.92) -- (546.7,164.92) .. controls (558.69,164.92) and (568.42,174.64) .. (568.42,186.63) -- (568.42,251.78) .. controls (568.42,263.78) and (558.69,273.5) .. (546.7,273.5) -- (410.13,273.5) .. controls (398.14,273.5) and (388.42,263.78) .. (388.42,251.78) -- cycle ;

%Rounded Rect [id:dp4878092184443583] 
\draw  [color={rgb, 255:red, 0; green, 0; blue, 0 }  ,draw opacity=1 ] (37.17,-159.42) .. controls (37.17,-175.34) and (50.08,-188.25) .. (66,-188.25) -- (632.83,-188.25) .. controls (648.76,-188.25) and (661.67,-175.34) .. (661.67,-159.42) -- (661.67,252.67) .. controls (661.67,268.59) and (648.76,281.5) .. (632.83,281.5) -- (66,281.5) .. controls (50.08,281.5) and (37.17,268.59) .. (37.17,252.67) -- cycle ;
%Rounded Rect [id:dp5424224977775751] 
\draw  [dash pattern={on 4.5pt off 4.5pt}] (50.42,-139.98) .. controls (50.42,-150.79) and (59.18,-159.56) .. (70,-159.56) -- (628.83,-159.56) .. controls (639.65,-159.56) and (648.42,-150.79) .. (648.42,-139.98) -- (648.42,-2.08) .. controls (648.42,8.73) and (639.65,17.5) .. (628.83,17.5) -- (70,17.5) .. controls (59.18,17.5) and (50.42,8.73) .. (50.42,-2.08) -- cycle ;
%Rounded Rect [id:dp02920837002057952] 
\draw  [draw opacity=0][fill={rgb, 255:red, 230; green, 226; blue, 226 }  ,fill opacity=1 ] (163.92,-132.68) .. controls (163.92,-140.57) and (170.32,-146.97) .. (178.21,-146.97) -- (520.62,-146.97) .. controls (528.52,-146.97) and (534.92,-140.57) .. (534.92,-132.68) -- (534.92,-89.79) .. controls (534.92,-81.9) and (528.52,-75.5) .. (520.62,-75.5) -- (178.21,-75.5) .. controls (170.32,-75.5) and (163.92,-81.9) .. (163.92,-89.79) -- cycle ;

% Text Node
\draw (349.42,-99.28) node  [font=\small] [align=left] {
 \citealt{bertsimas2010finite}, \citealt{hanasusanto2015k}\\ \citealt{subramanyam2020k} 
};
% Text Node
\draw  [draw opacity=0]  (272.92,-147.83) -- (425.92,-147.83) -- (425.92,-120.83) -- (272.92,-120.83) -- cycle  ;
\draw (349.42,-134.33) node  [font=\small] [align=left] {\textbf{Primal Approximation}};
% Text Node
\draw (349.42,-22.32) node  [font=\small] [align=left] {
\citealt{bertsimas2012adaptive}, \citealt{zeng2013solving}, \citealt{zhen2018adjustable}, \citealt{georghiou2020primal},\\\citealt{hashemi2020exploiting}, \citealt{arslan2022decomposition}
};
% Text Node
\draw  [draw opacity=0]  (326.92,-72.2) -- (371.92,-72.2) -- (371.92,-45.2) -- (326.92,-45.2) -- cycle  ;
\draw (349.42,-58.7) node  [font=\small] [align=left] {\textbf{Exact}};
% Text Node
\draw (349.42,105.51) node  [font=\small] [align=left] {\citealt{ben2004adjustable}, \citealt{chen2008linear}, \citealt{chen2009uncertain}, \citealt{goh2010distributionally},\\\citealt{see2010robust}, \citealt{bertsimas2011hierarchy}, \citealt{bertsimas2015design}, \\ \citealt{bertsimas2016multistage}, \citealt{postek2016multistage}, \citealt{bertsimas2018binary}, \\\citealt{ben2020tractable}, \citealt{romeijnders2020piecewise}, \citealt{xu2018improved}};
% Text Node
\draw  [draw opacity=0]  (272.92,29.82) -- (425.92,29.82) -- (425.92,56.82) -- (272.92,56.82) -- cycle  ;
\draw (349.42,43.32) node  [font=\small] [align=left] {\textbf{Primal Approximation}};
% Text Node
\draw  [draw opacity=0][fill={rgb, 255:red, 255; green, 255; blue, 255 }  ,fill opacity=1 ]  (318.92,-177.07) -- (379.92,-177.07) -- (379.92,-148.07) -- (318.92,-148.07) -- cycle  ;
\draw (349.42,-162.57) node   [align=left] { \ 2ARO };
% Text Node
\draw  [draw opacity=0][fill={rgb, 255:red, 255; green, 255; blue, 255 }  ,fill opacity=1 ]  (311.42,-204.26) -- (387.42,-204.26) -- (387.42,-175.26) -- (311.42,-175.26) -- cycle  ;
\draw (349.42,-189.76) node   [align=left] { \ MSARO };
% Text Node
\draw  [draw opacity=0]  (455.92,167.4) -- (500.92,167.4) -- (500.92,194.4) -- (455.92,194.4) -- cycle  ;
\draw (478.42,180.9) node  [font=\small] [align=left] {\textbf{Exact}};
% Text Node
\draw  [draw opacity=0]  (176.92,167.4) -- (263.92,167.4) -- (263.92,194.4) -- (176.92,194.4) -- cycle  ;
\draw (220.42,180.9) node  [font=\small] [align=left] {\textbf{Dual Bound}};
% Text Node
\draw (478.42,221.23) node  [font=\small] [align=left] {%\citealt{ben2004adjustable}\\
%\citealt{marandi2018static}\\
\citealt{zhen2018adjustable}\\\citealt{georghiou2019robust}};
% Text Node
\draw (220.42,232.73) node  [font=\small] [align=left] {\citealt{kuhn2011primal}\\ 
\citealt{hadjiyiannis2011scenario} \\
%\citealt{bertsimas2016duality},
 \citealt{georghiou2019robust} };

\end{tikzpicture}
    }
    \caption{Solution methods for \MSARO{}}
    \label{fig:literature}
\end{figure}

\emph{Exact solution methods} are scarce in the \MSARO{} literature and the existing studies  mostly focus on 2ARO. 
{\crev For 2ARO problems with fixed 
recourse and finite or polyhedral uncertainty set, \cite{zeng2013solving} developed a constraint-and-column generation algorithm. For the same type of problems restricted to continuous recourse, \cite{bertsimas2012adaptive} designed a Benders decomposition-type algorithm and applied it to a unit commitment problem, whereas \cite{georghiou2020primal} proposed a convergent method based on enumeration of the extreme points of the uncertainty set combined with affine decision rules to provide gradually improving primal and dual bounds.
For 2ARO problems with continuous fixed recourse, \cite{zhen2018adjustable} used Fourier-Motzkin elimination iteratively to remove the second-stage decisions, eventually forming an equivalent static robust optimization problem. This computationally expensive approach is also extended to multistage problems.
In the case of mixed-binary recourse and only objective uncertainty, \cite{arslan2022decomposition} proposed an exact method based on a Dantzig-Wolfe reformulation of the recourse problem based on a technical assumption on the structure of the linking constraints. Similarly, using Dantzig-Wolfe reformulation, for a subclass of 2ARO problems with fixed and mixed-integer recourse, block diagonal recourse matrix and a finite uncertainty set, \cite{hashemi2020exploiting} derived a static formulation which is amenable to Benders decomposition. For continuous \MSARO{} problems with a stage-wise rectangular uncertainty set, \cite{georghiou2019robust} developed robust dual dynamic programming (RDDP) and proved finite/asymptomatic convergence for various problem sub-classes. RDDP is an adaptation of the stochastic dual dynamic programming algorithm from the MSP literature \citep{pereira1991multi} to MSARO. 

Approximate solution methods are more common in the MSARO literature, with the central idea of restricting adaptable/adjustable decisions to follow a certain functional form, known as \emph{decision rules}. \cite{ben2004adjustable} proposed the first decision rule for \MSARO{} problems with \emph{continuous recourse}, LDRs, where recourse decisions are expressed as affine functions of uncertain parameters where the parameters of this function are to be optimized. The resulting LDR-restricted problem being a static optimization problem, it can be reformulated as a linear optimization problem in certain cases. Nonlinear decision rules were also proposed, such as deflected and segregated affine \citep{chen2008linear}, extended affine \citep{chen2009uncertain}, piecewise affine \citep{goh2010distributionally}, truncated linear \citep{see2010robust}, piecewise affine with exponentially many pieces \citep{ben2020tractable}, quadratic \citep{xu2018improved} and polynomial \citep{bertsimas2011hierarchy} decision rules. However, the resulting reformulations when using non-linear decision rules are often nonlinear, e.g, semidefinite or copositive programs. For a comprehensive list of nonlinear decision rules, interested reader may refer to the survey by \cite{yanikouglu2019survey}.  

In the case of \emph{mixed-integer recourse}, LDRs and most of its aforementioned extensions lead to non-adjustable decisions for the integer variables. Thus, alternative approaches have been proposed, with the key idea of (implicitly or explicitly) partitioning the uncertainty set and determining a constant recourse solution corresponding to each subset. 
A popular approach for 2ARO problems uses the notion of finite adaptability, first introduced by \cite{bertsimas2010finite}. In finite or $K$-adaptability, the decision-maker a priori commits to $K$ recourse decisions (while making the first-stage decisions), and then chooses among them after observing the uncertainty realization which leads to an implicit $K$-partition of the uncertainty set. 
While \cite{bertsimas2010finite} presented an exact formulation for the 2-adaptability problem, for the general $K$-adaptability case, 
\cite{hanasusanto2015k} proposed a monolithic formulation for problems with binary recourse, and \cite{subramanyam2020k} developed a branch-and-bound algorithm for problems with mixed-integer recourse. 
For MSARO problems on the other hand, explicit uncertainty set partitioning is considered in an iterative heuristic framework, with the aim of obtaining a sequence of improving approximations \citep{bertsimas2016multistage,postek2016multistage,romeijnders2020piecewise}. 
Lastly, for MSARO problems with pure-binary recourse, \cite{bertsimas2018binary} introduced binary decision rules, whereas for the mixed-binary recourse case, \cite{bertsimas2015design} implicitly designed piecewise linear/constant decision rules. 

While the aforementioned primal approximations can be shown to be exact in some special cases \citep{bertsimas2010optimality,bertsimas2012power,iancu2013supermodularity, hanasusanto2015k, zhen2018adjustable}, in general they do not provide optimal solutions. In order to assess the quality of their feasible solutions, dual bounds can be used. To this end, \cite{kuhn2011primal} presented the idea of deriving a dual problem for MSARO with only continuous variables and applying LDRs on the dual variables. {\crev Since their approach was originally derived for stochastic programs its application to MSAROs requires assigning a probability distribution to the uncertainty set. The impact of the chosen distribution on the quality of the obtained dual bound was observed by \cite{kuhn2011primal}, as such a distribution optimization problem was mentioned. This problem was later formalized by \cite{hadjiyiannis2011scenario} for a 2ARO with continuous variables and shown to be of the same theoretical difficulty as the original problem}. \cite{hadjiyiannis2011scenario} proposed to solve instead a 2ARO problem for a finite set of scenarios from the uncertainty set, selected based on a primal decision rule restriction, to reach a dual bound. This procedure can also be extended to obtain dual bounds for general \MSARO{}s.  
Finally, for special cases of 2ARO problems with continuous recourse, 
\cite{georghiou2020primal} proposed a framework to derive  progressive dual bounds, by considering the linear programming dual of their primal extreme point reformulation. 

}
%%%%%%%%%%%%%%%%%%%%%%%%%%%%%%%%
%%%%%%%%%%%%%%%%%%%%%%%%%%%%%%%%
%%%%%%%%%%%%%%%%%%%%%%%%%%%%%%%%
\section{Primal Bounding}\label{sec:primal}
{\crev The MSARO problem given in a nested form in  \eqref{eq:MSAROnested} can be reformulated as a monolithic optimization problem by explicitly introducing the functional form of the decision variables, $x_t(\xisupt): \xisuppt \rightarrow \R^{n_t-n^\texttt{i}_t}\times\Z^{n^\texttt{i}_t}$ for all $t\in \set{T}$, along with a deterministic variable, $z$, representing the worst-case objective value:}
\bsub
\label{eqs:msaro-mono}
\begin{align}
\exactB := \min\ \ &  z \label{eq:obj-msaro-mono} \\
\text{s.t.} \ \ &\sum_{t\in \set{T}}\Ct^\top x_t(\xisupt) \leq z &&  \xiT\in\xisupp \label{eq:obj-const-msaro-mono}\\
& \At x_t(\xisupt) + \Bt x_{t-1}(\xisuptminusone)\leq b_t(\xisupt) &&  t\in\set{2,T},\ \xisupt\in\xisuppt  \label{eq:const-state-equation-msaro-mono} \\
& x_t(\xisupt) \in X_{t}(\xisupt) &&   t\in \set{T},\ \xisupt\in\xisuppt. \label{eq:const-recourse-constraint-masaro-mono} 
\end{align}
\esub

{\crev 
Together with constraints \eqref{eq:obj-const-msaro-mono}, the objective function \eqref{eq:obj-msaro-mono} minimizes the worst outcome. Constraints \eqref{eq:const-state-equation-msaro-mono} and \eqref{eq:const-recourse-constraint-masaro-mono} are \emph{state} and \emph{recourse} constraints, respectively: while the former link different stages, the latter are local restrictions for a specific stage.

Throughout the paper, we make the following assumptions:
\begin{assumption}
\label{ass:nonemptyfirst}
$X_1(\xi^1)$ is non-empty.
\end{assumption}
\begin{assumption}
\label{ass:feasibility}
The problem has relatively complete recourse, i.e., for all $t\in[2,T]$, $\xisupt\in\xisuppt$ and a history of feasible decisions made up to $t$, $\{x_{t'}(\xi^{t'})\}_{t' \in [t-1]}$, there  always exists a feasible decision at stage $t$, $x_t(\xi^t)$.
\end{assumption}
\begin{assumption}\label{ass:bounded}
For $t\in\set{T}$ and $\xisupt\in\xisuppt$, the maximum diameter of $X_t(\xisupt)$ is finite, i.e., the feasibility sets are bounded. 
\end{assumption}
{We note that combined with the compactness assumption of the uncertainty set and its projections, these assumptions imply that the studied MSARO problem has a finite optimal objective value.}

To derive feasible policies to the \MSARO{} problem, it is quite common in the literature to restrict \emph{all} the decisions $x_t(\xisupt)$ to follow a simple functional form, such as an affine or piecewise constant decision rule. By breaking the temporal dependencies between stages, this approach approximates problem \eqref{eqs:msaro-mono} with a static robust optimization problem. 
Our goal in this section is to employ a new paradigm where a specific \emph{subset} of the decision variables $x_t(\xisupt)$ are enforced to follow a structured decision rule, leading to a restriction in the form of a 
2ARO
problem. We introduce this approach, two-stage decision rules for MSARO, in Section \ref{subsec:GeneralDecisionsRule}, then present its specific instantiation and possible solution methodologies in Sections \ref{sec:two-stage-LDR} and \ref{subsec:DecisionRulesDiscrete}.
}

\subsection{Two-stage Decision Rules}
\label{subsec:GeneralDecisionsRule}
In a similar manner to constraints \eqref{eq:const-state-equation-msaro-mono} and \eqref{eq:const-recourse-constraint-masaro-mono}, we partition the decision variables $x_t(\xisupt), t\in{\crev[T]}$ into $\statevart\in\R^{q_t}$ and $\recvart\in\R^{p_t}$, \emph{state} and \emph{recourse} variables, {\crev as those that appear in the state constraints of subsequent stages and the others, respectively. We have}  $q_t+p_t = n_t$, with integrality restrictions on the variables, if any, embedded in the set $X_t(\xisupt)$. Then the \MSARO{} {\crev \eqref{eqs:msaro-mono} can be written more explicitly} as follows:
\bsub
\label{eq:msaro-xt}
\begin{align}
\exactB = \min\ & z   \label{eq:obj-msaro-mono-xs} \\
\text{s.t.}\ & \sum_{t\in \set{T}}{\stateCt}^\top \statevart+{\recCt}^\top \recvart \leq z&& \xiT\in\xisupp \\
&\stateAt \statevart + \stateBt \statevartminusone + \recAt \recvart \leq b_{t}(\xisupt) && t\in \set{2,T}, \ \xisupt\in\xisuppt  \label{eq:const-state-equation-msaro-mono-xs} \\
&(\statevart,  \recvart)\in X_{t}(\xisupt)&&   t\in {\crev \set{T}},\ \xisupt\in\xisuppt \label{eq:const-recourse-constraint-masaro-mono-xs} 
 \end{align}
\esub
where $\stateCt,\stateAt,\stateBt, \crev{D_t^\texttt{s}(\xi^t)}$ are sub-arrays/sub-matrices of $\Ct,\At,\Bt, \crev{\Dt}$ associated with the state variables with appropriate dimensions, while $\recCt,\recAt, D_t^\texttt{r}(\xi^t)$ have the same role for the recourse variables, and {\crev $X_t(\xisupt) = \big\{\statevar_t \in \R^{q_t-q^\texttt{i}_t}\times\Z^{q^\texttt{i}_t}, \recvar_t \in \R^{p_t-p^\texttt{i}_t}\times\Z^{p^\texttt{i}_t}:\ D_t^\texttt{s}(\xi^t) \statevar_t + D_t^\texttt{r}(\xi^t) \recvar_t \leq d_t(\xisupt)\big\}$}. {\crev For notational convenience, we drop the parametrization for the first-stage variables as well as their feasible set and the objective vector, \emph{i.e.}, use $x_1 = (\statevar_1,\recvar_1), X_1,$ and $c_1 = (c_1^\texttt{s},c_1^\texttt{r})$.}

{\crev For $t \in [2,T],$} let $\statevart $ be approximated by a decision rule, i.e., $\statevart=\twostageBFt$, where $\twostageBF_t: \R^{\ell^{t}} \times \R^{K_t}  \rightarrow \R^{q_t}$ represents the rule, and $\twostageLDRVart \in \R^{K_t}$ {\crev represents} its vector of design parameters. By substituting this rule in problem \eqref{eq:msaro-xt}, 
we obtain an approximation that can be reformulated as:
\bsub
\label{eq:2RO}
\begin{align}
	\twostageUB := \min \ \ & c_1^\top x_1 + 
	\max_{\xiT\in\xisupp}\  
	\min_{\recvar\in\mathcal{X}({\crev \statevar_1}, \twostageLDR, \xiT)} \sum_{t\in \set{2,T}} {\stateCt}^\top \twostageBFt + {\recCt}^\top \recvar_t \label{eq:2RO-obj}\\
	\text{s.t.}\quad & x_1\in X_1 & \label{eq:2RO-const1}\\
	&\twostageLDRVart \in \R^{K_t}\qquad\ \ t\in\set{2,T},\label{eq:2RO-const2}
\end{align}
\esub
where:
\begin{align*}
	\mathcal{X}({\crev \statevar_1}, \twostageLDR,\xiT) := \Big\{&\big({\crev\recvar_t}\big)_{t\in\set{2,T}}\in\R^{p_2}\times\R^{p_3}\times \dots\times \R^{p_T}:\\
    &\ {\crev \recAt \recvar_t \leq b_{t}(\xisupt)-\Big(\stateAt \twostageBFt + \stateBt \statevar_1 \Big)} && {\crev t = 2} \\
	&\ \recAt \recvar_t \leq b_{t}(\xisupt)-\Big(\stateAt \twostageBFt + \stateBt \twostageBFtminusone \Big) && t\in \set{{\crev 3},T}   \\
	&\ (\twostageBFt, \recvar_t)\in X_{t}(\xisupt)&&   t\in \set{2,T}\ \Big\}.
\end{align*}
Note that the decision rules are solely applied to the state variables, whereas the recourse variables remain fully adjustable to the uncertain parameters (see Figure \ref{fig:two-stage-DR}). Problem \eqref{eq:2RO} is a 2ARO {\crev since} the temporal dependency between stages is removed thanks to the application of the two-stage decision rules {\crev (see Appendix~\ref{app:2AROreduction} for a detailed proof). We remark that the relatively complete recourse assumption stated in Assumption \ref{ass:feasibility} for the MSARO problem, does not guarantee that the 2ARO model has relatively complete recourse, but it can be ensured, for instance, by following the techniques mentioned in \citep{bodur2018two}.}

\begin{figure}[t] 
    \centering
    \scalebox{0.87}{
    \tikzset{every picture/.style={line width=0.75pt},
dotted_block/.style={draw=black!30!white, line width=1pt, dash pattern={on 4.5pt off 4.5pt},
            inner ysep=2.5mm,inner xsep=7mm, rectangle, rounded corners=9}
}

\begin{tikzpicture}[x=0.75pt,y=0.75pt,yscale=-1,xscale=1]
%uncomment if require: \path (0,300); %set diagram left start at 0, and has height of 300

%Rounded Rect [id:dp9477744974437592] 
% Stage t
\draw  [color={rgb, 255:red, 74; green, 74; blue, 74 }  ,draw opacity=1 ][fill={rgb, 255:red, 233; green, 229; blue, 229 }  ,fill opacity=1 ] (243,54.5) .. controls (243,51.19) and (245.69,48.5) .. (249,48.5) -- (327,48.5) .. controls (330.31,48.5) and (333,51.19) .. (333,54.5) -- (333,72.5) .. controls (333,75.81) and (330.31,78.5) .. (327,78.5) -- (249,78.5) .. controls (245.69,78.5) and (243,75.81) .. (243,72.5) -- cycle ;
%Rounded Rect [id:dp44992993718491825] 
% Stage 1
\draw  [color={rgb, 255:red, 74; green, 74; blue, 74 }  ,draw opacity=1 ][fill={rgb, 255:red, 233; green, 229; blue, 229 }  ,fill opacity=1 ] (122.5,54.5) .. controls (122.5,51.19) and (125.19,48.5) .. (128.5,48.5) -- (206.5,48.5) .. controls (209.81,48.5) and (212.5,51.19) .. (212.5,54.5) -- (212.5,72.5) .. controls (212.5,75.81) and (209.81,78.5) .. (206.5,78.5) -- (128.5,78.5) .. controls (125.19,78.5) and (122.5,75.81) .. (122.5,72.5) -- cycle ;
%Rounded Rect [id:dp2765950063531245] 
% Stage T
\draw  [color={rgb, 255:red, 74; green, 74; blue, 74 }  ,draw opacity=1 ][fill={rgb, 255:red, 233; green, 229; blue, 229 }  ,fill opacity=1 ] (361,54.5) .. controls (361,51.19) and (363.69,48.5) .. (367,48.5) -- (445,48.5) .. controls (448.31,48.5) and (451,51.19) .. (451,54.5) -- (451,72.5) .. controls (451,75.81) and (448.31,78.5) .. (445,78.5) -- (367,78.5) .. controls (363.69,78.5) and (361,75.81) .. (361,72.5) -- cycle ;
% dots
\draw [line width=1.5]  [dash pattern={on 1.69pt off 2.76pt}]  (220.5,63.25) -- (235,63.75) ;
% dots
\draw [line width=1.5]  [dash pattern={on 1.69pt off 2.76pt}]  (340.5,63.25) -- (355,63.75) ;
%Rounded Rect [id:dp3103912505033375] 
% Stage 1
\draw  [color={rgb, 255:red, 74; green, 74; blue, 74 }  ,draw opacity=1 ][fill={rgb, 255:red, 233; green, 229; blue, 229 }  ,fill opacity=1 ] (572.5,54.5) .. controls (572.5,51.19) and (575.19,48.5) .. (578.5,48.5) -- (656.5,48.5) .. controls (659.81,48.5) and (662.5,51.19) .. (662.5,54.5) -- (662.5,72.5) .. controls (662.5,75.81) and (659.81,78.5) .. (656.5,78.5) -- (578.5,78.5) .. controls (575.19,78.5) and (572.5,75.81) .. (572.5,72.5) -- cycle ;
%Straight Lines [id:da7036448025452298]
% Stage 2 
\draw  [color={rgb, 255:red, 74; green, 74; blue, 74 }  ,draw opacity=1 ][fill={rgb, 255:red, 233; green, 229; blue, 229 }  ,fill opacity=1 ] (690,54.5) .. controls (690,51.19) and (692.69,48.5) .. (696,48.5) -- (774,48.5) .. controls (777.31,48.5) and (780,51.19) .. (780,54.5) -- (780,72.5) .. controls (780,75.81) and (777.31,78.5) .. (774,78.5) -- (696,78.5) .. controls (692.69,78.5) and (690,75.81) .. (690,72.5) -- cycle ;

% Text Node
\draw (167.5,97) node   {$x_{1}$};
% Text Node
\draw (290,127.5) node    {$\statevart$};
% Text Node
\draw (290,97) node    {$\recvart$};
% Text Node
\draw (409.5,127.5) node  (C)  {$\statevarT$};
% Text Node
\draw (409.5,97) node    {$\recvarT$};
% Text Node
\draw (619,97) node    {$x_{1},\{\twostageLDRVart\}_{t\in\set{2,T}}$};
% Text Node
\draw (735,97) node  (E)  {$\left\{\recvart\right\}_{t\in\set{2,T}}$};
% Text Node
\draw (620,63.5) node   [align=left] (D) {Stage $\displaystyle 1$};
% Text Node
\draw (737,63.5) node   [align=left] {Stage $\displaystyle 2$};
% Text Node
\draw (167.5,63.5) node   [align=left] (A) {Stage 1};
% Text Node
\draw (290,63.5) node   [align=left] {Stage $\displaystyle t$};
% Text Node
\draw (408,63.5) node   [align=left] (B) {Stage $\displaystyle T$};
% Text Node
\draw (512,20) node  {$\statevart=\twostageBFt$};

% MSARO 
\node [dotted_block, xshift = -0.05cm, yshift = 0.1cm, fit = (A) (B) (C)] {};
\node [dotted_block, xshift = -0.25cm, yshift = 0.1cm, inner xsep=5mm, fit = (D) (E)] {};

\draw[->,line width=2pt] (474,67) -- (550,67) ;
\draw[->,line width=1pt] (503,30) -- (515,63) ;

% Text Node
\draw (287,25) node   [align=left] {MSARO};
\draw (675,25) node   [align=left] {2ARO};

\end{tikzpicture}
    }
    \caption{Two-stage decision rules}
    \label{fig:two-stage-DR}
\end{figure}

A practical result of such an approximation is that the resulting problem can be solved using the existing solution methods for 2ARO. {\crev In the following sections, we present two possible choices for the decision rule $\twostageBFt$, respectively applicable to continuous and integer state variables. Together, they permit the approximation of an \MSARO{} problem with mixed-integer state variables via a 2ARO model. We remark that the nature of the recourse variables does not impact the reduction of an MSARO to a 2ARO, but plays an important role in the choice of an appropriate solution method for the resulting 2ARO model. In what follows, we illustrate the application of these decision rules and the algorithmic solution of ensuing 2ARO models. For ease of exposition, we present MSARO problems with only continuous and only integer state variables separately.}

\subsection{Two-stage Linear Decision Rules {\crev for MSAROs with Continuous State Variables}}
\label{sec:two-stage-LDR}
If the state variables are continuous,
 we can approximate them {\crev via a decision rule with} an affine {\crev form}. {\crev For $t \in [2,T]$, letting} $\twostageLDRBFt = \big(\twostageLDRBF_{t1}(\xisupt),\dots,\twostageLDRBF_{tK_t}(\xisupt)\big) : \R^{\ell^{t}} \rightarrow \R^{{\crev q_t \times} K_t}$ be a vector of {\crev chosen} basis functions, the \emph{two-stage LDR} is enforced by using 
\begin{equation}
\label{eq:Theta_LDR}  
\twostageBFt = {\crev \twostageLDRBFt \twostageLDRVart}
\end{equation}
 in \eqref{eq:2RO} {\crev where we use a compact matrix representation\footnote{\crev This representation is obtained, without loss of generality, by concatenating individual LDR restrictions applied to each state variable $x^{\texttt{s}}_{ti}(\xisupt)$ for $i\in [q_t]$. For example, consider an instance where there are two state variables at stage $t=2$ and the history consists of two components $\xi^2 = (\xi_1,\xi_2)\in \R^2$ (with the convention that $\xi_1 = 1$). Then the  decisions rules $x_{21}(\xi^2)=\beta^1_{21}\xi_1+\beta^1_{22}\xi_2$ and $x_{22}(\xi^2)=\beta^2_{21}\xi_1+\beta^2_{22}\xi_2$ can be represented in the more compact matrix form with $K_2 = 4$ using the concatenated decision vector $\beta\in \R^4$ and the basis function matrix $\Phi_2(\xi^2) = 
 \begin{bmatrix}
\xi_1 & \xi_2 & 0 & 0 \\
0 & 0 & \xi_1 & \xi_2
\end{bmatrix}.$
} for notational convenience.}
 The resulting {\crev 2ARO} is {\crev written as:} 
 \bsub
 \label{eq:2RO-linear}
\begin{align}
	\twostageLDRUB := \min \ \ & c_1^\top x_1 + 
	\max_{\xiT\in\xisupp}\  
	\min_{\recvar\in\mathcal{X}({\crev \statevar_1}, \twostageLDR, \xiT)} \sum_{t\in \set{2,T}} {\stateCt}^\top {\crev\twostageLDRBFt\twostageLDRVart} + {\recCt}^\top \recvar_t \label{eq:2RO-linear-obj}\\
	\text{s.t.}\quad & x_1\in X_1 & \\
	&\twostageLDRVart \in \R^{K_t}\qquad\ \ t\in\set{2,T}.
\end{align}
\esub
{\crev As reviewed in Section \ref{sec:intro-lit},} the 2ARO problem \eqref{eq:2RO-linear} {\crev can either be solved by means of approximation (e.g., $K$-adaptability \citep{subramanyam2020k}, uncertainty set partitioning \citep{postek2016multistage,bertsimas2016multistage}, Neur2RO \citep{dumouchelle2023neur2ro}) }{\crev or exactly, most notably via the commonly used constraint-and-column generation (C\&CG) method, initially proposed by {\crev \cite{zeng2013solving}}, which we detail next.}

The C\&CG method draws on the fact that not all realizations in $\xisupp$ contribute to the {\crev worst-case} objective value. It then strives to identify {\it necessary {\crev realizations}} by starting from a smaller uncertainty set and gradually expanding it. {\crev This} leads to the generation of new columns and constraints, respectively corresponding to recourse variables and second-stage constraints for the newly identified {\crev realization}. More specifically, consider a relaxation of  problem \eqref{eq:2RO-linear} where instead of  
the uncertainty set $\xisupp$, a potentially empty subset 
$\XiMP \subseteq \xisupp$ is used {\crev to obtain} 
the following master problem, {\crev which we denote by $\MP$}: 
\bsub
\label{eq:2RO-MP}
\begin{align}
\min\ \ & c_1^\top x_1 +  \eta  \label{eq:2RO-MP-obj}  \\
\text{s.t.}\ \ & \eta \geq \sum_{t\in \set{2,T}} \Big({\stateCtFixed}{\crev (\xisuptval)}^\top {\crev\twostageLDRBFtVal\twostageLDRVart }+ {\recCtFixed}{\crev (\xisuptval)}^\top \recvartval\Big) && \xiTval\in\XiOpt \label{eq:2RO-MP-obj-const}\\
&{\crev {\crev \recAt} \recvartval + {\crev \stateAt \twostageLDRBFtVal\twostageLDRVart }+ {\crev \stateBt} \statevar_1 \leq b_{t}(\xisuptval)} && {\crev t = 2, {\crev \xiTval} \in\XiMP} \label{eq:2RO-MP-cut1-1}  \\
&{\crev \recAt} \recvartval + {\crev \stateAt \twostageLDRBFtVal\twostageLDRVart} + {\crev \stateBt} {\crev\twostageLDRBFtminusoneVal\twostageLDRVartminusone} \leq b_{t}(\xisuptval) && t\in \set{{\crev 3},T}, {\crev \xiTval}\in\XiMP \label{eq:2RO-MP-cut1-2}  \\
&({\crev\twostageLDRBFtVal\twostageLDRVart}, \recvartval)\in X_{t}(\xisuptval)&&   t\in \set{2,T}, {\crev \xiTval}\in\XiMP\label{eq:2RO-MP-cut2}\\
& x_1\in X_1,\quad\eta\in\R\\
&\twostageLDRVart \in \R^{K_t} && t\in\set{2,T},
\end{align}
\esub
where  
$\XiOpt \subseteq\XiMP$ is the set of identified {\crev necessary realizations} $\xiTval$ for which there exists a feasible first-stage solution $\hat{\twostageLDR}$ such that the feasibility space $\mathcal{X}({\crev \statevar_1},\hat{\twostageLDR},\xiTval)$ is nonempty. {\crev If the recourse variables are all continuous (i.e., $p^\texttt{i}_t=0, t\in\set{2,T}$), then  \eqref{eq:2RO-MP} is a linear program, otherwise it is a mixed-integer linear program. We remark that if $\XiOpt$ is empty at initialization, a valid lower bound on the $\eta$ variable can be added to the model to avoid unboundedness. If $\MP$ is infeasible for any $\XiMP \subseteq \xisupp$, the 2ARO model \eqref{eq:2RO-linear} is proven to be infeasible. Next, we consider the cases where $\MP$ is feasible and bounded.} 

The optimal objective value of {\crev $\MP$} is a lower bound on $\twostageLDRUB$, {\crev the optimal objective value of \eqref{eq:2RO-linear}}. 
{\crev To obtain} an exact solution to problem \eqref{eq:2RO-linear}, we {\crev may} need to {\crev gradually} expand $\XiMP$ with {\crev necessary realizations} (and consequently the set of recourse variable copies $\recvartval$).
At convergence, solving $\MP$ {\crev should} return an optimal solution $(\hat{x}_1,\hat{\twostageLDR},\hat{\eta})$ such that $\hat{\eta}$ accurately measures the worst-case second-stage cost over the complete uncertainty set $\xisupp$ ({\crev or conclude infeasibility of \eqref{eq:2RO-linear}). To check whether th{\crev is} convergence criterion is satisfied, we solve the adversarial problem for a given first-stage solution $(\hat{\twostageLDR},\hat{x}^\texttt{s}_1)$, which results in the following subproblem:
\begin{equation}
\SP := \max_{\xiT\in\xisupp}\ \Bigg\lbrace \sum_{t\in \set{2,T}} {\crev \stateCt}^\top {\crev\twostageLDRBFt\hat{\twostageLDR}_t}+   
\min_{\recvar\in\mathcal{X}({\crev \hat{x}^\texttt{s}_1},\hat{\twostageLDR}, \xiT)} \sum_{t\in \set{2,T}} {\crev \recCt}^\top \recvar_t\Bigg\rbrace. \label{eq:2RO-SP} 
\end{equation}
If $\hat{\eta} = \SP$,  then $\hat{\eta}$ exactly measures the worst-case cost of the second-stage problem and $(\hat{x}_1,\hat{\twostageLDR},\hat{\eta})$ is an optimal solution to problem \eqref{eq:2RO-linear}, i.e., $\twostageLDRUB =  c_1^\top\hat{x}_1 +  \hat{\eta} $. Otherwise, let $\hat{\xiparval}^T$ be the optimal solution of subproblem \eqref{eq:2RO-SP} if it is feasible, or a scenario such that $\mathcal{X}({\crev \hat{x}^\texttt{s}_1},\hat{\twostageLDR},\hat{\xiparval}^T)$ is an empty set. We create new variables $\recvar_{t,\hat{\xiparval}^{\crev T}}, t\in\set{2,T}$ and update the set of  {\crev necessary realizations} with $\XiMP = \XiMP \cup \big\{\hat{\xiparval}^T\big\}$, and accordingly update {\crev \eqref{eq:2RO-MP-cut1-1}-\eqref{eq:2RO-MP-cut2}}.
If  subproblem \eqref{eq:2RO-SP} is feasible, then we update $\XiOpt = \XiOpt \cup \big\{\hat{\xiparval}^{\crev T}\big\}$ and  \eqref{eq:2RO-MP-obj-const} as well.
In this case, constraints \eqref{eq:2RO-MP-obj-const}-\eqref{eq:2RO-MP-cut2} make up the optimality cuts. If subproblem \eqref{eq:2RO-SP} is infeasible, then {\crev \eqref{eq:2RO-MP-cut1-1}-\eqref{eq:2RO-MP-cut2}} act as feasibility cuts. We repeat solving the master problem and the subproblem in a cutting-plane fashion and add the appropriate cuts until $\XiMP$ includes all the {\crev necessary realizations} and $\hat{\eta} = \SP$, {\crev or $\MP$  becomes infeasible}. 

{\crev 
In general, subproblem \eqref{eq:2RO-SP} can be numerically challenging to solve especially if the inner minimization problem contains integer variables and/or non-linearities in $\xiT$. However, there are certain practical cases in which \eqref{eq:2RO-SP} can be reformulated as a mixed-integer linear program and then directly be given to an optimization solver.  In the context of our subproblem, one  notable example of this is provided in Remark \ref{rem:monoloithicSP}. Further, extending this case to allow integer recourse variables ($p^\texttt{i}_t\neq 0, t\in\set{2,T}$), while keeping the other assumptions in Remark \ref{rem:monoloithicSP}, \cite{zhao2012exact} proposed a nested constraint-and-column generation algorithm. 
In addition, we note that it is not necessary to solve the subproblem exactly at each iteration; the solution process can be stopped as soon as a violated cut for the master problem is identified. 
Recent studies have also considered the use of neural networks in order to approximate the inner minimization problem \citep{dumouchelle2023neur2ro}, especially to handle integer recourse variables and non-linearities.
\begin{remark}\label{rem:monoloithicSP}
 Consider an MSARO where the uncertainty set is a polytope, the basis functions $\twostageLDRBFt$ are chosen to be affine in $\xisupt$ for all $t\in\set{2,T}$, all the recourse variables are continuous  
(i.e., $p^\texttt{i}_t=0, t\in\set{2,T}$) and we have fixed parameters associated with the state variables, i.e., $\stateCt= c^{\texttt{s}}_t,\stateAt= A^{\texttt{s}}_t,\stateBt=B^{\texttt{s}}_t, D_t^\texttt{s}(\xi^t)=D^{\texttt{s}}_t$ for all $t\in\set{2,T}$ and $\xi^t \in \xisuppt$. In this case, using linear programming duality or the KKT optimality conditions (under a further relatively complete recourse assumption for the 2ARO problem \eqref{eq:2RO-linear}) yields a {\crev monolithic} bilinear continuous optimization problem, which can further be linearized using big-M constraints to obtain a mixed-integer linear program (see Appendix \ref{sec:monolithic-model} for details).  
\end{remark} 
}

{\crev
We further remark that the presence of decision-rule design variables $\beta$ may make the master problem $\MP$ numerically challenging to solve, compared to traditional use cases of C\&CG. In that regard, recently proposed enhancement ideas can be employed, such as the inexact C\&CG algorithm proposed by \cite{tsang2023inexact}. In their framework, master problems are solved to a given relative optimality gap which gradually reduces to zero over the course of the algorithm. 
Additionally, it is possible to design the basis functions $\Phi_t(\cdot)$ to use a small information basis rather than the entire history $\xi^t$ so that the number of design variables is reduced. This would be especially helpful for problems with larger number of decision stages. 
}

Another consideration concerning the exact solution of \eqref{eq:2RO-linear} using the C\&CG algorithm is its convergence. Indeed, the C\&CG  may not have finite convergence in general. {\crev If $\xisupp$ is a finite set, then the finite convergence is straightforward, otherwise more conditions are needed to ensure this property.}  
\cite{zeng2013solving} proved  finite convergence, considering a (bounded) polyhedral uncertainty set, for problems with only right-hand-side uncertainty represented as an affine function of uncertain parameters. However, in problem \eqref{eq:2RO-linear}, basis functions $\twostageLDRBFt$ appear as coefficients of $\twostageLDRVart$, both in the objective function \eqref{eq:2RO-linear-obj} (in the term ${\stateCt}^\top \twostageLDRBFt$) and the second-stage constraints (terms $A_t^{\texttt{s}}(\xisupt)\twostageLDRBFt, B_t^{\texttt{s}}(\xisupt)\twostageLDRBFtminusone$ and $D_t^{\texttt{s}}(\xisupt) \twostageLDRBFt$). The following proposition provides a sufficient condition for finite convergence of the C\&CG algorithm (proof is given in Appendix \ref{sec:proofs}) in this more general context.}

\begin{proposition}\label{prop:two-stage-ccg}
	{\crev Consider an \MSARO{} with only right-hand-side uncertainty,  continuous recourse, and (bounded) polyhedral uncertainty set. If the basis functions $\twostageLDRBFt$ are chosen to be affine in $\xisupt$ for all $t \in \set{2,T}$, the C\&CG algorithm converges to  $\twostageLDRUB$    
 in a finite number of iterations.}
\end{proposition}

In case the conditions of Proposition \ref{prop:two-stage-ccg} are not satisfied, the C\&CG algorithm still converge{\crev s but} asymptotically to an optimal solution of problem \eqref{eq:2RO-linear} {\crev if it is feasible (since it is bounded under the boundedness assumption imposed on the original MSARO problem)}.

{\crev 
Lastly, we note that in the case of continuous recourse, a linear decision rule can be applied to the recourse decision variables as well resulting in the LDR approach proposed by \cite{ben2004adjustable}.}
Since recourse variables $\recvart$ in \eqref{eq:2RO-linear} are fully adjustable, it immediately follows that
$$\exactB \leq \twostageLDRUB \leq \LDRUB,$$ 
{\crev where $\LDRUB$ refers to the bound obtained from the commonly used LDR approach. 
}

\subsection{Two-stage Piecewise-constant Decision Rules {\crev for MSAROs with Integer State Variables}}
\label{subsec:DecisionRulesDiscrete}

In this section, we study the application of two-stage decision rules in  another special case of \MSARO{}s, where each state variable $x^\texttt{s}_{ti}(\xisupt) \in \Z, i\in[q_t], t\in\set{2,T}$ is a bounded integer with a given domain $[\underline{\kappa}_{ti},\overline{\kappa}_{ti}]$ {\crev where boundedness follows from Assumption \ref{ass:bounded}}. 
We {\crev enforce} the \emph{two-stage piecewise-constant decision rule} (PCDR) {\crev by using $x^\texttt{s}_{ti}(\xisupt)=\Theta_{ti}(\xi^t,\beta^i_t)$ where} 

\begin{equation}
\label{eq:Theta_PCDR}
{\crev \Theta_{ti}(\xi^t,\beta^i_t)} = \left\lbrace
\begin{array}{lll}
\underline{\kappa}_{ti}&\qquad\quad& \twostageKBFti\in {\crev \K^i_{t1}} \\
\underline{\kappa}_{ti} + 1 && \twostageKBFti\in {\crev \K^i_{t2}} \\
\vdots&& \vdots \\
\overline{\kappa}_{ti}&&\twostageKBFti\in{\crev\K^i_{tJ_i}},
\end{array}
\right.
\end{equation}
where ${\crev \K^i_{tj}}\subset {\crev[-1,1]}, j\in\set{J_{\crev i}}$ are disjoint sets with $\bigcup_{j\in\set{J_{\crev i}}}{\crev\K^i_{tj}}={\crev[-1,1]}$, and  $\twostageKBFti:\R^{\ell^{t}} \times \R^{K_t}  \rightarrow {\crev[-1,1]}$ are functions defining the policy {\crev for $t \in [2,T]$ and $i\in[q_t]$}. Semantically, the PCDR partitions the {\crev interval} ${\crev[-1,1]}$ into subsets, and then assigns an integer value to each partition.

A special case of PCDRs can be defined by restricting the form of all mappings $\twostageKBFti $ to be a linear function of the {\crev decision rule design variables $\beta$}. In the following, we show that such a decision rule results in a model that is structurally very similar to \eqref{eq:2RO-linear}, thus it {\crev is amenable to} the C\&CG method. {\crev Let, without loss of generality, ${\crev \K^{i}_{t1}}=[{\crev a^i_{t1},b^{i}_{t1}}]$, and ${\crev \K^{i}_{tj}}=({\crev a^{i}_{tj},b^{i}_{tj}}], j\in[J_{\crev i}]\setminus \{1\}$ be intervals {\crev with $a^i_{t1} = -1$,  $b^{i}_{t,J_i} = 1$, and $a^{i}_{tj} = b^{i}_{t,j-1}$ for $j \in [2,J_i]$}. Let further, {\crev for $t \in [2,T],$} } $\twostageKBFti = \twostageKBFAffti^\top\twostageLDR_{t}^i$ {\crev be chosen as} an affine function of basis functions $\twostageKBFAffti$. 
This implies, in particular, that {\crev$-1\leq \twostageKBFAffti^\top\twostageLDR_{t}^i\leq 1$ which is enforced through robust constraints}. Then problem \eqref{eq:2RO} becomes:
\bsub
\label{eq:2RO-PLDR-deteq}
\begin{alignat}{2}
\twostagePLDRUB := \min \ & c_1^\top x_1 + \SPPLDR \\
\text{s.t.} \ & x_1\in X_1 \\
& \twostageLDRVart \in \R^{K_t} && t\in\set{2,T} \\
& {\crev -1 \leq \twostageKBFAffti^\top\twostageLDR_{t}^i \leq 1} \qquad && {\crev t\in\set{2,T}, i\in [q_t], \xiT\in \xisupp} \label{eq:2RO-PLDR-masterSI}
\end{alignat}
\esub
where:
\bsub
\label{eq:2RO-PLDR-SP}
\begin{align}
	\SPPLDR := \max_{\xiT\in\xisupp}\ & \sum_{t\in {\crev [2,T]} }{\crev{c^\texttt{s}_{t}(\xisupt)}^\top x^{\texttt{s}}_{t} }+ 
	\min_{\recvar\in\mathcal{X}({\crev \statevar_1},\twostageLDR, \xiT)} \sum_{t\in \set{2,T}} && \hspace*{-0.7cm}{\crev {\recCtFixed}(\xisupt)^\top \recvar} \label{eq:2RO-PLDR-obj}\\
	\text{s.t.}\ & \sum_{j\in\set{J_{\crev i}}} (\underline{\kappa}_{ti} + j-1)\upsilon_{tij} = \twostagePLDRVart && t\in\set{2,T}, i\in\set{q_t}\label{eq:2RO-PLDR-const1}\\
	& \sum_{j\in\set{J_{\crev i}}} \omega_{tij} = \twostageKBFAffti^\top\twostageLDR_{t}^i&&t\in\set{2,T}, i\in\set{q_t}\label{eq:2RO-PLDR-const2}\\
	& ( {\crev a^{i}_{tj}}+\epsilon_{\crev j})\upsilon_{tij} \leq \omega_{tij} \leq {\crev b^{i}_{tj}}\upsilon_{tij}&&t\in\set{2,T}, i\in\set{q_t}, j\in\set{J_{\crev i}}\label{eq:2RO-PLDR-const3}\\
	& \sum_{j\in\set{J_{\crev i}}} \upsilon_{tij} = 1 && t\in\set{2,T}, i\in\set{q_t}\label{eq:2RO-PLDR-const4}\\
	& \upsilon_{tij} \in\{0,1\}&&t\in\set{2,T}, i\in\set{q_t}, j\in\set{J_{\crev i}}.\label{eq:2RO-PLDR-const5}
\end{align}
\esub
{\crev with $\epsilon_1 =0$.}
The PCDR is modeled using the auxiliary variables $\upsilon_{tij}, \omega_{tij}$ and constraints \eqref{eq:2RO-PLDR-const1}-\eqref{eq:2RO-PLDR-const5}. {\crev Variables $\upsilon$ determine in which interval the quantity $\twostageKBFAffti^\top\twostageLDR_{t}^i$ falls, in accordance with the variables $\omega$, and the integer values assigned to the state variables.   Here, $\epsilon_{\crev j}$ is added to the lower bound on $\omega_{tij}$ to ensure that partitions $\K_{tj}$ are disjoint. We remark that, it is possible to choose $\epsilon_{\crev j}=0$ for $j\in [J_i]$, in which case the intervals would intersect at their boundaries. In this case, whenever the quantity $\twostageKBFAffti^\top\twostageLDR_{t}^i$ is a boundary point, the model allows assigning either one of the corresponding integer values to the associated state variable. The objective function then dictates that the solution leading to the worst objective value is chosen.}

{\crev Similar to Remark \ref{rem:monoloithicSP}, in certain special cases, we can reformulate problem \eqref{eq:2RO-PLDR-SP} as a monolithic mixed-integer linear program. 
\begin{remark}
\label{rem:monoloithicSPinPCDR}
 Consider an MSARO where the uncertainty set is a polytope, the basis functions $\twostageKBFAffti$ are chosen to be affine in $\xisupt$ for all $t\in\set{2,T}$ and $i\in [q_t]$, all the recourse variables are continuous,  and we have fixed parameters associated with the state variables, i.e., $\stateCt= c^{\texttt{s}}_t,\stateAt= A^{\texttt{s}}_t,\stateBt=B^{\texttt{s}}_t, D_t^\texttt{s}(\xi^t)=D^{\texttt{s}}_t$ for all $t\in\set{2,T}$ and $\xi^t \in \xisuppt$. In this case, using linear programming duality yields a {\crev monolithic} bilinear continuous optimization problem whose objective function involves products between variables $x^{\texttt{s}}_{it}$ and {\crev linear programming} dual variables of the inner minimization problem. Thanks to constraints \eqref{eq:2RO-PLDR-const1}, variables $x^{\texttt{s}}_{it}$ can be substituted for a weighted sum of binary variables $\upsilon$, as such these bilinear terms can be linearized
 using big-M constraints to obtain a mixed-integer linear program.  
\end{remark} 

Further, in applying the C\&CG method to model \eqref{eq:2RO-PLDR-deteq}, robust constraints \eqref{eq:2RO-PLDR-masterSI} will appear in the master problem. If the uncertainty set is a polytope and the basis functions $\twostageKBFAffti$ are chosen to be affine in $\xisupt$ for all $t\in\set{2,T}$ and $i\in [q_t]$ then these semi-infinite constraints can be reformulated as a finite set of linear constraints using classical robust optimization techniques based on {\crev linear programming} duality. 
} We also note that the arguments {\crev presented} in Section \ref{sec:two-stage-LDR} imply similarly that the C\&CG algorithm converges asymptotically to the optimal solution of $\twostagePLDRUB$ {\crev if it is feasible}.

{\crev Lastly, the methods presented in Sections \ref{sec:two-stage-LDR} and \ref{subsec:DecisionRulesDiscrete} can be combined to address MSAROs with mixed-integer state variables.

\begin{remark} 
\label{rem:2AROprimal_MixedIntegerState}
For an MSARO with mixed-integer state variables, the application of linear and piecewise constant decision rules, given by equations 
\eqref{eq:Theta_LDR} and \eqref{eq:Theta_PCDR}, to the continuous and integer state variables, respectively, yields a 2ARO approximation. The resulting model is presented in detail in Appendix \ref{app:sec:2AROprimal}. This model is similarly amenable to the C\&CG method for exact solution but can also benefit from other 2ARO solution methods from the literature.
\end{remark}
}
%%%%%%%%%%%%%%%%%%%%%%%%%%%%%%%%
%%%%%%%%%%%%%%%%%%%%%%%%%%%%%%%%
%%%%%%%%%%%%%%%%%%%%%%%%%%%%%%%%
%%%%%%%%%%%%%%%%%%%%%%%%%%%%%%%%
\section{Dual Bounding}
\label{sec:NA-dual}
In this section, we introduce a new dual problem that provides a lower bound for \MSARO{} problems with mixed-integer recourse. {\crev Due to the existence of integer variables, we rely on Lagrangian duality techniques, where we create a Lagrangian relaxation and optimize over the Lagrangian dual multipliers. As the MSARO involves constraints corresponding to every realization of uncertainty, Lagrangian multipliers are functions of uncertainty, as such they are usually high (possibly infinite) dimensional. To overcome the difficulty in their optimization, we propose to apply decision rule restrictions to Lagrangian multipliers, leveraging ideas rooted
in the MSP literature \citep{kuhn2011primal, daryalal2020lagrangian}. 
In deriving
a Lagrangian dual of the MSARO problem, we \emph{choose} a probability distribution with the support
as the uncertainty set, and use the associated density function to scale the constraints to be dualized, resulting in expectation terms in the objective function of the relaxation. Accordingly, we obtain a dual approximation of MSARO in the form of a two-stage stochastic program. This probability distribution-based approach has several benefits, most notably the possibility of leveraging state-of-the-art stochastic programming techniques to solve the dual problem. However, the quality of the resulting dual bound depends on the probability distribution used while developing the dual formulation, as previously observed by \cite{kuhn2011primal} for MSAROs with continuous recourse. With the aim of identifying the strongest such dual bound, we formally pose a distribution optimization problem {(\crev akin to what was developed in \citep{hadjiyiannis2011scenario})} and develop appropriate solution methods (tailored to the nature of the recourse variables) for the
resulting {\crev distribution optimization} problem. To the best of our knowledge, numerical solution of such a bounding problem and the quality of the obtained bounds have not been studied before. 

In what follows, in Section \ref{subsec:NAdual}, we introduce the nonanticipative reformulation of the MSARO problem and its Lagrangian dual. In Section \ref{subsec:LDDR}, we present the restricted Lagrangian dual problem and define the associated distribution optimization problem for which we develop solution methods in Section \ref{sec:solution}. Lastly, in Section \ref{subsec:DecomposableNA}, we propose an alternative dual problem, which can be weaker in terms of the quality of the obtained bound but has computational advantages thanks to its decomposable structure.
}

\subsection{Nonanticipative Dual of the \MSARO{}}
\label{subsec:NAdual}

The {\crev nonanticipative} (NA) dual is based on a reformulation of the \MSARO{} problem where we create a copy of decision variables {\crev for every stage and} every realization, and explicitly enforce nonanticipativity constraints. To this end, we introduce the copy variables $ y(\xiT)=(y_1(\xiT),\dots,y_T(\xiT))$ {\crev for all $\xiT \in \xisupp$} as perfect information variables depending on the entire realization $\xiT= (\xiparval_1,\dots,\xiparrand_T) $. {\crev We denote the decision variables by $y$ instead of $x$ used in previous sections to emphasize the fact that they are perfect information variables. We also note that we do not need to distinguish state and recourse variables in this section, thus vectors $y$ involve all the decisions.}
We can {\crev then obtain} the {\crev NA reformulation of the} \MSARO{} problem \eqref{eqs:msaro-mono} as:
\bsub
\label{eq:NA-reform}
\begin{align}
\min  \ \ & z  \label{eq:obj-NA-reform} \\
\text{s.t.} \ \ & \sum_{t\in \set{T}}c_t(\xisupt)^\top y_t(\xiT) \leq z &&  \xiT\in\xisupp \label{eq:const-obj-NA-reform}\\
& \At y_t(\xiT) + \Bt y_{t-1}(\xiT)\leq b_t(\xisupt) &&  t\in {\crev [2,T]},\ \xiT\in\xisupp \label{eq:const-state-NA-reform}\\
& \Dt y_t(\xiT)\leq d_t(\xisupt) &&   t\in \set{T},\ \xiT\in\xisupp \label{eq:const-rec-NA-reform}\\
& y_t(\xiT) = y_t(\xiprimeT) &&  t\in \set{T},\ \xiT,\xiprimeT\in\xisupp \text{ with }  \xisupt=\xiprimesupt \label{eq:const-reform-nonanticipative-lddro}\\
& y_t(\xiT)\in\R^{n_t-n^\texttt{i}_t}\times\Z^{n^\texttt{i}_t} &&   t\in \set{T},\ \xiT\in\xisupp. \label{eq:const-domain-NA-reform} 
\end{align}
\esub
Constraints  \eqref{eq:const-reform-nonanticipative-lddro} are \emph{nonanticipativity constraints} which ensure that at stage $t$ for every partial realization of $\xiT$, the decisions made are consistent (i.e., the decisions made in all realizations sharing the history $\xisupt$ are the same).

{\crev 

\begin{remark}
\label{remark:RedundantNAconsts}
The nonanticipativity constraints are redundant for stage $T$, however, we include them in model \eqref{eq:NA-reform} for notational convenience. In our implementation for the numerical results presented in Section \ref{sec:experiments}, we exclude those redundant constraints.
\end{remark}

In deriving a Lagrangian relaxation, we will first scale the constraints to be relaxed. As the constraints correspond to uncertainty realizations, we \emph{choose} the scaling factors in such a way that they induce a probability measure over the support $\xisupp$. To this end, we let $\dist$ denote this probability measure such that $\dist(\xisupp)=1$, which we interchangeably refer to as the probability distribution. 
We define $p^{\dist}: \xisupp \rightarrow \R^+$ as the associated density function. Lastly, we let $\mathcal{P}^{>}:=\{\dist\ |\ p^{\dist}(\xiT)>0, \xiT\in\xisupp\}$, i.e., every $\dist\in\mathcal{P}^{>}$ has a density function assigning a strictly positive value to all $\xiT\in\xisupp$. 
}

Before introducing the NA dual derived from \eqref{eq:NA-reform}, we present the following lemma, which we use to reformulate the nonanticipativity constraints (proof in Appendix \ref{sec:proofs}).

\begin{lemma}\label{lem:NA-const-alt}
For {\crev any $ \dist \in \mathcal{P}^{>}$}, constraints \eqref{eq:const-reform-nonanticipative-lddro} are equivalent to the following:
\begin{equation}\label{eq:NA-const-Exp-form}
    y_t(\xiT) = \Exp_{\xiprimeT\sim\dist}\left[y_t(\xiprimeT) \ \Big|\ \xiprimesupt = \xisupt \right],\quad  t\in \set{T},\ \xiT\in\xisupp.
\end{equation}
\end{lemma}

{\crev One advantage of this reformulation is the reduction in the number of constraints and in turn in the number of dual multipliers to be introduced. Let, to this end, $ \lambda_t(\cdot):\R^{\ell^t}\rightarrow\R^{n_t}$ for $t \in [T]$, where $\Exp_{\xiT\sim\dist} [\lambda_t(\xiT)]<+\infty $, be the dual functionals to be used in relaxing the nonanticipativity constraints \eqref{eq:NA-const-Exp-form}. Let further the feasibility space for $\xiT\in \xisupp$ be $Y(\xiT) := \big\{ \big(z,y_1(\xiT),\dots,y_T(\xiT)\big) : \eqref{eq:const-obj-NA-reform}-\eqref{eq:const-rec-NA-reform} \text{ and } \eqref{eq:const-domain-NA-reform}\big\} $.
}
{\crev Then, after scaling \eqref{eq:NA-const-Exp-form} with the probability densities associated with the distribution $\dist$, we obtain} the following NA Lagrangian relaxation problem:
\bsub
\label{eq:msaro--LR-obj}
\begin{align}
\mathcal{L}^{\textsc{NA}}_{\textsc{LR}}(\dist,\lambda_1(\cdot),\dots,\lambda_T(\cdot)) := \ \min  \ & z + \sum_{t\in[T]} \Exp_{\xiT\sim\dist} \Big[ {\lambda_t(\xiT)}^\top \Big( y_t(\xiT)-\Exp_{\xiprimeT\sim\dist}\left[y_t(\xiprimeT) \big| \xiprimesupt = \xisupt \right] \Big) \Big] \label{eq:obj-reform-stochastic-ip-LR} \\
\text{s.t.} \ & \big(z,y_1(\xiT),\dots,y_T(\xiT)\big) \in Y(\xiT)\qquad \xiT\in\xisupp.
\end{align}
\esub
{\crev whose optimal objective value} yields 
a lower bound for the original \MSARO{} problem {\crev for any $\dist \in \mathcal{P}^{>}$}. 
The NA Lagrangian dual problem aims to find the best bound among all such Lagrangian relaxation bounds:
\begin{equation}\label{eq:NA_bound}
\mathcal{L}^{\textsc{NA}}(\dist) := \max_{\lambda_1(\cdot),\dots,\lambda_T(\cdot)}\ \ \mathcal{L}^{\textsc{NA}}_{\textsc{LR}}(\dist,\lambda_1(\cdot),\dots,\lambda_T(\cdot)).
\end{equation}
The following proposition shows that regardless of the choice of $\dist$, $ \mathcal{L}^{\textsc{NA}}(\dist) $ is an exact dual bound for \MSARO{} problems {\crev with \emph{continuous recourse}} (proof in Appendix \ref{sec:proofs}).
{\crev
\begin{proposition}\label{prop:exactness-NA}
Let $\dist$ be any probability measure in {\crev $\mathcal{P}^{>}$}. For \MSARO{} problems with \emph{continuous recourse}, \eqref{eq:NA_bound} is a strong dual of \eqref{eqs:msaro-mono}, i.e., $\mathcal{L}^{\textsc{NA}}(\dist)=\exactB$.
\end{proposition}
}

We remark that our construction of the dual postulates that we multiply the constraints with a density function whose support is $\xisupp$. Therefore, the strictly positive density property of {\crev $\mathcal{P}^{>}$} is necessary for the exactness of our formulation.
Proposition \ref{prop:exactness-NA} suggests that the solution of problem \eqref{eq:NA_bound} gives an exact dual bound for \eqref{eqs:msaro-mono} (hence for \eqref{eq:NA-reform}) if all decision variables are continuous. Although this bound is not necessarily exact in the case of mixed-integer recourse, the potential of leveraging the literature of multistage stochastic programming in achieving a dual bound for \MSARO{} is quite appealing.

The objective function of the {\crev NA} Lagrangian relaxation problem \eqref{eq:msaro--LR-obj} contains (conditional) expectations of decision variables $y_t(\xiT)$, which is computationally challenging. {\crev Because of our initial assumption that $X_t(\xisupt)$ are bounded (Assumption \ref{ass:bounded}) we have, by letting $\tilde{Y}(\xiT):={\rm{proj}}_{y_1(\cdot),\ldots,y_T(\cdot)}Y(\xiT)$, that $\Exp_{\xiT\sim\dist} [\text{diam}(\tilde{Y}(\xiT))]<+\infty$. Since we further have that $\Exp_{\xiT\sim\dist} [\lambda_t(\xiT)]<+\infty $, we can apply Lemma 1 of \citep{daryalal2020lagrangian} and} replace the expectation term for $t$ in the objective function \eqref{eq:obj-reform-stochastic-ip-LR} with:
$$ \Exp_{\xiT\sim\dist} \Big[ \Big( \lambda_t(\xiT)-\Exp_{\xiprimeT\sim\dist}\left[\lambda_t(\xiprimeT) \big| \xiprimesupt = \xisupt \right] \Big)^\top y_t(\xiT) \Big]. $$

For given $\lambda_t(\xiT)$, this exchange allows us to compute the coefficients of $y_t(\xiT)$ in the Lagrangian relaxation problem. Still, the optimal form of the dual functionals $\lambda_t(\cdot)$ need to be determined, making the problem \eqref{eq:NA_bound} computationally intractable. 
In the next section, we restrict these Lagrangian multipliers to follow LDRs and obtain a  restricted dual problem with decision variables of {\crev smaller (finite)} dimension. Furthermore, this new dual problem is amenable to well-known solution techniques from the literature of two-stage stochastic programming which are designed to approximately solve a problem with expectation in the objective function.

\subsection{Lagrangian Dual Decision Rules}
\label{subsec:LDDR}
\noindent 
{\crev We restrict the NA Lagrangian dual problem \eqref{eq:NA_bound} for a given $\dist\in \mathcal{P}^{>}$ by enforcing LDRs on the Lagrangian multipliers,} referred to as Lagrangian dual decision rules (LDDRs). For a set of pre-determined basis functions $\Psi_t:\xisupp\rightarrow \R^{n_t\times K_t}$ and LDR decision variables $\naldr_t\in\R^{K_t}$, we restrict the form of $\lambda_t(\xiT)$ at stage {\crev $t\in [T]$} as follows:
$$ \lambda_{t}(\xiT) = \naBFt \naldr_t, $$
which gives us a restricted NA Lagrangian dual problem with respect to $\dist$:
\begin{equation}\label{eq:NA_Lag_restricted}
\mathcal{L}^{\textsc{NA}}_R(\dist) := \max_{\naldr_1,\dots,\naldr_T}\ \ \mathcal{L}^{\textsc{NA}}_{\textsc{LR}}(\dist,\naBF_{1}(\xiT)\naldr_{1},\dots, \naBF_{T}(\xiT) \naldr_{T}).
\end{equation}
Since problem \eqref{eq:NA_Lag_restricted} is a restriction of \eqref{eq:NA_bound}, we have $\mathcal{L}^{\textsc{NA}}_R(\dist)\leq\mathcal{L}^{\textsc{NA}}(\dist)$.

Using Lemma 2 in \citep{daryalal2020lagrangian}, the primal characterization of $\mathcal{L}^{\textsc{NA}}_R(\dist)$ is:
    \bsub
    \label{eq:restricted-NA-reform-primal-char}
    \begin{align}
    \min  \ \ & z  \label{eq:restricted-obj-NA-reform-primal-char} \\
    \text{s.t.} \ \ & \big(z,y_1(\xiT),\dots,y_T(\xiT)\big) \in \conv{Y(\xiT)} && \xiT\in\xisupp   \label{eq:restricted-const-remain-primal-char}\\
    & \Exp_{\xiT\sim\dist} \Bigg[ \naBFt^\top \Big(y_t(\xiT)-\Exp_{\xiprimeT\sim\dist}\big[y_t(\xiprimeT)\ |\ \xiprimesupt = \xisupt\big]\Big) \Bigg]  = \boldsymbol{0}&&  t\in\set{T}  \label{eq:restricted-const-NA-primal-char}
    \end{align}
    \esub
{\crev For a given $\dist$, comparing \eqref{eq:restricted-NA-reform-primal-char} to the primal characterization of the (unrestricted) NA Lagrangian dual problem \eqref{eq:NA_bound} (provided in \eqref{eq:NA-reform-primal-char} of Appendix~\ref{sec:proofs}), it is clear that, the former is a relaxation of the latter since constraints \eqref{eq:restricted-const-NA-primal-char} are an aggregation of their counterpart \eqref{eq:const-NA-primal-char}.} Consequently, unlike $\mathcal{L}^{\textsc{NA}}(\dist)$, even for \MSARO{} {\crev with continuous recourse},  $\mathcal{L}^{\textsc{NA}}_R(\dist)$ is not {\crev necessarily} a strong bound. Furthermore, 
due to constraints \eqref{eq:restricted-const-NA-primal-char}, the strength of the restricted NA Lagrangian dual bound depends on the choice of probability measure $\dist$.  {\crev A similar observation was made by \cite{kuhn2011primal} and \cite{hadjiyiannis2011scenario} concerning their dual bound for MSARO problems \emph{with only continuous variables}}. This observation motivates us to optimize over the probability distribution $\dist$ and LDR variables $\naldr$ to find the best such dual bound. As such, we propose to solve a {\crev distribution optimization (DO) problem over the set of probability distributions $\mathcal{P}^{>}$, defined as follows:}
\begin{equation}\label{eq:dro}
    {\crev \nu^{\textsc{NA-DO}}_{R}} := \sup_{\dist\in \mathcal{P}^{>}}\ \mathcal{L}^{\textsc{NA}}_R(\dist)
\end{equation}
{\crev where, as before, $\mathcal{P}^{>}=\{\dist\ |\ p^{\dist}(\xiT)>0, \xiT\in\xisupp\}$.}

In linear and mixed-integer programming, strict inequalities such as the ones required for $\dist\in {\crev \mathcal{P}^{>}}$ (the strictly positive density property for all $\xiT\in\xisupp$) often cause numerical and theoretical difficulties, thus 
are not desirable. To avoid these inequalities, in the following discussions, we modify the {\crev DO} model to improve its numerical behaviour. Denote by $ \mathcal{P}^{\geq}$
a superset of $\mathcal{P}^{>}$ that also admits distributions {\crev that} allow $p(\xiT)=0$ for some $\xiT\in\xisupp$.  Consider the problem:
\begin{equation}\label{eq:dro-p-bar}
    {\crev \bar{\nu}^{\textsc{NA-DO}}_R} := \max_{\dist\in \mathcal{P}^{\geq}} \ \mathcal{L}^{\textsc{NA}}_R(\dist).
\end{equation}
Because \eqref{eq:dro-p-bar} is a relaxation of \eqref{eq:dro}, it yields an upper bound for {\crev $\nu^{\textsc{NA-DO}}_R$}. Thus, it does not immediately follow that such a bound is a valid lower bound for the {\crev optimal value of the} MSARO problem, $\exactB$. 
The following proposition shows that \eqref{eq:dro-p-bar} indeed leads to a valid dual (lower) bound (see Appendix \ref{sec:proofs} for proof).
\begin{proposition}\label{eq:Pbar-valid}
{\crev $\bar{\nu}^{\textsc{NA-DO}}_R$} is a lower bound for $\exactB$.
\end{proposition}

{\crev Hereafter,} we refer to \eqref{eq:dro-p-bar} {\crev or its equivalent explicit form}
\bsub
\label{eq:restricted-NA-dual}
\begin{align}
{\crev \bar{\nu}^{\textsc{NA-DO}}_R} = \max_{\dist, \naldr}\ \ & \mathcal{Q}(\dist, \naldr)\\
\text{s.t.} \ \ & \naldr_t\in\R^{K_t}\quad t\in\set{T}\\
& \dist \in \mathcal{P}^{\geq},
\end{align}
\esub
{\crev as the {\crev DO} problem,} where
\bsub
\label{eq:Q-general}
\begin{align}
    \mathcal{Q}(\dist, \naldr) := \min\ \ & z + \sum_{t\in[T]} \Exp_{\xiT\sim\dist} \Bigg[ \bigg(\Big( \naBFt-\Exp_{\xiprimeT\sim\dist}\left[\naBF_{t}(\xiprimeT) \big| \xiprimesupt = \xisupt \right] \Big)\naldr_{t}\bigg)^\top y_{t}(\xiT)\Bigg]\label{eq:Q-general-obj}\\ 
    \text{s.t.}\ \ & \big(z,y_1(\xiT),\dots,y_T(\xiT)\big) \in {\crev Y(\xiT)}\qquad \xiT\in\xisupp.
\end{align}
\esub
Problem \eqref{eq:Q-general} is a \emph{two-stage 
stochastic program} (2SP) and can benefit from its rich literature. In the next section, we build on well-known stochastic programming techniques to design a decomposition method {\crev to solve the DO problem \eqref{eq:restricted-NA-dual}.}

\subsection{Solving the {\crev DO} Problem}\label{sec:solution}
{\crev There are two main challenges associated with the solution of the DO problem: the expectation terms in the objective of \eqref{eq:Q-general-obj} and the max-min structure in \eqref{eq:restricted-NA-dual}.}

If the uncertainty set of {\crev the \MSARO{} problem}, $\xisupp$, is not discrete, objective function \eqref{eq:Q-general-obj} includes the expectation of a nonsmooth concave function. {\crev Further, even when $\xisupp$ is discrete calculating the expectation term exactly can be prohibitive from a computational point of view.} The literature of two-stage stochastic programming addresses such a difficulty by means of sampling-based approaches that replace the expectation in the objective function with the average of a sample {\crev drawn from the underlying distribution} and has favourable theoretical convergence results (see e.g., \cite{shapiro2009lectures}). {\crev We follow the sample average approximation (SAA) approach in overcoming the first challenge and show in Section~\ref{subsubsec:SAA} that it leads to a valid dual bound for the MSARO problem.}

{\crev With regards to the second challenge, we propose a cutting plane algorithm that iteratively constructs improving approximations of $\mathcal{Q}(\dist, \naldr)$ through its supporting hyperplanes obtained by solving the SAA approximation of \eqref{eq:Q-general}. We present this general algorithm in Section~\ref{sec:disc-uncertainty-sets} and propose an alternative monolithic formulation in the special case of MSARO with continuous recourse in Section~\ref{sec:cont-uncertainty-sets}. Figure \ref{fig:NA_Dual} 
summarizes the methods presented in this section 
for solving the DO model.

}

\begin{figure}[htbp]
    \centering
    \scalebox{0.82}{
    \tikzset{every picture/.style={line width=0.75pt}} %set default line width to 0.75pt        

\begin{tikzpicture}[x=0.75pt,y=0.75pt,yscale=-1,xscale=1]
%uncomment if require: \path (0,347); %set diagram left start at 0, and has height of 347

%Straight Lines [id:da4261270415339775] 
\draw    (152.5,93) -- (138.13,122.21) ;
\draw [shift={(137.25,124)}, rotate = 296.19] [color={rgb, 255:red, 0; green, 0; blue, 0 }  ][line width=0.75]    (10.93,-3.29) .. controls (6.95,-1.4) and (3.31,-0.3) .. (0,0) .. controls (3.31,0.3) and (6.95,1.4) .. (10.93,3.29)   ;
%Straight Lines [id:da25569303454227377] 
\draw    (317.5,93) -- (335.68,121.81) ;
\draw [shift={(336.75,123.5)}, rotate = 237.74] [color={rgb, 255:red, 0; green, 0; blue, 0 }  ][line width=0.75]    (10.93,-3.29) .. controls (6.95,-1.4) and (3.31,-0.3) .. (0,0) .. controls (3.31,0.3) and (6.95,1.4) .. (10.93,3.29)   ;
%Curve Lines [id:da8321581431465329] 
\draw  [dash pattern={on 4.5pt off 4.5pt}]  (479.75,134.5) .. controls (455.75,133.5) and (407.75,84.5) .. (417.25,108) .. controls (426.42,130.68) and (349.45,89.57) .. (332.37,80.37) ;
\draw [shift={(330.75,79.5)}, rotate = 28.3] [color={rgb, 255:red, 0; green, 0; blue, 0 }  ][line width=0.75]    (10.93,-3.29) .. controls (6.95,-1.4) and (3.31,-0.3) .. (0,0) .. controls (3.31,0.3) and (6.95,1.4) .. (10.93,3.29)   ;
%Rounded Rect [id:dp6330273773861236] 
\draw  [color={rgb, 255:red, 74; green, 74; blue, 74 }  ,draw opacity=1 ][fill={rgb, 255:red, 233; green, 229; blue, 229 }  ,fill opacity=1 ] (50.5,132.5) .. controls (50.5,129.19) and (53.19,126.5) .. (56.5,126.5) -- (221.5,126.5) .. controls (224.81,126.5) and (227.5,129.19) .. (227.5,132.5) -- (227.5,150.5) .. controls (227.5,153.81) and (224.81,156.5) .. (221.5,156.5) -- (56.5,156.5) .. controls (53.19,156.5) and (50.5,153.81) .. (50.5,150.5) -- cycle ;
%Rounded Rect [id:dp24359293897994405] 
\draw  [color={rgb, 255:red, 74; green, 74; blue, 74 }  ,draw opacity=1 ][fill={rgb, 255:red, 233; green, 229; blue, 229 }  ,fill opacity=1 ] (246.75,132.5) .. controls (246.75,129.19) and (249.44,126.5) .. (252.75,126.5) -- (417.75,126.5) .. controls (421.06,126.5) and (423.75,129.19) .. (423.75,132.5) -- (423.75,150.5) .. controls (423.75,153.81) and (421.06,156.5) .. (417.75,156.5) -- (252.75,156.5) .. controls (249.44,156.5) and (246.75,153.81) .. (246.75,150.5) -- cycle ;
%Curve Lines [id:da9992466083193056] 
\draw    (360.25,295.37) .. controls (377.62,274.35) and (372.59,261.68) .. (361.87,249) ;
\draw [shift={(358.87,297)}, rotate = 310.82] [color={rgb, 255:red, 0; green, 0; blue, 0 }  ][line width=0.75]    (10.93,-3.29) .. controls (6.95,-1.4) and (3.31,-0.3) .. (0,0) .. controls (3.31,0.3) and (6.95,1.4) .. (10.93,3.29)   ;
%Curve Lines [id:da6175686459694207] 
\draw    (310.51,297) .. controls (295.07,280.59) and (299.18,266.05) .. (306.68,250.68) ;
\draw [shift={(307.51,249)}, rotate = 116.57] [color={rgb, 255:red, 0; green, 0; blue, 0 }  ][line width=0.75]    (10.93,-3.29) .. controls (6.95,-1.4) and (3.31,-0.3) .. (0,0) .. controls (3.31,0.3) and (6.95,1.4) .. (10.93,3.29)   ;
%Rounded Rect [id:dp23548429787653002] 
\draw  [dash pattern={on 4.5pt off 4.5pt}] (227,217.4) .. controls (227,201.72) and (239.72,189) .. (255.4,189) -- (411.35,189) .. controls (427.03,189) and (439.75,201.72) .. (439.75,217.4) -- (439.75,302.6) .. controls (439.75,318.28) and (427.03,331) .. (411.35,331) -- (255.4,331) .. controls (239.72,331) and (227,318.28) .. (227,302.6) -- cycle ;
%Straight Lines [id:da3469544917226628] 
\draw  [dash pattern={on 4.5pt off 4.5pt}]  (137.25,156.5) -- (137.72,183.5) ;
\draw [shift={(137.75,185.5)}, rotate = 269.01] [color={rgb, 255:red, 0; green, 0; blue, 0 }  ][line width=0.75]    (10.93,-3.29) .. controls (6.95,-1.4) and (3.31,-0.3) .. (0,0) .. controls (3.31,0.3) and (6.95,1.4) .. (10.93,3.29)   ;
%Straight Lines [id:da7652900117058608] 
\draw  [dash pattern={on 4.5pt off 4.5pt}]  (337.25,156.5) -- (337.72,183.5) ;
\draw [shift={(337.75,185.5)}, rotate = 269.01] [color={rgb, 255:red, 0; green, 0; blue, 0 }  ][line width=0.75]    (10.93,-3.29) .. controls (6.95,-1.4) and (3.31,-0.3) .. (0,0) .. controls (3.31,0.3) and (6.95,1.4) .. (10.93,3.29)   ;
%Straight Lines [id:da4512346598936602] 
\draw  [dash pattern={on 4.5pt off 4.5pt}]  (523.25,93.5) -- (523.71,120) ;
\draw [shift={(523.75,122)}, rotate = 268.99] [color={rgb, 255:red, 0; green, 0; blue, 0 }  ][line width=0.75]    (10.93,-3.29) .. controls (6.95,-1.4) and (3.31,-0.3) .. (0,0) .. controls (3.31,0.3) and (6.95,1.4) .. (10.93,3.29)   ;
%Rounded Rect [id:dp8246133830079057] 
\draw  [color={rgb, 255:red, 74; green, 74; blue, 74 }  ,draw opacity=1 ][fill={rgb, 255:red, 233; green, 229; blue, 229 }  ,fill opacity=1 ] (144.5,57) .. controls (144.5,52.03) and (148.53,48) .. (153.5,48) -- (312.5,48) .. controls (317.47,48) and (321.5,52.03) .. (321.5,57) -- (321.5,84) .. controls (321.5,88.97) and (317.47,93) .. (312.5,93) -- (153.5,93) .. controls (148.53,93) and (144.5,88.97) .. (144.5,84) -- cycle ;
%Rounded Rect [id:dp5844969701319849] 
\draw  [color={rgb, 255:red, 74; green, 74; blue, 74 }  ,draw opacity=1 ][fill={rgb, 255:red, 233; green, 229; blue, 229 }  ,fill opacity=1 ] (419.5,57) .. controls (419.5,52.03) and (423.53,48) .. (428.5,48) -- (610.5,48) .. controls (615.47,48) and (619.5,52.03) .. (619.5,57) -- (619.5,84) .. controls (619.5,88.97) and (615.47,93) .. (610.5,93) -- (428.5,93) .. controls (423.53,93) and (419.5,88.97) .. (419.5,84) -- cycle ;

% Text Node
\draw  [dash pattern={on 4.5pt off 4.5pt}]  (62.5,199.5) .. controls (62.5,194.53) and (66.53,190.5) .. (71.5,190.5) -- (206.5,190.5) .. controls (211.47,190.5) and (215.5,194.53) .. (215.5,199.5) -- (215.5,210.5) .. controls (215.5,215.47) and (211.47,219.5) .. (206.5,219.5) -- (71.5,219.5) .. controls (66.53,219.5) and (62.5,215.47) .. (62.5,210.5) -- cycle  ;
\draw (143,205) node  [color={rgb, 255:red, 0; green, 0; blue, 0 }  ,opacity=1 ] [align=center] {
{\crev Bilinear program} \eqref{eq:NA-Dual-Mono} 
};
% Text Node
\draw  [dash pattern={on 4.5pt off 4.5pt}]  (484.5,132) .. controls (484.5,127.03) and (488.53,123) .. (493.5,123) -- (554.5,123) .. controls (559.47,123) and (563.5,127.03) .. (563.5,132) -- (563.5,143) .. controls (563.5,147.97) and (559.47,152) .. (554.5,152) -- (493.5,152) .. controls (488.53,152) and (484.5,147.97) .. (484.5,143) -- cycle  ;
\draw (524,137.5) node  [color={rgb, 255:red, 0; green, 0; blue, 0 }  ,opacity=1 ] [align=center] {
 Sampling 
};
% Text Node
\draw (385.22,271) node    {$\hat{\beta }$};
% Text Node
\draw (335.25,206) node  [color={rgb, 255:red, 0; green, 0; blue, 0 }  ,opacity=1 ] [align=left] {
Cutting-plane method
};
% Text Node
\draw (139,141.5) node   [align=left] {Continuous {\crev recourse}};
% Text Node
\draw (335.25,141.5) node   [align=left] {Mixed-integer {\crev recourse}};
% Text Node
\draw (266.25,273) node  [color={rgb, 255:red, 0; green, 0; blue, 0 }  ,opacity=1 ] [align=left] {
Cut \eqref{eq:NA-Dual-CP-cut}
};
% Text Node
\draw    (232.75,229) .. controls (232.75,224.58) and (236.33,221) .. (240.75,221) -- (425.75,221) .. controls (430.17,221) and (433.75,224.58) .. (433.75,229) -- (433.75,242) .. controls (433.75,246.42) and (430.17,250) .. (425.75,250) -- (240.75,250) .. controls (236.33,250) and (232.75,246.42) .. (232.75,242) -- cycle  ;
\draw (333.25,235.5) node   [align=center] {
 \ Solve master problem \eqref{eq:NA-Dual-CP-MP} \ 
};
% Text Node
\draw    (246.25,305.5) .. controls (246.25,300.53) and (250.28,296.5) .. (255.25,296.5) -- (413.25,296.5) .. controls (418.22,296.5) and (422.25,300.53) .. (422.25,305.5) -- (422.25,316.5) .. controls (422.25,321.47) and (418.22,325.5) .. (413.25,325.5) -- (255.25,325.5) .. controls (250.28,325.5) and (246.25,321.47) .. (246.25,316.5) -- cycle  ;
\draw (334.25,311) node   [align=center] {
 \ Solve subproblem \eqref{eq:NA-Dual-CP-SP} \ 
};
% Text Node
\draw (233,70) node   [align=left] {Discrete uncertainty set};
% Text Node
\draw (519.5,60.5) node   [align=left] {Continuous/{\crev Large discrete}};
% Text Node
\draw (519.5,80.5) node   [align=left] {uncertainty set};

\end{tikzpicture}
    }
    \caption{Solution methods for the {\crev DO} problem 
    }
    \label{fig:NA_Dual}
\end{figure}

\subsubsection{{\crev Sample average approximation (SAA)}}
\label{subsubsec:SAA}

{\crev Let $\Omega \subseteq \xisupp$ be a finite subset of the uncertainty set $\xisupp$. We define the set of probability measures $\mathcal{P}^{\geq}_{\Omega}$ such that $\dist_\Omega\in \mathcal{P}^{\geq}_{\Omega}$ implies that $\dist_\Omega(\Omega)=1$ and the associated density function  has value zero for any {\crev realization} not in $\Omega$: 
$$p^{\dist_\Omega}(\xiT) =
    0,\ \xiT\in\xisupp\setminus\Omega,\text{ and }
    p^{\dist_\Omega}(\xiT)\geq 0, \ \xiT\in\Omega.$$

Since $\mathcal{P}^{\geq}_{\Omega} \subseteq \mathcal{P}^{\geq}$ and $\bar{\nu}^{\textsc{NA-DO}}_R$ is obtained by maximizing $\mathcal{L}^{\textsc{NA}}_R(\dist)$ over $\dist \in \mathcal{P}^{\geq}$, we have, for any $\dist_\Omega\in \mathcal{P}^{\geq}_{\Omega}$, that $\mathcal{L}^{\textsc{NA}}_R(\dist_\Omega)
\leq
{\crev \bar{\nu}^{\textsc{NA-DO}}_R}
\leq
\exactB$. It then follows that 
\begin{equation}
\label{eq:DO_Omega}
\max_{\dist_{\Omega}\in \mathcal{P}_{\Omega} } \mathcal{L}^{\textsc{NA}}_R(\dist_\Omega)
\leq
{\crev \bar{\nu}^{\textsc{NA-DO}}_R}
\leq
\exactB.
\end{equation}
We finally have that
\begin{align}
\label{eq:discrete_DO}
\max_{\dist_{\Omega}\in \mathcal{P}^{\geq}_{\Omega} } \mathcal{L}^{\textsc{NA}}_R(\dist_\Omega) = \max_{\dist_{\Omega}, \naldr}\ \ & \mathcal{Q}(\dist_{\Omega}, \naldr)\\
\nonumber\text{s.t.} \ \ & \naldr_t\in\R^{K_t}\quad t\in\set{T}\\
\nonumber& \dist_{\Omega} \in \mathcal{P}^{\geq}_{\Omega},
\end{align}
where for any $\dist_{\Omega}\in \mathcal{P}^{\geq}_{\Omega}$ the expectation terms in the objective function \eqref{eq:Q-general-obj}, used in calculating $\mathcal{Q}(\dist_{\Omega}, \naldr)$, are replaced by their sample average over $\Omega$.

In the remainder of this section, we omit the notation $\Omega$ and write our models with a finite discrete uncertainty set $\xisupp$ which can either be the full uncertainty set of the MSARO problem or a set of realizations sampled from it. 
}

\subsubsection{{\crev MSARO with mixed-integer recourse}}\label{sec:disc-uncertainty-sets}
In a discrete uncertainty set with realizations $\xiTval\in\xisupp$,  $\mathcal{P}^{\geq}$ can be modeled by a set of (in)equalities, that is $ \mathcal{P}^{\geq} = \big\{\pDistVar\in \R_+^{|\xisupp|}\ |\ \boldsymbol{1}^\top \pDistVar = 1\big\}$, where a vector $\pDistVar^{\dist}\in\mathcal{P}^{\geq}$ characterizes the probability measure $\dist$ such that for all $\xiTval\in\xisupp$,  $\pDist$ is the probability of {\crev realization} $\xiTval$ with respect to $\dist$.  
{\crev Define $\pNormalization:=\sum_{\substack{\xiprimeTval\in\xisupp:\\ \xiprimesuptval = \xisuptval}} \pDistprime$ as the sum of the probabilities of all realizations sharing the same history  $\xisuptval$ up to stage $t$ and let $\rho^{\dist}_{\xiprimeTval\mid\xi^t}$ be the conditional probability of realization $\xiprimeTval\in \xisupp$ given history $\xi^t$. More precisely, given $\xisuptval$,  for $\xiprimeTval\in\xisupp$ with $\xiprimesuptval = \xisuptval$ we have that $\rho^{\dist}_{\xiprimeTval\mid\xi^t}=\frac{\pDistprime}{\pNormalization}$ if $\pNormalization>0$ and that $\rho^{\dist}_{\xiprimeTval\mid\xi^t}=0$ otherwise. Then,} the expectation term in the objective function \eqref{eq:Q-general-obj} can be written as:
$$\sum_{t\in[T]}\sum_{\xiTval\in\xisupp} \pDist\bigg(\Big( \Psi_{t}(\xiTval)-\sum_{\substack{\xiprimeTval\in\xisupp:\\ \xiprimesuptval = \xisuptval}} {\crev \rho^{\dist}_{\xiprimeTval\mid\xi^t}} \Psi_{t}(\xiprimeTval) \Big)\naldr_{t}\bigg)^\top y_{t}(\xiTval).$$
Let $\Psi_{tk}(\xiTval)$ be the $k^{\text{th}}$ column of the matrix $\Psi_{t}(\xiTval)$.
Then $\mathcal{Q}(\dist, \naldr)$ can be expressed as:
\begin{align*}
	\min \ \ & \displaystyle z + \sum_{t\in[T]} \sum_{\xiTval\in\xisupp} \pDist\bigg(\sum_{k\in[K_t]}\Big( \Psi_{tk}(\xiTval)-\sum_{\substack{\xiprimeTval\in\xisupp:\\ \xiprimesuptval = \xisuptval}} {\crev \rho^{\dist}_{\xiprimeTval\mid\xi^t}}\Psi_{tk}(\xiprimeTval)  \Big)\alpha_{tk}\bigg)^\top y_{t}(\xiTval)\\
	\text{s.t.} \ \ & \big(z,y_1(\xiTval),\dots,y_T(\xiTval)\big) \in Y(\xiTval)\quad \xiTval\in\xisupp.
\end{align*}
Let  $\gamma_{tk\xiTval} := \pDist\alpha_{tk}$ and {\crev define }$\displaystyle\beta_{t\xiTval} := \sum_{k\in[K_t]}\Big( \gamma_{tk\xiTval}\Psi_{tk}(\xiTval)-\sum_{\substack{\xiprimeTval\in\xisupp:\\ \xiprimesuptval = \xisuptval}} {\crev \rho^{\dist}_{\xiprimeTval\mid\xi^t}}\gamma_{tk\xiTval}\Psi_{tk}(\xiprimeTval)  \Big)$ {\crev as the coefficient vector of variables $y_{t}(\xiTval)$}. Then with a change of variables in the {\crev DO problem} \eqref{eq:restricted-NA-dual} we have:
\bsub
\label{eq:DRO-reform}
\begin{alignat}{2}
	{\crev \bar{\nu}^{\textsc{NA-DO}}_R} = 
 \max_{\pDistVar, \alpha, \gamma, \beta}\ & \mathcal{Q}(\beta)\\
	\text{s.t.}& \sum_{\xiTval\in\xisupp}\pDistVar_{\xiTval} = 1\label{eq:probs}\\
	& \gamma_{tk\xiTval} = \pDistVar_{\xiTval}\alpha_{tk}&& t\in\set{T}, k\in\set{K_t},\xiTval\in\xisupp \label{eq:NADO_GammaConsts} \\
	&\displaystyle{\crev\sum_{\substack{\xiprimeTval\in\xisupp:\\ \xiprimesuptval = \xisuptval}} \pDistprime}\beta_{t\xiTval} - \hspace*{-0.1cm} \sum_{k\in[K_t]} \hspace*{-0.1cm} \Big( {\crev \sum_{\substack{\xiprimeTval\in\xisupp:\\ \xiprimesuptval = \xisuptval}} \hspace*{-0.15cm} \pDistprime}\gamma_{tk\xiTval}\Psi_{tk}(\xiTval) && - \hspace*{-0.1cm}\sum_{\substack{\xiprimeTval\in\xisupp:\\ \xiprimesuptval = \xisuptval}} \hspace*{-0.1cm} \pDistVar_{\xiTval}\gamma_{tk\xiTval}\Psi_{tk}(\xiprimeTval)\Big) = \boldsymbol{0} \nonumber \\
 & && t\in\set{T}, \xiTval\in\xisupp \label{eq:NADOBetaConsts}\\
    &{\crev -\boldsymbol{M} \pDistVar_{\xiTval}\leq \beta_{t\xiTval} \leq \boldsymbol{M} \pDistVar_{\xiTval}} &&  {\crev t\in\set{T}, \xiTval\in\xisupp} \label{eq:NADOLogicalConsts} \\
	& \pDistVar \in \R_+^{|\xisupp|} \\
	& \naldr_t\in\R^{K_t},\ {\crev \gamma_t\in\R^{K_t \times |\xisupp|},\ \beta_t \in\R^{|\xisupp|\times n_t}}&&  t\in\set{T}\label{eq:bounds}
\end{alignat}
\esub
where
\begin{equation}\label{eq:NA-Dual-CP-SP}
\mathcal{Q}(\beta) = \min\Big\{ z + \sum_{t\in[T]} \sum_{\xiTval\in\xisupp} \beta_{t\xiTval}^\top y_{t}(\xiTval)\ \Big|\ \big(z,y_1(\xiTval),\dots,y_T(\xiTval)\big) \in Y(\xiTval),\ \xiTval\in\xisupp \Big\},
\end{equation}
{\crev and $\boldsymbol{M}\geq \boldsymbol{0}$ is a vector of sufficiently large numbers. }

{\crev In model \eqref{eq:DRO-reform}, constraints \eqref{eq:probs} along with the nonnegativity restrictions imposed on variables $\rho$ induce a probability distribution over the realizations $\xisupp$ with some realizations potentially assigned zero probability. For given $\xiTval$, if its associated probability $\pDistVar_{\xiTval}$ is zero, then constraints \eqref{eq:NADO_GammaConsts} and \eqref{eq:NADOLogicalConsts} imply, respectively, that all variables $\gamma$ and $\beta$ indexed by $\xiTval$ are zero, as such constraints \eqref{eq:NADOBetaConsts} trivially hold. Otherwise, since probability $\pDistVar_{\xiTval}>0$, we have that $\sum_{\substack{\xiprimeTval\in\xisupp:\\ \xiprimesuptval = \xisuptval}} \pDistprime>0$, as such constraints \eqref{eq:NADOBetaConsts} correctly impose the definition of $\beta$ by dividing all terms by $\sum_{\substack{\xiprimeTval\in\xisupp:\\ \xiprimesuptval = \xisuptval}} \pDistprime$ and plugging in the definition of $\gamma$ from constraints \eqref{eq:NADO_GammaConsts}. Then given coefficients $\beta$, \eqref{eq:NA-Dual-CP-SP} evaluates the expected value of the optimal decisions $y(\cdot)$. 

\begin{remark}
The choice of values for the vector $\boldsymbol{M}$ impacts the quality of the bound obtained from model \eqref{eq:DRO-reform}. In particular, in the limit case where $\boldsymbol{M}=\boldsymbol{0}$, one obtains the perfect information bound, that is, all nonanticipativity constraints are relaxed from the NA reformulation of the MSARO problem. Otherwise, for larger values of $\boldsymbol{M}$, the optimal value of \eqref{eq:DRO-reform} is lower bounded by the perfect information bound since choosing $\beta = \boldsymbol{0}$ is always feasible. One can therefore expect to obtain a better bound from \eqref{eq:DRO-reform}. 
\end{remark}
}

Model \eqref{eq:DRO-reform} can be solved via a cutting-plane method in which $\mathcal{Q}(\beta)$ is approximated by a set of linear inequalities. At each iteration, we solve the following bilinear program as the master problem:
\begin{equation}\label{eq:NA-Dual-CP-MP}
    \max_{{\rho}, {\alpha}, {\gamma}, {\beta}}\ \Big\lbrace \eta\ \Big|\ \eqref{eq:probs}-\eqref{eq:bounds},\ (\eta,{\beta}) \in \mathcal{H}\Big\rbrace,
\end{equation}
where $\eta$ is an auxiliary variable representing $\mathcal{Q}(\beta)$, and $\mathcal{H}$ is a set described by optimality cuts approximating $\mathcal{Q}({\beta})$. Note that, as $\beta$ only parameterizes the objective function of $\mathcal{Q}({\beta})$, i.e., it does not impact the feasibility space, there is no need for feasibility cuts. With $(\hat{\eta},\hat{{\beta}})$ returned from solving the master problem \eqref{eq:NA-Dual-CP-MP}, we solve the subproblem \eqref{eq:NA-Dual-CP-SP} to compute $\mathcal{Q}(\hat{{\beta}})$, resulting in $\hat{y}_{t}(\xiTval)$ as the optimal solution. If $\hat{\eta} \leq \mathcal{Q}(\hat{{\beta}}) $, we have found the optimal solution of the {\crev DO problem}. Otherwise we add the following optimality cut to the master problem:
\begin{equation}\label{eq:NA-Dual-CP-cut}
\eta \leq \mathcal{Q}(\hat{{\beta}}) + \sum_{t\in[T]} \sum_{\xiTval\in\xisupp} \Big(\beta_{t\xiTval}-\hat{\beta}_{t\xiTval}\Big)^\top \hat{y}_{t}(\xiTval).
\end{equation}
This procedure continues until no more optimality cuts are found. The objective function of the subproblem,  $\mathcal{Q}(\beta)$, is a concave function in $\beta$ (pointwise minimum of linear  functions with respect to $\beta$). Each cut \eqref{eq:NA-Dual-CP-cut} is a hyperplane approximating the subproblem from above. The cutting-plane method iteratively finds improving  approximations of $\mathcal{Q}(\beta)$. {\crev We note that it is not necessary to execute the cutting-plane procedure until convergence in order to obtain a valid dual bound. Indeed at each iteration of the algorithm the value $\mathcal{Q}(\hat{\beta})$ provides a valid dual bound for $\exactB$. }

{
\crev
\begin{remark}
In our implementation for the numerical results presented in Section \ref{sec:experiments}, rather than creating the copies of first-stage variables $y_1(\cdot)$ and relaxing their nonanticipativity constraints, we keep them as static variables in \eqref{eq:NA-Dual-CP-SP}, same as variable $z$.  
\end{remark}

}

\subsubsection{{\crev \MSARO{} with continuous recourse}}
\label{sec:cont-uncertainty-sets}
As a special case, we study \MSARO{} {\crev with continuous recourse.} Our goal here is to use this particular structure and derive a monolithic formulation as an alternative to the cutting-plane algorithm, to leverage off-the-shelf solvers. Denote by $\nadualobj$, $\nadualstate$ and $\nadualrec$, the dual variables associated with the set of constraints described by $Y(\xiTval)$ for given $\xiTval\in\xisupp$ (corresponding to  \eqref{eq:const-obj-NA-reform}, \eqref{eq:const-state-NA-reform} and \eqref{eq:const-rec-NA-reform}, respectively). The {\crev linear programming} dual of the inner minimization problem \eqref{eq:NA-Dual-CP-SP}, i.e., the subproblem of the cutting-plane algorithm, is:
\bsub
\label{eqs:Q-dual}
\begin{align}
	\mathcal{Q}^D(\beta) = \max\ \ &  \sum_{t\in[T]}\sum_{\xiTval\in\xisupp}b_t(\xiTval)^\top\nadualstate + \sum_{t\in[T]}\sum_{\xiTval\in\xisupp}d_t(\xisuptval)^\top\nadualrec \label{eq:obj-Q-dual} \\
	\text{s.t.} \ \ &\sum_{\xiTval\in\xisupp}\nadualobj = 1 && \label{eq:const-Q-dual-z}\\
	&-c_T(\xiTval)\nadualobj + A_T(\xiTval)^\top v_{T\xiTval} + D_{T\xiTval}^\top w_{T\xiTval}-\beta_{T\xiTval}=\boldsymbol{0} &&  \xiTval\in\xisupp \label{eq:const-Q-dual-y-T} \\
	& -\Ctval\nadualobj + \Atval^\top\nadualstate + \Dtval^\top\nadualrec +\nonumber\\ 
	& \qquad\qquad\Btplusoneval^\top\nadualstateplusone - \beta_{t\xiTval} =\boldsymbol{0} &&  t\in \set{T-1} , \xiTval\in\xisupp \label{eq:const-Q-dual-y} \\
	& \nadualobj\geq 0 && \xiTval\in\xisupp \\
	& \nadualstate,\nadualrec \leq \boldsymbol{0} &&   t\in \set{T}, \xiTval\in\xisupp. \label{eq:QP-dual-bounds}
\end{align}
\esub
Merging the two maximization problems in \eqref{eq:DRO-reform}, we get the monolithic \emph{bilinear} program:
\bsub
\label{eq:NA-Dual-Mono}
\begin{align}
	{\crev \bar{\nu}^{\textsc{NA-DO}}_R} = \max\ \ &  \sum_{t\in[T]}b_t(\xiTval)^\top\nadualstate + \sum_{t\in[T]}\sum_{\xiTval\in\xisupp}d_t(\xisuptval)^\top\nadualrec \\
	\text{s.t.} \ \
	&\eqref{eq:probs}-\eqref{eq:bounds}\\
	&\eqref{eq:const-Q-dual-z}-\eqref{eq:QP-dual-bounds}.
\end{align}
\esub

There is a large body of research on solution methods for bilinear problems that can be used in solving model \eqref{eq:NA-Dual-Mono}. Furthermore, many optimization solvers, such as MOSEK \citep{aps2020mosek} and  Gurobi \citep{gurobi}, offer off-the-shelf alternatives to solve problems of type \eqref{eq:NA-Dual-Mono}. {\crev Further, \eqref{eq:NA-Dual-Mono} can be solved heuristically to obtain a valid dual bound. For instance, one could alternate between optimizing over variables $\rho$ and $(\alpha,\beta,\gamma)$ (given that all bilinear terms involve variables $\rho$) to converge towards a local optimal solution.}

\subsection{Restricted Decomposable NA Dual}
\label{subsec:DecomposableNA}
In solving problem \eqref{eq:DRO-reform} using the cutting-plane method, we frequently optimize \eqref{eq:NA-Dual-CP-SP} to compute $\mathcal{Q}(\beta)$. This can become computationally demanding when there is a large  number of {\crev realizations}. Further, we cannot decompose $\mathcal{Q}(\beta)$ by  {\crev realizations}, as they are linked through the $z$ variables. In this section, we present an alternative NA {\crev reformulation of the MSARO} problem that can yield, {\crev through its associated DO problem,} a potentially weaker bound than the one provided by \eqref{eq:DRO-reform}. {\crev However, this alternative reformulation offers a computational advantage since in the framework of the cutting-plane method it leads to decomposable subproblems when calculating $\mathcal{Q}(\beta)$. To this end,
 in addition to the} decision variable copies $y_t(\xiT)$, we introduce copy variables $\dnavar(\xiT)$ and explicitly enforce them to be equal {\crev via nonanticipativity constraints.} For an assigned probability measure $\dist\in \mathcal{P}^{>}$, the alternative NA reformulation of the MSARO problem \eqref{eqs:msaro-mono} is:
\bsub
\label{eq:NA-reform-2}
\begin{align}
\min  \ \ & \Exp_{\xiT\sim\dist}\left[\dnavar(\xiT)\right]  \label{eq:obj-NA-reform-2} \\
\text{s.t.} \ \ & \sum_{t\in \set{T}}c_t(\xisupt)^\top y_t(\xiT) \leq \dnavar(\xiT) &&  \xiT\in\xisupp \label{eq:const-obj-NA-reform-2}\\
& \dnavar(\xiT) = \Exp_{\xiprimeT\sim\dist}\left[\dnavar(\xiprimeT)\right] && \xiT\in\xisupp \label{eq:const-z-NA-reform-2}\\
& \eqref{eq:const-state-NA-reform},  \eqref{eq:const-rec-NA-reform},  \eqref{eq:const-domain-NA-reform}, \eqref{eq:NA-const-Exp-form}.
\end{align}
\esub
Together, 
\eqref{eq:const-obj-NA-reform-2} and \eqref{eq:const-z-NA-reform-2} capture the semantics of the worst-case outcome, which is minimized in the objective function \eqref{eq:obj-NA-reform-2}. {\crev We remark that constraints \eqref{eq:const-z-NA-reform-2} are obtained from individual nonanticipativity constraints similarly to the derivation of constraints \eqref{eq:NA-const-Exp-form} provided in Lemma~\ref{lem:NA-const-alt}.}

Relaxation of the nonanticipativity constraints \eqref{eq:NA-const-Exp-form} and \eqref{eq:const-z-NA-reform-2} with the Lagrangian multipliers $ \lambda^y_t(\cdot)$ and $\lambda^z(\cdot)$ {\crev such that $\Exp_{\xiT\sim\dist} [\lambda^{y}_t(\xiT)]<+\infty$ and $\Exp_{\xiT\sim\dist} [\lambda^{z}_t(\xiT)]<+\infty$}, leads to the {\crev decomposable} NA Lagrangian dual problem 
$\mathcal{L}^{\textsc{DNA}}(\dist) = \displaystyle \max_{\lambda^y(\cdot), \lambda^z(\cdot)} \mathcal{L}^{\textsc{DNA}}_\textsc{LR}(\dist,\lambda^y_1(\cdot),\dots,\lambda^y_T(\cdot),\lambda^z(\cdot)),$ where
\bsub
\begin{align*}
\mathcal{L}^{\textsc{DNA}}_\textsc{LR}(\dist,\lambda^y_1(\cdot),\dots,\lambda^y_T(\cdot),\lambda^z(\cdot)) = 
\min  \ \ & \Exp_{\xiT\sim\dist} \Big[  \Big( 1 + \lambda^z(\xiT)-\Exp_{\xiprimeT\sim\dist}\left[\lambda^z(\xiprimeT)\right] \Big)\dnavar(\xiT) \Big] + \nonumber \\
& \sum_{t\in[T]} \Exp_{\xiT\sim\dist} \Big[ \Big( \lambda^y_t(\xiT)-\Exp_{\xiprimeT\sim\dist}\left[\lambda^y_t(\xiprimeT) \big| \xiprimesupt = \xisupt \right] \Big)^\top y_t(\xiT)  \Big]  \\
\text{s.t.} \ \ & \big(\dnavar(\xiT), y_1(\xiT),\dots,y_T(\xiT)\big) \in Y(\xiT)\qquad \xiT\in\xisupp,
\end{align*}
\esub
with $Y(\xiT)$ the scenario feasibility space described by constraints \eqref{eq:const-obj-NA-reform-2}, \eqref{eq:const-state-NA-reform}, \eqref{eq:const-rec-NA-reform}, and \eqref{eq:const-domain-NA-reform}. We note that, in the first expectation term in the objective function the inner expectation is not conditional, since constraints \eqref{eq:const-z-NA-reform-2} (and accordingly their associated dual functions) are not defined at every stage $t\in\set{T}$ as decision variable $\dnavar(\xiT)\in\R$ captures the cost of an entire {\crev realization $\xiT$}. {\crev Since in $\mathcal{L}^{\textsc{DNA}}_\textsc{LR}(\cdot)$, the objective function, constraints and variables are decomposable in $\xi^T$ they can be optimized individually.  }

{\crev In deriving the {\crev decomposable NA} Lagrangian relaxation problem $\mathcal{L}^{\textsc{DNA}}_\textsc{LR}(\cdot)$ we apply Lemma 1 of \citep{daryalal2020lagrangian} to obtain both expectation terms. We remark that the lemma requires the condition $\Exp_{\xiT\sim\dist} [\text{diam}({\rm{proj}}_{z(\cdot)}Y(\xiT))]<+\infty$ which is not naturally satisfied. However, as a result of Assumptions~\ref{ass:nonemptyfirst}-\ref{ass:bounded} and the compactness of the uncertainty set $\xisupp$ the optimal value of the MSARO problem is bounded. As such the functionals $z(\cdot)$ can be artificially bounded without changing the optimal value of $\mathcal{L}^{\textsc{DNA}}_\textsc{LR}(\cdot)$.   }

After substituting the decision rules $ \lambda^y_{t}(\xiT) = \naBFt \naldr^y_t,  t\in\set{T} $ and $ \lambda^z(\xiT) = \naBF_{T}(\xiT) \naldr^z $ in $\mathcal{L}^{\textsc{DNA}}(\dist)$, where $\naldr^z\in\R^{K_T}$, {\crev and merging with the optimization over the probability distributions $\dist \in \mathcal{P}^{\geq}$} the {\crev decomposable DO} problem is:
$$ {\crev \bar{\nu}^{\text{DNA-DO}}_R} := \max_{\dist \in \mathcal{P}^{\geq},\naldr^y_1,\ldots,\naldr^y_T, \naldr^z} \mathcal{L}^{\textsc{DNA}}_\textsc{LR}(\dist,\naBF_{1}(\xiT)\naldr^y_1,\dots,\naBF_{T}(\xiT) \naldr^y_T,\naBF_{T}(\xiT) \naldr^z).$$
This {\crev DO problem} can be solved using the same methods developed in the previous sections for the non-decomposable {\crev DO} problem. Due to the relaxation of the nonanticipativity constraints on $z(\cdot)$ variables, for given $\dist$ and $\lambda^y(\cdot)$, subproblem $\mathcal{L}^{\textsc{DNA}}_\textsc{LR}(\cdot)$ is a relaxation of $\mathcal{L}^{\textsc{NA}}_\textsc{LR}(\cdot)$. Therefore ${\crev \bar{\nu}^{\text{DNA-DO}}_R} \leq {\crev \bar{\nu}^{\text{NA-DO}}_R}$, i.e., {\crev $\bar{\nu}^{\text{NA-DO}}_R$} is a potentially stronger bound. However, the fact that $\mathcal{L}^{\textsc{DNA}}_\textsc{LR}(\cdot)$ is decomposable is highly desirable. Particularly for continuous {\crev or large discrete} uncertainty sets where we rely on sampling, and the quality of the bound varies based on the selected sample. Since the decomposable model can afford samples of larger sizes it can potentially yield better bounds compared to the bound obtained from non-decomposable model over a smaller sample. We explore this trade-off in our numerical section.

%%%%%%%%%%%%%%%%%%%%%%%%%%%%%%%%
%%%%%%%%%%%%%%%%%%%%%%%%%%%%%%%%
%%%%%%%%%%%%%%%%%%%%%%%%%%%%%%%%
%%%%%%%%%%%%%%%%%%%%%%%%%%%%%%%%
\section{Numerical Experiments}\label{sec:experiments}
We evaluate the performance of the proposed bounding framework over multistage versions of three classical decision-making problems under uncertainty: ($i$) the newsvendor problem, ($ii$) the location-transportation problem, and ($iii$) the capital budgeting problem. Depending on their structure, each problem is solved by using the appropriate models and methods described in Sections \ref{sec:primal} and \ref{sec:NA-dual}, illustrating the applicability of the developed concepts to a large array of problem classes.

\subsection{Benchmarks and Implementation Details}
To assess the quality of the primal and dual bounds, we measure the relative distance of the bounds from the true optimal value when an exact solution of \MSARO{} is available (in small-size instances). Otherwise, we report the optimality gap between the bounds obtained from the proposed methods, and compare it against a  gap {\crev from} traditional bounding methods if one exists. In {\crev problems with} continuous {\crev recourse}, we consider LDRs (i.e., their application to all decision variables) as the benchmark for the primal decision rules. On the dual side, we use the perfect information (PI) bound (denoted by $\PILB$) for comparison, which often can be conveniently evaluated for a general \MSARO{} problem{\crev, as well as the bound obtained by solving model \eqref{eqs:msaro-mono} using only the binding realizations identified from the primal decision rule solution}. The PI bound corresponds to the optimal objective value of the \MSARO{} problem reformulated as in \eqref{eq:NA-reform} without the nonanticipativity constraints \eqref{eq:const-reform-nonanticipative-lddro}, i.e., it finds the cost of every {\crev realization} in the uncertainty set individually, and then selects the one with the worst-case cost.

The algorithms are implemented in Python and use the Gurobi Optimizer 9.5.1 \citep{gurobi} as the mixed-integer/bilinear programming solver. The computational experiments are carried out on the Niagara supercomputer servers \citep{loken2010scinet, ponce2019deploying}. {\crev For instances with discrete uncertainty sets, the programs have a time limit of 1 hour. We report the valid lower/upper bound at the point of termination. For instances with continuous uncertainty sets, this time limit is extended to 10 hours.} As a common design choice for the basis functions of the LDRs, we use the uncertain parameters themselves, i.e., the standard basis functions. Any implementation nuances and enhancements used for improving the performance of the algorithms are discussed for each problem class in a dedicated section, along with the characteristics of the studied instances.

\subsection{Robust Newsvendor Problem}
In this section, we extend the two-stage newsvendor problem studied in \citep{xu2018improved} to the multistage setting. In this problem, a decision-maker (the newsvendor) needs to order from a set of items to be sold (only) at the next decision stage, with the objective of maximizing the worst-case profit over the planning horizon.
Let $d_{it}(\xisupt)$ be the uncertain demand of item $i\in\set{I}$ at stage $t\in\set{2,T}$, $c_i $ and $s_i$ the purchase and shortage costs of item $i$, respectively, and $r_i$ its sale price. To meet the customers' demands of stage $t$, at stage $t-1$ the decision-maker decides on the amounts to be ordered from each item, such that the total spending over the $T$ stages does not exceed a predetermined budget of $B$. Denote by $x_{it}(\xisupt)$ the decision variable for the amount of item $i$ ordered at stage $t\in\set{T-1}$. The multistage multi-item budgeted newsvendor problem is:
\bsub
\label{eq:newsvendor}
\begin{flalign}
	\max \ \ & z  \label{eq:obj-newsvendor} \\
	\text{s.t.} \ \ & z\leq \sum_{i\in \set{I}}\sum_{t\in \set{2,T}}y_{it}(\xisupt) &&  \xiT\in\xisupp \label{eq:obj-newsvendor-const}\\
	& y_{it}(\xisupt) \leq \big(r_i-c_i\big) x_{i,t-1}(\xisuptminusone) - r_i \big(x_{i,t-1}(\xisuptminusone) - d_i(\xisupt)\big) &&  i\in\set{I}, t\in\set{2,T}, \xisupt\in\xisuppt \label{eq:const-newsvendor-profit1}\\
	& y_{it}(\xisupt) \leq \big(r_i-c_i\big) x_{i,t-1}(\xisuptminusone) - s_i \big(d_i(\xisupt) - x_{i,t-1}(\xisuptminusone)\big) &&  i\in\set{I}, t\in\set{2,T},  \xisupt\in\xisuppt\label{eq:const-newsvendor-profit2}\\
	& \sum_{i\in\set{I}}\sum_{t\in \set{T-1}} x_{it}(\xisupt) \leq B&&   \xiTminusone\in\Xi^{T-1} \label{eq:const-newsvendor-budget}\\
	&x_t(\xisupt)\in\R_+^{I} &&t\in \set{T-1}, \xisupt\in\xisuppt.&&
\end{flalign}
\esub
where auxiliary variable $y_{it}(\xisupt)$ captures the profit from item $i$ at stage $t$, by means of 
\eqref{eq:const-newsvendor-profit1}- \eqref{eq:const-newsvendor-profit2}. Constraints \eqref{eq:const-newsvendor-budget} impose a budget of $B$ over the order amounts throughout the planning horizon. The objective function is the worst-case profit of the newsvendor, modeled via \eqref{eq:obj-newsvendor} and \eqref{eq:obj-newsvendor-const}. 

\subsubsection{Problem Instances}
Our instance generation loosely follows the procedure described in \citep{ardestani2021linearized} for the two-stage robust newsvendor problem. Parameters  
$r_i$, $s_i$ and $c_i$ are drawn uniformly from the intervals $[140,160]$, $[80,90]$ and $[50,70]$, respectively.  We consider a discrete uncertainty set modeled as a stagewise-dependent scenario tree with branching factor $\textsc{br}$ (i.e., every node of the tree prior to the leaves has $\textsc{br}$ many child nodes). Demand realizations $d_{it}(\xisupt), i\in\set{I}, t\in\set{2,T}$ at a child node are drawn uniformly from $[\mu_{it} - \sigma_{it},\mu_{it} + \sigma_{it}]$, where $\mu_{it}$ and $\sigma_{it}$ are  
uniformly drawn from the intervals $[20,40]$ and $[10,20]$. We have generated 26 small-size instances with $T\in\set{3,5}$, $I\in\set{2,5}$, $B\in\{100,150,200,250,300\}$, and $\textsc{br}\in\{2,3,4,5,10\}$, such that the number of {\crev realizations} $|\xisupp|=\textsc{br}^{T-1}$ is less than 150. Additionally, we have generated 18 large-size instances with $T\in\set{4,8}$. For $T=4$, the number of items  $I$ lies in the set $\set{3,5}$, with a budget $B\in\{200,300\}$ and $\textsc{br}\in\{10,15,20\}$. For $T\in\set{5,8}$, our instances have $I\in\{3,4\}$ items, budget of $B\in\{300,400,500,600\}$, and $\textsc{br}\in\set{3,6}$,  restricted to the cases with $|\xisupp|\leq 3000$.

\subsubsection{Quality of the Bounds}\label{sec:newsvendor-quality}
{\crev The small-size instances are easily optimized by solving model \eqref{eq:NA-reform} over all {\crev realizations} in the uncertainty set. In our case, this computation takes less than 3 seconds.} Using {\crev these optimal values}, we can examine the quality of a primal/dual bound by measuring its relative distance to the optimal objective value $\exactB$.
For this problem, all primal and dual problems are solved by the extensive form (i.e., monolithic) of their respective models. Figures \ref{fig:newsvendor_primal_Small} and \ref{fig:newsvendor_dual_Small} present the gap between the bound and the optimal value of the exact solution, defined as $100\big(\frac{\exactB-\nu^{(\cdot)}}{\exactB}\big)$ and $100\big(\frac{\nu^{(\cdot)} - \exactB}{\exactB}\big)$ for primal and dual bounds, respectively, and presented as a percentage (detailed results {\crev along with solution times} are given in {\crev Appendix~\ref{app:results-newsvendor}}). 
In each figure, the solid bars depict the performance of the newly proposed bounds, while the hatched bars represent the benchmarks. The results show that $\twostageLDRUB$ and  {\crev $\bar{\nu}^{\textsc{NA-DO}}_R $} outperform the benchmark bounds by orders of magnitude. 
\begin{figure}[tbp]
    \centering
    \scalebox{0.9}{
    \footnotesize
    \input{Figures/Primal_bounds_newsvendor.tex}
    }
    \caption{Quality of the {\crev primal} bounds from LDRs and two-stage LDRs for small-size newsvendor instances}
    \label{fig:newsvendor_primal_Small}
\end{figure}
\begin{figure}[b]
    \centering
    \scalebox{0.9}{
    \footnotesize
    \input{Figures/Dual_bounds_newsvendor.tex}
    }
    \caption{
    Quality of the {\crev dual} bounds from the PI and LDDRs for small-size newsvendor instances}
    \label{fig:newsvendor_dual_Small}
\end{figure}
More precisely, $\twostageLDRUB$ on average achieves 84\% improvement over $\LDRUB$, with reductions in relative distance ranging from 64\% to 98\%. Interestingly, the quality of the bound $\twostageLDRUB$ remains rather stable with changes in the number of items $I$ and budget $B$, compared to the drastic changes of $\LDRUB$ with variations in the inputs. The notable performance of the two-stage LDRs for the newsvendor problem contrasted with the LDRs can be explained by the nature of the  recourse variables. 
In model \eqref{eq:newsvendor}, $y_{it}(\xisupt)$ 
determines the  profit of a given realization at stage $t$ for item $i$, which  for the newsvendor problem is by definition nonlinear. In fact, $y_{it}(\xisupt)$ is an auxiliary variable, defined to linearize the following net profit at stage $t$  from item $i$: 
$$ r_i\min\big\{x_{i,t-1}(\xiparrand^{t-1}),d_{it}(\xiparrand^{t})\big\} - c_ix_{i,t-1}(\xiparrand^{t-1}) - s_i\max\big\{d_{it}(\xiparrand^{t})-x_{i,t-1}(\xiparrand^{t-1}),0\big\}.$$
Therefore, for the newsvendor problem LDRs \emph{always} return suboptimal decisions as they restrict the form of the nonlinear profit function to be affine, while two-stage LDRs allow them to take any form, giving them an immediate advantage over LDRs.

From the dual perspective, {\crev $\bar{\nu}^{\textsc{NA-DO}}_R $} achieves an average improvement of 55\% compared to $\PILB$, and in 9 instances fully closes the gap. From the results of Figure \ref{fig:newsvendor_dual_Small}, a common observation is that LDDRs return a bound of higher quality for smaller values of the ratio $\frac{B}{(T-1)\times I}$, which is an estimate of the average budget available per product at every stage. This trend suggests that the restricted NA dual bound performs better with tighter budget, when there is a higher dependency between stages, in which case the importance of intermediate decisions becomes more pronounced.

\subsubsection{Optimality Gap}
For instances with larger number of {\crev realizations}, we compare the optimality gap from the benchmark methods ($\optT$) with the gap obtained by applying the newly proposed methods, namely two-stage LDRs and LDDRs ($\optN$): 
$$\displaystyle\optT = 100\Big(\frac{\PILB - \LDRUB}{\LDRUB}\Big),\quad \optN = 100\Big(\frac{ {\crev \bar{\nu}^{\textsc{NA-DO}}_R} - \twostageLDRUB}{\twostageLDRUB}\Big).$$
Table \ref{tab:improvement-newsvendor} presents the optimality gaps from the benchmark and proposed models  (whose running times are provided in Appendix \ref{app:results-newsvendor}). 
\begin{table}[htbp]
\crev
\small
  \centering
  \caption{Optimality gaps for larger instances of the newsvendor problem}
  \scalebox{0.88}{
\begin{tabular}{lcrrccrrrrrrc}
\toprule
\multicolumn{1}{c}{\multirow{2}{*}{Instance}} & \multirow{2}{*}{$\ T\ $} & \multirow{2}{*}{\textsc{br}} & \multirow{2}{*}{$|\xisupp|$} & \multirow{2}{*}{$\ I\ $} & \multirow{2}{*}{$B$} & \multicolumn{2}{c}{Primal Bounds} & \multicolumn{2}{c}{Dual Bounds} & \multicolumn{2}{c}{Optimality Gap} & \multirow{2}{*}{Gap reduction} \\
\cmidrule(lr){7-8} \cmidrule(rl){9-10} \cmidrule(lr){11-12}      &       &       &       &       &       & \multicolumn{1}{c}{$\LDRUB$} & \multicolumn{1}{c}{$\twostageLDRUB$} & \multicolumn{1}{c}{$\bar{\nu}^{\textsc{NA-DO}}_R$} & \multicolumn{1}{c}{$\PILB$} & \multicolumn{1}{c}{$\optN$} & \multicolumn{1}{c}{$\optT$} &  \\
\midrule
1     & \multirow{5}{*}{4} & \multirow{5}{*}{10} & \multirow{5}{*}{1000} & 3     & 200   & 5648.8 & 8142.9 & 9002.8 & 9353.0 & 10.6\% & 65.6\% & 83.9\% \\
    2     &       &       &       & 3     & 300   & 9143.5 & 13687.6 & 17853.0 & 17853.0 & 30.4\% & 95.3\% & 68.1\% \\
    3     &       &       &       & 4     & 200   & -104.4 & 440.0 & 642.4 & 919.0 & 46.0\% & 980.4\% & 95.3\% \\
    4     &       &       &       & 4     & 300   & 11029.5 & 15854.7 & 18432.0 & 18432.0 & 16.3\% & 67.1\% & 75.8\% \\
    5     &       &       &       & 5     & 300   & 724.9 & 6125.5 & 6614.0 & 7368.0 & 8.0\% & 916.5\% & 99.1\% \\
\cmidrule{2-13}6     & 4     & 15    & 3375  & 3     & 200   & 5111.3 & 7072.0 & 8222.0 & 8222.0 & 16.3\% & 60.9\% & 73.3\% \\
\cmidrule{2-13}7     & 4     & 20    & 8000  & 3     & 300   & 9030.9 & 13040.7 & 17615.0 & 17615.0 & 35.1\% & 95.1\% & 63.1\% \\
\midrule
8     & \multirow{2}{*}{5} & \multirow{2}{*}{5} & \multirow{2}{*}{625} & 3     & 300   & 11134.5 & 15053.4 & 16216.4 & 16680.0 & 7.7\% & 49.8\% & 84.5\% \\
    9     &       &       &       & 4     & 300   & 3651.5 & 9494.8 & 10403.1 & 10628.0 & 9.6\% & 191.1\% & 95.0\% \\
\cmidrule{2-13}10    & 5     & 6     & 1296  & 3     & 400   & 15013.4 & 20271.7 & 25157.0 & 25157.0 & 24.1\% & 67.6\% & 64.3\% \\
\midrule
11    & \multirow{3}{*}{6} & \multirow{3}{*}{4} & \multirow{3}{*}{1024} & 3     & 400   & 15492.3 & 21457.9 & 28163.0 & 28163.0 & 31.2\% & 81.8\% & 61.8\% \\
    12    &       &       &       & 4     & 400   & 7124.3 & 14805.2 & 15473.0 & 15473.0 & 4.5\% & 117.2\% & 96.2\% \\
    13    &       &       &       & 4     & 500   & 14445.4 & 24887.7 & 32973.0 & 32973.0 & 32.5\% & 128.3\% & 74.7\% \\
\midrule
14    & \multirow{3}{*}{7} & \multirow{3}{*}{3} & \multirow{3}{*}{729} & 3     & 300   & 24.3  & 1495.9 & 2150.0 & 2383.0 & 43.7\% & 9692.2\% & 99.5\% \\
    15    &       &       &       & 3     & 400   & 12994.7 & 17774.1 & 19554.5 & 19983.0 & 10.0\% & 53.8\% & 81.4\% \\
    16    &       &       &       & 4     & 400   & 3267.5 & 4303.1 & 5118.1 & 5427.0 & 18.9\% & 66.1\% & 71.3\% \\
\midrule
17    & \multirow{2}{*}{8} & \multirow{2}{*}{3} & \multirow{2}{*}{2187} & 3     & 500   & 16149.1 & 25321.4 & 30259.0 & 30259.0 & 19.5\% & 87.4\% & 77.7\% \\
    18    &       &       &       & 4     & 600   & 16892.3 & 27608.1 & 28747.0 & 28747.0 & 4.1\% & 70.2\% & 94.1\% \\
\bottomrule
\end{tabular}%
}
  \label{tab:improvement-newsvendor}%
\end{table}%
{Achieving an average gap reduction of {\crev 81.6\%}, there is a considerable value in using the two-stage LDRs and LDDRs in devising policies for the multistage newsvendor problem. The optimality gaps $\optN$ range from 4\% to 46\%. There are some interesting cases such as instance 3 where  LDR policies can even lead to a profit loss in the worst case. 
In the majority of the instances, both two-stage LDRs and LDDRs contribute to the improvement of the gap, although the primal side clearly has the larger impact. For instance, in five instances, the dual bound {\crev $\bar{\nu}^{\textsc{NA-DO}}_R$} and the benchmark $\PILB$ are the same. In some instances, such as instances 8, 12 and 18, the PI bound might already be strong enough so that using LDDRs does not make a tangible difference in strengthening it. 

{\crev Lastly, for the primal bounds, we observe that the solution times for obtaining the traditional and two-stage LDR bounds are comparable. For the dual bounds, the  times for obtaining the PI bound are quite small whereas for our dual problem the solution time significantly increases with problem size with several instances not being solved to optimality within the 1-hour time limit.    }

\subsection{Robust Location-Transportation Problem}
The two-stage robust location-transportation problem studied by \cite{zeng2013solving} is as follows. Given a set of $I$ potential facilities with
building cost $f_i$ and unit capacity cost $a_i, i\in\set{I}$, we have to meet the uncertain demand of a set of customers $J$ with unit transportation cost $c_{ij}, i\in\set{I}, j\in\set{J}$ . The goal is to decide which facilities to open and their initial capacities, such that the worst-case total cost of facility deployment and future transportation is minimized. Letting $d_{jt}(\xisupt)$ be the demand  of customer $j$ at stage $t$, we define the \MSARO{} extension of the problem:
\bsub
\label{eq:location-transportation}
\begin{align}
	\min  \ \ & z  \label{eq:obj-LT} \\
	\text{s.t.} \ \ & z\geq \sum_{i\in \set{I}} (f_iy_{i} + {\crev a_i s_{i1}) + } \sum_{t\in {\crev [2,T]}}\sum_{i\in \set{I}} \Big( a_is_{it}(\xisupt) + \sum_{j\in \set{J}}c_{ij}x_{ijt}(\xisupt) \Big) \quad  &&  \hspace{-5mm} \xiT\in\xisupp \label{eq:obj-LT-worst} \\
	& s_{i1} \leq K_iy_i &&  \hspace{-5mm} i\in\set{I}\label{eq:const-cap-ub} \\
 	& {\crev s_{it}(\xisupt) = s_{i1} - \sum_{j\in \set{J}}x_{ijt}(\xisupt)}&&  \hspace{-5mm}i\in\set{I}, {\crev t =2},  \xisupt\in\xisuppt\label{eq:const-state-1}\\
	& s_{it}(\xisupt) = s_{i,t-1}(\xisuptminusone) - \sum_{j\in \set{J}}x_{ijt}(\xisupt)&&  \hspace{-5mm}i\in\set{I}, t\in\set{{\crev 3},T},  \xisupt\in\xisuppt\label{eq:const-state-2}\\
	& \sum_{i\in \set{I}}x_{ijt}(\xisupt) \geq d_{jt}(\xisupt)&&  \hspace{-5mm}j\in\set{J}, t\in\set{2,T}, \xisupt\in\xisuppt\label{eq:const-demand}\\
	&y\in\{0,1\}^I,\ s_{1}\in\R_+^I\\
	&s_t(\xisupt)\in\R_+^{I},\ x_t(\xisupt)\in\R_+^{I\times J}&&\hspace{-5mm} t\in\set{2,T},  \xisupt\in\xisuppt,
\end{align}
\esub
where $y_i$ is a binary variable equal to 1 if facility $i$ is built,
$s_{i1}$ determines the initial capacity of facility $i$, $s_{it}(\xisupt)$ is the state variable calculating the remaining capacity of facility $i$ at stage $t$, and $x_{ijt}(\xisupt)$ is the amount of goods transported from facility $i$ to customer $j$ at stage $t$. Objective function \eqref{eq:obj-LT} together with constraint 
\eqref{eq:obj-LT-worst} measures the worst-case cost. Constraints \eqref{eq:const-cap-ub} bound the initial capacity of the facilities, whereas constraints  \eqref{eq:const-state-1}-{\crev \eqref{eq:const-state-2}} are the state equations calculating the remaining capacities. Constraints \eqref{eq:const-demand} ensure that the customer demands are met.

\subsubsection{Problem Instances}
Our instances are generated using the instance parameters described in \citep{zeng2013solving} for the two-stage problem. For number of stages $T\in\set{3,5}$, we have five combinations for the number of facilities and customers, $(I,J)\in\{(5,5),(5,7),(5,10),(10,10),(20,20)\}$. Fixed installation, unit capacity, and unit transportation costs are drawn from $f_i \in[100,1000],
a_i \in[10,100], c_{ij} \in[1,1000]$, respectively. Maximal capacity is set to $K_i = 2\times 10^4$ based on preliminary experiments to make sure the instances are feasible. Customer demands are random parameters with support $\big[\mu_{jt},(1+\alpha^d) \mu_{jt}\big]$, where $\mu_{jt} \in[10,500]$ and $\alpha^d\in\{0.1,0.3,0.5\}$ are given, with $\alpha^d$ a parameter controlling the variation among demand realizations of customer $j$. For this problem, we have generated 32 instances with $\alpha^d=0.5$ over small scenario trees to compare the bounds with the optimal objective value. In addition, we have  generated 83 instances {\crev with demands $d_{jt}=\mu_{jt}+\xiparrand_{jt}\sigma_{jt},\    j\in\set{J},\ t\in\set{2,T}$, where $\xiparrand$ belongs to the following budgeted uncertainty set:
$$\xisupp=\bigg\{\xiparrand\in \R_+^{J\times T-1}\ \Big|\   \xiparrand_{jt}\in[0,1],\ j\in\set{J},\ t\in\set{2,T},\ \sum_{t\in\set{2,T}}\sum_{j\in\set{J}}\xiparrand_{jt}\leq\Gamma\bigg\}.$$}
Parameter $\Gamma$ correlates the demands of all customers and stages together, which results in stagewise (temporal) dependence between the decision stages. 
We use $\Gamma =\alpha^uI$, with $\alpha^u\in\{0.1,0.4,0.7,1\}$. 

\subsubsection{Scenario-Tree Instances}
In solving the primal (2S-LDR) and exact models of the scenario-tree instances, we have used their respective extensive forms, while for computing {\crev our} dual bound we implemented the cutting-plane method described in Section \ref{sec:solution}. 

\begin{table}[b] 
\crev
\caption{Quality of the 
bounds for the scenario-tree instances of the location-transportation problem}
  \small
  \centering 
  \begin{minipage}[t]{0.48\textwidth}
  \vspace{0pt}
  \scalebox{0.88}{
  \begin{tabular}{crcrrrcrr}
\toprule
$T$   & $I$   & $J$   & $\textsc{br}$ & $|\xisupp|$ & $\LDRUB$ & $\twostageLDRUB$ & $\bar{\nu}^{\textsc{NA-DO}}_R$ & $\PILB$ \\
\midrule
\multirow{8}{*}{3} & \multirow{8}{*}{5} & \multirow{8}{*}{10} & 3     & 9     & 2.9\% & 0.0\% & 1.4\% & 4.7\% \\
          &       &       & 4     & 16    & 0.8\% & 0.0\% & 3.3\% & 3.4\% \\
          &       &       & 5     & 25    & 4.6\% & 0.0\% & 2.5\% & 13.2\% \\
          &       &       & 6     & 36    & 5.6\% & 0.0\% & 17.1\% & 27.5\% \\
          &       &       & 7     & 49    & 8.2\% & 0.0\% & 9.4\% & 9.4\% \\
          &       &       & 8     & 64    & 8.8\% & 0.0\% & 4.8\% & 14.6\% \\
          &       &       & 9     & 81    & 6.5\% & 0.0\% & 6.1\% & 6.1\% \\
          &       &       & 10    & 100   & 6.6\% & 0.0\% & 8.0\% & 10.2\% \\
\midrule
\multirow{8}{*}{3} & \multirow{8}{*}{10} & \multirow{8}{*}{10} & 3     & 9     & 2.3\% & 0.0\% & 2.5\% & 13.2\% \\
          &       &       & 4     & 16    & 2.7\% & 0.1\% & 1.2\% & 8.4\% \\
          &       &       & 5     & 25    & 6.1\% & 0.0\% & 12.0\% & 19.3\% \\
          &       &       & 6     & 36    & 9.8\% & 0.0\% & 10.4\% & 22.6\% \\
          &       &       & 7     & 49    & 13.7\% & 0.1\% & 6.7\% & 14.1\% \\
          &       &       & 8     & 64    & 13.9\% & 0.3\% & 5.5\% & 5.5\% \\
          &       &       & 9     & 81    & 12.2\% & 0.2\% & 15.5\% & 15.5\% \\
          &       &       & 10    & 100   & 12.3\% & 0.2\% & 10.4\% & 19.0\% \\
\midrule
\multirow{5}{*}{3} & \multirow{5}{*}{20} & \multirow{5}{*}{20} & 3     & 9     & 2.5\% & 0.7\% & 9.5\% & 11.3\% \\
          &       &       & 4     & 16    & 3.9\% & 0.5\% & 7.0\% & 7.0\% \\
          &       &       & 5     & 25    & 5.8\% & 0.5\% & 28.4\% & 28.4\% \\
          &       &       & 6     & 36    & 6.6\% & 0.5\% & 13.0\% & 19.2\% \\
          &       &       & 7     & 49    & 9.1\% & 0.5\% & 13.3\% & 31.3\% \\
\bottomrule
\end{tabular}
}
\end{minipage}
\begin{minipage}[t]{0.48\textwidth}
\vspace{0pt}
\scalebox{0.88}{
  \begin{tabular}{crcrrrcrr} %p{1mm}
\toprule
$T$   & $I$   & $J$   & $\textsc{br}$ & $|\xisupp|$ & $\LDRUB$ & $\twostageLDRUB$ & $\bar{\nu}^{\textsc{NA-DO}}_R$ & $\PILB$ \\
\midrule
\multirow{3}{*}{4} & \multirow{3}{*}{5} & \multirow{3}{*}{10} & 3     & 27    & 8.0\% & 0.0\% & 2.6\% & 8.0\% \\
          &       &       & 4     & 64    & 6.8\% & 0.0\% & 2.0\% & 2.0\% \\
          &       &       & 5     & 125   & 13.5\% & 0.1\% & 5.0\% & 5.0\% \\
\midrule
\multirow{3}{*}{4} & \multirow{3}{*}{10} & \multirow{3}{*}{10} & 3     & 27    & 11.9\% & 0.1\% & 17.3\% & 17.3\% \\
          &       &       & 4     & 64    & 13.3\% & 0.0\% & 10.2\% & 21.0\% \\
          &       &       & 5     & 125   & 13.1\% & 0.4\% & 10.0\% & 10.0\% \\
\midrule
\multirow{2}{*}{4} & \multirow{2}{*}{20} & \multirow{2}{*}{20} & 3     & 27    & 9.0\% & 0.9\% & 11.7\% & 18.1\% \\
          &       &       & 4     & 64    & 11.9\% & 0.5\% & 9.8\% & 21.7\% \\
\midrule
\multirow{2}{*}{5} & \multirow{2}{*}{5} & \multirow{2}{*}{10} & 3     & 81    & 14.7\% & 1.0\% & 4.4\% & 6.3\% \\
          &       &       & 4     & 256   & 13.9\% & 0.5\% & 11.7\% & 19.1\% \\
\midrule
5     & 10    & 10    & 3     & 81    & 16.8\% & 0.7\% & 8.2\% & 8.2\% \\
\bottomrule
\end{tabular}
}
\end{minipage}
\label{tab:LT-ST-results}
\end{table}

Table \ref{tab:LT-ST-results} presents the gap between the bound and the optimal value of the exact solution, $100\big(\frac{\nu^{(\cdot)}-\exactB}{\exactB}\big)$ and $100\big(\frac{\exactB-\nu^{(\cdot)}}{\exactB}\big)$ for the primal and dual problems, respectively.  Results show that, among the group of instances with similar characteristics, as the branching factor or the number of stages increases, the $\LDRUB$ bound gets noticeably worse. In contrast, $\twostageLDRUB$ stays very close to the optimal value, with an average relative distance of $0.2\%$ among all the instances (compared to $8.7\%$ for $\LDRUB$). We do not observe the same trend for the dual bounds, and their relative distance to the optimal value fluctuates even between two instances that only differ in the number of {\crev realizations}. Nevertheless, {\crev $\bar{\nu}^{\textsc{NA-DO}}_R$} considerably outperforms 
PI, with an average improvement of {\crev $36.2\%$}.

\subsubsection{Budgeted-Uncertainty Instances} 
{\crev For the instances with the budgeted uncertainty set, we use, as benchmarks, $\LDRUB$ as an upper bound and $\nu^{\Omega(\text{LDR})}$ as a lower bound, obtained by solving problem \eqref{eq:location-transportation} using only the binding {\crev realizations} identified from the benchmark primal solution.
We calculate $\twostageLDRUB$  using the C\&CG method described in Section \ref{sec:two-stage-LDR} (detailed models are given in Appendix \ref{sec:CCG-LT-detail}).  Similar to the scenario-tree instances, we use the cutting-plane method to obtain the dual bound  {\crev $\bar{\nu}^{\textsc{NA-DO}}_R$}. To do so, we use a sample of size at least $ 50(T-2)$, which includes both the binding {\crev realizations} from the two-stage LDR solution and additional randomly generated realizations from the uncertainty set. }
For both algorithms, we stop when the optimality gap of the method falls below 5\% or we reach the 10-hour time limit. 
{\crev Table \ref{tab:LT_Budgeted} presents the  optimality gaps between the benchmark bounds and the proposed bounds, respectively, for the 83 instances considered. Detailed results for each instance are provided in Appendix \ref{app:LT-results}.} 
\begin{table}[htbp]
\crev
  \caption{\crev Optimality gaps for the budgeted-uncertainty instances, calculated using both benchmark and proposed bounding methods. For this  sampling instance, $\optT$ is defined as the gap between $\LDRUB$ and $\nu^{\Omega(\text{LDR})}$. As before, $\optN$ represents the gap between $\twostageLDRUB$ and $\bar{\nu}^{\textsc{NA-DO}}_R$.\\}
    
    \centering
    
    \scalebox{0.82}{
    \setlength{\tabcolsep}{10pt}
\begin{tabular}{cccccccrrrrrrr}
\cmidrule{1-6}\cmidrule{9-14}\multicolumn{2}{c}{$(T, I,J)$} & $\alpha^d$ & $\alpha^u$ & {$\optT$} & {$\optN$} & $\ $  & \multicolumn{1}{c}{$\ $} & \multicolumn{2}{c}{$(T,I,J)$} & \multicolumn{1}{c}{$\alpha^d$} & \multicolumn{1}{c}{$\alpha^u$} & \multicolumn{1}{c}{$\optT$} & \multicolumn{1}{c}{$\optN$} \\
\cmidrule{1-6}\cmidrule{9-14}\multicolumn{2}{c}{\multirow{12}[6]{*}{(3,10,10)}} & \multirow{4}[2]{*}{0.1} & 0.1   & 9.0\% & 7.4\% &       &       & \multicolumn{2}{c}{\multirow{12}[6]{*}{(4,5,10)}} & \multicolumn{1}{c}{\multirow{4}[2]{*}{0.1}} & \multicolumn{1}{c}{0.1} & \multicolumn{1}{c}{22.4\%} & \multicolumn{1}{c}{16.8\%} \\
\multicolumn{2}{c}{} &       & 0.4   & 20.8\% & 12.4\% &       &       & \multicolumn{2}{c}{} &       & \multicolumn{1}{c}{0.4} & \multicolumn{1}{c}{32.1\%} & \multicolumn{1}{c}{28.1\%} \\
\multicolumn{2}{c}{} &       & 0.7   & 32.0\% & 22.5\% &       &       & \multicolumn{2}{c}{} &       & \multicolumn{1}{c}{0.7} & \multicolumn{1}{c}{15.2\%} & \multicolumn{1}{c}{10.2\%} \\
\multicolumn{2}{c}{} &       & 1     & 30.6\% & 22.2\% &       &       & \multicolumn{2}{c}{} &       & \multicolumn{1}{c}{1} & \multicolumn{1}{c}{31.2\%} & \multicolumn{1}{c}{19.0\%} \\
\cmidrule{3-6}\cmidrule{11-14}\multicolumn{2}{c}{} & \multirow{4}[2]{*}{0.3} & 0.1   & 15.3\% & 12.0\% &       &       & \multicolumn{2}{c}{} & \multicolumn{1}{c}{\multirow{4}[2]{*}{0.3}} & \multicolumn{1}{c}{0.1} & \multicolumn{1}{c}{28.5\%} & \multicolumn{1}{c}{20.8\%} \\
\multicolumn{2}{c}{} &       & 0.4   & 27.0\% & 16.7\% &       &       & \multicolumn{2}{c}{} &       & \multicolumn{1}{c}{0.4} & \multicolumn{1}{c}{37.7\%} & \multicolumn{1}{c}{30.6\%} \\
\multicolumn{2}{c}{} &       & 0.7   & 30.3\% & 21.5\% &       &       & \multicolumn{2}{c}{} &       & \multicolumn{1}{c}{0.7} & \multicolumn{1}{c}{42.9\%} & \multicolumn{1}{c}{27.6\%} \\
\multicolumn{2}{c}{} &       & 1     & 29.1\% & 21.3\% &       &       & \multicolumn{2}{c}{} &       & \multicolumn{1}{c}{1} & \multicolumn{1}{c}{24.3\%} & \multicolumn{1}{c}{14.2\%} \\
\cmidrule{3-6}\cmidrule{11-14}\multicolumn{2}{c}{} & \multirow{4}[2]{*}{0.5} & 0.1   & 17.8\% & 12.9\% &       &       & \multicolumn{2}{c}{} & \multicolumn{1}{c}{\multirow{4}[2]{*}{0.5}} & \multicolumn{1}{c}{0.1} & \multicolumn{1}{c}{33.0\%} & \multicolumn{1}{c}{27.7\%} \\
\multicolumn{2}{c}{} &       & 0.4   & 28.3\% & 15.1\% &       &       & \multicolumn{2}{c}{} &       & \multicolumn{1}{c}{0.4} & \multicolumn{1}{c}{22.2\%} & \multicolumn{1}{c}{11.4\%} \\
\multicolumn{2}{c}{} &       & 0.7   & 32.7\% & 28.4\% &       &       & \multicolumn{2}{c}{} &       & \multicolumn{1}{c}{0.7} & \multicolumn{1}{c}{27.7\%} & \multicolumn{1}{c}{17.9\%} \\
\multicolumn{2}{c}{} &       & 1     & 26.8\% & 7.4\% &       &       & \multicolumn{2}{c}{} &       & \multicolumn{1}{c}{1} & \multicolumn{1}{c}{33.1\%} & \multicolumn{1}{c}{23.5\%} \\
\cmidrule{1-6}\cmidrule{9-14}\multicolumn{2}{c}{\multirow{12}[6]{*}{(3,10,15)}} & \multirow{4}[2]{*}{0.1} & 0.1   & 32.8\% & 24.4\% &       &       & \multicolumn{2}{c}{\multirow{12}[6]{*}{(4,10,10)}} & \multicolumn{1}{c}{\multirow{4}[2]{*}{0.1}} & \multicolumn{1}{c}{0.1} & \multicolumn{1}{c}{19.4\%} & \multicolumn{1}{c}{17.5\%} \\
\multicolumn{2}{c}{} &       & 0.4   & 31.5\% & 27.1\% &       &       & \multicolumn{2}{c}{} &       & \multicolumn{1}{c}{0.4} & \multicolumn{1}{c}{32.1\%} & \multicolumn{1}{c}{26.6\%} \\
\multicolumn{2}{c}{} &       & 0.7   & 12.1\% & 7.2\% &       &       & \multicolumn{2}{c}{} &       & \multicolumn{1}{c}{0.7} & \multicolumn{1}{c}{38.7\%} & \multicolumn{1}{c}{32.7\%} \\
\multicolumn{2}{c}{} &       & 1     & 31.5\% & 17.3\% &       &       & \multicolumn{2}{c}{} &       & \multicolumn{1}{c}{1} & \multicolumn{1}{c}{21.7\%} & \multicolumn{1}{c}{10.8\%} \\
\cmidrule{3-6}\cmidrule{11-14}\multicolumn{2}{c}{} & \multirow{4}[2]{*}{0.3} & 0.1   & 12.8\% & 8.3\% &       &       & \multicolumn{2}{c}{} & \multicolumn{1}{c}{\multirow{4}[2]{*}{0.3}} & \multicolumn{1}{c}{0.1} & \multicolumn{1}{c}{11.1\%} & \multicolumn{1}{c}{6.6\%} \\
\multicolumn{2}{c}{} &       & 0.4   & 35.3\% & 25.6\% &       &       & \multicolumn{2}{c}{} &       & \multicolumn{1}{c}{0.4} & \multicolumn{1}{c}{19.8\%} & \multicolumn{1}{c}{16.5\%} \\
\multicolumn{2}{c}{} &       & 0.7   & 17.2\% & 5.7\% &       &       & \multicolumn{2}{c}{} &       & \multicolumn{1}{c}{0.7} & \multicolumn{1}{c}{22.6\%} & \multicolumn{1}{c}{17.4\%} \\
\multicolumn{2}{c}{} &       & 1     & 26.2\% & 13.6\% &       &       & \multicolumn{2}{c}{} &       & \multicolumn{1}{c}{1} & \multicolumn{1}{c}{38.1\%} & \multicolumn{1}{c}{28.0\%} \\
\cmidrule{3-6}\cmidrule{11-14}\multicolumn{2}{c}{} & \multirow{4}[2]{*}{0.5} & 0.1   & 32.4\% & 28.2\% &       &       & \multicolumn{2}{c}{} & \multicolumn{1}{c}{\multirow{4}[2]{*}{0.5}} & \multicolumn{1}{c}{0.1} & \multicolumn{1}{c}{24.9\%} & \multicolumn{1}{c}{12.4\%} \\
\multicolumn{2}{c}{} &       & 0.4   & 25.8\% & 15.3\% &       &       & \multicolumn{2}{c}{} &       & \multicolumn{1}{c}{0.4} & \multicolumn{1}{c}{23.3\%} & \multicolumn{1}{c}{11.0\%} \\
\multicolumn{2}{c}{} &       & 0.7   & 29.4\% & 18.2\% &       &       & \multicolumn{2}{c}{} &       & \multicolumn{1}{c}{0.7} & \multicolumn{1}{c}{31.4\%} & \multicolumn{1}{c}{18.2\%} \\
\multicolumn{2}{c}{} &       & 1     & 29.7\% & 20.7\% &       &       & \multicolumn{2}{c}{} &       & \multicolumn{1}{c}{1} & \multicolumn{1}{c}{31.2\%} & \multicolumn{1}{c}{25.7\%} \\
\cmidrule{1-6}\cmidrule{9-14}\multicolumn{2}{c}{\multirow{12}[6]{*}{(4,5,5)}} & \multirow{4}[2]{*}{0.1} & 0.1   & 14.3\% & 6.0\% &       &       & \multicolumn{2}{c}{\multirow{12}[6]{*}{(5,5,10)}} & \multicolumn{1}{c}{\multirow{4}[2]{*}{0.1}} & \multicolumn{1}{c}{0.1} & \multicolumn{1}{c}{17.2\%} & \multicolumn{1}{c}{13.7\%} \\
\multicolumn{2}{c}{} &       & 0.4   & 33.2\% & 25.7\% &       &       & \multicolumn{2}{c}{} &       & \multicolumn{1}{c}{0.4} & \multicolumn{1}{c}{24.5\%} & \multicolumn{1}{c}{20.6\%} \\
\multicolumn{2}{c}{} &       & 0.7   & 27.7\% & 10.8\% &       &       & \multicolumn{2}{c}{} &       & \multicolumn{1}{c}{0.7} & \multicolumn{1}{c}{34.2\%} & \multicolumn{1}{c}{21.8\%} \\
\multicolumn{2}{c}{} &       & 1     & 12.7\% & 5.6\% &       &       & \multicolumn{2}{c}{} &       & \multicolumn{1}{c}{1} & \multicolumn{1}{c}{20.3\%} & \multicolumn{1}{c}{8.6\%} \\
\cmidrule{3-6}\cmidrule{11-14}\multicolumn{2}{c}{} & \multirow{4}[2]{*}{0.3} & 0.1   & 13.5\% & 7.8\% &       &       & \multicolumn{2}{c}{} & \multicolumn{1}{c}{\multirow{4}[2]{*}{0.3}} & \multicolumn{1}{c}{0.1} & \multicolumn{1}{c}{28.9\%} & \multicolumn{1}{c}{25.6\%} \\
\multicolumn{2}{c}{} &       & 0.4   & 25.3\% & 18.1\% &       &       & \multicolumn{2}{c}{} &       & \multicolumn{1}{c}{0.4} & \multicolumn{1}{c}{32.5\%} & \multicolumn{1}{c}{19.3\%} \\
\multicolumn{2}{c}{} &       & 0.7   & 23.6\% & 13.8\% &       &       & \multicolumn{2}{c}{} &       & \multicolumn{1}{c}{0.7} & \multicolumn{1}{c}{39.1\%} & \multicolumn{1}{c}{31.9\%} \\
\multicolumn{2}{c}{} &       & 1     & 26.2\% & 9.6\% &       &       & \multicolumn{2}{c}{} &       & \multicolumn{1}{c}{1} & \multicolumn{1}{c}{37.5\%} & \multicolumn{1}{c}{29.4\%} \\
\cmidrule{3-6}\cmidrule{11-14}\multicolumn{2}{c}{} & \multirow{4}[2]{*}{0.5} & 0.1   & 31.8\% & 27.1\% &       &       & \multicolumn{2}{c}{} & \multicolumn{1}{c}{\multirow{4}[2]{*}{0.5}} & \multicolumn{1}{c}{0.1} & \multicolumn{1}{c}{36.6\%} & \multicolumn{1}{c}{31.8\%} \\
\multicolumn{2}{c}{} &       & 0.4   & 22.6\% & 14.1\% &       &       & \multicolumn{2}{c}{} &       & \multicolumn{1}{c}{0.4} & \multicolumn{1}{c}{39.8\%} & \multicolumn{1}{c}{26.6\%} \\
\multicolumn{2}{c}{} &       & 0.7   & 29.8\% & 19.6\% &       &       & \multicolumn{2}{c}{} &       & \multicolumn{1}{c}{0.7} & \multicolumn{1}{c}{39.2\%} & \multicolumn{1}{c}{17.6\%} \\
\multicolumn{2}{c}{} &       & 1     & 29.8\% & 5.4\% &       &       & \multicolumn{2}{c}{} &       & \multicolumn{1}{c}{1} & \multicolumn{1}{c}{41.6\%} & \multicolumn{1}{c}{27.1\%} \\
\cmidrule{1-6}\cmidrule{9-14}\multicolumn{2}{c}{\multirow{11}[6]{*}{(4,5,7)}} & \multirow{3}[2]{*}{0.1} & 0.1   & 15.7\% & 13.6\% &       &       &       &       &       &       &       &  \\
\multicolumn{2}{c}{} &       & 0.4   & 24.0\% & 19.4\% &       &       &       &       &       &       &       &  \\
\multicolumn{2}{c}{} &       & 1     & 38.3\% & 24.5\% &       &       &       &       &       &       &       &  \\
\cmidrule{3-6}\multicolumn{2}{c}{} & \multirow{4}[2]{*}{0.3} & 0.1   & 23.2\% & 19.3\% &       &       &       &       &       &       &       &  \\
\multicolumn{2}{c}{} &       & 0.4   & 30.8\% & 21.6\% &       &       &       &       &       &       &       &  \\
\multicolumn{2}{c}{} &       & 0.7   & 10.2\% & 6.4\% &       &       &       &       &       &       &       &  \\
\multicolumn{2}{c}{} &       & 1     & 30.7\% & 18.7\% &       &       &       &       &       &       &       &  \\
\cmidrule{3-6}\multicolumn{2}{c}{} & \multirow{4}[2]{*}{0.5} & 0.1   & 24.2\% & 21.5\% &       &       &       &       &       &       &       &  \\
\multicolumn{2}{c}{} &       & 0.4   & 39.5\% & 27.1\% &       &       &       &       &       &       &       &  \\
\multicolumn{2}{c}{} &       & 0.7   & 32.7\% & 17.9\% &       &       &       &       &       &       &       &  \\
\multicolumn{2}{c}{} &       & 1     & 41.0\% & 23.3\% &       &       &       &       &       &       &       &  \\
\cmidrule{1-6}\end{tabular}%
}
  \label{tab:LT_Budgeted}%
\end{table}%
Our methods return an average optimality gap of {\crev $17.3\%$ across all instances, compared to an average gap of $27\%$ from the benchmarks. Given the strength of the two-stage LDR bound observed in the scenario-tree instances, and its resilience to increases in the size of the tree, it is likely that the dual bounds are further from the optimal value. 

Figure \ref{fig:LT-dual-bound} illustrates the improvement achieved from $\bar{\nu}^{\textsc{NA-DO}}_R$ over the benchmark dual bound, $\nu^{\Omega(\text{LDR})}$. Each bar in the figure represents an individual instance. Instances are grouped into blocks based on their shared parameters $T$, $(I,J)$, and $\alpha^d$, and within each block, instances are arranged in ascending order of their $\alpha^u$ values. Across all instances, $\bar{\nu}^{\textsc{NA-DO}}_R$ demonstrates an average improvement of 7\%. This figure also highlights that improved identification of critical {\crev realizations}, guided by the two-stage LDR solution, and the subsequent solution of the discretized relaxation of the problem using these {\crev realizations}, independently contribute to an average 3.6\% enhancement of the dual bound. A key consideration here is that} 
the value of {\crev $\bar{\nu}^{\textsc{NA-DO}}_R$ depends on the sample used for its computation. An insufficiently large sample can lead to a poor bound.}  However, because of the bilinear form of the cutting-plane master problem, we were not able to solve the model with large samples using  off-the-shelf commercial solvers. By employing specialized algorithms developed for bilinear problems, it could be possible to increase the size of the sample and improve the {\crev $\bar{\nu}^{\textsc{NA-DO}}_R$} bound, which we leave for future research.

\begin{figure}
    \centering
    \scalebox{0.95}{\input{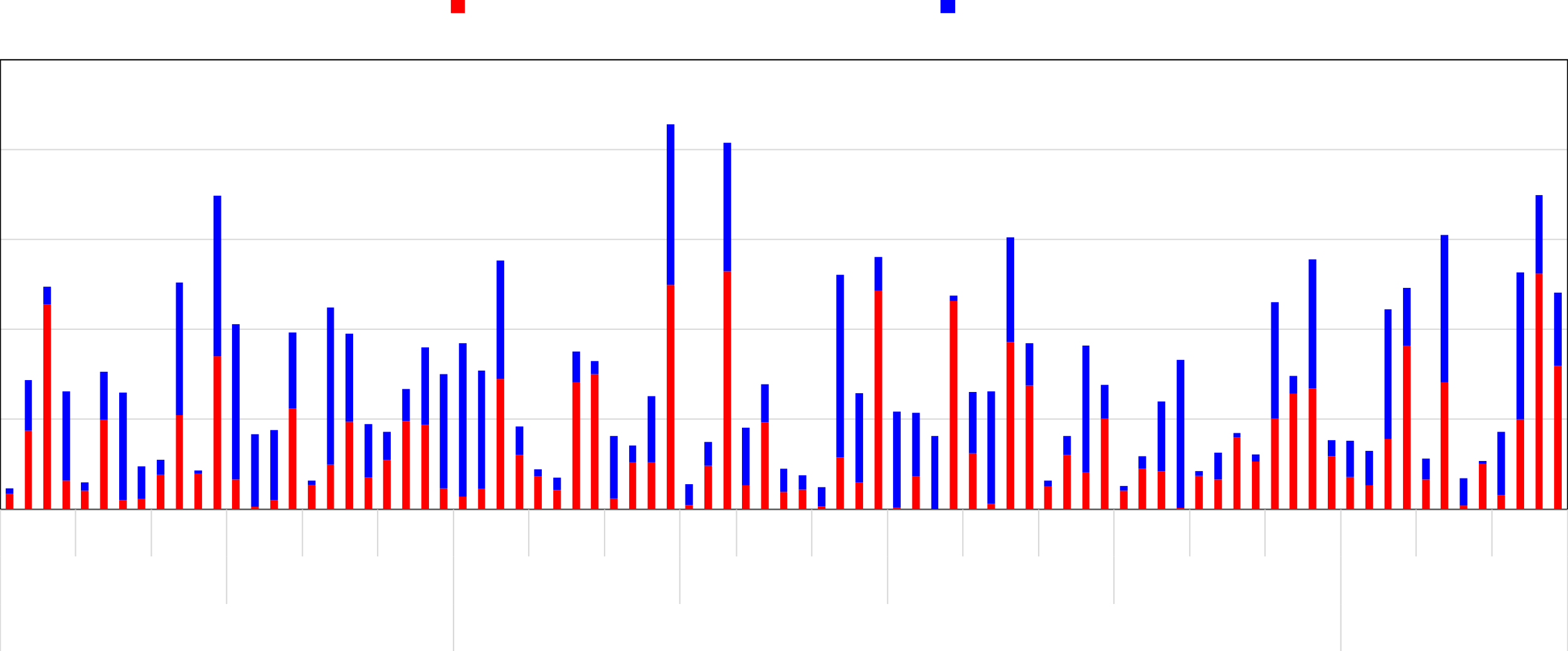}}
    \caption{\crev Improvement over the benchmark bound $\nu^{\Omega(\text{LDR})}$ in the the  budgeted-uncertainty instances. Here, $\nu^{\Omega(\text{LDR})}$ and $\nu^{\Omega(\text{2S-LDR})}$ are defined as the bounds obtained by solving problem \eqref{eq:location-transportation} using only the binding {\crev realizations} of $\LDRUB$ and $\twostageLDRUB$, respectively. }
    \label{fig:LT-dual-bound}
\end{figure}

{\crev Analyzing the performance of the C\&CG and cutting-plane algorithms in solving our instances provides further insight into the quality of the bounds.} Solution times, number of iterations{\crev, and final optimality gaps at the time limit (namely the termination gap)} are provided in Appendix \ref{app:LT-results}. {\crev Our findings show that while our methods do require extra computational effort, they yield stronger bounds as a result. Another key aspect is that a large termination gap can negatively impact the quality of the bounds.} {\crev Our results show that} in all instances, the C\&CG stops at a solution of the 2ARO model with less than 5\% termination gap. In fact, the method proves to be quite powerful in detecting the significant {\crev realizations} for the 2ARO approximation, such that in 69 out of 83 instances it achieves this gap after only two iterations. On the other hand, in 39 instances the cutting-plane method is not able to reach the optimality gap of 5\% within the time limit of 10 hours {\crev albeit having less than 10\% optimality gap in all instances with one exception of 13.4\%.}  {\crev Consequently, for calculating $\optN$, we} use the best lower bound on {\crev $\bar{\nu}^{\textsc{NA-DO}}_R$} obtained at the end of 10 hours. This can further contribute to an increased optimality gap. This suggests that, in addition to having a difficult nonlinear master problem, the cutting-plane method itself requires algorithmic enhancements, such as the design of stronger cuts.

\subsection{Robust Capital Budgeting Problem}
 In the capital budgeting problem, a company  wants to invest in a subset of $I$ projects with uncertain cost and profit, subject to an initial budget of $B$ that can be increased by getting a loan. A variant of the two-stage problem is studied by \cite{subramanyam2020k}. In the following, we formulate the multistage capital budgeting problem as an \MSARO{}. Over a planning horizon of $T$ stages, let $x_{it}(\xisupt)$ be a binary decision variable taking the value of 1 if the company decides to invest in the project $i\in\set{I}$ at stage $t\in\set{T}$, with a cost of $c_{it}(\xisupt)$ and profit of $r_{it}(\xisupt)$. Further, let $L_t(\xisupt)$ be a continuous decision variable determining the amount of loan the company decides to get at stage $t\in\set{T}$ with a unit cost of $c_{\textsc{L}}\mu^{t-1}, \mu > 1$. The \MSARO{} model is as follows:
 \bsub
 \label{eqs:capital-budget}
 \begin{align}
 	\max  \ \ & z  \label{eq:capital-budgeting} \\
 	\text{s.t.} \ \ & z\leq \sum_{t\in \set{T}}\sum_{i\in \set{I}}r_{it}(\xisupt)x_{it}(\xisupt) - \sum_{t\in \set{T}} c_{\textsc{l}}\mu^{t-1}L_t(\xisupt) &&  \xiT\in\xisupp\label{eq:capital-budgeting-obj} \\
 	& B_t(\xisupt) - B_{t-1}(\xisuptminusone) + C_{t-1}(\xisuptminusone) - L_t(\xisupt) = 0&& t\in\set{T}, \xisupt\in\xisuppt\label{eq:capital-budgeting-funds}\\
 	& \sum_{i\in \set{I}}c_{it}(\xisupt)x_{it}(\xisupt) = C_t(\xisupt) && t\in\set{T}, \xisupt\in\xisuppt\label{eq:capital-budgeting-cost-const}\\
 	& B_t(\xisupt) - C_t(\xisupt) \geq 0 && t\in\set{T}, \xisupt\in\xisuppt\label{eq:capital-budgeting-budget-const}\\
 	&x_{t}(\xisupt)\in\{0,1\}^{I}, L_t(\xisupt)\in\R_+ &&t\in\set{T}, \xisupt\in\xisuppt,\label{eq:capital-budgeting-domain}
 \end{align}
 \esub
 where $B_t(\xisupt)$ is the amount of available funds at stage $t\in\set{T}$,  determined by constraints \eqref{eq:capital-budgeting-funds}, while $C_t(\xisupt)$ is the expenditure calculated through constraints \eqref{eq:capital-budgeting-cost-const}, with initial values of $B_0(\boldsymbol{\xi}^0) = B$ and $C_0(\boldsymbol{\xi}^0)=0$. Constraints \eqref{eq:capital-budgeting-budget-const} bound the expenditure amount by  available funds. The objective is to maximize the worst-case profit over the planning horizon, measured by constraints \eqref{eq:capital-budgeting-obj}.

\subsubsection{Problem Instances}
{\crev Our instance generation  follows the procedure of \cite{subramanyam2020k} for the two-stage problem. To incorporate the dynamic nature of multistage capital budgeting, project costs and profits are modeled as affine functions of evolving risk factors $\boldsymbol{\xi}_t, t\in\set{2,T}$:
$$ c_{it}\left(\boldsymbol{\xi}_t\right)=c_{i t}^0\left(1+\boldsymbol{\Phi}_{it}^{\top} \boldsymbol{\xi}_t / 2\right),
\qquad r_{it}\left(\boldsymbol{\xi}_t\right)=r_{i t}^0\left(1+\boldsymbol{\Psi}_{it}^{\top} \boldsymbol{\xi}_t / 2\right),
$$
where $c_{it}^0$ and $r_{it}^0$ represent the baseline cost and profit, respectively, assuming all risk factors are held at their neutral value of zero; $\boldsymbol{\Phi}_{it}$ and $\boldsymbol{\Psi}_{it}$ are factor loading vectors governing the sensitivity of cost and profit to deviations from the neutral state; and $c_{i1}$ and $ r_{i1}$ are set to the nominal values. 
We consider four risk factors at each stage $t\in\set{2,T}$ such that $\boldsymbol{\xi}_t \in [-1,1]^4$. Thus, $\boldsymbol{\Phi}_{it}\in\mathbb{R}^4$ and $\boldsymbol{\Psi}_{it}\in\mathbb{R}^4$, quantify the influence of each of the four risk factors on the project's financial outcome. When sampling from $\mathbb{R}^4$, we ensure that $\boldsymbol{\Phi}_{it}^{\top} \mathbf{e} = \boldsymbol{\Psi}_{it}^{\top} \mathbf{e} = 1$ for all $i\in\mathcal{I}$ and $t\in\set{2,T}$, where $\mathbf{e}$ is a vector of ones.

For $T\in\{3,4,5\}$ and $I\in\{5,8,10\}$, the nominal cost vector $c_{it}^0, t\in\set{2,T}, i\in\set{I}$ is drawn uniformly from  $[0,10]^I$, and nominal profits are set as $r_{it}^0=\frac{c_{it}^0}{5}$. Loan purchase cost is $c_L\mu^{t-1}=0.12(1.2)^{t-1}$ per unit of loan. For each combination of $T$ and $I$ we consider different levels of initial budget which impacts the dependence between stages. 
}

\subsubsection{Optimality Gap}
Due to presence of binary variables $x_{t}(\xisupt)$, with the existing methods in the literature of 2ARO we cannot calculate  $\twostageLDRUB$ exactly.  Therefore, we solve an approximation of it using the $K$-adaptability method  of \cite{subramanyam2020k} with $K=2$ and use the obtained bound $\nu^{\texttt{K}}$ in measuring the optimality gap. On the dual side, we study three options: ($i$) $\nu^{\Omega}$, the upper bound from solving the model \eqref{eqs:capital-budget} with a sample $|\Omega| = 250$,  ($ii$) the NA bound {\crev $\bar{\nu}^{\textsc{NA-DO}}_R$} with the same sample $\Omega$, ($iii$) the decomposable NA bound {\crev $\bar{\nu}^{\textsc{DNA-DO}}_R$} with a sample of size 500{\crev, which includes the sample $\Omega$}. For each option, when the model is not solved to optimality within the time limit, the best valid bound is used in the calculations. Note that, if we can solve both options ($i$) and ($ii$) to optimality, we should expect a better bound from option ($i$).
{\crev Figure \ref{fig:capital-budget-results}} presents the optimality gaps of the capital budgeting instances between the bound $\nu^{\texttt{K}}$ and the three choices of upper bound (detailed results are given in Appendix \ref{sec:CB-results}). 
\begin{figure}[htbp]
\footnotesize
    \centering
    \scalebox{0.65}{
    % This file was created with tikzplotlib v0.10.1.
\begin{tikzpicture}

\definecolor{darkgray176}{RGB}{176,176,176}
\definecolor{forestgreen015351}{RGB}{0,153,51}
\definecolor{lightgray204}{RGB}{204,204,204}
\definecolor{magenta}{RGB}{255,0,255}
\definecolor{peru23712549}{RGB}{237,125,49}

\begin{groupplot}[group style={
    group size=1 by 4, 
    vertical sep=2cm  
  },width=1.2\textwidth,height=6cm,]
\nextgroupplot[
tick align=outside,
tick pos=left,
title={$T = 3$},
x grid style={darkgray176},
xlabel={},
xmin=1, xmax=24,
xtick=data,
xtick style={color=black},
xticklabel style={rotate=90.0,anchor=east},
y grid style={darkgray176},
ylabel={Optimality gap},
ymin=5, ymax=37,
ytick style={color=black},yticklabel=\pgfmathprintnumber{\tick},
yticklabel=\pgfmathprintnumber{\tick}\%
]
\addplot [semithick, magenta, mark=*, mark size=3, mark options={solid}]
table {%
1 19.7
2 18
3 16.2
4 13.1
5 18.7
6 18.9
7 18.5
8 18
9 16.5
10 19.4
11 18.7
12 20.4
13 16
14 21.1
15 19
16 17.5
17 14.2
18 16.5
19 21.6
20 20.5
21 18.8547380939758
22 22.3
23 22.3
24 19.7
};
\addplot [semithick, forestgreen015351, mark=asterisk, mark size=3, mark options={solid}]
table {%
1 16.07
2 14.1
3 12.77
4 9.12000000000001
5 13.8
6 14.47
7 15.18
8 14.09
9 13.26
10 15.82
11 14.26
12 16.29
13 11.18
14 17.61
15 15.79
16 14.44
17 8.25000000000001
18 12.65
19 17.83
20 16.12
21 15.3547380939759
22 17.94
23 19.13
24 15.76
};
\addplot [semithick, peru23712549, mark=triangle*, mark size=3, mark options={solid}]
table {%
1 14.5
2 12.55
3 10.35
4 7.03
5 12.27
6 11.79
7 12.56
8 11.81
9 10.98
10 12.97
11 11.8
12 13.87
13 8.88999999999999
14 15.55
15 13.78
16 12.41
17 6.71999999999999
18 10.9
19 16.37
20 14.31
21 13.9247380939758
22 16.47
23 17.61
24 14.26
};

\nextgroupplot[
tick align=outside,
tick pos=left,
title={$T = 4$},
x grid style={darkgray176},
xlabel={},
xmin=25, xmax=60,
xtick=data,
xtick style={color=black},
xticklabel style={rotate=90.0,anchor=east},
y grid style={darkgray176},
ylabel={Optimality gap},
ymin=5, ymax=37,
ytick style={color=black},
yticklabel=\pgfmathprintnumber{\tick}\%
]
\addplot [semithick, magenta, mark=*, mark size=3, mark options={solid}]
table {%
25 17.1
26 19.2
27 17.2
28 18.8
29 17.4
30 16.6992028905506
31 18.3
32 19.2
33 22.1
34 19.7
35 22.2
36 18.1
37 18
38 18.2
39 18.2
40 17.8
41 20.9
42 18.1
43 20.6
44 18.9
45 21.6
46 19.2
47 18.3
48 21.7
49 22.2
50 25.6
51 22.5
52 20.5
53 22.6
54 23.3
55 21.6
56 29.7
57 22.3
58 26.1
59 25.1
60 25.7
};
\addplot [semithick, forestgreen015351, mark=asterisk, mark size=3, mark options={solid}]
table {%
25 10.68
26 13.65
27 12
28 13.05
29 13.39
30 9.70000000000002
31 12.68
32 12.34
33 17.69
34 14.71
35 18.04
36 13.53
37 11.72
38 11.01
39 13.98
40 13.69
41 16.64
42 11.96
43 17
44 11.95
45 17.82
46 12.39
47 12.87
48 17.79
49 18.18
50 18.46
51 16.36
52 15.27
53 18.74
54 18.83
55 14.55
56 25.45
57 16.41
58 20.12
59 19.61
60 20.48
};
\addplot [semithick, peru23712549, mark=triangle*, mark size=3, mark options={solid}]
table {%
25 7.99
26 11.11
27 9.53999999999999
28 10.25
29 10.73
30 6.94000000000001
31 9.88999999999999
32 9.37
33 15.41
34 11.74
35 15.2
36 10.67
37 9.32999999999999
38 8.17999999999999
39 11.85
40 11.29
41 13.81
42 9.67999999999999
43 15.09
44 9.62
45 15.41
46 10.11
47 10.97
48 15.64
49 16
50 16.65
51 14.48
52 13.09
53 16.86
54 16.67
55 12.64
56 24.24
57 14.85
58 18.37
59 17.8
60 18.99
};

\nextgroupplot[
tick align=outside,
tick pos=left,
title={$T = 5$},
x grid style={darkgray176},
xlabel={},
xmin=61, xmax=108,
xtick=data,
xtick style={color=black},
xticklabel style={rotate=90.0,anchor=east},
y grid style={darkgray176},
ylabel={Optimality gap},
ymin=5, ymax=37,
ytick style={color=black},
yticklabel=\pgfmathprintnumber{\tick}\%
]
\addplot [semithick, magenta, mark=*, mark size=3, mark options={solid}]
table {%
61 23.6
62 23.5
63 22.4
64 24.2
65 21.7
66 25.6
67 22.6
68 22
69 19.5
70 29
71 21.9
72 24.9
73 21.9
74 29.9
75 27.7
76 31.9
77 28.1
78 26.2
79 29.9
80 30.5
81 30.7
82 28.1
83 32.2
84 29.6
85 23.1
86 27.9
87 26.6
88 23.1
89 24.5
90 34.2
91 25.9
92 28.4
93 25.9
94 24.3
95 24.6
96 30.5
97 25.9
98 26.4
99 27
100 27.6
101 26.1
102 31.3
103 28.7
104 32.8
105 33.9
106 31.9
107 25.5
108 22.4
};
\addplot [semithick, forestgreen015351, mark=asterisk, mark size=3, mark options={solid}]
table {%
61 17.36
62 16.88
63 16.44
64 19.98
65 15.55
66 19.28
67 17.66
68 15.62
69 12.55
70 23.29
71 17.59
72 18.78
73 17.13
74 24.77
75 22.39
76 25.06
77 20.82
78 22.28
79 22.48
80 23.69
81 26.34
82 23.3
83 27.78
84 23.98
85 19.59
86 23.72
87 22.55
88 17.7
89 21.76
90 31.2
91 21.72
92 25.71
93 23.41
94 21.26
95 19.68
96 25.88
97 21.4
98 23.44
99 23.02
100 24.27
101 21.79
102 27.44
103 26.02
104 28.69
105 30.23
106 27.84
107 22.1
108 18.93
};
\addplot [semithick, peru23712549, mark=triangle*, mark size=3, mark options={solid}]
table {%
61 14.54
62 13.32
63 12.99
64 16.55
65 12.12
66 16.48
67 14.46
68 12.32
69 9.61
70 20.89
71 14.61
72 15.35
73 13.83
74 22.34
75 19.31
76 22.35
77 17.45
78 19.59
79 18.91
80 19.96
81 23.79
82 20.79
83 25.38
84 21.61
85 16.41
86 20.67
87 19.98
88 15.27
89 18.92
90 28
91 19.11
92 23.27
93 20.86
94 18.31
95 17.63
96 23.13
97 18.97
98 20.62
99 21.19
100 21.82
101 19.78
102 25.31
103 24.1
104 27.22
105 27.37
106 26.23
107 20.04
108 16.54
};

\nextgroupplot[
legend cell align={left},
legend style={
  fill opacity=0.8,
  draw opacity=1,
  text opacity=1,
  at={(0.97,0.03)},
  anchor=south east,
  draw=lightgray204
},
tick align=outside,
tick pos=left,
title={$T = 6$},
x grid style={darkgray176},
xlabel={Instance},
xmin=109, xmax=162,
xtick=data,
xtick style={color=black},
xticklabel style={rotate=90.0,anchor=east},
y grid style={darkgray176},
ylabel={Optimality gap},
ymin=5, ymax=37,
ytick style={color=black},
yticklabel=\pgfmathprintnumber{\tick}\%
]
\addplot [semithick, magenta, mark=*, mark size=3, mark options={solid}]
table {%
109 20.2
110 20.5
111 24
112 31.4
113 27.7
114 30.5
115 30.1
116 27.2
117 30.7
118 27.8
119 25.8
120 25.4
121 28.3
122 29.2
123 31.5
124 30.4
125 28.3
126 29.9
127 24.7
128 30.7
129 25.8
130 31.2
131 30.8
132 34.3
133 34.8
134 29.1
135 31.2
136 27.6
137 25.4
138 27.9
139 26.4
140 28.4
141 29.8
142 30
143 32.6
144 31.6
145 28.2
146 31.5
147 31.3
148 29.9
149 32.5
150 29.4
151 30.4
152 29.6
153 32
154 35.2
155 29.5
156 34.4
157 33.9
158 35.5
159 29.7
160 32.4
161 30.2
162 32.2
};
\addlegendentry{${(\nu^{\Omega}-\nu^{\texttt{K}})}/{\nu^{\texttt{K}}}$}
\addplot [semithick, forestgreen015351, mark=asterisk, mark size=3, mark options={solid}]
table {%
109 17.98
110 18.04
111 19.83
112 28.35
113 24.89
114 27.67
115 27.2
116 24.24
117 26.92
118 26.04
119 23.49
120 22.35
121 25.8
122 26.48
123 28.74
124 28.08
125 26.26
126 26.97
127 22.47
128 28.74
129 23.91
130 29.27
131 28.34
132 31.72
133 32.9
134 26.83
135 29
136 26.46
137 24.26
138 26.32
139 25.17
140 27.23
141 28.83
142 28.32
143 31.49
144 30.59
145 27.23
146 30.3
147 29.97
148 28.5
149 31.04
150 28.46
151 28.85
152 28.59
153 30.91
154 34.05
155 28.34
156 33.13
157 32.64
158 34.36
159 28.46
160 31.41
161 28.98
162 31.13
};
\addlegendentry{${({\crev \bar{\nu}^{\textsc{NA-DO}}_R}-\nu^{\texttt{K}})}/{\nu^{\texttt{K}}}$}
\addplot [semithick, peru23712549, mark=triangle*, mark size=3, mark options={solid}]
table {%
109 14.08
110 14.19
111 16.35
112 24.81
113 21.15
114 23.07
115 22.76
116 20.27
117 23.05
118 22.51
119 19.04
120 18.22
121 21.52
122 23.06
123 24.65
124 23.99
125 22.19
126 23.3
127 19.33
128 25.8
129 20.94
130 26.11
131 25.5
132 28.74
133 29.71
134 23.76
135 25.9
136 23.9
137 21.6
138 23.73
139 22.83
140 24.37
141 26.42
142 25.99
143 29.26
144 28.08
145 24.93
146 28.21
147 28.03
148 26.31
149 29.21
150 26.55
151 26.84
152 26.68
153 28.69
154 32.15
155 26.57
156 31.53
157 30.78
158 32.49
159 26.87
160 29.43
161 27.12
162 29.17
};
\addlegendentry{${({\crev \bar{\nu}^{\textsc{DNA-DO}}_R}-\nu^{\texttt{K}})}/{\nu^{\texttt{K}}}$}
\end{groupplot}

\end{tikzpicture}
    }
    \caption{\crev Optimality gap improvements for capital budgeting problems when using LDDR-based methods. In the legend, $\nu^{\Omega}$ denotes the bound obtained by solving model \eqref{eqs:capital-budget} over a sample set.}
    \label{fig:capital-budget-results}
\end{figure}
Results show that, even though the bound $\nu^{\Omega}$ should theoretically be at least as good as the bound {\crev $\bar{\nu}^{\textsc{NA-DO}}_R$}, on average the NA bound returns a better upper bound within the same time limit. This is a testament to the difficulty of the multistage problem even when it is solved for a discrete set of {\crev realizations}. The decomposable NA bound {\crev $\bar{\nu}^{\textsc{DNA-DO}}_R$} further improves the results by using a larger sample which is viable because of its superior computational performance. {\crev Note that, the best gaps in Figure \ref{fig:capital-budget-results}, ranging between 7\% to 33\%, are obtained from approximations over approximations on both primal and dual side.}  Accordingly, these rather large gaps can be attributed to both bounding methods. 

\paragraph{\crev General Integer Recourse Variables.}
{\crev To assess our algorithms' performance with general integer recourse, we conducted a set of experiments where loan amounts, $L_t(\xisupt)$, were restricted to integer values. Figure \ref{fig:capital-budget-int-results} in Appendix~\ref{app:results-generalinteger} presents the same analysis as in Figure \ref{fig:capital-budget-results} for the capital budgeting problem, but with the added constraint of $L_t(\xisupt)\in\Z_+$.
Interestingly, requiring integer loan amounts did not significantly alter the optimal investment decisions compared to the continuous case. This explains the visual similarity between the two figures, with both exhibiting similar patterns in optimal decisions. While there are slight differences in objective function values, the overall investment strategies remain largely unaffected by the integrality constraints.

It is important to note that although the advantage of our dual methods appears to diminish with an increasing number of stages, even the initial gap  relies on the strength of our primal solution, $\nu^{\texttt{K}}$. Therefore, these figures demonstrate our ability to achieve further improvements beyond the initial strong primal bound.

Appendix \ref{sec:CB-results} provides a detailed analysis of the computational requirements for bounding the capital budgeting problem, considering both fractional and integral loans. For the primal side, we present the solution times of our approach, as no alternative primal method is available for comparison. On the dual side, we compare the solution times of the three previously discussed bounding approaches. The solution times demonstrate that the extra computational effort associated with the decomposable NA model yields a demonstrably improved dual bound.
}

%%%%%%%%%%%%%%%%%%%%%%%%%%%%%%%%
%%%%%%%%%%%%%%%%%%%%%%%%%%%%%%%%
%%%%%%%%%%%%%%%%%%%%%%%%%%%%%%%%
%%%%%%%%%%%%%%%%%%%%%%%%%%%%%%%%
\section{Conclusion}

Robust optimization models are built on a different premise than stochastic programming in the sense that they do not assume any knowledge about the probability distribution, focusing instead on optimizing the worst-case outcomes. {\crev In this paper, we study general MSAROs for which we develop primal and dual bounding methods by adapting two decision rule approximations from the stochastic programming literature (namely two-stage LDRs and LDDRs). These approximations allow us to reduce MSARO to a 2ARO from the primal side and a two-stage stochastic programming problem from the dual side. As such, the resulting approximations drastically reduce the theoretical complexity of the studied problems. Since our dual bounds are dependent on the choice of a probability distribution while deriving the dual model, we propose to solve a distribution optimization problem to obtain a stronger bound. We develop various solution methods for our proposed bounding problems where we also leverage existing  methods from both the robust optimization and the stochastic programming literature. Our extensive numerical study demonstrates that our methods considerably improve both primal and dual bounds compared to the commonly used approaches in the literature. Our work opens the door to the direct application of existing two-stage robust optimization and stochastic programming algorithms and other future algorithmic developments in these areas to MSAROs. For instance, as the algorithms such as C\&CG and K-adaptability for 2ARO improve, our models can be solved more efficiently.

We believe that our work can initiate additional methodological and numerical developments. Both from the primal and the dual side the question of solving the problems we pose in a more numerically efficient manner definitely merits more attention. Further methodological work may also explore the synergies between primal and dual decision rules.  Finally, the following directions can be the subject of future research: exploiting problem structure in order to approximate MSAROs with numerically more favorable, e.g, decomposable models,  identifying special cases of MSAROs in which the proposed approximations can be proven to be exact, and application of similar approaches in related fields such as distributionally robust optimization. 
}

\ACKNOWLEDGMENT{This work was supported by Natural Sciences and Engineering Research Council of Canada [Grants RGPIN-2018-04984 and RGPIN-2024-05908] and Agence Nationale de la Recherche of France [Grant ANR-22-CE048-0018]. Computations were performed on the Niagara supercomputer at the SciNet HPC Consortium. SciNet is funded by: the Canada Foundation for Innovation; the Government of Ontario; Ontario Research Fund - Research Excellence; and the University of Toronto.

For the purpose of Open Access, a CC-BY public copyright licence
has been applied by the authors to the present document and will
be applied to all subsequent versions up to the Author Accepted
Manuscript arising from this submission.\\[3mm]
\includegraphics[width=2cm]{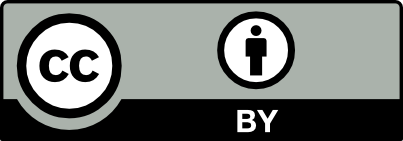}}

% ------------------------------
\bibliographystyle{informs2014}
\bibliography{References/MSARO,References/MSMIP}
% ------------------------------
\newpage
\begin{appendices}
% \input{discussions}
% \newpage
{\crev
\section{Obtaining 2ARO Model via Two-stage Decision Rules} \label{app:2AROreduction}
In this section, we detail the transition from the 
MSARO model \eqref{eq:msaro-xt} to the 2ARO model \eqref{eq:2RO} as a result of applying two-stage decision rules proposed in Section \ref{subsec:GeneralDecisionsRule}. Substituting the state variables in the MSARO problem with the decision rules $\statevart=\twostageBFt, t \in [2,T]$, we obtain the following model: 
\bsub
\label{eq:msaro_2SDR}
\begin{align}
\min\ & c_1^\top x_1 + z^\text{rest}   \label{eq:obj-msaro_2SDR-mono-xs} \\
\text{s.t.}\ 
& x_1\in X_1 & \\
&\twostageLDRVart \in \R^{K_t}\qquad\ \ && t\in\set{2,T} \\
& z^\text{rest} \geq
 \sum_{t\in \set{2,T}}{\stateCt}^\top \twostageBFt + \sum_{t\in \set{2,T}} {\recCt}^\top \recvart && \xiT\in\xisupp \\
&\ \recAt \recvart \leq b_{t}(\xisupt)-\Big(\stateAt \twostageBFt + \stateBt \statevar_1 \Big) && t = 2, \xisupt \in \xisuppt \label{eq:msaro_2SDR_t2} \\
	&\ \recAt \recvart \leq b_{t}(\xisupt)-\Big(\stateAt \twostageBFt + \stateBt \twostageBFtminusone \Big) && t\in \set{3,T}, \xisupt \in \xisuppt   \\
	&\ (\twostageBFt, \recvart)\in X_{t}(\xisupt)&&   t\in \set{2,T}, \xisupt \in \xisuppt  \label{eq:msaro_2SDR_Xt}
\end{align}
\esub
Given the first-stage decisions $x_1$ and $\beta$, observe that the feasible set of the recourse variables, defined by \eqref{eq:msaro_2SDR_t2}-\eqref{eq:msaro_2SDR_Xt} decomposes by stage $t$ and history $\xi^t$ (as both the decision variables, $\recvart$, and the constraints, \eqref{eq:msaro_2SDR_t2}-\eqref{eq:msaro_2SDR_Xt}, are separately defined for each stage and history, i.e., there is no link between them). For convenience, let us represent this feasible space in a decomposed form via $\recvart \in \mathcal{X}^\texttt{r}_t(\statevar_1, \beta,\xi^t), t \in [2,T], \xisupt \in \xisuppt$, and write \eqref{eq:msaro_2SDR} more compactly as follows:
\bsub
\begin{align}
\min_{x_1\in X_1, \beta} \ & c_1^\top x_1 + z^\text{rest} \\
\text{s.t.}\ 
& z^\text{rest} \geq
 \sum_{t\in \set{2,T}}{\stateCt}^\top \twostageBFt + \sum_{t\in \set{2,T}} {\recCt}^\top \recvart \label{eq:msaro_2SDRcompact_zRest} && \xiT\in\xisupp \\
	&\recvart \in \mathcal{X}^\texttt{r}_t(\statevar_1, \beta,\xi^t)&&   t\in \set{2,T}, \xisupt \in \xisuppt 
\end{align}
\esub
In order to minimize $z^\text{rest}$, the second summation term in \eqref{eq:msaro_2SDRcompact_zRest} should be minimized over the recourse decisions $\recvart$. Since this term is additively separable over stages $t\in [2,T]$ and their associated history $\xisupt \in \xisuppt $, the recourse decisions of each stage can be optimized separately, yielding the following equivalent model:   
\bsub
\begin{align}
\min_{x_1\in X_1, \beta} \ & c_1^\top x_1 + z^\text{rest} \\
\text{s.t.}\ 
& z^\text{rest} \geq
 \sum_{t\in \set{2,T}}{\stateCt}^\top \twostageBFt + \sum_{t\in \set{2,T}} \min_{\recvart \in \mathcal{X}^\texttt{r}_t(\statevar_1, \beta,\xi^t)} {\recCt}^\top \recvart  && \xiT\in\xisupp  \label{eq:msaro_2SDRcompact_InsideMin}
\end{align}
\esub
Inspecting the minimization problems in \eqref{eq:msaro_2SDRcompact_InsideMin}, we observe that the objective function coefficients and the feasible set are parametrized by the history, as such changing the parametrization of the recourse decisions from the history $\xi^t$, to a full uncertainty realization $\xi^T$ should not change the optimal objective value. Thus, instead, we can consider the following equivalent model:
\bsub
\begin{align}
\min_{x_1\in X_1, \beta} \ & c_1^\top x_1 + z^\text{rest} \\
\text{s.t.}\ 
& z^\text{rest} \geq
 \sum_{t\in \set{2,T}}{\stateCt}^\top \twostageBFt + \sum_{t\in \set{2,T}} \min_{\recvar_t(\xi^T) \in \mathcal{X}^\texttt{r}_t(\statevar_1, \beta,\xi^t)} {\recCt}^\top \recvar_t(\xi^T)  && \xiT\in\xisupp  \label{eq:msaro_2SDRcompact_InsideMinFullHistory} 
\end{align}
\esub
{\crev More formally, we can show that given an optimal solution to the inner minimization problem in \eqref{eq:msaro_2SDRcompact_InsideMinFullHistory}, $\hat{x}^\texttt{r}_t: \xi^T \rightarrow \R^{p_t}$, we can construct a feasible solution $\tilde{x}^\texttt{r}_t: \xi^t \rightarrow \R^{p_t}$ to the inner minimization problem in \eqref{eq:msaro_2SDRcompact_InsideMin} which attains the same objective value. Since the minimization problem in \eqref{eq:msaro_2SDRcompact_InsideMinFullHistory} is a relaxation of that in \eqref{eq:msaro_2SDRcompact_InsideMin}, having more flexible decision variables, the aforementioned construction is sufficient to conclude the proof of our equivalence claim. Given stage $t\in[2,T]$ and history  $\Breve{\xi}^t \in \xisuppt$, 
we let $\tilde{x}^\texttt{r}_t(\Breve{\xi}^t) := \hat{x}^\texttt{r}_t(\ring{\xi}^T)$ where $\ring{\xi}^T$ is an arbitrarily selected element from the set of realizations with the same history $\{ \xi^T \in \xisupp : \xi^t = \Breve{\xi}^t \}$. The feasibility of the constructed policy is straightforward since 
the two inner minimization problems have the same feasibility set, defined by $ \mathcal{X}^\texttt{r}_t(\statevar_1, \beta,\xi^t), t \in [2,T]$. Next, we observe that ${\recCt}^\top \hat{x}^\texttt{r}_t(\xi^T)$ is the same for all realizations in $\{ \xi^T \in \xisupp : \xi^t = \Breve{\xi}^t \}$ since the minimization problem in \eqref{eq:msaro_2SDRcompact_InsideMinFullHistory} has the same feasible set and objective coefficient vector for any of those realizations and $\hat{x}^\texttt{r}_t(\cdot)$ is chosen to be an optimal policy. Lastly, by construction, the policy $\tilde{x}^\texttt{r}_t(\cdot)$ achieves the same objective value.

}

As there is no link between the different stage optimization problems in \eqref{eq:msaro_2SDRcompact_InsideMinFullHistory}, we can swap the summation and minimization operators and optimize over the recourse variables associated with all the stages together given an uncertainty realization: 
\bsub
\begin{align}
\min_{x_1\in X_1, \beta} \ & c_1^\top x_1 + z^\text{rest} \\
\text{s.t.}\ 
& z^\text{rest} \geq
 \sum_{t\in \set{2,T}}{\stateCt}^\top \twostageBFt + \min_{(\recvar_t(\xi^T))_{t \in [2,T]} \in \mathcal{X}({\crev \statevar_1}, \twostageLDR,\xiT)} \sum_{t\in \set{2,T}}
 {\recCt}^\top \recvar_t(\xi^T)  && \xiT\in\xisupp  \label{eq:msaro_2SDRcompact_InsideMinFullHistoryTwo} 
\end{align}
\esub
where $\mathcal{X}({\crev \statevar_1}, \twostageLDR,\xiT) = \mathop{\scalebox{1.5}{$\times$}}\limits_{t \in [2,T]}
 \mathcal{X}^\texttt{r}_t(\statevar_1, \beta,  \xi^t)$. Since $z^\text{rest}$ is equal to the maximum of the right-hand side of \eqref{eq:msaro_2SDRcompact_InsideMinFullHistory} over the realizations $\xiT\in\xisupp$ at an optimal solution, we can turn this problem into the following nested formulation:
\bsub
\begin{align}
	\min \ \ & c_1^\top x_1 + 
	\max_{\xiT\in\xisupp}\  
	\min_{\recvar\in\mathcal{X}({\crev \statevar_1}, \twostageLDR, \xiT)} \sum_{t\in \set{2,T}} {\stateCt}^\top \twostageBFt + {\recCt}^\top \recvar_t 
\end{align}
\esub
where the parametrization of the recourse variables is omitted since $\xi^T$ is given as an input to the inner minimization problem. 
}

\section{Proofs} \label{sec:proofs}
In this section, we present the proofs of the three propositions and one lemma mentioned in the body of the paper, for which we also restate the claims for convenience. 
\begin{repeattheorem}[Proposition \ref{prop:two-stage-ccg}.]
	{\crev Consider an \MSARO{} with only right-hand-side uncertainty,  continuous recourse, and (bounded) polyhedral uncertainty set. If the basis functions $\twostageLDRBFt$ are chosen to be affine in $\xisupt$ for all $t \in \set{2,T}$, the C\&CG algorithm converges to  $\twostageLDRUB$   
 in a finite number of iterations.}
\end{repeattheorem}
\begin{proof}{Proof}
	Denote by $f({\crev \statevar_1},\twostageLDR,\xiT)$, the objective function of the inner m{\crev in}imization in the two-stage problem \eqref{eq:2RO-linear}. Then, we can rewrite  problem \eqref{eq:2RO-linear} as: $$\twostageLDRUB = \min \big\{ c_1^\top x_1 + 
	\max_{\xiT\in\xisupp} f({\crev \statevar_1},\twostageLDR,\xiT)\ |\ x_1\in X_1,\  \twostageLDRVart \in \R^{K_t},\  t\in\set{2,T}\big\}.$$
	First, we show that $f({\crev \statevar_1},\twostageLDR,\xiT)$ is convex in $\xiT$ for given $x_1 {\crev = ({\crev \statevar_1}, {\crev \recvar_1}) }\in X_1$ {\crev and $\twostageLDR$}. For $\hat{\xiparrand}^T,\tilde{\xiparrand}^T\in\xisupp$ and  $\lambda\in\set{0,1}$,  the following (in)equalities hold: 
	\bsub
	\begin{align}
    &\lambda f({\crev \statevar_1},\twostageLDR,\hat{\xiparrand}^T) + (1-\lambda)f({\crev \statevar_1},\twostageLDR,\tilde{\xiparrand}^T) &&= \nonumber\\
    &\sum_{t\in \set{2,T}} {\stateCtFixed}^\top \big(\lambda {\crev \twostageLDRBF_{\crev t}(\hat{\xiparrand}^T)} + (1-\lambda) {\crev \twostageLDRBF_{\crev t}(\tilde{\xiparrand}^t)} \big) {\crev \twostageLDRVart} + \nonumber\\
    & \qquad  \lambda\ \big(\min_{\recvar\in\mathcal{X}({\crev \statevar_1},\twostageLDR, \hat{\xiparrand}^T)} \sum_{t\in \set{2,T}} {\recCtFixed}^\top \recvar_t\big) + (1-\lambda) \big(\min_{\recvar\in\mathcal{X}({\crev \statevar_1},\twostageLDR, \tilde{\xiparrand}^T)} \sum_{t\in \set{2,T}} {\recCtFixed}^\top \recvar_t\big) && = \label{eq:BF-affine}\\ 
    &\sum_{t\in \set{2,T}} {\stateCtFixed}^\top {\crev \twostageLDRBF_{\crev t}(\lambda\hat{\xiparrand}^t+(1-\lambda)\tilde{\xiparrand}^t)  \twostageLDRVart} + \nonumber\\
    & \qquad  \lambda\ \big(\min_{\recvar\in\mathcal{X}({\crev \statevar_1},\twostageLDR, \hat{\xiparrand}^T)} \sum_{t\in \set{2,T}} {\recCtFixed}^\top \recvar_t\big) + (1-\lambda) \big(\min_{\recvar\in\mathcal{X}({\crev \statevar_1},\twostageLDR, \tilde{\xiparrand}^T)} \sum_{t\in \set{2,T}} {\recCtFixed}^\top \recvar_t\big) &&\geq\label{eq:LP-convex}\\ 
    &\sum_{t\in \set{2,T}} {\stateCtFixed}^\top {\crev \twostageLDRBF_{\crev t}(\lambda\hat{\xiparrand}^t+(1-\lambda)\tilde{\xiparrand}^t) \twostageLDRVart} + \min_{\recvar\in\mathcal{X}({\crev \statevar_1},\twostageLDR, \lambda\hat{\xiparrand}^T {\crev +} (1-\lambda)\tilde{\xiparrand}^T)} \sum_{t\in \set{2,T}} {\recCtFixed}^\top \recvar_t && = f({\crev \statevar_1},\twostageLDR,\lambda\hat{\xiparrand}^T+ (1-\lambda)\tilde{\xiparrand}^T).\nonumber
	\end{align}
	\esub
    Equality \eqref{eq:BF-affine} holds because $\twostageLDRBF_{\crev t}(\xisupt)$ is affine in $\xisupt$. Inequality \eqref{eq:LP-convex} follows from the convexity of the optimal value of the inner minimization problem {\crev as a function of $\xiT$, since uncertainty} appears only on the right-hand-sides of constraints and {\crev recourse variables} are continuous {\crev by assumption}.
    
    Now, let $\XiExt$ be the set of extreme points of $\xisupp$. 
    In maximization of a convex function over a compact polyhedral set, there is an optimal solution that is an extreme point \citep{hendrix2010introduction}.  Then, the two-stage problem becomes:
	$$\twostageLDRUB = \min \big\{ c_1^\top x_1 + \eta
	\ |\ \eta \geq f({\crev \statevar_1},\xiT,\twostageLDR),\  \xiT\in\XiExt,\  x_1\in X_1,\  \twostageLDRVart \in \R^{K_t},\  t\in\set{2,T}\big\}.$$
	C\&CG is then the process of gradually adding constraints for each extreme point. Because $|\XiExt|<+\infty$, iterating over all extreme points   takes finitely many steps, which concludes the proof. 
	\hfill\Halmos
\end{proof}
\begin{repeattheorem}[Lemma \ref{lem:NA-const-alt}.]
For {\crev any $ \dist \in \mathcal{P}^{>} $}, constraints \eqref{eq:const-reform-nonanticipative-lddro} are equivalent to the following:
\begin{equation}
\tag{\ref{eq:NA-const-Exp-form}}
    y_t(\xiT) = \Exp_{\xiprimeT\sim\dist}\left[y_t(\xiprimeT) \ \Big|\ \xiprimesupt = \xisupt \right],\quad  t\in \set{T},\ \xiT\in\xisupp.
\end{equation}
\end{repeattheorem}
\begin{proof}{Proof}
{\crev 
For a fixed $t$ and {\crev$\xi^T$}, we start by multiplying both sides of \eqref{eq:const-reform-nonanticipative-lddro} by the density $p^{\dist}(\xiprimeT)$ for every $\xiprimeT$ (including $\xiT$) that shares the same history with $\xiT$ up to $t$: 
$$p^{\dist}(\xiprimeT) y_t(\xiT) = p^{\dist}(\xiprimeT) y_t(\xiprimeT)\qquad  t\in \set{T},\ \xiT,\xiprimeT\in\xisupp \text{ with }  \xisupt=\xiprimesupt.$$
Since $p^{\dist}(\xiprimeT) >0 $ for every realization $\xiprimeT$ under any probability distribution $\dist \in \mathcal{P}^{>}$, the feasible set of model \eqref{eq:NA-reform} remains the same. Then, we integrate (or sum if $\xisupp$ is discrete) the scaled constraints over all such $\xiprimeT$ realizations to obtain
$$y_t(\xiT) 
\int\limits_{\xiprimeT : \xiprimesupt = \xisupt} \hspace*{-0.3cm} p^{\dist}(\xiprimeT) d\xiprimeT = \int\limits_{\xiprimeT : \xiprimesupt = \xisupt} \hspace*{-0.3cm} p^{\dist}(\xiprimeT) y_t(\xiprimeT) d\xiprimeT \qquad  t\in \set{T},\ \xiT\in\xisupp.$$
Let $\delta := \int\limits_{\xiprimeT : \xiprimesupt = \xisupt} \hspace*{-0.3cm} p^{\dist}(\xiprimeT) d\xiprimeT$ where $\delta > 0$ since $\xi^T\in \{\xiprimeT : \xiprimesupt = \xisupt\}$ with $p^{\dist}(\xiT) > 0$. Then, we obtain \eqref{eq:NA-const-Exp-form} via  
$$y_t(\xiT) = \int\limits_{\xiprimeT : \xiprimesupt = \xisupt} \hspace*{-0.3cm} \frac{p^{\dist}(\xiprimeT)}{\delta} y_t(\xiprimeT) d\xiprimeT = \Exp_{\xiprimeT\sim\dist}\left[y_t(\xiprimeT) \ \Big|\ \xiprimesupt = \xisupt \right] \qquad  t\in \set{T},\ \xiT\in\xisupp.$$
Since constraints \eqref{eq:NA-const-Exp-form} are obtained as an aggregation of the original constraints \eqref{eq:const-reform-nonanticipative-lddro}, they are valid for model \eqref{eq:NA-reform}. Further, they imply the nonanticipativity constraints \eqref{eq:const-reform-nonanticipative-lddro} since for any given $t$ along with $\xiT$ and $\hat{\xi}^T$ such that $ \xisupt = \hat{\xi}^t$, the right-hand side of \eqref{eq:NA-const-Exp-form} is the same, i.e., $\Exp_{\xiprimeT\sim\dist}\left[y_t(\xiprimeT) \ \Big|\ \xiprimesupt = \xisupt \right] = \Exp_{\xiprimeT\sim\dist}\left[y_t(\xiprimeT) \ \Big|\ \xiprimesupt = \hat{\xi}^t \right]$, enforcing $y_t(\xiT)=y_t(\hat{\xi}^T)$. \hfill\Halmos
}
\end{proof}
\begin{repeattheorem}[Proposition \ref{prop:exactness-NA}.]
{\crev
Let $\dist$ be any probability measure in {\crev $\mathcal{P}^{>}$}. For \MSARO{} problems with \emph{continuous recourse}, \eqref{eq:NA_bound} is a strong dual of \eqref{eqs:msaro-mono}, i.e., $\mathcal{L}^{\textsc{NA}}(\dist)=\exactB$.
}
\end{repeattheorem}
\begin{proof}{Proof}
    The strength of a Lagrangian dual problem can be studied by its primal characterization, derived  by \cite{Geoffrion1974} for a mixed-integer linear optimization problem.
    Then, using standard Lagrangian duality theory (see, for instance, \cite{wolsey1999integer}), the primal characterization of \eqref{eq:NA_bound} is:
    \bsub
    \label{eq:NA-reform-primal-char}
    \begin{align}
    \min  \ \ & z  \label{eq:obj-NA-reform-primal-char} \\
    \text{s.t.} \ \ & \big(z,y_1(\xiT),\dots,y_T(\xiT)\big) \in \conv{Y(\xiT)} &&  \xiT\in\xisupp \label{eq:const-remain-primal-char}\\
    & y_t(\xiT) = \Exp_{\xiprimeT\sim\dist}\left[y_t(\xiprimeT) \ \Big|\ \xiprimesupt = \xisupt \right]&&  t\in \set{T}, \xiT\in\xisupp. \label{eq:const-NA-primal-char}
    \end{align}
    \esub
    The result follows from the equivalence of \eqref{eq:NA-const-Exp-form} and \eqref{eq:const-reform-nonanticipative-lddro}, as well as the fact that for \MSARO{} problems {\crev with continuous recourse}, we have that $\conv{Y(\xiT)} = Y(\xiT)$ for all $\xiT\in\xisupp$.
    \hfill\Halmos
\end{proof}
\begin{repeattheorem}[Proposition \ref{eq:Pbar-valid}.]
{\crev $\bar{\nu}^{\textsc{NA-DO}}_R$} is a lower bound for $\exactB$.
\end{repeattheorem}

\begin{proof}{Proof}
We have shown that $\mathcal{L}^{\textsc{NA}}_R({\dist})\leq\exactB$ for any ${\dist} \in \mathcal{P}^{>}$.
Now we claim that, for any $\hat{\dist}\in \mathcal{P}^{\geq}\setminus\mathcal{P}^{>} = \big\{\texttt{P}\in\mathcal{P}^{\geq}\ |\ \exists \xiT\in\xisupp:\ p^{\texttt{P}}(\xiT) = 0\big\}$ where $p^{\hat{\dist}}:\xisupp\rightarrow\R_+$ is the density function of  $\hat{\dist}$, the inequality $\mathcal{L}^{\textsc{NA}}_R(\hat{\dist})\leq\exactB$ also holds. 
{\crev
Given $\hat{\dist}$, let $\hat{\Xi}^T := \{ \xiT\in\xisupp : p^{\hat{\dist}}(\xiT) > 0\}$. 
Using $\hat{\Xi}^T$, we create a relaxation of \eqref{eq:NA-reform}, the NA reformulation of the MSARO problem, as 
\bsub
\label{eq:NA-reform_GivenPmiddle}
\begin{align}
\hat{\nu}(\hat{\mathbb{P}}) := \min  \ \ & z  \label{eq:obj-NA-reformmiddle} \\
\text{s.t.} \ \ & \sum_{t\in \set{T}}c_t(\xisupt)^\top y_t(\xiT) \leq z &&  \xiT\in\xisupp \label{eq:const-obj-NA-reform_GivenPmiddle}\\
& \At y_t(\xiT) + \Bt y_{t-1}(\xiT)\leq b_t(\xisupt) &&  t\in {\crev [2,T]},\ \xiT\in\xisupp \label{eq:const-state-NA-reform_GivenPmiddle}\\
& \Dt y_t(\xiT)\leq d_t(\xisupt) &&   t\in \set{T},\ \xiT\in\xisupp \label{eq:const-rec-NA-reform_GivenPmiddle}\\
& y_t(\xiT) = y_t(\xiprimeT) &&  t\in \set{T},\ \xiT \in\hat{\Xi}^T,\xiprimeT\in\hat{\Xi}^T \text{ with }  \xisupt=\xiprimesupt \label{eq:const-reform-nonanticipative-lddro_GivenPmiddle}\\
& y_t(\xiT)\in\R^{n_t-n^\texttt{i}_t}\times\Z^{n^\texttt{i}_t} &&   t\in \set{T},\ \xiT\in\xisupp \label{eq:const-domain-NA-reform_GivenPmiddle} 
\end{align}
\esub
where the NA constraints are only imposed for the pairs of realizations in $\hat{\Xi}^T \subset \Xi^T$. Therefore, we have $\hat{\nu}(\hat{\mathbb{P}}) \leq \nu^*$. 

Now that $\hat{\dist}$ assigns a positive density to all the realizations from $\hat{\Xi}^T$, via Lemma \ref{lem:NA-const-alt}, the NA constraints in \eqref{eq:const-reform-nonanticipative-lddro_GivenPmiddle} can be equivalently reformulated as 
$$y_t(\xiT) = \Exp_{\xiprimeT\sim\hat{\mathbb{P}}}\left[y_t(\xiprimeT) \ \Big|\ \xiprimesupt = \xisupt \right],\quad   t\in \set{T},\ \xiT\in \hat{\Xi}^T.$$ 
}
Subsequently, relaxation of {\crev these reformulated NA constraints} and construction of the Lagrangian dual problem \eqref{eq:NA_bound} with respect to $\hat{\dist}$ yields a relaxation of a relaxation of a minimization problem, where for each fixed $\lambda_t(\cdot)$ we have $\mathcal{L}^{\textsc{NA}}_{\textsc{LR}}(\hat{\dist},\lambda_1(\cdot).\dots,\lambda_T(\cdot)) \leq {\crev \hat{\nu}(\hat{\mathbb{P}})}$,  and consequently $\mathcal{L}^{\textsc{NA}}(\hat{\dist})\leq {\crev \hat{\nu}(\hat{\mathbb{P}})}$. {\crev Furthermore, restricting the Lagrangian duals to follow LDRs, we obtain}   $\mathcal{L}^{\textsc{NA}}_R(\hat{\dist})\leq\mathcal{L}^{\textsc{NA}}(\hat{\dist})$.

{\crev Since we showed  $\mathcal{L}^{\textsc{NA}}_R(\dist)\leq \nu^*$ for all $\dist \in \mathcal{P}^{\geq}$, we have $\bar{\nu}^{\textsc{NA-DO}}_R := \max_{\dist\in \mathcal{P}^{\geq}} \ \mathcal{L}^{\textsc{NA}}_R(\dist) \leq \nu^*$,} which completes the proof.
\hfill\Halmos
\end{proof}
\section{Detailed Models}

\subsection{Monolithic Form of Model \eqref{eq:2RO-SP}}\label{sec:monolithic-model}
{\crev
Consider the case mentioned in Remark \ref{rem:monoloithicSP}, namely an MSARO where the uncertainty set is a polytope, the basis functions $\twostageLDRBFt$ are chosen to be affine in $\xisupt$ for all $t\in\set{2,T}$, all the recourse variables are continuous, we have fixed parameters associated with the state variables, and the 2ARO problem \eqref{eq:2RO-linear} has relatively complete recourse. We next detail how a monolithic mixed-integer linear programming formulation of the inner minimization problem \eqref{eq:2RO-SP} can be derived as in \citep{ayoub2016decomposition, zeng2013solving}. 
}

Let $\dualSPstate$ and $\dualSPrec$ be the {\crev linear programming} dual variables associated with the state and recourse constraints in $\mathcal{X}({\crev \hat{x}^\texttt{s}_1}, \hat{\beta}, \xiT)$, in the inner minimization problem of  \eqref{eq:2RO-SP}, respectively. Then, {\crev using KKT conditions,} subproblem \eqref{eq:2RO-SP} can be modelled as follows:
\bsub
 \label{eq:SP-mono}
\begin{align}
	\max\ \ & \sum_{t\in \set{2,T}} {\stateCtFixed}^\top {\crev \twostageLDRBFt \hat{\twostageLDR}_t} + {\recCtFixed}^\top \recvar_t \\
	\text{s.t.}\ 
    & {\crev \recAtFixed \recvar_t  + \stateAtFixed {\crev \twostageLDRBFt \hat{\twostageLDR}_t} + \stateBtFixed \hat{x}^\texttt{s}_1 -b_{t}(\xisupt) \leq 0} && {\crev t = 2} \label{eq:two-stage-LDR-cont-primal1-Part1} \\
    &\recAtFixed \recvar_t  + \stateAtFixed {\crev \twostageLDRBFt \hat{\twostageLDR}_t} + \stateBtFixed{\crev \twostageLDRBFtminusone \hat{\twostageLDR}_{t-1}} -b_{t}(\xisupt) \leq 0 && t\in \set{{\crev 3},T}  \label{eq:two-stage-LDR-cont-primal1-Part2} \\
	& \stateDtFixed {\crev \twostageLDRBFt \hat{\twostageLDR}_t} + \recDtFixed\recvar_t-d_t(\xisupt) \leq 0 &&   t\in \set{2,T}\label{eq:two-stage-LDR-cont-primal2}\\
	& {\recAtFixed}^\top\dualSPstatet + {\recDtFixed}^\top\dualSPrect = \recCtFixed && t\in \set{2,T}\label{eq:two-stage-LDR-cont-dual_feasibility}\\
 	& {\crev \Big(\recAtFixed \recvar_t  + \stateAtFixed {\crev \twostageLDRBFt \hat{\twostageLDR}_t} + \stateBtFixed \hat{x}^\texttt{s}_1 -b_{t}(\xisupt)\Big)^\top \dualSPstatet = 0}  && {\crev t = 2} \label{eq:complementary-slack-state-Part1}\\
	& \Big(\recAtFixed \recvar_t  + \stateAtFixed {\crev \twostageLDRBFt \hat{\twostageLDR}_t} + \stateBtFixed {\crev \twostageLDRBFtminusone \hat{\twostageLDR}_{t-1}} -b_{t}(\xisupt)\Big)^\top \dualSPstatet = 0  && t\in \set{{\crev 3},T} \label{eq:complementary-slack-state-Part2}\\
	& \Big(\stateDtFixed {\crev \twostageLDRBFt \hat{\twostageLDR}_t} + \recDtFixed\recvar_t-d_t(\xisupt)\Big)^\top  \dualSPrect = 0 && t\in \set{2,T} \label{eq:complementary-slack-recourse}\\
	& \recvar_t\in\R^{p_t},\quad \dualSPstatet\in\R_-^{{\mstatet}},\quad \dualSPrect\in\R_-^{{\mrecourset}}&& t\in \set{2,T} \\
 & \xiT \in\xisupp.
\end{align}
\esub
Inequalities \eqref{eq:two-stage-LDR-cont-primal1-Part1}-\eqref{eq:two-stage-LDR-cont-primal2} are the primal feasibility constraints at $\hat{\beta}_t$, \eqref{eq:two-stage-LDR-cont-dual_feasibility} the dual feasibility constraints, and \eqref{eq:complementary-slack-state-Part1}-\eqref{eq:complementary-slack-recourse} are the complementary slackness constraints. The latter include bilinear terms, as basis functions $\twostageLDRBFt$ are functions of $\xiT$ which are decision variables in \eqref{eq:SP-mono}. They can be linearized with the addition of binary decision variables via the so-called big-$M$ constraints. {\crev Since the basis functions are chosen to be affine, the resulting model is a mixed-integer linear program.} A detailed example of building a mixed-integer linear model for a multistage location-allocation problem is provided in Appendix \ref{sec:LT-models}.

{\crev
\subsection{Two-stage Decision Rules for MSAROs with Mixed-integer State Variables}
\label{app:sec:2AROprimal}
As mentioned in Remark \ref{rem:2AROprimal_MixedIntegerState}, for MSAROs with mixed-integer state variables, linear and piecewise constant decision rules can be combined to obtain a 2ARO approximation. To this end, consider the partition of the index set of the state variables into sets $\mathcal{I}_t^{\texttt{i}}$ and $ \mathcal{I}_t^{\texttt{c}}$ for integer and continuous variables, respectively, i.e.,  
$ \mathcal{I}_t^{\texttt{i}} \cup \mathcal{I}_t^{\texttt{c}} = [q_t], \mathcal{I}_t^{\texttt{i}} \cap \mathcal{I}_t^{\texttt{c}} = \emptyset,   |\mathcal{I}_t^{\texttt{i}}| = q_t^\texttt{i}$. Similarly, let 
$x^\texttt{s}_t = (x^\texttt{s,i}_t,x^\texttt{s,c}_t)$ and $\beta = (\beta^\texttt{i},\beta^\texttt{c})$ be the vectors of state variables and decision rule design variables with sub-vectors corresponding to the integer and continuous state variables. Then, the application of LDRs \eqref{eq:Theta_LDR} to the continuous state variables and PCDRs \eqref{eq:Theta_PCDR} to the integer state variables  yields the following 2ARO model:
\bsub
\begin{alignat}{2}
\nu^\text{2S-LDR-PCDR} := \min \ & c_1^\top x_1 + \mathcal{SP}^\text{2S-LDR-PCDR}(\beta,x_1^\texttt{s}) \\
\text{s.t.} \ & x_1\in X_1 \\
& \twostageLDRVart^{\texttt{c}} \in \R^{K^{\texttt{c}}_t} && t\in\set{2,T} \\
& \twostageLDRVart^{\texttt{i}} \in \R^{K^{\texttt{i}}_t} && t\in\set{2,T} \\
& {\crev -1 \leq \twostageKBFAffti^\top\twostageLDR_{ti}^{\texttt{i}} \leq 1} \qquad && {\crev t\in\set{2,T}, i\in \mathcal{I}_{t}^{\texttt{i}}, \xiT\in \xisupp} 
\end{alignat}
\esub
where 
\bsub
\label{eq:2RO-LDR-PCDR-SP}
\begin{align}
	\mathcal{SP}^\text{2S-LDR-PCDR}(\beta,x_1^\texttt{s}) := \max_{\xiT\in\xisupp} \quad & \sum_{t\in {\crev [2,T]} }{\crev{c^\texttt{s}_{t}(\xisupt)}^\top x^{\texttt{s}}_{t} }+ 
	\min_{\recvar\in\mathcal{X}({\crev \statevar}, \xiT)} && \hspace*{-0.5cm} \sum_{t\in \set{2,T}} {\crev {\recCtFixed}(\xisupt)^\top \recvar} 
 \\*[0.12cm]
	\text{s.t.}\quad & \ \twostageLDRBFt \twostageLDR^\texttt{c}_t = x^{\texttt{s,c}}_{t} && t\in\set{2,T}
 \\ 
 &  \sum_{j\in\set{J_{\crev i}}} (\underline{\kappa}_{ti} + j-1)\upsilon_{tij} = x^{\texttt{s,i}}_{ti} && t\in\set{2,T}, i\in\mathcal{I}_t^{\texttt{i}}
 \\
	& \sum_{j\in\set{J_{\crev i}}} \omega_{tij} = \twostageKBFAffti^\top\twostageLDR_{t}^i&&t\in\set{2,T}, i\in\mathcal{I}_t^{\texttt{i}}
 \\
	& ( {\crev a^{i}_{tj}}+\epsilon_{\crev j})\upsilon_{tij} \leq \omega_{tij} \leq {\crev b^{i}_{tj}}\upsilon_{tij}&&t\in\set{2,T}, i\in\mathcal{I}_t^{\texttt{i}}, j\in\set{J_{\crev i}}
 \\
	& \sum_{j\in\set{J_{\crev i}}} \upsilon_{tij} = 1 && t\in\set{2,T}, i\in\mathcal{I}_t^{\texttt{i}}
 \\
	& \upsilon_{tij} \in\{0,1\}&&t\in\set{2,T}, i\in\mathcal{I}_t^{\texttt{i}}, j\in\set{J_{\crev i}}.
\end{align}
\esub
Given the state variables and an uncertainty realization, the recourse feasible set is defined as 
\begin{align*}
	\mathcal{X}(\statevar,\xi^T) := \Big\{&\big({\crev\recvar_t}\big)_{t\in\set{2,T}}\in\R^{p_2}\times\R^{p_3}\times \dots\times \R^{p_T}:\\
    &\ {\crev \recAt \recvar_t \leq b_{t}(\xisupt)-\Big(\stateAt \statevar_t + \stateBt \statevar_1 \Big)} && {\crev t = 2} \\
	&\ \recAt \recvar_t \leq b_{t}(\xisupt)-\Big(\stateAt \statevar_t + \stateBt \statevar_{t-1} \Big) && t\in \set{{\crev 3},T}   \\
	&\ (\statevar_t, \recvar_t)\in X_{t}(\xisupt)&&   t\in \set{2,T}\ \Big\}.
\end{align*}

Lastly, we note that in the special case where the assumptions of 
Remarks \ref{rem:monoloithicSP} and  
\ref{rem:monoloithicSPinPCDR} are satisfied, the subproblem \eqref{eq:2RO-LDR-PCDR-SP} can be similarly reformulated as a monolithic mixed-integer linear program.
}

\subsection{Models for the Location-Transportation}\label{sec:LT-models}
\subsubsection{Column-and-constraint Generation with Two-stage Linear Decision Rules}\label{sec:CCG-LT-detail}

Applying LDRs on the state variables of model \eqref{eq:location-transportation},
$\displaystyle s_{it}(\xisupt) = \beta^0_{it} + \sum_{t'\in \set{2,t}}\sum_{j\in \set{J}}d_{jt'}(\xisupt)\beta_{it}^{jt'}, i\in\set{I}, t\in\set{2,T} $, results in the following 2ARO problem in monolithic form:
\begin{align*}
	\min  \ \ &z \\
	\text{s.t.} \ \ & z\geq \sum_{i\in \set{I}}f_iy_{i}  + \sum_{i\in \set{I}}a_is_{i1} + \sum_{t\in \set{2,T}}\sum_{i\in \set{I}}a_i\big(\beta^0_{it} + \sum_{t'\in \set{2,t}}\sum_{j\in \set{J}}d_{jt'}(\xiprimesupt)\beta_{it}^{jt'}\big) + \nonumber\\
	& \quad\sum_{t\in \set{2,T}}\sum_{i\in \set{I}}\sum_{j\in \set{J}}c_{ij}x_{ijt}(\xisupt) & \hspace{-3cm} \xiT\in\xisupp  \\
	& s_{i1} \leq K_iy_i &  \hspace{-3cm}i\in\set{I} \\
	& \beta^0_{i2} + \sum_{j\in \set{J}}d_{j2}(\xiparrand^2)\beta_{i2}^{j2} = s_{i1}-\sum_{j\in \set{J}}x_{ij2}(\xiparrand^2)& \hspace{-3cm}i\in\set{I}, \xiparrand^2\in\Xi^2\\
	& \beta^0_{it} + \sum_{t'\in \set{2,t}}\sum_{j\in \set{J}}d_{jt'}(\xisupt)\beta_{it}^{jt'} = \beta^0_{i,t-1} + \sum_{t'\in \set{2,t-1}}\sum_{j\in \set{J}}d_{jt'}(\xisuptminusone)\beta_{i,t-1}^{jt'}-\sum_{j\in \set{J}}x_{ijt}(\xisupt)&  \nonumber\\
	&& \hspace{-3cm}i\in\set{I}, t\in\set{3,T},  \xisupt\in\xisuppt\\
	& \sum_{i\in \set{I}}x_{ijt}(\xisupt) \geq d_{jt}(\xisupt)& \hspace{-3cm} j\in\set{J}, t\in\set{2,T}, \xisupt\in\xisuppt\\
	& \beta^0_{it} + \sum_{t'\in \set{2,t}}\sum_{j\in \set{J}}d_{jt'}(\xisupt)\beta_{it}^{jt'}\geq 0&\hspace{-3cm} i\in\set{I}, t\in\set{2,T}, \xisupt\in\xisuppt\\
	&y\in\{0,1\}^I,\ z,x\geq 0.
\end{align*}
For better clarity, let us re-write the above formulation in the more common min-max-min form of 2ARO:\\

\noindent
\scalebox{0.88}{
\begin{tikzpicture}[x=0.75pt,y=0.75pt,yscale=-1,xscale=1]
% Text Node
\draw (-35,1.7) node [anchor=north west][inner sep=0.75pt][align=left] {$\begin{aligned}
	\min \ & \sum_{i\in \set{I}}f_iy_{i} + \sum_{i\in \set{I}}a_i\big(s_{i1} + \sum_{t\in \set{2,T}}\beta^0_{it}\big)\qquad + 
	\qquad\max_{\xiT\in\xisupp} \\
	\text{s.t.}\ & s_{i1} \leq K_iy_i \quad i\in\set{I} &  \\
 & {\crev\beta^0_{it} + \sum_{t'\in \set{2,t}}\sum_{j\in \set{J}}d_{jt'}(\xisupt)\beta_{it}^{jt'}\geq 0 \quad i\in\set{I}, t\in\set{2,T}, \xisupt\in\xisuppt}\\
	& y\in\{0,1\}^I. 
\end{aligned}$};
% Text Node
\draw (340,0) node [anchor=north west][inner sep=0.75pt][align=left]  {$\begin{aligned}
\min\ & \sum_{t'\in \set{2,t}}a_id_{jt'}(\xisupt)\beta_{it}^{jt'}+\sum_{t\in \set{2,T}}\sum_{i\in \set{I}}\sum_{j\in \set{J}} c_{ij}x_{ijt} \\
\text{s.t.} \ 
& \sum_{j\in \set{J}}x_{ij2} = s_{i1} -\beta^0_{i2}-\sum_{j\in \set{J}}d_{j2}(\xiparrand^2)\beta_{i2}^{j2}  &\hspace{-2cm}i\in\set{I}\\
& \sum_{j\in \set{J}}x_{ijt} = \beta^0_{i,t-1}-\beta^0_{it}-\sum_{j\in \set{J}} \sum_{t'\in \set{2,t}}d_{jt'}(\xisupt)\beta_{it}^{jt'} \\
&\qquad\qquad+\sum_{j\in \set{J}}\sum_{t'\in \set{2,t-1}}d_{jt'}(\xisuptminusone)\beta_{i,t-1}^{jt'}  &  \nonumber\\
&&\hspace{-2cm}i\in\set{I}, t\in\set{3,T}\\
& \sum_{i\in \set{I}}x_{ijt} \geq d_{jt}(\xisupt)&  \hspace{-2cm}j\in\set{J}, t\in\set{2,T}\\
& x\geq 0
\end{aligned}$};
\end{tikzpicture}
}

{\crev \noindent Note that we have strengthened the outer minimization problem by adding the non-negativity constraints from the inner minimization problem as robust constraints. This constraint can be rewritten as follows:
$$\begin{aligned}
& \beta^0_{it} + \sum_{t'\in \set{2,t}}\sum_{j\in \set{J}}d_{jt'}(\xisupt)\beta_{it}^{jt'}\geq 0&  i\in\set{I}, t\in\set{2,T}, \xisupt\in\xisuppt\rightarrow\\
 & \beta^0_{it} + \sum_{t'\in \set{2,t}}\sum_{j\in \set{J}} \mu_{jt'}\beta_{it}^{jt'} + \min \Big\{\sum_{t'\in \set{2,t}}\sum_{j\in \set{J}}\xiparrand_{jt'}\sigma_{jt'}\beta_{it}^{jt'}:\\
 &\hspace{49mm} \sum_{t'\in\set{2,T}}\sum_{j\in\set{J}}\xiparrand_{jt'}\leq \Gamma,\\
 &\hspace{49mm} 0\leq\xiparrand_{jt'}\leq 1, j\in\set{J}, t'\in\set{2,T}\Big\}\geq 0&  i\in\set{I}, t\in\set{2,T} \rightarrow\\
 & \beta^0_{it} + \sum_{t'\in \set{2,t}}\sum_{j\in \set{J}} \mu_{jt'}\beta_{it}^{jt'} + \max \Big\{u\Gamma + \sum_{j\in \set{J}} \sum_{t'\in \set{2,T}}\omega_{jt'}:\\
 &\hspace{52mm} u + \omega_{jt'} \leq \sigma_{jt'}\beta_{it}^{jt'},  j\in\set{J}, t'\in\set{2,t},\\
 & \hspace{52mm} u\leq 0, \omega_{jt'} \leq 0, j\in\set{J}, t'\in\set{2,T}
  \Big\}\geq 0&  i\in\set{I}, t\in\set{2,T} \rightarrow
\end{aligned}$$
$$\left\lbrace \begin{array}{lr}\displaystyle
      \beta^0_{it} + \sum_{t'\in \set{2,t}}\sum_{j\in \set{J}} \mu_{jt'}\beta_{it}^{jt'} +u\Gamma + \sum_{j\in \set{J}} \sum_{t'\in \set{2,T}}\omega_{jt'}\geq 0&  i\in\set{I}, t\in\set{2,T}  \\
      u + \omega_{jt'} \leq \sigma_{jt'}\beta_{it}^{jt'} & i\in\set{I}, j\in\set{J}, t\in\set{2,T}, t'\in\set{2,t}\\ 
      \omega_{jt'} \leq 0& i\in\set{I}, j\in\set{J}, t\in\set{2,T}, t'\in\set{2,T} \\ 
      u\leq 0\\[2mm]
 \end{array}\right.$$

 \noindent In order to establish relatively complete recourse, we add the following constraints to the outer minimization problem:

 $$ 
 \begin{aligned}
    & \sum_{i\in\set{I}}\bigg[s_{i1} -\beta^0_{i2}-\sum_{j\in \set{J}}d_{j2}(\xiparrand^2)\beta_{i2}^{j2}\bigg]
    \geq
    \sum_{j\in\set{J}}d_{j2}(\xiparrand^2)&   \xi^2\in\Xi^2\\
& \sum_{i\in\set{I}}\bigg[\beta^0_{i,t-1}-\beta^0_{it}-\sum_{j\in \set{J}} \sum_{t'\in \set{2,t}}d_{jt'}(\xisupt)\beta_{it}^{jt'}+\sum_{j\in \set{J}}\sum_{t'\in \set{2,t-1}}d_{jt'}(\xisuptminusone)\beta_{i,t-1}^{jt'}\bigg]
\geq 
\sum_{j\in\set{J}}d_{jt}(\xisupt)&  t\in\set{3,T}, \xisupt\in\xisuppt,
 \end{aligned}
 $$
which can be reformulated as follows: \\ 
 $$ 
 \begin{aligned}
    & \sum_{i\in\set{I}}\Big[s_{i1} -\beta^0_{i2}\Big]
    \geq
    \max_{\xi^2\in\Xi^2}\bigg\{\sum_{j\in\set{J}}d_{j2}(\xiparrand^2)\Big[1+\sum_{i\in\set{I}}\beta_{i2}^{j2}\Big]\bigg\}\\
& \sum_{i\in\set{I}}\Big[\beta^0_{i,t-1}-\beta^0_{it}\Big]
\geq 
\max_{\xisupt\in\xisuppt}\bigg\{
\sum_{i\in\set{I}}\bigg[\sum_{j\in \set{J}} \sum_{t'\in \set{2,t}}d_{jt'}(\xisupt)\beta_{it}^{jt'}-\sum_{j\in \set{J}}\sum_{t'\in \set{2,t-1}}d_{jt'}(\xisuptminusone)\beta_{i,t-1}^{jt'}\bigg]+\sum_{j\in\set{J}}d_{jt}(\xisupt)\big\}&  t\in\set{3,T}
 \end{aligned} 
 $$
The linearization process is similar to the one previously discussed for non-negativity constraints.
 }

\noindent Denote by $\pi^1, \pi^2, \pi^3$ the dual variables associated with three constraint sets of the inner minimization problem. After taking its {\crev linear programming} dual and the adding the KKT conditions, we get the following as the subproblem of the column-and-constraint generation:
\begin{align*}
	\max  \ \ & \sum_{t\in \set{2,T}}\sum_{i\in \set{I}}\sum_{j\in \set{J}}  \sum_{t'\in \set{2,t}}a_id_{jt'}\beta_{it}^{jt'} + \sum_{t\in \set{2,T}}\sum_{i\in \set{I}}\sum_{j\in \set{J}}c_{ij}x_{ijt}\\
	\text{s.t.} \ \ 
	% Primal
	& \sum_{j\in \set{J}}x_{ij2} = s_{i1} -\beta^0_{i2}-\sum_{j\in \set{J}}d_{j2}\beta_{i2}^{j2}  &\hspace{-1cm}i\in\set{I}\\
    & \sum_{j\in \set{J}}x_{ijt} = \beta^0_{i,t-1}-\beta^0_{it}-\sum_{j\in \set{J}} \sum_{t'\in \set{2,t}}d_{jt'}\beta_{it}^{jt'} +\sum_{j\in \set{J}}\sum_{t'\in \set{2,t-1}}d_{jt'}\beta_{i,t-1}^{jt'}  &\hspace{-1cm}  i\in\set{I}, t\in\set{3,T}\\
    & \sum_{i\in \set{I}}x_{ijt} \geq d_{jt}&\hspace{-1cm}  j\in\set{J}, t\in\set{2,T}\\[3mm]
    % Dual
	& \pi^1_i + \pi^3_{j2} \leq c_{ij} &\hspace{-1cm} i\in\set{I}, j\in\set{J}	\\
	& \pi^2_{it} + \pi^3_{jt} \leq c_{ij} &\hspace{-1cm} i\in\set{I}, j\in\set{J}, t\in\set{3,T}\\[3mm]
	% Complementary slackness
	& \Big(d_{jt}-\sum_{i\in \set{I}}x_{ijt}\Big)\pi^3_{jt} = 0 \xrightarrow{\text{linearization}}\\
	&\qquad\qquad\pi^3_{jt}\leq M(1-\ell^R_{jt}),\ \  \sum_{i\in \set{I}}x_{ijt}-d_{jt} \leq M\ell^R_{jt} &\hspace{-1cm}j\in\set{J}, t\in\set{2,T} \\
	& \big(\pi^1_i + \pi^3_{j2}-c_{ij}\big) x_{ij2} = 0\xrightarrow{\text{linearization}}\\
	&\qquad\qquad x_{ij2}\leq M(1-\ell^D_{ij2}),\ \ c_{ij}-\pi^1_i-\pi^3_{j2} \leq M\ell^D_{ij2}
	&\hspace{-1cm}i\in\set{I}, j\in\set{J}\\
	& \big(\pi^2_{it} + \pi^3_{jt}-c_{ij}\big)x_{ijt} = 0 \xrightarrow{\text{linearization}}\\
	&\qquad\qquad x_{ijt}\leq M(1-\ell^D_{ijt}),\ \ c_{ij}-\pi^2_{it}-\pi^3_{jt} \leq M\ell^D_{ijt}
	&\hspace{-1cm} i\in\set{I}, j\in\set{J}, t\in\set{3,T} \\[3mm]
	% Uncertainty set
	& d_{jt} = \mu_{jt} + \xiparrand_{jt}\sigma_{jt}&\hspace{-1cm} j\in\set{J}, t\in\set{2,T}\\
	& \sum_{t\in\set{2,T}}\sum_{j\in\set{J}}\xiparrand_{jt}\leq \Gamma\\
	& x_{ijt}\geq 0,\ \ell^D_{ijt}\in\{0,1\}&\hspace{-1cm}i\in\set{I}, j\in\set{J}, t\in\set{2,T}\\
	& 0\leq\xiparrand_{jt}\leq 1,\ 0\leq\pi^3_{jt},\ \ell^R_{jt}\in\{0,1\} &\hspace{-1cm}j\in\set{J}, t\in\set{2,T}.
\end{align*}

\subsubsection{Restricted NA Dual}
The $\mathcal{Q}(\beta^s,\beta^x)$ function for the location-transportation problem: 
\begin{align*}
\mathcal{Q}(\beta^s,\beta^x) = \\
\min  \ \ &z+ \sum_{t\in \set{2,T}}\sum_{i\in \set{I}}\beta^s_{it}(\xiT)s_{it}(\xiT) + \sum_{t\in \set{2,T}}\sum_{i\in \set{I}}\sum_{j\in \set{J}}\beta^x_{ijt}(\xiT)x_{ijt}(\xiT)  \\
\text{s.t.} \ \ & z\geq \sum_{i\in \set{I}}y_{i} + \sum_{t\in \set{T}}\sum_{i\in \set{I}}a_is_{it}(\xiT) + \sum_{t\in \set{2,T}}\sum_{i\in \set{I}}\sum_{j\in \set{J}}c_{ij}x_{ijt}(\xiT) & \hspace{-3cm} \xiT\in\xisupp  \\
& s_{i1} \leq K_iy_i, &  i\in\set{I}\\
& s_{it}(\xiT) = s_{i,t-1}(\xiTminusone)-\sum_{j\in \set{J}}x_{ijt}(\xiT)&\hspace{-3cm}  i\in\set{I}, t\in\set{2,T},  \xiT\in\xisupp\\
& \sum_{i\in \set{I}}x_{ijt}(\xiT) \geq d_{jt}(\xisupt)&\hspace{-3cm}  i\in\set{I}, t\in\set{2,T}, \xiT\in\xisupp\\
&y\in\{0,1\}^I, z\in\R_+^{I\times T},x\in\R_+^{I\times J\times (T-1)}.
\end{align*}
For this mixed-integer subproblem, the cutting-plane method presented in Section \ref{sec:solution} is applicable.

\newpage
\section{Detailed Results}\label{app:results}
\subsection{Newsvendor Problem}
\label{app:results-newsvendor}
\begin{table}[htbp]
    \caption{Detailed results for the newsvendor problem}
    \label{tab:detailed-newsvendor}
    \crev
    \centering
    \scalebox{0.88}{
    \begin{tabular}{p{5mm}rrccrrrrr}
       \toprule
$T$   & \textsc{br} & $|\xisupp|$ & $\ \ I\ $   & $B$   & $\LDRUB$ & $\twostageLDRUB$ & $\exactB$ & $\bar{\nu}^{\textsc{NA-DO}}_R$ & $\PILB$ \\
\midrule
         \multirow{6}{*}{3} & \multirow{6}{*}{5} & \multirow{6}{*}{25} & 3     & 100   & 656.6 & 970.2 & 975.3 & 975.3 & 1457.0 \\
          &       &       & 3     & 150   & 6668.5 & 8410.4 & 8647.2 & 9802.6 & 10036.0 \\
          &       &       & 3     & 200   & 8662.0 & 11446.2 & 11569.7 & 13709.0 & 13709.0 \\
          &       &       & 4     & 150   & 5154.8 & 5384.4 & 5466.8 & 5466.8 & 5781.0 \\
          &       &       & 4     & 200   & 10049.4 & 12452.0 & 12759.9 & 14501.9 & 14525.0 \\
          &       &       & 5     & 200   & 5944.7 & 7120.6 & 7402.4 & 7402.4 & 7713.0 \\
\midrule\multirow{5}{*}{3} & \multirow{5}{*}{10} & \multirow{5}{*}{100} & 3     & 150   & 5752.1 & 7620.9 & 7984.6 & 8590.3 & 8895.0 \\
          &       &       & 3     & 200   & 6279.6 & 9199.0 & 9700.9 & 11833.0 & 11833.0 \\
          &       &       & 4     & 150   & 3502.7 & 4886.5 & 4910.2 & 5110.1 & 5323.0 \\
          &       &       & 4     & 200   & 7922.9 & 10961.1 & 11273.4 & 14036.5 & 14060.0 \\
          &       &       & 5     & 200   & 2865.5 & 5877.7 & 6063.1 & 6063.1 & 6737.0 \\
\midrule
\multirow{5}{*}{4} & \multirow{5}{*}{4} & \multirow{5}{*}{64} & 2     & 100   & 669   & 1047  & 1214  & 1214  & 1337 \\
          &       &       & 3     & 200   & 7807  & 9198  & 9491  & 9492  & 10157 \\
          &       &       & 3     & 300   & 13421 & 16414 & 17317 & 19105 & 19188 \\
          &       &       & 4     & 200   & 2034  & 2297  & 2447  & 2447  & 2640 \\
          &       &       & 4     & 300   & 13645 & 17860 & 18721 & 20139 & 20139 \\
\midrule\multirow{5}{*}{4} & \multirow{5}{*}{5} & \multirow{5}{*}{125} & 3     & 200   & 6643.8 & 8070.1 & 8295.6 & 8785.8 & 9249.0 \\
          &       &       & 3     & 300   & 11854.3 & 15035.4 & 16606.8 & 18240.0 & 18240.0 \\
          &       &       & 4     & 200   & 1443.7 & 1700.8 & 1739.8 & 1917.7 & 2122.0 \\
          &       &       & 4     & 300   & 13186.4 & 17579.4 & 18450.8 & 19618.0 & 19635.0 \\
          &       &       & 5     & 300   & 5142.0 & 10494.4 & 11299.5 & 11804.9 & 12343.0 \\
\midrule\multirow{5}{*}{4} & \multirow{5}{*}{10} & \multirow{5}{*}{1000} & 3     & 200   & 5648.8 & 8142.9 &       & 9353.0 & 9353.0 \\
          &       &       & 3     & 300   & 9143.5 & 13687.6 &       & 17853.0 & 17853.0 \\
          &       &       & 4     & 200   & -104.4 & 440.0 &       & 642.4 & 919.0 \\
          &       &       & 4     & 300   & 11029.5 & 15854.7 &       & 18508.0 & 18432.0 \\
          &       &       & 5     & 300   & 724.9 & 6125.5 &       & 6740.1 & 7368.0 \\
\midrule\multirow{2}{*}{4} & \multirow{2}{*}{15} & \multirow{2}{*}{3375} & 3     & 200   & 5111.3 & 7072.0 &       & 7646.6 & 8222.0 \\
          &       &       & 3     & 300   & 9030.9 & 13040.7 &       & 17615.0 & 17615.0 \\
          \midrule
\multirow{4}{*}{5} & \multirow{4}{*}{3} & \multirow{4}{*}{81} & 2     & 150   & 2490.9 & 2673.5 & 2763.4 & 2763.3 & 2975.4 \\
          &       &       & 2     & 200   & 8888.4 & 11192.5 & 11313.7 & 11525.6 & 11525.6 \\
          &       &       & 3     & 250   & 7809.0 & 8800.3 & 9065.4 & 9104.8 & 9382.3 \\
          &       &       & 3     & 300   & 14702.3 & 16898.8 & 17214.3 & 17308.7 & 18156.1 \\
\midrule5     & 4     & 256   & 3     & 300   & 12367.7 & 15591.0 & 15989.0 & 16511.1 & 17192.6 \\
\midrule\multirow{2}{*}{5} & \multirow{2}{*}{5} & \multirow{2}{*}{625} & 3     & 300   & 11134.5 & 15053.4 &       & 16191.6 & 16680.0 \\
          &       &       & 4     & 300   & 3651.5 & 9494.8 &       & 10628.0 & 10628.0 \\
\midrule
5     & 6     & 1296  & 3     & 400   & 15013.4 & 20271.7 &       & 25157.0 & 25157.0 \\
\midrule
\multirow{3}{*}{6} & \multirow{3}{*}{4} & \multirow{3}{*}{1024} & 3     & 400   & 15492.3 & 21457.9 &       & 28163.0 & 28163.0 \\
          &       &       & 4     & 400   & 7124.3 & 14805.2 &       & 15400.6 & 15473.0 \\
          &       &       & 4     & 500   & 14445.4 & 24887.7 &       & 32973.0 & 32973.0 \\
\midrule
\multirow{3}{*}{7} & \multirow{3}{*}{3} & \multirow{3}{*}{729} & 3     & 300   & 24.3  & 1495.9 &       & 2137.3 & 2383.0 \\
          &       &       & 3     & 400   & 12994.7 & 17774.1 &       & 19554.5 & 19983.0 \\
          &       &       & 4     & 400   & 3267.5 & 4303.1 &       & 5388.1 & 5427.0 \\
\midrule
\multirow{2}{*}{8} & \multirow{2}{*}{3} & \multirow{2}{*}{2187} & 3     & 500   & 16149.1 & 25321.4 &       & 30259.0 & 30259.0 \\
          &       &       & 4     & 600   & 16892.3 & 27608.1 &       & 28717.7 & 28747.0 \\
          \bottomrule
    \end{tabular}
    }
\end{table}

\begin{table}[htbp]
\crev
\setlength{\tabcolsep}{3mm}
  \centering
  \caption{Running times for the larger-size newsvendor problem instances\\}
    \begin{tabular}{crrccrrrc}
    \toprule
    \multirow{2}{*}{$\ T\ $} & \multicolumn{1}{c}{\multirow{2}[0]{*}{\textsc{br}}} & \multicolumn{1}{c}{\multirow{2}[0]{*}{$|\xisupp|$}} & \multirow{2}[1]{*}{$\ I\ $} & \multirow{2}[1]{*}{$B$} & \multicolumn{4}{c}{Time (s)} \\
\cmidrule{6-9}          &       &       &       &       & $\LDRUB$ & $\twostageLDRUB$ & {\crev $\bar{\nu}^{\textsc{NA-DO}}_R$} & $\text{PI}$ \\
    \midrule
    \multirow{7}[6]{*}{4} & \multirow{5}[2]{*}{10} & \multirow{5}[2]{*}{1000} & 3     & 200   & 38    & 31    & 89    & 1 \\
      &       &       & 3     & 300   & 50    & 31    & 20    & 1 \\
      &       &       & 4     & 200   & 47    & 28    & 637   & 1 \\
      &       &       & 4     & 300   & 58    & 30    & 2712  & 1 \\
      &       &       & 5     & 300   & 59    & 33    & 3088  & 2 \\
\cmidrule{2-9}      & 15    & 3375  & 3     & 200   & 158   & 140   & 558   & 3 \\
\cmidrule{2-9}      & 20    & 8000  & 3     & 300   & 726   & 675   & 2313  & 7 \\
\midrule
\multirow{3}[4]{*}{5} & \multirow{2}[2]{*}{5} & \multirow{2}[2]{*}{625} & 3     & 300   & 16    & 15    & 1423  & 1 \\
      &       &       & 4     & 300   & 22    & 19    & $>1h$  & 1 \\
\cmidrule{2-9}      & 6     & 1296  & 3     & 400   & 54    & 54    & 64    & 1 \\
\midrule
\multirow{3}[2]{*}{6} & \multirow{3}[2]{*}{4} & \multirow{3}[2]{*}{1024} & 3     & 400   & 57    & 53    & 99   & 1 \\
      &       &       & 4     & 400   & 70    & 67    & 12    & 2 \\
      &       &       & 4     & 500   & 74    & 69    & 3405  & 2 \\
\midrule
\multirow{3}[2]{*}{7} & \multirow{3}[2]{*}{3} & \multirow{3}[2]{*}{729} & 3     & 300   & 37    & 38    & 280   & 1 \\
      &       &       & 3     & 400   & 37    & 36    & $>1h$ & 1 \\
      &       &       & 4     & 400   & 50    & 50    & $>1h$  & 1 \\
\midrule
\multirow{2}[2]{*}{8} & \multirow{2}[2]{*}{3} & \multirow{2}[2]{*}{2187} & 3     & 500   & 364   & 335   & 408   & 4 \\
      &       &       & 4     & 600   & 462   & 464   & $>1h$ & 5 \\
\bottomrule
    \end{tabular}%
  \label{tab:newsvendor-time}%
\end{table}%

\newpage
\subsection{Location-Transportation Problem}\label{app:LT-results}
\begin{table}[htbp]
\crev
  \centering
  \caption{Detailed results for the location-transportation problem over a budgeted uncertainty set}
  \label{tab:detailed-LT-Bounds}
  \setlength{\tabcolsep}{3mm}
\scalebox{0.9}{
\begin{tabular}{crrrp{10pt}rrrr}
\cmidrule{1-4}\cmidrule{6-9}    \multicolumn{1}{c}{$(T,I,J,\alpha^d)$} & \multicolumn{1}{c}{$\alpha^u$} & \multicolumn{1}{c}{$\twostageLDRUB$} & \multicolumn{1}{c}{$\bar{\nu}^{\textsc{NA-DO}}_R$} &       & \multicolumn{1}{c}{$(T,I,J,\alpha^d)$} & \multicolumn{1}{c}{$\alpha^u$} & \multicolumn{1}{c}{$\twostageLDRUB$} & \multicolumn{1}{c}{$\bar{\nu}^{\textsc{NA-DO}}_R$} \\
\cmidrule{1-4}\cmidrule{6-9}    \multirow{4}[2]{*}{$(3,10,10,0.1)$} & 0.1   & 738969.9 & 684460.5 &       & \multicolumn{1}{c}{\multirow{4}[2]{*}{$(4,5,7,0.5)$}} & 0.1   & 1052392.2 & 826648.8 \\
          & 0.4   & 738967.8 & 647069.1 &       &       & 0.4   & 1052393.0 & 767446.4 \\
          & 0.7   & 738967.0 & 573052.3 &       &       & 0.7   & 1052392.4 & 863507.3 \\
          & 1     & 738969.3 & 574580.2 &       &       & 1     & 1052392.9 & 806706.2 \\
\cmidrule{1-4}\cmidrule{6-9}    \multirow{4}[2]{*}{$(3,10,10,0.3)$} & 0.1   & 1657362.2 & 1457789.5 &       & \multicolumn{1}{c}{\multirow{4}[2]{*}{$(4,5,10,0.1)$}} & 0.1   & 2414351.2 & 2008281.8 \\
          & 0.4   & 1637485.0 & 1364705.6 &       &       & 0.4   & 2414351.2 & 1736011.9 \\
          & 0.7   & 1535360.3 & 1205974.5 &       &       & 0.7   & 2414351.2 & 2168851.3 \\
          & 1     & 1554797.6 & 1223643.1 &       &       & 1     & 2414351.2 & 1954758.7 \\
\cmidrule{1-4}\cmidrule{6-9}    \multirow{4}[2]{*}{$(3,10,10,0.5)$} & 0.1   & 1641883.1 & 1430838.7 &       & \multicolumn{1}{c}{\multirow{4}[2]{*}{$(4,5,10,0.3)$}} & 0.1   & 2414351.2 & 1912871.0 \\
          & 0.4   & 1605346.2 & 1363127.9 &       &       & 0.4   & 2414351.2 & 1675507.4 \\
          & 0.7   & 1646129.3 & 1178537.7 &       &       & 0.7   & 2414351.2 & 1747647.4 \\
          & 1     & 1592763.4 & 1474506.6 &       &       & 1     & 2414351.5 & 2071113.7 \\
\cmidrule{1-4}\cmidrule{6-9}    \multirow{4}[2]{*}{$(3,10,15,0.1)$} & 0.1   & 1133251.5 & 857198.3 &       & \multicolumn{1}{c}{\multirow{4}[2]{*}{$(4,5,10,0.5)$}} & 0.1   & 2414351.1 & 1745320.2 \\
          & 0.4   & 1133251.5 & 825580.3 &       &       & 0.4   & 2342846.9 & 2076574.6 \\
          & 0.7   & 1157049.6 & 1073536.8 &       &       & 0.7   & 2386394.3 & 1959412.9 \\
          & 1     & 1156385.8 & 956570.9 &       &       & 1     & 2386394.3 & 1824441.4 \\
\cmidrule{1-4}\cmidrule{6-9}    \multirow{4}[2]{*}{$(3,10,10,0.3)$} & 0.1   & 1157948.4 & 1062129.6 &       & \multicolumn{1}{c}{\multirow{4}[2]{*}{$(4,10,10,0.1)$}} & 0.1   & 2502997.3 & 2064097.4 \\
          & 0.4   & 1157948.4 & 861882.4 &       &       & 0.4   & 2508960.6 & 1842514.1 \\
          & 0.7   & 1157948.4 & 1092108.1 &       &       & 0.7   & 2508960.5 & 1688369.1 \\
          & 1     & 1157948.4 & 1000064.7 &       &       & 1     & 2508960.8 & 2237510.8 \\
\cmidrule{1-4}\cmidrule{6-9}    \multirow{4}[2]{*}{$(3,10,10,0.5)$} & 0.1   & 1157948.4 & 831435.5 &       & \multicolumn{1}{c}{\multirow{4}[2]{*}{$(4,10,10,0.3)$}} & 0.1   & 2344765.1 & 2189919.0 \\
          & 0.4   & 1157948.4 & 981328.1 &       &       & 0.4   & 2508960.6 & 2095408.1 \\
          & 0.7   & 1156079.7 & 946148.4 &       &       & 0.7   & 2508961.0 & 2071323.3 \\
          & 1     & 1116748.8 & 885303.1 &       &       & 1     & 2401193.3 & 1728094.5 \\
\cmidrule{1-4}\cmidrule{6-9}    \multirow{4}[2]{*}{$(4,5,5,0.1)$} & 0.1   & 770943.7 & 724852.8 &       & \multicolumn{1}{c}{\multirow{4}[2]{*}{$(4,10,10,0.5)$}} & 0.1   & 2432431.8 & 2130923.4 \\
          & 0.4   & 770943.7 & 573157.6 &       &       & 0.4   & 2474962.8 & 2203068.7 \\
          & 0.7   & 770943.7 & 687320.0 &       &       & 0.7   & 2347843.6 & 1921134.8 \\
          & 1     & 770943.7 & 727980.6 &       &       & 1     & 2347843.6 & 1743911.8 \\
\cmidrule{1-4}\cmidrule{6-9}    \multirow{4}[2]{*}{$(4,5,5,0.3)$} & 0.1   & 770943.7 & 710778.0 &       & \multicolumn{1}{c}{\multirow{4}[2]{*}{$(5,5,10,0.1)$}} & 0.1   & 3562840.8 & 3075701.1 \\
          & 0.4   & 770943.7 & 631020.7 &       &       & 0.4   & 3595989.4 & 2855492.7 \\
          & 0.7   & 770943.7 & 664643.9 &       &       & 0.7   & 3563956.9 & 2788612.9 \\
          & 1     & 770943.7 & 696703.2 &       &       & 1     & 3563956.4 & 3256824.0 \\
\cmidrule{1-4}\cmidrule{6-9}    \multirow{4}[2]{*}{$(4,5,5,0.5)$} & 0.1   & 770943.7 & 562210.5 &       & \multicolumn{1}{c}{\multirow{4}[2]{*}{$(5,5,10,0.3)$}} & 0.1   & 3623768.3 & 2694805.6 \\
          & 0.4   & 770943.7 & 662298.1 &       &       & 0.4   & 3501499.4 & 2826482.5 \\
          & 0.7   & 770943.7 & 620073.7 &       &       & 0.7   & 3649157.4 & 2483339.2 \\
          & 1     & 770943.7 & 729544.4 &       &       & 1     & 3470084.5 & 2450441.3 \\
\cmidrule{1-4}\cmidrule{6-9}    \multirow{3}[2]{*}{$(4,5,7,0.1)$} & 0.1   & 1052392.2 & 908927.4 &       & \multicolumn{1}{c}{\multirow{4}[4]{*}{$(5,5,10,0.5)$}} & 0.1   & 3630095.7 & 2474291.5 \\
          & 0.4   & 1052392.2 & 847845.0 &       &       & 0.4   & 3573873.7 & 2622461.9 \\
          & 1     & 1052392.2 & 794593.6 &       &       & 0.7   & 3610736.5 & 2973602.4 \\
\cmidrule{1-4}    \multirow{4}[4]{*}{$(4,5,7,0.3)$} & 0.1   & 1052392.2 & 848889.1 &       &       & 1     & 3489038.0 & 2543037.8 \\
\cmidrule{6-9}          & 0.4   & 1052392.2 & 824873.8 &       &       &       &       &  \\
          & 0.7   & 1052392.2 & 984627.8 &       &       &       &       &  \\
          & 1     & 1052392.7 & 855154.4 &       &       &       &       &  \\
\cmidrule{1-4}    
\end{tabular}%
}
\end{table}%

\begin{table}[htbp]
\crev
\small
  \centering
  \caption{Algorithmic details for the location-transportation problem over a budgeted uncertainty set}
\label{tab:detailed-LT-algorithms}
\scalebox{0.8}{
    \begin{tabular}{crrrrrrrrrr}
\toprule    \multirow{2}[2]{*}{$(T,I,J,\alpha^d)$} & \multicolumn{1}{c}{\multirow{2}[2]{*}{$\alpha^u$}} & \multicolumn{1}{c}{$\LDRUB$} & \multicolumn{3}{c}{$\twostageLDRUB$, C\&CG} & \multicolumn{1}{c}{$\nu^{\Omega(\text{LDR})}$} & \multicolumn{1}{c}{$\nu^{\Omega(\text{2S-LDR})}$} & \multicolumn{3}{c}{$\bar{\nu}^{\textsc{NA-DO}}_R$} \\
\cmidrule(rl){3-3}  \cmidrule(rl){4-6} \cmidrule(rl){7-7} \cmidrule(rl){8-8} \cmidrule(rl){9-11}        &       & \multicolumn{1}{c}{Time (s)} & \multicolumn{1}{c}{Time (s)} & \multicolumn{1}{c}{\#iterations} & \multicolumn{1}{c}{Gap} & \multicolumn{1}{c}{Time (s)} & \multicolumn{1}{c}{Time (s)} & Time (s) & \multicolumn{1}{c}{\#iterations} & \multicolumn{1}{c}{Gap} \\
    \midrule
    \multirow{4}[2]{*}{$(3,10,10,0.1)$} & 0.1   & 0.5   & 168.3 & 2     & 0.0\% & 2.4   & 1.4   & 1376.9 & 617   & 4.9\% \\
          & 0.4   & 6.5   & 59.6  & 2     & 0.0\% & 9.2   & 11.5  & 1113.5 & 809   & 2.6\% \\
          & 0.7   & 10.8  & 87.0  & 2     & 0.0\% & 9.3   & 9.8   & 1088.4 & 726   & 4.3\% \\
          & 1     & 1.3   & 141.8 & 2     & 0.0\% & 16.7  & 10.9  & $>10h$ & 1315  & 7.9\% \\
    \midrule
    \multirow{4}[2]{*}{$(3,10,10,0.3)$} & 0.1   & 14.6  & 867.2 & 12    & 4.8\% & 2.2   & 0.1   & $>10h$ & 1013  & 8.3\% \\
          & 0.4   & 4.7   & 51.9  & 2     & 4.9\% & 8.1   & 16.5  & 1082.4 & 242   & 4.1\% \\
          & 0.7   & 12.6  & 66.1  & 2     & 1.2\% & 3.1   & 1.6   & $>10h$ & 730   & 7.1\% \\
          & 1     & 8.4   & 71.6  & 2     & 4.9\% & 2.6   & 0.9   & 588.0 & 984   & 3.0\% \\
    \midrule
    \multirow{4}[2]{*}{$(3,10,10,0.5)$} & 0.1   & 12.9  & 124.6 & 2     & 4.7\% & 7.9   & 10.5  & $>10h$ & 271   & 7.0\% \\
          & 0.4   & 16.3  & 111.8 & 3     & 4.9\% & 1.0   & 0.7   & $>10h$ & 665   & 9.1\% \\
          & 0.7   & 1.8   & 83.2  & 3     & 4.3\% & 1.7   & 0.7   & 1132.4 & 1064  & 3.2\% \\
          & 1     & 25.2  & 96.5  & 3     & 3.9\% & 2.1   & 2.2   & $>10h$ & 1405  & 4.1\% \\
    \midrule
    \multirow{4}[2]{*}{$(3,10,15,0.1)$} & 0.1   & 10.2  & 150.1 & 3     & 1.2\% & 2.1   & 1.3   & 2634.8 & 258   & 2.9\% \\
          & 0.4   & 3.5   & 112.6 & 2     & 1.2\% & 2.3   & 1.1   & $>10h$ & 337   & 6.9\% \\
          & 0.7   & 29.3  & 96.6  & 2     & 0.0\% & 1.7   & 2.0   & $>10h$ & 432   & 5.5\% \\
          & 1     & 8.3   & 29.3  & 2     & 0.0\% & 2.2   & 2.7   & $>10h$ & 824   & 9.9\% \\
    \midrule
    \multirow{4}[2]{*}{$(3,10,10,0.3)$} & 0.1   & 2.5   & 169.7 & 2     & 1.1\% & 3.0   & 2.4   & 4193.8 & 771   & 2.1\% \\
          & 0.4   & 13.4  & 95.6  & 2     & 1.2\% & 3.8   & 0.4   & 4706.4 & 974   & 4.9\% \\
          & 0.7   & 13.0  & 94.6  & 2     & 1.1\% & 2.3   & 1.2   & 5043.6 & 475   & 2.0\% \\
          & 1     & 10.2  & 119.9 & 2     & 1.1\% & 1.6   & 2.0   & 5015.7 & 1187  & 4.4\% \\
    \midrule
    \multirow{4}[2]{*}{$(3,10,10,0.5)$} & 0.1   & 5.7   & 97.1  & 6     & 1.3\% & 2.0   & 1.9   & 3757.0 & 1253  & 1.1\% \\
          & 0.4   & 2.8   & 155.3 & 2     & 1.2\% & 3.0   & 1.5   & 2902.3 & 803   & 4.5\% \\
          & 0.7   & 0.7   & 112.0 & 2     & 0.0\% & 8.5   & 1.6   & $>10h$ & 884   & 9.7\% \\
          & 1     & 2.3   & 89.9  & 2     & 1.2\% & 2.3   & 4.5   & 3898.8 & 1120  & 4.6\% \\
    \midrule
    \multirow{4}[2]{*}{$(4,5,5,0.1)$} & 0.1   & 0.2   & 0.8   & 2     & 0.0\% & 7.3   & 16.5  & 1198.6 & 951   & 4.8\% \\
          & 0.4   & 0.3   & 1.2   & 2     & 0.0\% & 7.1   & 21.8  & $>10h$ & 1483  & 7.2\% \\
          & 0.7   & 0.6   & 1.1   & 2     & 0.0\% & 2.5   & 5.4   & 1558.8 & 629   & 2.1\% \\
          & 1     & 0.5   & 1.2   & 2     & 0.0\% & 1.4   & 2.6   & $>10h$ & 1267  & 7.3\% \\
    \midrule
    \multirow{4}[2]{*}{$(4,5,5,0.3)$} & 0.1   & 0.7   & 1.3   & 2     & 0.0\% & 5.6   & 0.9   & 1649.9 & 1399  & 0.6\% \\
          & 0.4   & 0.4   & 0.9   & 2     & 0.0\% & 8.2   & 3.0   & $>10h$ & 569   & 9.0\% \\
          & 0.7   & 0.3   & 0.4   & 2     & 0.0\% & 3.5   & 1.1   & 905.6 & 1016  & 3.3\% \\
          & 1     & 0.9   & 0.4   & 2     & 0.0\% & 1.5   & 3.3   & $>10h$ & 548   & 8.4\% \\
    \midrule
    \multirow{4}[2]{*}{$(4,5,5,0.5)$} & 0.1   & 0.4   & 0.9   & 2     & 0.0\% & 2.2   & 0.0   & $>10h$ & 550   & 7.0\% \\
          & 0.4   & 0.5   & 1.3   & 2     & 0.0\% & 9.5   & 14.1  & $>10h$ & 284   & 7.0\% \\
          & 0.7   & 0.2   & 1.4   & 2     & 0.0\% & 3.9   & 2.0   & 1064.2 & 920   & 2.6\% \\
          & 1     & 0.2   & 1.2   & 2     & 0.0\% & 3.9   & 2.5   & $>10h$ & 912   & 5.6\% \\
    \midrule
    \multirow{3}[2]{*}{$(4,5,7,0.1)$} & 0.1   & 4.7   & 77.2  & 2     & 0.0\% & 1.4   & 0.9   & 1207.0 & 332   & 4.8\% \\
          & 0.4   & 9.7   & 133.8 & 2     & 0.0\% & 4.8   & 10.1  & $>10h$ & 433   & 8.7\% \\
          & 1     & 1.1   & 159.8 & 2     & 0.0\% & 2.7   & 5.5   & $>10h$ & 1073  & 6.8\% \\
    \midrule
    \multirow{4}[2]{*}{$(4,5,7,0.3)$} & 0.1   & 2.4   & 123.2 & 2     & 0.0\% & 2.9   & 4.6   & 886.7 & 864   & 4.9\% \\
          & 0.4   & 6.6   & 197.4 & 2     & 0.0\% & 4.9   & 2.1   & 1364.0 & 1361  & 4.4\% \\
          & 0.7   & 5.3   & 179.7 & 2     & 0.0\% & 2.2   & 2.1   & $>10h$ & 787   & 6.9\% \\
          & 1     & 3.2   & 156.8 & 2     & 0.0\% & 3.7   & 0.8   & $>10h$ & 966   & 9.3\% \\
    \midrule
    \multirow{4}[2]{*}{$(4,5,7,0.5)$} & 0.1   & 6.2   & 121.9 & 2     & 0.0\% & 1.5   & 1.2   & 1047.0 & 1296  & 4.5\% \\
          & 0.4   & 3.9   & 101.3 & 2     & 0.0\% & 10.3  & 14.3  & 988.4 & 392   & 4.7\% \\
          & 0.7   & 25.5  & 135.8 & 2     & 0.0\% & 2.2   & 6.6   & 970.9 & 562   & 4.8\% \\
          & 1     & 4.6   & 204.7 & 2     & 0.0\% & 4.1   & 3.1   & $>10h$ & 591   & 8.9\% \\
    \midrule
    \multirow{4}[2]{*}{$(4,5,10,0.1)$} & 0.1   & 4.8   & 83.7  & 2     & 1.7\% & 3.0   & 7.2   & $>10h$ & 1508  & 9.4\% \\
          & 0.4   & 13.6  & 88.8  & 2     & 1.6\% & 1.1   & 1.0   & $>10h$ & 723   & 6.4\% \\
          & 0.7   & 6.8   & 102.2 & 2     & 1.6\% & 2.1   & 5.0   & 4502.0 & 479   & 4.9\% \\
          & 1     & 6.4   & 118.4 & 2     & 1.8\% & 12.6  & 4.2   & $>10h$ & 566   & 7.4\% \\
    \midrule
    \multirow{4}[2]{*}{$(4,5,10,0.3)$} & 0.1   & 17.9  & 199.3 & 2     & 1.8\% & 3.8   & 8.2   & 7782.2 & 1059  & 4.2\% \\
          & 0.4   & 25.3  & 103.5 & 2     & 1.8\% & 6.7   & 6.3   & 12134.1 & 667   & 2.6\% \\
          & 0.7   & 13.8  & 50.2  & 2     & 1.7\% & 9.0   & 10.3  & 11770.6 & 586   & 4.5\% \\
          & 1     & 8.6   & 107.1 & 2     & 1.7\% & 2.0   & 3.8   & 9143.7 & 421   & 2.2\% \\
    \midrule
    \multirow{4}[2]{*}{$(4,5,10,0.5)$} & 0.1   & 10.4  & 189.7 & 2     & 1.7\% & 1.2   & 0.6   & $>10h$ & 366   & 9.0\% \\
          & 0.4   & 6.2   & 38.2  & 2     & 1.8\% & 5.7   & 5.5   & $>10h$ & 1429  & 9.3\% \\
          & 0.7   & 4.0   & 72.6  & 2     & 1.8\% & 3.1   & 3.1   & 12449.6 & 617   & 4.9\% \\
          & 1     & 1.0   & 41.1  & 2     & 1.8\% & 17.4  & 40.7  & 11130.6 & 739   & 4.3\% \\
    \bottomrule
    \end{tabular}%
}
\end{table}%

\begin{table}[htbp]
\ContinuedFloat
\crev
\small
  \centering
  \caption{Algorithmic details for the location-transportation problem over a budgeted uncertainty set (continued)}
\label{tab:detailed-LT-algorithms}
\scalebox{0.8}{
    \begin{tabular}{crrrrrrrrrr}
\toprule    \multirow{2}[2]{*}{$(T,I,J,\alpha^d)$} & \multicolumn{1}{c}{\multirow{2}[2]{*}{$\alpha^u$}} & \multicolumn{1}{c}{$\LDRUB$} & \multicolumn{3}{c}{$\twostageLDRUB$, C\&CG} & \multicolumn{1}{c}{$\nu^{\Omega(\text{LDR})}$} & \multicolumn{1}{c}{$\nu^{\Omega(\text{2S-LDR})}$} & \multicolumn{3}{c}{$\bar{\nu}^{\textsc{NA-DO}}_R$} \\
\cmidrule(rl){3-3}  \cmidrule(rl){4-6} \cmidrule(rl){7-7} \cmidrule(rl){8-8} \cmidrule(rl){9-11}        &       & \multicolumn{1}{c}{Time (s)} & \multicolumn{1}{c}{Time (s)} & \multicolumn{1}{c}{\#iterations} & \multicolumn{1}{c}{Gap} & \multicolumn{1}{c}{Time (s)} & \multicolumn{1}{c}{Time (s)} & Time (s) & \multicolumn{1}{c}{\#iterations} & \multicolumn{1}{c}{Gap} \\
    \midrule
    \multirow{4}[2]{*}{$(4,10,10,0.1)$} & 0.1   & 10.7  & 13318.8 & 1     & 4.1\% & 21.0  & 43.0  & $>10h$ & 1880  & 5.2\% \\
          & 0.4   & 18.9  & 13171.6 & 1     & 4.8\% & 14.8  & 9.2   & $>10h$ & 2154  & 8.7\% \\
          & 0.7   & 33.7  & 4419.5 & 2     & 4.6\% & 18.0  & 49.5  & $>10h$ & 776   & 8.2\% \\
          & 1     & 3.4   & 7159.0 & 2     & 4.8\% & 14.6  & 6.9   & 14098.5 & 1275  & 4.1\% \\
    \midrule
    \multirow{4}[2]{*}{$(4,10,10,0.3)$} & 0.1   & 35.8  & 10763.5 & 1     & 2.4\% & 4.6   & 9.1   & 9702.9 & 885   & 4.8\% \\
          & 0.4   & 12.7  & 7688.3 & 2     & 4.8\% & 9.3   & 8.2   & $>10h$ & 755   & 8.6\% \\
          & 0.7   & 24.5  & 12108.4 & 1     & 4.6\% & 8.5   & 5.2   & $>10h$ & 1343  & 13.4\% \\
          & 1     & 33.8  & 7580.7 & 2     & 3.2\% & 8.8   & 1.2   & 13045.5 & 862   & 4.9\% \\
    \midrule
    \multirow{4}[2]{*}{$(4,10,10,0.5)$} & 0.1   & 23.1  & 8361.0 & 2     & 2.4\% & 36.8  & 7.7   & 10972.8 & 765   & 1.5\% \\
          & 0.4   & 17.5  & 9441.9 & 2     & 4.4\% & 36.9  & 69.7  & $>10h$ & 791   & 5.7\% \\
          & 0.7   & 16.8  & 5638.9 & 2     & 2.3\% & 38.5  & 50.6  & 14347.7 & 1511  & 4.4\% \\
          & 1     & 26.7  & 3497.6 & 2     & 2.4\% & 29.5  & 20.6  & $>10h$ & 696   & 8.7\% \\
    \midrule
    \multirow{4}[2]{*}{$(5,5,10,0.1)$} & 0.1   & 2.3   & 372.6 & 4     & 1.6\% & 34.3  & 17.0  & $>10h$ & 1439  & 7.5\% \\
          & 0.4   & 5.8   & 989.5 & 4     & 2.8\% & 28.3  & 38.4  & $>10h$ & 267   & 8.2\% \\
          & 0.7   & 5.4   & 657.6 & 2     & 1.9\% & 21.4  & 31.2  & 15167.1 & 605   & 4.6\% \\
          & 1     & 3.6   & 784.6 & 3     & 1.7\% & 31.3  & 52.0  & $>10h$ & 904   & 7.0\% \\
    \midrule
    \multirow{4}[2]{*}{$(5,5,10,0.3)$} & 0.1   & 1.8   & 883.5 & 2     & 2.8\% & 39.3  & 84.0  & $>10h$ & 789   & 5.2\% \\
          & 0.4   & 5.5   & 607.0 & 2     & 1.9\% & 24.5  & 23.9  & 16175.2 & 1053  & 4.8\% \\
          & 0.7   & 4.5   & 414.9 & 10    & 2.8\% & 27.4  & 33.6  & 13106.8 & 762   & 4.8\% \\
          & 1     & 3.5   & 959.7 & 2     & 1.3\% & 19.5  & 32.5  & 14376.6 & 1016  & 2.0\% \\
    \midrule
    \multirow{4}[2]{*}{$(5,5,10,0.5)$} & 0.1   & 5.0   & 1009.4 & 2     & 1.8\% & 35.5  & 94.5  & $>10h$ & 418   & 5.6\ \\
          & 0.4   & 6.2   & 723.7 & 2     & 3.0\% & 38.5  & 52.0  & 13067.1 & 1329  & 4.5\% \\
          & 0.7   & 7.9   & 383.6 & 2     & 2.5\% & 39.5  & 53.6  & 15368.2 & 911   & 4.0\% \\
          & 1     & 9.9   & 724.3 & 2     & 1.6\% & 39.5  & 66.5  & 14150.1 & 521   & 4.7\% \\
    \bottomrule
    \end{tabular}%
}
\end{table}%

\newpage

\subsection{Capital Budgeting with Loan}\label{sec:CB-results}

\begin{table}[htbp]
    \crev
  \small
  \centering 
  \caption{Results for the capital budgeting problem with unrestricted loans}
  \label{tab:detailed-CB}
  \setlength{\tabcolsep}{2mm}
\scalebox{0.8}{
    \begin{tabular}{llrrrrrrrr}
         \toprule
\multirow{2}[2]{*}{Instance} & \multirow{2}[2]{*}{$(T,I)$} & \multirow{2}[2]{*}{$B$} & \multirow{2}[2]{*}{$\nu^{\texttt{K}}$} & \multicolumn{3}{c}{UB} & \multicolumn{3}{c}{Optimality Gap ($\%$)} \\
\cmidrule(rl){5-7} \cmidrule(rl){8-10}   &    &       &       & $\nu^{\Omega}$ & $\bar{\nu}^{\textsc{NA-DO}}_R$ & $\bar{\nu}^{\textsc{DNA-DO}}_R$ & $\big(\frac{\nu^{\Omega}-\nu^{\texttt{K}}}{\nu^{\texttt{K}}}\big)$ & $\big(\frac{\bar{\nu}^{\textsc{NA-DO}}_R-\nu^{\texttt{K}}}{\nu^{\texttt{K}}}\big)$ & $\big(\frac{\bar{\nu}^{\textsc{DNA-DO}}_R-\nu^{\texttt{K}}}{\nu^{\texttt{K}}}\big)$ \\
         \midrule
1     & \multirow{4}[0]{*}{(3,5)} & 0     & 1.8   & 2.2   & 2.1   & 2.1   & 19.7\% & 16.1\% & 14.5\% \\
2     &       & 50    & 6.7   & 7.9   & 7.6   & 7.5   & 18.0\% & 14.1\% & 12.6\% \\
3     &       & 100   & 7.6   & 8.8   & 8.5   & 8.4   & 16.2\% & 12.8\% & 10.4\% \\
4     &       & 150   & 7.8   & 8.8   & 8.5   & 8.3   & 13.1\% & 9.1\% & 7.0\% \\
\midrule
5     & \multirow{4}[0]{*}{(3,10)} & 0     & 6.7   & 7.9   & 7.6   & 7.5   & 18.7\% & 13.8\% & 12.3\% \\
6     &       & 50    & 11.7  & 13.9  & 13.4  & 13.1  & 18.9\% & 14.5\% & 11.8\% \\
7     &       & 100   & 16.7  & 19.7  & 19.2  & 18.8  & 18.5\% & 15.2\% & 12.6\% \\
8     &       & 150   & 16.7  & 19.7  & 19.1  & 18.7  & 18.0\% & 14.1\% & 11.8\% \\
\midrule
9     & \multirow{4}[0]{*}{(3,15)} & 0     & 8.9   & 10.4  & 10.1  & 9.9   & 16.5\% & 13.3\% & 11.0\% \\
10    &       & 50    & 13.7  & 16.4  & 15.9  & 15.5  & 19.4\% & 15.8\% & 13.0\% \\
11    &       & 100   & 18.9  & 22.4  & 21.6  & 21.1  & 18.7\% & 14.3\% & 11.8\% \\
12    &       & 150   & 21.6  & 26.0  & 25.1  & 24.6  & 20.4\% & 16.3\% & 13.9\% \\
\midrule
13    & \multirow{4}[0]{*}{(3,20)} & 0     & 11.9  & 13.8  & 13.2  & 13.0  & 16.0\% & 11.2\% & 8.9\% \\
14    &       & 50    & 16.4  & 19.8  & 19.2  & 18.9  & 21.1\% & 17.6\% & 15.6\% \\
15    &       & 100   & 21.7  & 25.8  & 25.1  & 24.7  & 19.0\% & 15.8\% & 13.8\% \\
16    &       & 150   & 27.1  & 31.8  & 31.0  & 30.4  & 17.5\% & 14.4\% & 12.4\% \\
\midrule
17    & \multirow{4}[0]{*}{(3,25)} & 0     & 17.7  & 20.2  & 19.1  & 18.8  & 14.2\% & 8.3\% & 6.7\% \\
18    &       & 50    & 22.5  & 26.2  & 25.3  & 24.9  & 16.5\% & 12.7\% & 10.9\% \\
19    &       & 100   & 26.5  & 32.2  & 31.2  & 30.8  & 21.6\% & 17.8\% & 16.4\% \\
20    &       & 150   & 31.7  & 38.2  & 36.8  & 36.2  & 20.5\% & 16.1\% & 14.3\% \\
\midrule
21    & \multirow{4}[0]{*}{(3,30)} & 0     & 22.2  & 26.4  & 25.6  & 25.3  & 18.9\% & 15.4\% & 13.9\% \\
22    &       & 50    & 26.5  & 32.4  & 31.2  & 30.8  & 22.3\% & 17.9\% & 16.5\% \\
23    &       & 100   & 31.4  & 38.4  & 37.4  & 36.9  & 22.3\% & 19.1\% & 17.6\% \\
24    &       & 150   & 37.1  & 44.4  & 42.9  & 42.4  & 19.7\% & 15.8\% & 14.3\% \\
\bottomrule
    \end{tabular}
}
\end{table}

\begin{table}[htbp]
\ContinuedFloat
    \crev
  \small
  \centering 
  \caption{Results for the capital budgeting problem with unrestricted loans (continued)}
  \label{tab:detailed-CB}
  \setlength{\tabcolsep}{2mm}
\scalebox{0.88}{
    \begin{tabular}{llrrrrrrrr}
         \toprule
\multirow{2}[2]{*}{Instance} & \multirow{2}[2]{*}{$(T,I)$} & \multirow{2}[2]{*}{$B$} & \multirow{2}[2]{*}{$\nu^{\texttt{K}}$} & \multicolumn{3}{c}{UB} & \multicolumn{3}{c}{Optimality Gap ($\%$)} \\
\cmidrule(rl){5-7} \cmidrule(rl){8-10}   &    &       &       & $\nu^{\Omega}$ & $\bar{\nu}^{\textsc{NA-DO}}_R$ & $\bar{\nu}^{\textsc{DNA-DO}}_R$ & $\big(\frac{\nu^{\Omega}-\nu^{\texttt{K}}}{\nu^{\texttt{K}}}\big)$ & $\big(\frac{\bar{\nu}^{\textsc{NA-DO}}_R-\nu^{\texttt{K}}}{\nu^{\texttt{K}}}\big)$ & $\big(\frac{\bar{\nu}^{\textsc{DNA-DO}}_R-\nu^{\texttt{K}}}{\nu^{\texttt{K}}}\big)$ \\
\midrule
25    & \multirow{6}[0]{*}{(4,5)} & 0     & 3.3   & 3.8   & 3.6   & 3.5   & 17.1\% & 10.7\% & 8.0\% \\
26    &       & 50    & 8.3   & 9.9   & 9.4   & 9.2   & 19.2\% & 13.7\% & 11.1\% \\
27    &       & 100   & 8.7   & 10.2  & 9.7   & 9.5   & 17.2\% & 12.0\% & 9.5\% \\
28    &       & 150   & 7.6   & 9.1   & 8.6   & 8.4   & 18.8\% & 13.1\% & 10.3\% \\
29    &       & 200   & 8.7   & 10.2  & 9.9   & 9.6   & 17.4\% & 13.4\% & 10.7\% \\
30    &       & 250   & 8.7   & 10.2  & 9.6   & 9.3   & 16.7\% & 9.7\% & 6.9\% \\
\midrule
31    & \multirow{6}[0]{*}{(4,10)} & 0     & 8.4   & 9.9   & 9.5   & 9.2   & 18.3\% & 12.7\% & 9.9\% \\
32    &       & 50    & 13.4  & 16.0  & 15.1  & 14.7  & 19.2\% & 12.3\% & 9.4\% \\
33    &       & 100   & 18.0  & 22.0  & 21.2  & 20.8  & 22.1\% & 17.7\% & 15.4\% \\
34    &       & 150   & 20.8  & 25.0  & 23.9  & 23.3  & 19.7\% & 14.7\% & 11.7\% \\
35    &       & 200   & 20.4  & 25.0  & 24.1  & 23.5  & 22.2\% & 18.0\% & 15.2\% \\
36    &       & 250   & 21.1  & 25.0  & 24.0  & 23.4  & 18.1\% & 13.5\% & 10.7\% \\
\midrule
37    & \multirow{6}[0]{*}{(4,15)} & 0     & 12.9  & 15.3  & 14.5  & 14.2  & 18.0\% & 11.7\% & 9.3\% \\
38    &       & 50    & 18.0  & 21.3  & 20.0  & 19.5  & 18.2\% & 11.0\% & 8.2\% \\
39    &       & 100   & 23.1  & 27.3  & 26.3  & 25.8  & 18.2\% & 14.0\% & 11.9\% \\
40    &       & 150   & 28.3  & 33.3  & 32.1  & 31.4  & 17.8\% & 13.7\% & 11.3\% \\
41    &       & 200   & 31.6  & 38.2  & 36.9  & 36.0  & 20.9\% & 16.6\% & 13.8\% \\
42    &       & 250   & 32.3  & 38.2  & 36.2  & 35.5  & 18.1\% & 12.0\% & 9.7\% \\
\midrule
43    & \multirow{6}[0]{*}{(4,20)} & 0     & 16.0  & 19.3  & 18.8  & 18.4  & 20.6\% & 17.0\% & 15.1\% \\
44    &       & 50    & 21.3  & 25.3  & 23.9  & 23.4  & 18.9\% & 12.0\% & 9.6\% \\
45    &       & 100   & 25.8  & 31.3  & 30.4  & 29.7  & 21.6\% & 17.8\% & 15.4\% \\
46    &       & 150   & 31.3  & 37.3  & 35.2  & 34.5  & 19.2\% & 12.4\% & 10.1\% \\
47    &       & 200   & 36.6  & 43.3  & 41.3  & 40.6  & 18.3\% & 12.9\% & 11.0\% \\
48    &       & 250   & 39.7  & 48.3  & 46.8  & 45.9  & 21.7\% & 17.8\% & 15.6\% \\
\midrule
49    & \multirow{6}[0]{*}{(4,25)} & 0     & 22.7  & 27.7  & 26.8  & 26.3  & 22.2\% & 18.2\% & 16.0\% \\
50    &       & 50    & 26.8  & 33.7  & 31.8  & 31.3  & 25.6\% & 18.5\% & 16.7\% \\
51    &       & 100   & 32.4  & 39.7  & 37.7  & 37.1  & 22.5\% & 16.4\% & 14.5\% \\
52    &       & 150   & 37.9  & 45.7  & 43.7  & 42.9  & 20.5\% & 15.3\% & 13.1\% \\
53    &       & 200   & 42.2  & 51.7  & 50.1  & 49.3  & 22.6\% & 18.7\% & 16.9\% \\
54    &       & 250   & 46.8  & 57.7  & 55.6  & 54.6  & 23.3\% & 18.8\% & 16.7\% \\
\midrule
55    & \multirow{6}[0]{*}{(4,30)} & 0     & 23.0  & 27.9  & 26.3  & 25.9  & 21.6\% & 14.6\% & 12.6\% \\
56    &       & 50    & 26.1  & 33.9  & 32.8  & 32.5  & 29.7\% & 25.5\% & 24.2\% \\
57    &       & 100   & 32.6  & 39.9  & 38.0  & 37.5  & 22.3\% & 16.4\% & 14.9\% \\
58    &       & 150   & 36.4  & 45.9  & 43.7  & 43.1  & 26.1\% & 20.1\% & 18.4\% \\
59    &       & 200   & 41.5  & 51.9  & 49.6  & 48.9  & 25.1\% & 19.6\% & 17.8\% \\
60    &       & 250   & 46.1  & 57.9  & 55.5  & 54.8  & 25.7\% & 20.5\% & 19.0\% \\
\specialrule{.1em}{.05em}{.05em} \\[-3mm]
61    & \multirow{8}[0]{*}{(5,5)} & 0     & 4.9   & 6.1   & 5.8   & 5.6   & 23.6\% & 17.4\% & 14.5\% \\
62    &       & 50    & 10.5  & 13.0  & 12.3  & 12.0  & 23.5\% & 16.9\% & 13.3\% \\
63    &       & 100   & 14.4  & 17.6  & 16.7  & 16.2  & 22.4\% & 16.4\% & 13.0\% \\
64    &       & 150   & 14.1  & 17.6  & 17.0  & 16.5  & 24.2\% & 20.0\% & 16.6\% \\
65    &       & 200   & 14.4  & 17.6  & 16.7  & 16.2  & 21.7\% & 15.6\% & 12.1\% \\
66    &       & 250   & 14.0  & 17.6  & 16.7  & 16.3  & 25.6\% & 19.3\% & 16.5\% \\
67    &       & 300   & 14.3  & 17.6  & 16.9  & 16.4  & 22.6\% & 17.7\% & 14.5\% \\
68    &       & 350   & 14.4  & 17.6  & 16.6  & 16.2  & 22.0\% & 15.6\% & 12.3\% \\
\midrule
69    & \multirow{8}[0]{*}{(5,10)} & 0     & 10.0  & 12.0  & 11.3  & 11.0  & 19.5\% & 12.6\% & 9.6\% \\
70    &       & 50    & 14.0  & 18.0  & 17.2  & 16.9  & 29.0\% & 23.3\% & 20.9\% \\
71    &       & 100   & 19.7  & 24.0  & 23.1  & 22.6  & 21.9\% & 17.6\% & 14.6\% \\
72    &       & 150   & 24.0  & 30.0  & 28.5  & 27.7  & 24.9\% & 18.8\% & 15.4\% \\
73    &       & 200   & 24.6  & 30.0  & 28.8  & 28.0  & 21.9\% & 17.1\% & 13.8\% \\
74    &       & 250   & 23.1  & 30.0  & 28.8  & 28.2  & 29.9\% & 24.8\% & 22.3\% \\
75    &       & 300   & 23.5  & 30.0  & 28.7  & 28.0  & 27.7\% & 22.4\% & 19.3\% \\
76    &       & 350   & 22.7  & 30.0  & 28.4  & 27.8  & 31.9\% & 25.1\% & 22.4\% \\
\bottomrule
    \end{tabular}
}
\end{table}

\begin{table}[htbp]
\ContinuedFloat
    \crev
  \small
  \centering 
  \caption{Results for the capital budgeting problem with unrestricted loans (continued)}
  \label{tab:detailed-CB}
  \setlength{\tabcolsep}{2mm}
\scalebox{0.88}{
    \begin{tabular}{llrrrrrrrr}
         \toprule
\multirow{2}[2]{*}{Instance} & \multirow{2}[2]{*}{$(T,I)$} & \multirow{2}[2]{*}{$B$} & \multirow{2}[2]{*}{$\nu^{\texttt{K}}$} & \multicolumn{3}{c}{UB} & \multicolumn{3}{c}{Optimality Gap ($\%$)} \\
\cmidrule(rl){5-7} \cmidrule(rl){8-10}   &    &       &       & $\nu^{\Omega}$ & $\bar{\nu}^{\textsc{NA-DO}}_R$ & $\bar{\nu}^{\textsc{DNA-DO}}_R$ & $\big(\frac{\nu^{\Omega}-\nu^{\texttt{K}}}{\nu^{\texttt{K}}}\big)$ & $\big(\frac{\bar{\nu}^{\textsc{NA-DO}}_R-\nu^{\texttt{K}}}{\nu^{\texttt{K}}}\big)$ & $\big(\frac{\bar{\nu}^{\textsc{DNA-DO}}_R-\nu^{\texttt{K}}}{\nu^{\texttt{K}}}\big)$ \\
\midrule
77    & \multirow{8}[0]{*}{(5,15)} & 0     & 12.7  & 16.2  & 15.3  & 14.9  & 28.1\% & 20.8\% & 17.5\% \\
78    &       & 50    & 17.6  & 22.2  & 21.5  & 21.1  & 26.2\% & 22.3\% & 19.6\% \\
79    &       & 100   & 21.7  & 28.2  & 26.6  & 25.8  & 29.9\% & 22.5\% & 18.9\% \\
80    &       & 150   & 26.2  & 34.2  & 32.4  & 31.5  & 30.5\% & 23.7\% & 20.0\% \\
81    &       & 200   & 30.8  & 40.2  & 38.9  & 38.1  & 30.7\% & 26.3\% & 23.8\% \\
82    &       & 250   & 31.6  & 40.5  & 39.0  & 38.2  & 28.1\% & 23.3\% & 20.8\% \\
83    &       & 300   & 30.7  & 40.5  & 39.2  & 38.4  & 32.2\% & 27.8\% & 25.4\% \\
84    &       & 350   & 31.3  & 40.5  & 38.8  & 38.0  & 29.6\% & 24.0\% & 21.6\% \\
\midrule
85    & \multirow{8}[0]{*}{(5,20)} & 0     & 19.1  & 23.6  & 22.9  & 22.3  & 23.1\% & 19.6\% & 16.4\% \\
86    &       & 50    & 23.1  & 29.6  & 28.6  & 27.9  & 27.9\% & 23.7\% & 20.7\% \\
87    &       & 100   & 28.1  & 35.6  & 34.4  & 33.7  & 26.6\% & 22.6\% & 20.0\% \\
88    &       & 150   & 33.8  & 41.6  & 39.7  & 38.9  & 23.1\% & 17.7\% & 15.3\% \\
89    &       & 200   & 38.2  & 47.6  & 46.5  & 45.4  & 24.5\% & 21.8\% & 18.9\% \\
90    &       & 250   & 39.9  & 53.6  & 52.4  & 51.1  & 34.2\% & 31.2\% & 28.0\% \\
91    &       & 300   & 46.8  & 58.9  & 57.0  & 55.7  & 25.9\% & 21.7\% & 19.1\% \\
92    &       & 350   & 45.9  & 58.9  & 57.7  & 56.6  & 28.4\% & 25.7\% & 23.3\% \\
\midrule
93    & \multirow{8}[0]{*}{(5,25)} & 0     & 23.5  & 29.6  & 29.0  & 28.4  & 25.9\% & 23.4\% & 20.9\% \\
94    &       & 50    & 28.6  & 35.6  & 34.7  & 33.9  & 24.3\% & 21.3\% & 18.3\% \\
95    &       & 100   & 33.4  & 41.6  & 39.9  & 39.3  & 24.6\% & 19.7\% & 17.6\% \\
96    &       & 150   & 36.5  & 47.6  & 45.9  & 44.9  & 30.5\% & 25.9\% & 23.1\% \\
97    &       & 200   & 42.6  & 53.6  & 51.7  & 50.6  & 25.9\% & 21.4\% & 19.0\% \\
98    &       & 250   & 47.1  & 59.6  & 58.2  & 56.9  & 26.4\% & 23.4\% & 20.6\% \\
99    &       & 300   & 51.6  & 65.6  & 63.5  & 62.6  & 27.0\% & 23.0\% & 21.2\% \\
100   &       & 350   & 56.1  & 71.6  & 69.7  & 68.3  & 27.6\% & 24.3\% & 21.8\% \\
\midrule
101   & \multirow{8}[0]{*}{(5,30)} & 0     & 25.6  & 32.2  & 31.1  & 30.6  & 26.1\% & 21.8\% & 19.8\% \\
102   &       & 50    & 29.1  & 38.2  & 37.1  & 36.5  & 31.3\% & 27.4\% & 25.3\% \\
103   &       & 100   & 34.4  & 44.2  & 43.3  & 42.6  & 28.7\% & 26.0\% & 24.1\% \\
104   &       & 150   & 37.8  & 50.2  & 48.7  & 48.1  & 32.8\% & 28.7\% & 27.2\% \\
105   &       & 200   & 42.0  & 56.2  & 54.7  & 53.5  & 33.9\% & 30.2\% & 27.4\% \\
106   &       & 250   & 47.2  & 62.2  & 60.3  & 59.5  & 31.9\% & 27.8\% & 26.2\% \\
107   &       & 300   & 54.4  & 68.2  & 66.4  & 65.3  & 25.5\% & 22.1\% & 20.0\% \\
108   &       & 350   & 60.6  & 74.2  & 72.1  & 70.7  & 22.4\% & 18.9\% & 16.5\% \\
\specialrule{.1em}{.05em}{.05em} \\[-3mm]
109   & \multirow{9}[0]{*}{(6,5)} & 0     & 5.6   & 6.8   & 6.6   & 6.4   & 20.2\% & 18.0\% & 14.1\% \\
110   &       & 50    & 10.6  & 12.8  & 12.5  & 12.1  & 20.5\% & 18.0\% & 14.2\% \\
111   &       & 100   & 13.6  & 16.9  & 16.3  & 15.8  & 24.0\% & 19.8\% & 16.4\% \\
112   &       & 150   & 12.8  & 16.9  & 16.5  & 16.0  & 31.4\% & 28.4\% & 24.8\% \\
113   &       & 200   & 13.2  & 16.9  & 16.5  & 16.0  & 27.7\% & 24.9\% & 21.2\% \\
114   &       & 250   & 12.9  & 16.9  & 16.5  & 15.9  & 30.5\% & 27.7\% & 23.1\% \\
115   &       & 300   & 13.0  & 16.9  & 16.5  & 15.9  & 30.1\% & 27.2\% & 22.8\% \\
116   &       & 350   & 13.3  & 16.9  & 16.5  & 16.0  & 27.2\% & 24.2\% & 20.3\% \\
117   &       & 400   & 12.9  & 16.9  & 16.4  & 15.9  & 30.7\% & 26.9\% & 23.1\% \\
\midrule
118   & \multirow{9}[0]{*}{(6,10)} & 0     & 10.9  & 14.0  & 13.8  & 13.4  & 27.8\% & 26.0\% & 22.5\% \\
119   &       & 50    & 15.9  & 20.0  & 19.6  & 18.9  & 25.8\% & 23.5\% & 19.0\% \\
120   &       & 100   & 20.7  & 26.0  & 25.4  & 24.5  & 25.4\% & 22.4\% & 18.2\% \\
121   &       & 150   & 24.9  & 32.0  & 31.4  & 30.3  & 28.3\% & 25.8\% & 21.5\% \\
122   &       & 200   & 27.1  & 35.0  & 34.2  & 33.3  & 29.2\% & 26.5\% & 23.1\% \\
123   &       & 250   & 26.6  & 35.0  & 34.2  & 33.2  & 31.5\% & 28.7\% & 24.7\% \\
124   &       & 300   & 26.8  & 35.0  & 34.4  & 33.3  & 30.4\% & 28.1\% & 24.0\% \\
125   &       & 350   & 27.3  & 35.0  & 34.4  & 33.3  & 28.3\% & 26.3\% & 22.2\% \\
126   &       & 400   & 26.9  & 35.0  & 34.2  & 33.2  & 29.9\% & 27.0\% & 23.3\% \\
\bottomrule
    \end{tabular}
}
\end{table}

\newpage
\ 
\begin{table}[H]
\ContinuedFloat
    \crev
  \small
  \centering 
  \caption{Results for the capital budgeting problem with unrestricted loans (continued)}
  \label{tab:detailed-CB}
  \setlength{\tabcolsep}{2mm}
\scalebox{0.88}{
    \begin{tabular}{llrrrrrrrr}
         \toprule
\multirow{2}[2]{*}{Instance} & \multirow{2}[2]{*}{$(T,I)$} & \multirow{2}[2]{*}{$B$} & \multirow{2}[2]{*}{$\nu^{\texttt{K}}$} & \multicolumn{3}{c}{UB} & \multicolumn{3}{c}{Optimality Gap ($\%$)} \\
\cmidrule(rl){5-7} \cmidrule(rl){8-10}   &    &       &       & $\nu^{\Omega}$ & $\bar{\nu}^{\textsc{NA-DO}}_R$ & $\bar{\nu}^{\textsc{DNA-DO}}_R$ & $\big(\frac{\nu^{\Omega}-\nu^{\texttt{K}}}{\nu^{\texttt{K}}}\big)$ & $\big(\frac{\bar{\nu}^{\textsc{NA-DO}}_R-\nu^{\texttt{K}}}{\nu^{\texttt{K}}}\big)$ & $\big(\frac{\bar{\nu}^{\textsc{DNA-DO}}_R-\nu^{\texttt{K}}}{\nu^{\texttt{K}}}\big)$ \\
\midrule
127   & \multirow{9}[0]{*}{(6,15)} & 0     & 15.9  & 19.9  & 19.5  & 19.0  & 24.7\% & 22.5\% & 19.3\% \\
128   &       & 50    & 19.8  & 25.9  & 25.5  & 24.9  & 30.7\% & 28.7\% & 25.8\% \\
129   &       & 100   & 25.3  & 31.9  & 31.4  & 30.6  & 25.8\% & 23.9\% & 20.9\% \\
130   &       & 150   & 28.9  & 37.9  & 37.3  & 36.4  & 31.2\% & 29.3\% & 26.1\% \\
131   &       & 200   & 33.5  & 43.9  & 43.1  & 42.1  & 30.8\% & 28.3\% & 25.5\% \\
132   &       & 250   & 37.0  & 49.7  & 48.7  & 47.6  & 34.3\% & 31.7\% & 28.7\% \\
133   &       & 300   & 36.9  & 49.7  & 49.0  & 47.8  & 34.8\% & 32.9\% & 29.7\% \\
134   &       & 350   & 38.5  & 49.7  & 48.8  & 47.6  & 29.1\% & 26.8\% & 23.8\% \\
135   &       & 400   & 37.9  & 49.7  & 48.9  & 47.7  & 31.2\% & 29.0\% & 25.9\% \\
\midrule
136   & \multirow{9}[0]{*}{(6,20)} & 0     & 22.4  & 28.6  & 28.3  & 27.8  & 27.6\% & 26.5\% & 23.9\% \\
137   &       & 50    & 27.6  & 34.6  & 34.3  & 33.5  & 25.4\% & 24.3\% & 21.6\% \\
138   &       & 100   & 31.7  & 40.6  & 40.1  & 39.3  & 27.9\% & 26.3\% & 23.7\% \\
139   &       & 150   & 36.9  & 46.6  & 46.1  & 45.3  & 26.4\% & 25.2\% & 22.8\% \\
140   &       & 200   & 41.0  & 52.6  & 52.1  & 50.9  & 28.4\% & 27.2\% & 24.4\% \\
141   &       & 250   & 45.1  & 58.6  & 58.2  & 57.1  & 29.8\% & 28.8\% & 26.4\% \\
142   &       & 300   & 49.7  & 64.6  & 63.8  & 62.6  & 30.0\% & 28.3\% & 26.0\% \\
143   &       & 350   & 53.2  & 70.6  & 70.0  & 68.8  & 32.6\% & 31.5\% & 29.3\% \\
144   &       & 400   & 54.3  & 71.5  & 70.9  & 69.6  & 31.6\% & 30.6\% & 28.1\% \\
\midrule
145   & \multirow{9}[0]{*}{(6,25)} & 0     & 29.0  & 37.2  & 36.9  & 36.2  & 28.2\% & 27.2\% & 24.9\% \\
146   &       & 50    & 32.8  & 43.2  & 42.8  & 42.1  & 31.5\% & 30.3\% & 28.2\% \\
147   &       & 100   & 37.5  & 49.2  & 48.7  & 48.0  & 31.3\% & 30.0\% & 28.0\% \\
148   &       & 150   & 42.5  & 55.2  & 54.6  & 53.7  & 29.9\% & 28.5\% & 26.3\% \\
149   &       & 200   & 46.2  & 61.2  & 60.5  & 59.7  & 32.5\% & 31.0\% & 29.2\% \\
150   &       & 250   & 51.9  & 67.2  & 66.7  & 65.7  & 29.4\% & 28.5\% & 26.6\% \\
151   &       & 300   & 56.1  & 73.2  & 72.3  & 71.2  & 30.4\% & 28.9\% & 26.8\% \\
152   &       & 350   & 61.1  & 79.2  & 78.6  & 77.4  & 29.6\% & 28.6\% & 26.7\% \\
153   &       & 400   & 64.5  & 85.2  & 84.5  & 83.0  & 32.0\% & 30.9\% & 28.7\% \\
\midrule
154   & \multirow{9}[0]{*}{(6,30)} & 0     & 31.5  & 42.6  & 42.2  & 41.6  & 35.2\% & 34.1\% & 32.2\% \\
155   &       & 50    & 37.5  & 48.6  & 48.1  & 47.5  & 29.5\% & 28.3\% & 26.6\% \\
156   &       & 100   & 40.6  & 54.6  & 54.1  & 53.4  & 34.4\% & 33.1\% & 31.5\% \\
157   &       & 150   & 45.2  & 60.6  & 60.0  & 59.2  & 33.9\% & 32.6\% & 30.8\% \\
158   &       & 200   & 49.1  & 66.6  & 66.0  & 65.1  & 35.5\% & 34.4\% & 32.5\% \\
159   &       & 250   & 56.0  & 72.6  & 71.9  & 71.0  & 29.7\% & 28.5\% & 26.9\% \\
160   &       & 300   & 59.4  & 78.6  & 78.0  & 76.8  & 32.4\% & 31.4\% & 29.4\% \\
161   &       & 350   & 65.0  & 84.6  & 83.8  & 82.6  & 30.2\% & 29.0\% & 27.1\% \\
162   &       & 400   & 68.5  & 90.6  & 89.8  & 88.5  & 32.2\% & 31.1\% & 29.2\% \\
\bottomrule
    \end{tabular}
}
\end{table}

\newpage

\subsection{\crev Capital  Budgeting with Integer Loan}
\label{app:results-generalinteger}
\begin{figure}[htbp]
\footnotesize
    \centering
    \scalebox{0.7}{
    % This file was created with tikzplotlib v0.10.1.
\begin{tikzpicture}

\definecolor{darkgray176}{RGB}{176,176,176}
\definecolor{forestgreen015351}{RGB}{0,153,51}
\definecolor{lightgray204}{RGB}{204,204,204}
\definecolor{magenta}{RGB}{255,0,255}
\definecolor{peru23712549}{RGB}{237,125,49}

\begin{groupplot}[group style={
    group size=1 by 4, 
    vertical sep=2cm  
  },width=1.2\textwidth,height=6cm,]
\nextgroupplot[
tick align=outside,
tick pos=left,
title={$T = 3$},
x grid style={darkgray176},
xlabel={},
xmin=1, xmax=24,
xtick=data,
xtick style={color=black},
xticklabel style={rotate=90.0,anchor=east},
y grid style={darkgray176},
ylabel={Optimality gap},
ymin=5, ymax=37,
ytick style={color=black},
yticklabel=\pgfmathprintnumber{\tick}\%
]
\addplot [semithick, magenta, mark=*, mark size=3, mark options={solid}]
table {%
1 19.7
2 18
3 16.2
4 13.1
5 18.7
6 18.9
7 18.5
8 18
9 16.5
10 19.4
11 18.7
12 20.4
13 16
14 21.1
15 19
16 17.5
17 14.2
18 16.5
19 21.6
20 20.5
21 18.8547380939758
22 22.3
23 22.3
24 19.7
};
\addplot [semithick, forestgreen015351, mark=asterisk, mark size=3, mark options={solid}]
table {%
1 16.07
2 14.1
3 12.77
4 9.12000000000001
5 13.8
6 14.47
7 15.18
8 14.09
9 13.26
10 15.82
11 14.26
12 16.29
13 11.18
14 17.61
15 15.79
16 14.44
17 8.25000000000001
18 12.65
19 17.83
20 16.12
21 15.3547380939759
22 17.94
23 19.13
24 15.76
};
\addplot [semithick, peru23712549, mark=triangle*, mark size=3, mark options={solid}]
table {%
1 14.5
2 12.55
3 10.35
4 7.03
5 12.27
6 11.79
7 12.56
8 11.81
9 10.98
10 12.97
11 11.8
12 13.87
13 8.88999999999999
14 15.55
15 13.78
16 12.41
17 6.71999999999999
18 10.9
19 16.37
20 14.31
21 13.9247380939758
22 16.47
23 17.61
24 14.26
};

\nextgroupplot[
tick align=outside,
tick pos=left,
title={$T = 4$},
x grid style={darkgray176},
xlabel={},
xmin=25, xmax=60,
xtick=data,
xtick style={color=black},
xticklabel style={rotate=90.0,anchor=east},
y grid style={darkgray176},
ylabel={Optimality gap},
ymin=5, ymax=37,
ytick style={color=black},
yticklabel=\pgfmathprintnumber{\tick}\%
]
\addplot [semithick, magenta, mark=*, mark size=3, mark options={solid}]
table {%
25 17.1
26 19.2
27 17.2
28 18.8
29 17.4
30 16.6992028905506
31 18.3
32 19.2
33 22.1
34 19.7
35 22.2
36 18.1
37 18
38 18.2
39 18.2
40 17.8
41 20.9
42 18.1
43 20.6
44 18.9
45 21.6
46 19.2
47 18.3
48 21.7
49 22.2
50 25.6
51 22.5
52 20.5
53 22.6
54 23.3
55 21.6
56 29.7
57 22.3
58 26.1
59 25.1
60 25.7
};
\addplot [semithick, forestgreen015351, mark=asterisk, mark size=3, mark options={solid}]
table {%
25 10.68
26 13.65
27 12
28 13.05
29 13.39
30 9.70000000000002
31 12.68
32 12.34
33 17.69
34 14.71
35 18.04
36 13.53
37 11.72
38 11.01
39 13.98
40 13.69
41 16.64
42 11.96
43 17
44 11.95
45 17.82
46 12.39
47 12.87
48 17.79
49 18.18
50 18.46
51 16.36
52 15.27
53 18.74
54 18.83
55 14.55
56 25.45
57 16.41
58 20.12
59 19.61
60 20.48
};
\addplot [semithick, peru23712549, mark=triangle*, mark size=3, mark options={solid}]
table {%
25 7.99
26 11.11
27 9.53999999999999
28 10.25
29 10.73
30 6.94000000000001
31 9.88999999999999
32 9.37
33 15.41
34 11.74
35 15.2
36 10.67
37 9.32999999999999
38 8.17999999999999
39 11.85
40 11.29
41 13.81
42 9.67999999999999
43 15.09
44 9.62
45 15.41
46 10.11
47 10.97
48 15.64
49 16
50 16.65
51 14.48
52 13.09
53 16.86
54 16.67
55 12.64
56 24.24
57 14.85
58 18.37
59 17.8
60 18.99
};

\nextgroupplot[
tick align=outside,
tick pos=left,
title={$T = 5$},
x grid style={darkgray176},
xlabel={},
xmin=61, xmax=108,
xtick=data,
xtick style={color=black},
xticklabel style={rotate=90.0,anchor=east},
y grid style={darkgray176},
ylabel={Optimality gap},
ymin=5, ymax=37,
ytick style={color=black},
yticklabel=\pgfmathprintnumber{\tick}\%
]
\addplot [semithick, magenta, mark=*, mark size=3, mark options={solid}]
table {%
61 23.6
62 23.5
63 22.4
64 24.2
65 21.7
66 25.6
67 22.6
68 22
69 19.5
70 29
71 21.9
72 24.9
73 21.9
74 29.9
75 27.7
76 31.9
77 28.1
78 26.2
79 29.9
80 30.5
81 30.7
82 28.1
83 32.2
84 29.6
85 23.1
86 27.9
87 26.6
88 23.1
89 24.5
90 34.2
91 25.9
92 28.4
93 25.9
94 24.3
95 24.6
96 30.5
97 25.9
98 26.4
99 27
100 27.6
101 26.1
102 31.3
103 28.7
104 32.8
105 33.9
106 31.9
107 25.5
108 22.4
};
\addplot [semithick, forestgreen015351, mark=asterisk, mark size=3, mark options={solid}]
table {%
61 17.36
62 16.88
63 16.44
64 19.98
65 15.55
66 19.28
67 17.66
68 15.62
69 12.55
70 23.29
71 17.59
72 18.78
73 17.13
74 24.77
75 22.39
76 25.06
77 20.82
78 22.28
79 22.48
80 23.69
81 26.34
82 23.3
83 27.78
84 23.98
85 19.59
86 23.72
87 22.55
88 17.7
89 21.76
90 31.2
91 21.72
92 25.71
93 23.41
94 21.26
95 19.68
96 25.88
97 21.4
98 23.44
99 23.02
100 24.27
101 21.79
102 27.44
103 26.02
104 28.69
105 30.23
106 27.84
107 22.1
108 18.93
};
\addplot [semithick, peru23712549, mark=triangle*, mark size=3, mark options={solid}]
table {%
61 14.54
62 13.32
63 12.99
64 16.55
65 12.12
66 16.48
67 14.46
68 12.32
69 9.61
70 20.89
71 14.61
72 15.35
73 13.83
74 22.34
75 19.31
76 22.35
77 17.45
78 19.59
79 18.91
80 19.96
81 23.79
82 20.79
83 25.38
84 21.61
85 16.41
86 20.67
87 19.98
88 15.27
89 18.92
90 28
91 19.11
92 23.27
93 20.86
94 18.31
95 17.63
96 23.13
97 18.97
98 20.62
99 21.19
100 21.82
101 19.78
102 25.31
103 24.1
104 27.22
105 27.37
106 26.23
107 20.04
108 16.54
};

\nextgroupplot[
legend cell align={left},
legend style={
  fill opacity=0.8,
  draw opacity=1,
  text opacity=1,
  at={(0.97,0.03)},
  anchor=south east,
  draw=lightgray204
},
tick align=outside,
tick pos=left,
title={$T = 6$},
xmin=109, xmax=162,
xtick=data,
xtick style={color=black},
xticklabel style={rotate=90.0,anchor=east},
y grid style={darkgray176},
ylabel={Optimality gap},
ymin=5, ymax=37,
ytick style={color=black},
yticklabel=\pgfmathprintnumber{\tick}\%
]
\addplot [semithick, magenta, mark=*, mark size=3, mark options={solid}]
table {%
109 20.2
110 20.5
111 24
112 31.4
113 27.7
114 30.5
115 30.1
116 27.2
117 30.7
118 27.8
119 25.8
120 25.4
121 28.3
122 29.2
123 31.5
124 30.4
125 28.3
126 29.9
127 24.7
128 30.7
129 25.8
130 31.2
131 30.8
132 34.3
133 34.8
134 29.1
135 31.2
136 27.6
137 25.4
138 27.9
139 26.4
140 28.4
141 29.8
142 30
143 32.6
144 31.6
145 28.2
146 31.5
147 31.3
148 29.9
149 32.5
150 29.4
151 30.4
152 29.6
153 32
154 35.2
155 29.5
156 34.4
157 33.9
158 35.5
159 29.7
160 32.4
161 30.2
162 32.2
};
\addlegendentry{${(\nu^{\Omega}-\nu^{\texttt{K}})}/{\nu^{\texttt{K}}}$}
\addplot [semithick, forestgreen015351, mark=asterisk, mark size=3, mark options={solid}]
table {%
109 17.98
110 18.04
111 19.83
112 28.35
113 24.89
114 27.67
115 27.2
116 24.24
117 26.92
118 26.04
119 23.49
120 22.35
121 25.8
122 26.48
123 28.74
124 28.08
125 26.26
126 26.97
127 22.47
128 28.74
129 23.91
130 29.27
131 28.34
132 31.72
133 32.9
134 26.83
135 29
136 26.46
137 24.26
138 26.32
139 25.17
140 27.23
141 28.83
142 28.32
143 31.49
144 30.59
145 27.23
146 30.3
147 29.97
148 28.5
149 31.04
150 28.46
151 28.85
152 28.59
153 30.91
154 34.05
155 28.34
156 33.13
157 32.64
158 34.36
159 28.46
160 31.41
161 28.98
162 31.13
};
\addlegendentry{${({\crev \bar{\nu}^{\textsc{NA-DO}}_R}-\nu^{\texttt{K}})}/{\nu^{\texttt{K}}}$}
\addplot [semithick, peru23712549, mark=triangle*, mark size=3, mark options={solid}]
table {%
109 14.08
110 14.19
111 16.35
112 24.81
113 21.15
114 23.07
115 22.76
116 20.27
117 23.05
118 22.51
119 19.04
120 18.22
121 21.52
122 23.06
123 24.65
124 23.99
125 22.19
126 23.3
127 19.33
128 25.8
129 20.94
130 26.11
131 25.5
132 28.74
133 29.71
134 23.76
135 25.9
136 23.9
137 21.6
138 23.73
139 22.83
140 24.37
141 26.42
142 25.99
143 29.26
144 28.08
145 24.93
146 28.21
147 28.03
148 26.31
149 29.21
150 26.55
151 26.84
152 26.68
153 28.69
154 32.15
155 26.57
156 31.53
157 30.78
158 32.49
159 26.87
160 29.43
161 27.12
162 29.17
};
\addlegendentry{${({\crev \bar{\nu}^{\textsc{DNA-DO}}_R}-\nu^{\texttt{K}})}/{\nu^{\texttt{K}}}$}
\end{groupplot}

\end{tikzpicture}
    }
    \caption{\crev Optimality gap improvements for capital budgeting instances with general integer loan decisions.}
    \label{fig:capital-budget-int-results}
\end{figure}

\newpage
\begin{table}[htbp]
\crev
  \centering
  \caption{\crev Solution times of different bounding methods in seconds for the capital budgeting problem}
    \label{tab:sol-time-CB}%
  \setlength{\tabcolsep}{3mm}
  \scalebox{0.8}{
    \begin{tabular}{rrrrrrrrrrr}
    \toprule
    \multirow{2}[4]{*}{Instance} & \multirow{2}[4]{*}{$(T,I)$} & \multirow{2}[4]{*}{$B$} & \multicolumn{4}{c}{Fractional loan} & \multicolumn{4}{c}{Integral loan} \\
\cmidrule(rl){4-7}  \cmidrule(rl){8-11}         &       &       & \cellcolor[rgb]{ .949,  .949,  .949}$\nu^{\texttt{K}}$ & $\nu^{\Omega}$ & \cellcolor[rgb]{ .949,  .949,  .949}$\bar{\nu}^{\textsc{NA}}_R$ & $\bar{\nu}^{\textsc{DNA}}_R$ & \cellcolor[rgb]{ .949,  .949,  .949}$\nu^{\texttt{K}}$ & $\nu^{\Omega}$ & \cellcolor[rgb]{ .949,  .949,  .949}$\bar{\nu}^{\textsc{NA}}_R$ & $\bar{\nu}^{\textsc{DNA}}_R$ \\
    \midrule
    1     & \multirow{4}[2]{*}{(3,5)} & 0     & \cellcolor[rgb]{ .949,  .949,  .949}27.0 & 63.0  & \cellcolor[rgb]{ .949,  .949,  .949}56.1 & 62.1  & \cellcolor[rgb]{ .949,  .949,  .949}67.2 & 72.7  & \cellcolor[rgb]{ .949,  .949,  .949}112.3 & 80.8 \\
    2     &       & 50    & \cellcolor[rgb]{ .949,  .949,  .949}16.2 & 85.2  & \cellcolor[rgb]{ .949,  .949,  .949}83.8 & 74.9  & \cellcolor[rgb]{ .949,  .949,  .949}52.6 & 99.3  & \cellcolor[rgb]{ .949,  .949,  .949}678.9 & 142.2 \\
    3     &       & 100   & \cellcolor[rgb]{ .949,  .949,  .949}25.8 & 60.8  & \cellcolor[rgb]{ .949,  .949,  .949}55.5 & 52.2  & \cellcolor[rgb]{ .949,  .949,  .949}43.5 & 95.6  & \cellcolor[rgb]{ .949,  .949,  .949}260.8 & 287.1 \\
    4     &       & 150   & \cellcolor[rgb]{ .949,  .949,  .949}25.2 & 56.2  & \cellcolor[rgb]{ .949,  .949,  .949}55.3 & 35.5  & \cellcolor[rgb]{ .949,  .949,  .949}64.0 & 59.3  & \cellcolor[rgb]{ .949,  .949,  .949}259.9 & 156.1 \\
    \midrule
    5     & \multirow{4}[2]{*}{(3,10)} & 0     & \cellcolor[rgb]{ .949,  .949,  .949}35.4 & 547.4 & \cellcolor[rgb]{ .949,  .949,  .949}547.4 & 521.2 & \cellcolor[rgb]{ .949,  .949,  .949}111.3 & 840.0 & \cellcolor[rgb]{ .949,  .949,  .949}1133.1 & 948.6 \\
    6     &       & 50    & \cellcolor[rgb]{ .949,  .949,  .949}40.5 & 384.5 & \cellcolor[rgb]{ .949,  .949,  .949}254.5 & 376.1 & \cellcolor[rgb]{ .949,  .949,  .949}77.5 & 771.1 & \cellcolor[rgb]{ .949,  .949,  .949}413.1 & 557.9 \\
    7     &       & 100   & \cellcolor[rgb]{ .949,  .949,  .949}34.5 & 340.5 & \cellcolor[rgb]{ .949,  .949,  .949}392.5 & 366.2 & \cellcolor[rgb]{ .949,  .949,  .949}96.5 & 373.2 & \cellcolor[rgb]{ .949,  .949,  .949}478.1 & 432.2 \\
    8     &       & 150   & \cellcolor[rgb]{ .949,  .949,  .949}44.4 & 514.4 & \cellcolor[rgb]{ .949,  .949,  .949}657.4 & 471.2 & \cellcolor[rgb]{ .949,  .949,  .949}107.7 & 1215.7 & \cellcolor[rgb]{ .949,  .949,  .949}1447.6 & 604.0 \\
    \midrule
    9     & \multirow{4}[2]{*}{(3,15)} & 0     & \cellcolor[rgb]{ .949,  .949,  .949}90.9 & 1502.9 & \cellcolor[rgb]{ .949,  .949,  .949}1178.9 & 1390.9 & \cellcolor[rgb]{ .949,  .949,  .949}218.5 & 2029.6 & \cellcolor[rgb]{ .949,  .949,  .949}3812.6 & 3569.0 \\
    10    &       & 50    & \cellcolor[rgb]{ .949,  .949,  .949}90.0 & 1169.0 & \cellcolor[rgb]{ .949,  .949,  .949}1007.0 & 1061.0 & \cellcolor[rgb]{ .949,  .949,  .949}201.8 & 1467.1 & \cellcolor[rgb]{ .949,  .949,  .949}2698.8 & 2868.9 \\
    11    &       & 100   & \cellcolor[rgb]{ .949,  .949,  .949}96.3 & 1438.3 & \cellcolor[rgb]{ .949,  .949,  .949}1168.3 & 1556.3 & \cellcolor[rgb]{ .949,  .949,  .949}169.1 & 2261.3 & \cellcolor[rgb]{ .949,  .949,  .949}2553.9 & 5394.1 \\
    12    &       & 150   & \cellcolor[rgb]{ .949,  .949,  .949}94.8 & 1535.8 & \cellcolor[rgb]{ .949,  .949,  .949}1247.8 & 1355.8 & \cellcolor[rgb]{ .949,  .949,  .949}200.5 & 1974.3 & \cellcolor[rgb]{ .949,  .949,  .949}3441.4 & 3313.6 \\
    \midrule
    13    & \multirow{4}[2]{*}{(3,20)} & 0     & \cellcolor[rgb]{ .949,  .949,  .949}611.4 & 3138.4 & \cellcolor[rgb]{ .949,  .949,  .949}1896.4 & 2568.4 & \cellcolor[rgb]{ .949,  .949,  .949}689.5 & 4305.1 & \cellcolor[rgb]{ .949,  .949,  .949}3531.1 & 5606.8 \\
    14    &       & 50    & \cellcolor[rgb]{ .949,  .949,  .949}653.1 & 3809.1 & \cellcolor[rgb]{ .949,  .949,  .949}2498.1 & 3005.1 & \cellcolor[rgb]{ .949,  .949,  .949}772.2 & 5003.3 & \cellcolor[rgb]{ .949,  .949,  .949}3962.0 & 7371.5 \\
    15    &       & 100   & \cellcolor[rgb]{ .949,  .949,  .949}667.2 & 3057.2 & \cellcolor[rgb]{ .949,  .949,  .949}2643.2 & 4295.2 & \cellcolor[rgb]{ .949,  .949,  .949}743.4 & 3777.1 & \cellcolor[rgb]{ .949,  .949,  .949}5651.2 & 12202.7 \\
    16    &       & 150   & \cellcolor[rgb]{ .949,  .949,  .949}609.6 & 3539.6 & \cellcolor[rgb]{ .949,  .949,  .949}2918.6 & 2557.6 & \cellcolor[rgb]{ .949,  .949,  .949}765.9 & 4162.9 & \cellcolor[rgb]{ .949,  .949,  .949}5034.6 & 5319.8 \\
    \midrule
    17    & \multirow{4}[2]{*}{(3,25)} & 0     & \cellcolor[rgb]{ .949,  .949,  .949}4292.4 & 8689.4 & \cellcolor[rgb]{ .949,  .949,  .949}4392.9 & 5521.4 & \cellcolor[rgb]{ .949,  .949,  .949}4561.3 & 9868.9 & \cellcolor[rgb]{ .949,  .949,  .949}7116.6 & 10188.8 \\
    18    &       & 50    & \cellcolor[rgb]{ .949,  .949,  .949}3762.0 & 6794.0 & \cellcolor[rgb]{ .949,  .949,  .949}3465.3 & 5684.5 & \cellcolor[rgb]{ .949,  .949,  .949}3953.5 & 7354.5 & \cellcolor[rgb]{ .949,  .949,  .949}4636.6 & 11134.0 \\
    19    &       & 100   & \cellcolor[rgb]{ .949,  .949,  .949}3759.0 & 6921.0 & \cellcolor[rgb]{ .949,  .949,  .949}3830.0 & 4160.5 & \cellcolor[rgb]{ .949,  .949,  .949}4054.2 & 8311.5 & \cellcolor[rgb]{ .949,  .949,  .949}6707.6 & 6271.3 \\
    20    &       & 150   & \cellcolor[rgb]{ .949,  .949,  .949}4012.8 & 7998.8 & \cellcolor[rgb]{ .949,  .949,  .949}4716.5 & 4608.3 & \cellcolor[rgb]{ .949,  .949,  .949}4336.0 & 8751.4 & \cellcolor[rgb]{ .949,  .949,  .949}6883.0 & 8159.8 \\
    \midrule
    21    & \multirow{4}[2]{*}{(3,30)} & 0     & \cellcolor[rgb]{ .949,  .949,  .949}1808.1 & 5929.1 & \cellcolor[rgb]{ .949,  .949,  .949}3754.7 & 8207.1 & \cellcolor[rgb]{ .949,  .949,  .949}2188.2 & 6894.2 & \cellcolor[rgb]{ .949,  .949,  .949}6115.2 & 20660.0 \\
    22    &       & 50    & \cellcolor[rgb]{ .949,  .949,  .949}2083.8 & 5512.8 & \cellcolor[rgb]{ .949,  .949,  .949}2817.2 & 8133.8 & \cellcolor[rgb]{ .949,  .949,  .949}2394.6 & 6036.5 & \cellcolor[rgb]{ .949,  .949,  .949}3818.2 & 15524.7 \\
    23    &       & 100   & \cellcolor[rgb]{ .949,  .949,  .949}2027.1 & 5422.1 & \cellcolor[rgb]{ .949,  .949,  .949}3152.7 & 9087.1 & \cellcolor[rgb]{ .949,  .949,  .949}2307.5 & 5763.8 & \cellcolor[rgb]{ .949,  .949,  .949}4756.4 & 23493.2 \\
    24    &       & 150   & \cellcolor[rgb]{ .949,  .949,  .949}1850.7 & 6203.7 & \cellcolor[rgb]{ .949,  .949,  .949}4069.8 & 7591.7 & \cellcolor[rgb]{ .949,  .949,  .949}2225.7 & 7524.8 & \cellcolor[rgb]{ .949,  .949,  .949}5719.4 & 15765.4 \\
    \midrule
    25    & \multirow{6}[2]{*}{(4,5)} & 0     & \cellcolor[rgb]{ .949,  .949,  .949}89.2 & 1012.2 & \cellcolor[rgb]{ .949,  .949,  .949}626.8 & 1185.2 & \cellcolor[rgb]{ .949,  .949,  .949}173.3 & 3395.6 & \cellcolor[rgb]{ .949,  .949,  .949}686.6 & 1483.1 \\
    26    &       & 50    & \cellcolor[rgb]{ .949,  .949,  .949}97.2 & 1141.2 & \cellcolor[rgb]{ .949,  .949,  .949}592.8 & 1236.2 & \cellcolor[rgb]{ .949,  .949,  .949}190.1 & 1536.5 & \cellcolor[rgb]{ .949,  .949,  .949}665.1 & 1478.5 \\
    27    &       & 100   & \cellcolor[rgb]{ .949,  .949,  .949}106.0 & 1189.0 & \cellcolor[rgb]{ .949,  .949,  .949}624.7 & 1259.0 & \cellcolor[rgb]{ .949,  .949,  .949}173.2 & 2337.4 & \cellcolor[rgb]{ .949,  .949,  .949}706.7 & 1488.1 \\
    28    &       & 150   & \cellcolor[rgb]{ .949,  .949,  .949}80.8 & 849.8 & \cellcolor[rgb]{ .949,  .949,  .949}470.5 & 932.8 & \cellcolor[rgb]{ .949,  .949,  .949}170.7 & 1740.1 & \cellcolor[rgb]{ .949,  .949,  .949}516.6 & 1071.5 \\
    29    &       & 200   & \cellcolor[rgb]{ .949,  .949,  .949}85.2 & 838.2 & \cellcolor[rgb]{ .949,  .949,  .949}390.8 & 952.2 & \cellcolor[rgb]{ .949,  .949,  .949}160.4 & 3339.3 & \cellcolor[rgb]{ .949,  .949,  .949}423.1 & 1060.8 \\
    30    &       & 250   & \cellcolor[rgb]{ .949,  .949,  .949}107.2 & 1039.2 & \cellcolor[rgb]{ .949,  .949,  .949}596.8 & 1173.2 & \cellcolor[rgb]{ .949,  .949,  .949}202.7 & 2826.5 & \cellcolor[rgb]{ .949,  .949,  .949}665.2 & 1332.0 \\
    \midrule
    31    & \multirow{6}[2]{*}{(4,10)} & 0     & \cellcolor[rgb]{ .949,  .949,  .949}590.0 & 2420.0 & \cellcolor[rgb]{ .949,  .949,  .949}1277.3 & 3150.0 & \cellcolor[rgb]{ .949,  .949,  .949}713.9 & 3186.4 & \cellcolor[rgb]{ .949,  .949,  .949}1544.7 & 5161.8 \\
    32    &       & 50    & \cellcolor[rgb]{ .949,  .949,  .949}593.2 & 2247.2 & \cellcolor[rgb]{ .949,  .949,  .949}1050.1 & 3178.2 & \cellcolor[rgb]{ .949,  .949,  .949}720.7 & 5845.5 & \cellcolor[rgb]{ .949,  .949,  .949}1295.2 & 4602.0 \\
    33    &       & 100   & \cellcolor[rgb]{ .949,  .949,  .949}580.8 & 2189.8 & \cellcolor[rgb]{ .949,  .949,  .949}1049.2 & 2941.8 & \cellcolor[rgb]{ .949,  .949,  .949}779.5 & 12174.6 & \cellcolor[rgb]{ .949,  .949,  .949}1236.0 & 4424.5 \\
    34    &       & 150   & \cellcolor[rgb]{ .949,  .949,  .949}470.0 & 2230.0 & \cellcolor[rgb]{ .949,  .949,  .949}1449.3 & 3485.0 & \cellcolor[rgb]{ .949,  .949,  .949}636.5 & 2278.7 & \cellcolor[rgb]{ .949,  .949,  .949}1820.4 & 5099.7 \\
    35    &       & 200   & \cellcolor[rgb]{ .949,  .949,  .949}417.2 & 2208.2 & \cellcolor[rgb]{ .949,  .949,  .949}1285.5 & 2833.2 & \cellcolor[rgb]{ .949,  .949,  .949}602.3 & 5185.9 & \cellcolor[rgb]{ .949,  .949,  .949}1572.6 & 4482.1 \\
    36    &       & 250   & \cellcolor[rgb]{ .949,  .949,  .949}431.6 & 2198.6 & \cellcolor[rgb]{ .949,  .949,  .949}1353.7 & 3231.6 & \cellcolor[rgb]{ .949,  .949,  .949}583.6 & 10838.1 & \cellcolor[rgb]{ .949,  .949,  .949}1576.6 & 4862.5 \\
    \midrule
    37    & \multirow{6}[2]{*}{(4,15)} & 0     & \cellcolor[rgb]{ .949,  .949,  .949}1372.8 & 5153.8 & \cellcolor[rgb]{ .949,  .949,  .949}3131.9 & 5853.8 & \cellcolor[rgb]{ .949,  .949,  .949}1679.6 & 7835.0 & \cellcolor[rgb]{ .949,  .949,  .949}4739.6 & 11098.8 \\
    38    &       & 50    & \cellcolor[rgb]{ .949,  .949,  .949}1370.4 & 4828.4 & \cellcolor[rgb]{ .949,  .949,  .949}2509.6 & 6153.4 & \cellcolor[rgb]{ .949,  .949,  .949}1724.7 & 8235.8 & \cellcolor[rgb]{ .949,  .949,  .949}3476.6 & 12602.2 \\
    39    &       & 100   & \cellcolor[rgb]{ .949,  .949,  .949}1268.0 & 5458.0 & \cellcolor[rgb]{ .949,  .949,  .949}3486.7 & 7861.0 & \cellcolor[rgb]{ .949,  .949,  .949}1525.9 & 15446.3 & \cellcolor[rgb]{ .949,  .949,  .949}5541.5 & 21591.5 \\
    40    &       & 150   & \cellcolor[rgb]{ .949,  .949,  .949}1507.6 & 5316.6 & \cellcolor[rgb]{ .949,  .949,  .949}2784.4 & 7977.6 & \cellcolor[rgb]{ .949,  .949,  .949}1750.2 & 17220.9 & \cellcolor[rgb]{ .949,  .949,  .949}3909.3 & 21449.1 \\
    41    &       & 200   & \cellcolor[rgb]{ .949,  .949,  .949}1426.4 & 5243.4 & \cellcolor[rgb]{ .949,  .949,  .949}3039.6 & 6315.4 & \cellcolor[rgb]{ .949,  .949,  .949}1613.6 & 12527.3 & \cellcolor[rgb]{ .949,  .949,  .949}4895.8 & 11460.3 \\
    42    &       & 250   & \cellcolor[rgb]{ .949,  .949,  .949}1589.2 & 5702.2 & \cellcolor[rgb]{ .949,  .949,  .949}2990.8 & 8046.2 & \cellcolor[rgb]{ .949,  .949,  .949}1911.2 & 19692.0 & \cellcolor[rgb]{ .949,  .949,  .949}4633.7 & 21338.5 \\
    \midrule
    43    & \multirow{6}[2]{*}{(4,20)} & 0     & \cellcolor[rgb]{ .949,  .949,  .949}4676.4 & 10535.4 & \cellcolor[rgb]{ .949,  .949,  .949}5935.6 & 17371.4 & \cellcolor[rgb]{ .949,  .949,  .949}5228.1 & 20173.7 & \cellcolor[rgb]{ .949,  .949,  .949}9827.4 &$>10h$\\
    44    &       & 50    & \cellcolor[rgb]{ .949,  .949,  .949}4988.0 & 10909.0 & \cellcolor[rgb]{ .949,  .949,  .949}5992.7 & 16613.0 & \cellcolor[rgb]{ .949,  .949,  .949}5402.2 & 14597.9 & \cellcolor[rgb]{ .949,  .949,  .949}7704.6 &$>10h$\\
    45    &       & 100   & \cellcolor[rgb]{ .949,  .949,  .949}4942.8 & 10226.8 & \cellcolor[rgb]{ .949,  .949,  .949}5601.9 & 17168.8 & \cellcolor[rgb]{ .949,  .949,  .949}5356.3 & 20126.0 & \cellcolor[rgb]{ .949,  .949,  .949}9011.5 &$>10h$\\
    46    &       & 150   & \cellcolor[rgb]{ .949,  .949,  .949}4111.2 & 9830.2 & \cellcolor[rgb]{ .949,  .949,  .949}5337.5 & 13249.2 & \cellcolor[rgb]{ .949,  .949,  .949}4642.0 & 27734.9 & \cellcolor[rgb]{ .949,  .949,  .949}7855.0 &$>10h$\\
    47    &       & 200   & \cellcolor[rgb]{ .949,  .949,  .949}4633.6 & 9934.6 & \cellcolor[rgb]{ .949,  .949,  .949}6239.1 & 16356.6 & \cellcolor[rgb]{ .949,  .949,  .949}5176.0 & 16374.7 & \cellcolor[rgb]{ .949,  .949,  .949}9543.7 &$>10h$\\
    48    &       & 250   & \cellcolor[rgb]{ .949,  .949,  .949}4344.4 & 9147.4 & \cellcolor[rgb]{ .949,  .949,  .949}5522.3 & 14321.4 & \cellcolor[rgb]{ .949,  .949,  .949}4749.8 & 18371.0 & \cellcolor[rgb]{ .949,  .949,  .949}7313.3 &$>10h$\\
    \midrule
    49    & \multirow{6}[2]{*}{(4,25)} & 0     & \cellcolor[rgb]{ .949,  .949,  .949}8113.2 & 16000.2 & \cellcolor[rgb]{ .949,  .949,  .949}10434.8 & 25348.2 & \cellcolor[rgb]{ .949,  .949,  .949}8769.9 &$>10h$ & \cellcolor[rgb]{ .949,  .949,  .949}13120.7 &$>10h$\\
        50    &       & 50    & \cellcolor[rgb]{ .949,  .949,  .949}7384.0 & 14978.0 & \cellcolor[rgb]{ .949,  .949,  .949}9753.3 & 21113.0 & \cellcolor[rgb]{ .949,  .949,  .949}7679.0 &$>10h$ & \cellcolor[rgb]{ .949,  .949,  .949}12535.0 &$>10h$\\
    51    &       & 100   & \cellcolor[rgb]{ .949,  .949,  .949}7619.2 & 14613.2 & \cellcolor[rgb]{ .949,  .949,  .949}9046.1 & 20856.2 & \cellcolor[rgb]{ .949,  .949,  .949}7865.7 &$>10h$ & \cellcolor[rgb]{ .949,  .949,  .949}12054.9 &$>10h$\\
    52    &       & 150   & \cellcolor[rgb]{ .949,  .949,  .949}7864.4 & 15145.4 & \cellcolor[rgb]{ .949,  .949,  .949}9942.3 & 22623.4 & \cellcolor[rgb]{ .949,  .949,  .949}8270.7 &$>10h$ & \cellcolor[rgb]{ .949,  .949,  .949}12352.3 &$>10h$\\
    53    &       & 200   & \cellcolor[rgb]{ .949,  .949,  .949}7440.4 & 14736.4 & \cellcolor[rgb]{ .949,  .949,  .949}8354.9 & 25575.4 & \cellcolor[rgb]{ .949,  .949,  .949}7731.2 &$>10h$ & \cellcolor[rgb]{ .949,  .949,  .949}9985.8 &$>10h$\\
    54    &       & 250   & \cellcolor[rgb]{ .949,  .949,  .949}7986.4 & 15734.4 & \cellcolor[rgb]{ .949,  .949,  .949}8246.9 & 27221.4 & \cellcolor[rgb]{ .949,  .949,  .949}8265.6 &$>10h$ & \cellcolor[rgb]{ .949,  .949,  .949}10005.2 &$>10h$\\
    \bottomrule
    \end{tabular}%
    }
\end{table}%

\begin{table}[htbp]
\ContinuedFloat
\crev
  \centering
  \caption{\crev Solution times of different bounding methods in seconds for the capital budgeting problem (continued)}
    \label{tab:sol-time-CB}%
  \setlength{\tabcolsep}{3mm}
  \scalebox{0.8}{
    \begin{tabular}{rrrrrrrrrrr}
    \toprule
    \multirow{2}[4]{*}{Instance} & \multirow{2}[4]{*}{$(T,I)$} & \multirow{2}[4]{*}{$B$} & \multicolumn{4}{c}{Fractional loan} & \multicolumn{4}{c}{Integral loan} \\
\cmidrule(rl){4-7}  \cmidrule(rl){8-11}         &       &       & \cellcolor[rgb]{ .949,  .949,  .949}$\nu^{\texttt{K}}$ & $\nu^{\Omega}$ & \cellcolor[rgb]{ .949,  .949,  .949}$\bar{\nu}^{\textsc{NA}}_R$ & $\bar{\nu}^{\textsc{DNA}}_R$ & \cellcolor[rgb]{ .949,  .949,  .949}$\nu^{\texttt{K}}$ & $\nu^{\Omega}$ & \cellcolor[rgb]{ .949,  .949,  .949}$\bar{\nu}^{\textsc{NA}}_R$ & $\bar{\nu}^{\textsc{DNA}}_R$ \\
    \midrule
    55    & \multirow{6}[2]{*}{(4,30)} & 0     & \cellcolor[rgb]{ .949,  .949,  .949}13898.8 & 20608.8 & \cellcolor[rgb]{ .949,  .949,  .949}12197.9 &$>10h$ & \cellcolor[rgb]{ .949,  .949,  .949}14146.4 &$>10h$ & \cellcolor[rgb]{ .949,  .949,  .949}15337.6 &$>10h$\\
    56    &       & 50    & \cellcolor[rgb]{ .949,  .949,  .949}13221.6 & 20082.6 & \cellcolor[rgb]{ .949,  .949,  .949}10940.4 & 33232.6 & \cellcolor[rgb]{ .949,  .949,  .949}13470.7 &$>10h$ & \cellcolor[rgb]{ .949,  .949,  .949}14137.2 &$>10h$\\
    57    &       & 100   & \cellcolor[rgb]{ .949,  .949,  .949}12854.0 & 20045.0 & \cellcolor[rgb]{ .949,  .949,  .949}13182.0 &$>10h$ & \cellcolor[rgb]{ .949,  .949,  .949}13435.0 &$>10h$ & \cellcolor[rgb]{ .949,  .949,  .949}19973.4 &$>10h$\\
    58    &       & 150   & \cellcolor[rgb]{ .949,  .949,  .949}13757.6 & 20898.6 & \cellcolor[rgb]{ .949,  .949,  .949}13569.7 &$>10h$ & \cellcolor[rgb]{ .949,  .949,  .949}14320.1 &$>10h$ & \cellcolor[rgb]{ .949,  .949,  .949}19708.7 &$>10h$\\
    59    &       & 200   & \cellcolor[rgb]{ .949,  .949,  .949}12260.8 & 19245.8 & \cellcolor[rgb]{ .949,  .949,  .949}10654.5 &$>10h$ & \cellcolor[rgb]{ .949,  .949,  .949}12773.9 &$>10h$ & \cellcolor[rgb]{ .949,  .949,  .949}14818.3 &$>10h$\\
    60    &       & 250   & \cellcolor[rgb]{ .949,  .949,  .949}13656.4 & 20471.4 & \cellcolor[rgb]{ .949,  .949,  .949}12559.6 &$>10h$ & \cellcolor[rgb]{ .949,  .949,  .949}14180.5 &$>10h$ & \cellcolor[rgb]{ .949,  .949,  .949}18532.9 &$>10h$\\
    \midrule
    61    & \multirow{8}[2]{*}{(5,5)} & 0     & \cellcolor[rgb]{ .949,  .949,  .949}1188.0 & 4539.0 & \cellcolor[rgb]{ .949,  .949,  .949}2692.7 & 6018.0 & \cellcolor[rgb]{ .949,  .949,  .949}1302.5 & 6041.9 & \cellcolor[rgb]{ .949,  .949,  .949}2993.6 & 7706.0 \\
    62    &       & 50    & \cellcolor[rgb]{ .949,  .949,  .949}1219.0 & 4453.0 & \cellcolor[rgb]{ .949,  .949,  .949}2802.0 & 5922.0 & \cellcolor[rgb]{ .949,  .949,  .949}1297.0 & 4489.5 & \cellcolor[rgb]{ .949,  .949,  .949}3164.9 & 6700.7 \\
    63    &       & 100   & \cellcolor[rgb]{ .949,  .949,  .949}1469.0 & 4346.0 & \cellcolor[rgb]{ .949,  .949,  .949}2864.0 & 5529.0 & \cellcolor[rgb]{ .949,  .949,  .949}1545.1 & 12850.6 & \cellcolor[rgb]{ .949,  .949,  .949}3207.0 & 6774.4 \\
    64    &       & 150   & \cellcolor[rgb]{ .949,  .949,  .949}1315.0 & 4093.0 & \cellcolor[rgb]{ .949,  .949,  .949}2428.7 & 5442.0 & \cellcolor[rgb]{ .949,  .949,  .949}1443.7 & 7996.8 & \cellcolor[rgb]{ .949,  .949,  .949}2647.2 & 6284.1 \\
    65    &       & 200   & \cellcolor[rgb]{ .949,  .949,  .949}1099.0 & 4026.0 & \cellcolor[rgb]{ .949,  .949,  .949}2384.0 & 5303.0 & \cellcolor[rgb]{ .949,  .949,  .949}1247.0 & 12676.3 & \cellcolor[rgb]{ .949,  .949,  .949}2573.5 & 6182.0 \\
    66    &       & 250   & \cellcolor[rgb]{ .949,  .949,  .949}1194.0 & 4438.0 & \cellcolor[rgb]{ .949,  .949,  .949}2825.3 & 6075.0 & \cellcolor[rgb]{ .949,  .949,  .949}1294.4 & 6630.7 & \cellcolor[rgb]{ .949,  .949,  .949}3110.0 & 7217.1 \\
    67    &       & 300   & \cellcolor[rgb]{ .949,  .949,  .949}812.5 & 3603.5 & \cellcolor[rgb]{ .949,  .949,  .949}2135.7 & 4741.5 & \cellcolor[rgb]{ .949,  .949,  .949}923.0 & 12057.9 & \cellcolor[rgb]{ .949,  .949,  .949}2253.7 & 5600.9 \\
    68    &       & 350   & \cellcolor[rgb]{ .949,  .949,  .949}1376.5 & 5019.5 & \cellcolor[rgb]{ .949,  .949,  .949}3046.3 & 6457.5 & \cellcolor[rgb]{ .949,  .949,  .949}1484.3 & 19850.1 & \cellcolor[rgb]{ .949,  .949,  .949}3340.3 & 7406.8 \\
    \midrule
    69    & \multirow{8}[2]{*}{(5,10)} & 0     & \cellcolor[rgb]{ .949,  .949,  .949}4275.5 & 10708.5 & \cellcolor[rgb]{ .949,  .949,  .949}9053.2 & 14882.5 & \cellcolor[rgb]{ .949,  .949,  .949}4512.4 & 14699.0 & \cellcolor[rgb]{ .949,  .949,  .949}11601.7 & 20366.7 \\
    70    &       & 50    & \cellcolor[rgb]{ .949,  .949,  .949}3354.0 & 11245.0 & \cellcolor[rgb]{ .949,  .949,  .949}9677.3 & 13469.0 & \cellcolor[rgb]{ .949,  .949,  .949}3547.8 & 13864.4 & \cellcolor[rgb]{ .949,  .949,  .949}11116.0 & 16596.1 \\
    71    &       & 100   & \cellcolor[rgb]{ .949,  .949,  .949}3764.5 & 12332.5 & \cellcolor[rgb]{ .949,  .949,  .949}10597.7 & 14540.5 & \cellcolor[rgb]{ .949,  .949,  .949}3959.8 & 20538.2 & \cellcolor[rgb]{ .949,  .949,  .949}12756.1 & 18039.9 \\
    72    &       & 150   & \cellcolor[rgb]{ .949,  .949,  .949}3769.0 & 10903.0 & \cellcolor[rgb]{ .949,  .949,  .949}10864.4 & 13991.0 & \cellcolor[rgb]{ .949,  .949,  .949}3991.5 & 22267.6 & \cellcolor[rgb]{ .949,  .949,  .949}13768.9 & 19160.7 \\
    73    &       & 200   & \cellcolor[rgb]{ .949,  .949,  .949}3507.5 & 10891.5 & \cellcolor[rgb]{ .949,  .949,  .949}10037.7 & 13358.5 & \cellcolor[rgb]{ .949,  .949,  .949}3679.3 & 21045.9 & \cellcolor[rgb]{ .949,  .949,  .949}12135.6 & 18143.1 \\
    74    &       & 250   & \cellcolor[rgb]{ .949,  .949,  .949}4005.5 & 12899.5 & \cellcolor[rgb]{ .949,  .949,  .949}11522.3 & 16746.5 & \cellcolor[rgb]{ .949,  .949,  .949}4244.5 & 17293.1 & \cellcolor[rgb]{ .949,  .949,  .949}13665.4 & 20972.2 \\
    75    &       & 300   & \cellcolor[rgb]{ .949,  .949,  .949}3507.0 & 10998.0 & \cellcolor[rgb]{ .949,  .949,  .949}13136.7 & 13005.0 & \cellcolor[rgb]{ .949,  .949,  .949}3724.3 & 25089.7 & \cellcolor[rgb]{ .949,  .949,  .949}18343.2 & 17303.2 \\
    76    &       & 350   & \cellcolor[rgb]{ .949,  .949,  .949}3803.0 & 9881.0 & \cellcolor[rgb]{ .949,  .949,  .949}8232.7 & 12413.0 & \cellcolor[rgb]{ .949,  .949,  .949}4006.9 & 21655.3 & \cellcolor[rgb]{ .949,  .949,  .949}9519.8 & 17092.7 \\
    \midrule
    77    & \multirow{8}[2]{*}{(5,15)} & 0     & \cellcolor[rgb]{ .949,  .949,  .949}7936.0 & 22038.0 & \cellcolor[rgb]{ .949,  .949,  .949}14625.3 & 28044.0 & \cellcolor[rgb]{ .949,  .949,  .949}8368.5 &$>10h$ & \cellcolor[rgb]{ .949,  .949,  .949}16465.2 &$>10h$\\
    78    &       & 50    & \cellcolor[rgb]{ .949,  .949,  .949}8516.5 & 22334.5 & \cellcolor[rgb]{ .949,  .949,  .949}14089.7 & 29802.5 & \cellcolor[rgb]{ .949,  .949,  .949}8893.6 &$>10h$ & \cellcolor[rgb]{ .949,  .949,  .949}15924.1 &$>10h$\\
    79    &       & 100   & \cellcolor[rgb]{ .949,  .949,  .949}9214.0 & 21936.0 & \cellcolor[rgb]{ .949,  .949,  .949}14357.3 & 29063.0 & \cellcolor[rgb]{ .949,  .949,  .949}9757.7 &$>10h$ & \cellcolor[rgb]{ .949,  .949,  .949}15971.1 &$>10h$\\
    80    &       & 150   & \cellcolor[rgb]{ .949,  .949,  .949}9150.0 & 22127.0 & \cellcolor[rgb]{ .949,  .949,  .949}13484.7 & 28673.0 & \cellcolor[rgb]{ .949,  .949,  .949}9541.0 &$>10h$ & \cellcolor[rgb]{ .949,  .949,  .949}14723.9 &$>10h$\\
    81    &       & 200   & \cellcolor[rgb]{ .949,  .949,  .949}8487.5 & 21143.5 & \cellcolor[rgb]{ .949,  .949,  .949}13829.0 & 27331.5 & \cellcolor[rgb]{ .949,  .949,  .949}8850.1 &$>10h$ & \cellcolor[rgb]{ .949,  .949,  .949}15523.1 &$>10h$\\
    82    &       & 250   & \cellcolor[rgb]{ .949,  .949,  .949}7664.5 & 21657.5 & \cellcolor[rgb]{ .949,  .949,  .949}13638.3 & 28499.5 & \cellcolor[rgb]{ .949,  .949,  .949}8033.2 &$>10h$ & \cellcolor[rgb]{ .949,  .949,  .949}14898.5 &$>10h$\\
    83    &       & 300   & \cellcolor[rgb]{ .949,  .949,  .949}8442.5 & 22287.5 & \cellcolor[rgb]{ .949,  .949,  .949}14258.3 & 27426.5 & \cellcolor[rgb]{ .949,  .949,  .949}8785.5 &$>10h$ & \cellcolor[rgb]{ .949,  .949,  .949}15829.6 &$>10h$\\
    84    &       & 350   & \cellcolor[rgb]{ .949,  .949,  .949}9361.5 & 19241.5 & \cellcolor[rgb]{ .949,  .949,  .949}12027.7 & 27496.5 & \cellcolor[rgb]{ .949,  .949,  .949}9747.9 &$>10h$ & \cellcolor[rgb]{ .949,  .949,  .949}13072.9 &$>10h$\\
    \midrule
    85    & \multirow{8}[2]{*}{(5,20)} & 0     & \cellcolor[rgb]{ .949,  .949,  .949}12474.0 & 31538.0 & \cellcolor[rgb]{ .949,  .949,  .949}24198.3 &$>10h$ & \cellcolor[rgb]{ .949,  .949,  .949}13006.9 &$>10h$ & \cellcolor[rgb]{ .949,  .949,  .949}28183.8 &$>10h$\\
    86    &       & 50    & \cellcolor[rgb]{ .949,  .949,  .949}13804.0 & 30383.0 & \cellcolor[rgb]{ .949,  .949,  .949}22819.2 &$>10h$ & \cellcolor[rgb]{ .949,  .949,  .949}14552.9 &$>10h$ & \cellcolor[rgb]{ .949,  .949,  .949}26801.1 &$>10h$\\
    87    &       & 100   & \cellcolor[rgb]{ .949,  .949,  .949}13531.5 & 27682.5 & \cellcolor[rgb]{ .949,  .949,  .949}21402.1 &$>10h$ & \cellcolor[rgb]{ .949,  .949,  .949}14425.2 &$>10h$ & \cellcolor[rgb]{ .949,  .949,  .949}25016.9 &$>10h$\\
    88    &       & 150   & \cellcolor[rgb]{ .949,  .949,  .949}12393.5 & 30610.5 & \cellcolor[rgb]{ .949,  .949,  .949}23529.6 &$>10h$ & \cellcolor[rgb]{ .949,  .949,  .949}13144.1 &$>10h$ & \cellcolor[rgb]{ .949,  .949,  .949}28047.3 &$>10h$\\
    89    &       & 200   & \cellcolor[rgb]{ .949,  .949,  .949}13991.5 & 32146.5 & \cellcolor[rgb]{ .949,  .949,  .949}25330.4 &$>10h$ & \cellcolor[rgb]{ .949,  .949,  .949}14732.0 &$>10h$ & \cellcolor[rgb]{ .949,  .949,  .949}31559.2 &$>10h$\\
    90    &       & 250   & \cellcolor[rgb]{ .949,  .949,  .949}12486.0 & 28563.0 & \cellcolor[rgb]{ .949,  .949,  .949}22135.8 &$>10h$ & \cellcolor[rgb]{ .949,  .949,  .949}13412.9 &$>10h$ & \cellcolor[rgb]{ .949,  .949,  .949}26835.3 &$>10h$\\
    91    &       & 300   & \cellcolor[rgb]{ .949,  .949,  .949}14566.5 & 32117.5 & \cellcolor[rgb]{ .949,  .949,  .949}24472.9 &$>10h$ & \cellcolor[rgb]{ .949,  .949,  .949}15184.6 &$>10h$ & \cellcolor[rgb]{ .949,  .949,  .949}28655.3 &$>10h$\\
    92    &       & 350   & \cellcolor[rgb]{ .949,  .949,  .949}14411.0 & 32416.0 & \cellcolor[rgb]{ .949,  .949,  .949}26284.2 &$>10h$ & \cellcolor[rgb]{ .949,  .949,  .949}15094.5 &$>10h$ & \cellcolor[rgb]{ .949,  .949,  .949}31906.3 &$>10h$\\
    \midrule
    93    & \multirow{8}[2]{*}{(5,25)} & 0     & \cellcolor[rgb]{ .949,  .949,  .949}19659.5 &$>10h$ & \cellcolor[rgb]{ .949,  .949,  .949}22510.6 &$>10h$ & \cellcolor[rgb]{ .949,  .949,  .949}20733.1 &$>10h$ & \cellcolor[rgb]{ .949,  .949,  .949}$>10h$ &$>10h$\\
    94    &       & 50    & \cellcolor[rgb]{ .949,  .949,  .949}19808.5 &$>10h$ & \cellcolor[rgb]{ .949,  .949,  .949}26676.1 &$>10h$ & \cellcolor[rgb]{ .949,  .949,  .949}20452.7 &$>10h$ & \cellcolor[rgb]{ .949,  .949,  .949}$>10h$ &$>10h$\\
    95    &       & 100   & \cellcolor[rgb]{ .949,  .949,  .949}21565.0 &$>10h$ & \cellcolor[rgb]{ .949,  .949,  .949}27794.4 &$>10h$ & \cellcolor[rgb]{ .949,  .949,  .949}22605.3 &$>10h$ & \cellcolor[rgb]{ .949,  .949,  .949}$>10h$ &$>10h$\\
    96    &       & 150   & \cellcolor[rgb]{ .949,  .949,  .949}20473.5 &$>10h$ & \cellcolor[rgb]{ .949,  .949,  .949}25248.3 &$>10h$ & \cellcolor[rgb]{ .949,  .949,  .949}21390.7 &$>10h$ & \cellcolor[rgb]{ .949,  .949,  .949}$>10h$ &$>10h$\\
    97    &       & 200   & \cellcolor[rgb]{ .949,  .949,  .949}20571.5 &$>10h$ & \cellcolor[rgb]{ .949,  .949,  .949}25523.9 &$>10h$ & \cellcolor[rgb]{ .949,  .949,  .949}21776.2 &$>10h$ & \cellcolor[rgb]{ .949,  .949,  .949}$>10h$ &$>10h$\\
    98    &       & 250   & \cellcolor[rgb]{ .949,  .949,  .949}21483.0 &$>10h$ & \cellcolor[rgb]{ .949,  .949,  .949}24536.7 &$>10h$ & \cellcolor[rgb]{ .949,  .949,  .949}22207.9 &$>10h$ & \cellcolor[rgb]{ .949,  .949,  .949}$>10h$ &$>10h$\\
    99    &       & 300   & \cellcolor[rgb]{ .949,  .949,  .949}22307.5 &$>10h$ & \cellcolor[rgb]{ .949,  .949,  .949}27119.4 &$>10h$ & \cellcolor[rgb]{ .949,  .949,  .949}23392.4 &$>10h$ & \cellcolor[rgb]{ .949,  .949,  .949}$>10h$ &$>10h$\\
    100   &       & 350   & \cellcolor[rgb]{ .949,  .949,  .949}18784.5 &$>10h$ & \cellcolor[rgb]{ .949,  .949,  .949}25705.0 &$>10h$ & \cellcolor[rgb]{ .949,  .949,  .949}19746.9 &$>10h$ & \cellcolor[rgb]{ .949,  .949,  .949}$>10h$ &$>10h$\\
    \midrule
    101   & \multirow{8}[2]{*}{(5,30)} & 0     & \cellcolor[rgb]{ .949,  .949,  .949}16882.5 &$>10h$ & \cellcolor[rgb]{ .949,  .949,  .949}21286.1 &$>10h$ & \cellcolor[rgb]{ .949,  .949,  .949}17358.9 &$>10h$ & \cellcolor[rgb]{ .949,  .949,  .949}$>10h$ &$>10h$\\
    102   &       & 50    & \cellcolor[rgb]{ .949,  .949,  .949}15126.5 &$>10h$ & \cellcolor[rgb]{ .949,  .949,  .949}17585.0 &$>10h$ & \cellcolor[rgb]{ .949,  .949,  .949}15856.9 &$>10h$ & \cellcolor[rgb]{ .949,  .949,  .949}$>10h$ &$>10h$\\
    103   &       & 100   & \cellcolor[rgb]{ .949,  .949,  .949}17118.0 &$>10h$ & \cellcolor[rgb]{ .949,  .949,  .949}25825.6 &$>10h$ & \cellcolor[rgb]{ .949,  .949,  .949}17614.8 &$>10h$ & \cellcolor[rgb]{ .949,  .949,  .949}$>10h$ &$>10h$\\
    104   &       & 150   & \cellcolor[rgb]{ .949,  .949,  .949}16249.5 &$>10h$ & \cellcolor[rgb]{ .949,  .949,  .949}20582.8 &$>10h$ & \cellcolor[rgb]{ .949,  .949,  .949}16861.9 &$>10h$ & \cellcolor[rgb]{ .949,  .949,  .949}$>10h$ &$>10h$\\
    105   &       & 200   & \cellcolor[rgb]{ .949,  .949,  .949}17442.5 &$>10h$ & \cellcolor[rgb]{ .949,  .949,  .949}22880.6 &$>10h$ & \cellcolor[rgb]{ .949,  .949,  .949}18291.0 &$>10h$ & \cellcolor[rgb]{ .949,  .949,  .949}$>10h$ &$>10h$\\
    106   &       & 250   & \cellcolor[rgb]{ .949,  .949,  .949}16167.0 &$>10h$ & \cellcolor[rgb]{ .949,  .949,  .949}21657.8 &$>10h$ & \cellcolor[rgb]{ .949,  .949,  .949}16757.5 &$>10h$ & \cellcolor[rgb]{ .949,  .949,  .949}$>10h$ &$>10h$\\
    107   &       & 300   & \cellcolor[rgb]{ .949,  .949,  .949}16427.5 &$>10h$ & \cellcolor[rgb]{ .949,  .949,  .949}19808.3 &$>10h$ & \cellcolor[rgb]{ .949,  .949,  .949}17289.3 &$>10h$ & \cellcolor[rgb]{ .949,  .949,  .949}$>10h$ &$>10h$\\
    108   &       & 350   & \cellcolor[rgb]{ .949,  .949,  .949}16505.0 &$>10h$ & \cellcolor[rgb]{ .949,  .949,  .949}21061.1 &$>10h$ & \cellcolor[rgb]{ .949,  .949,  .949}17336.6 &$>10h$ & \cellcolor[rgb]{ .949,  .949,  .949}$>10h$ &$>10h$\\
    \bottomrule
    \end{tabular}%
    }
\end{table}%

\begin{table}[htbp]
\ContinuedFloat
\crev
  \centering
  \caption{\crev Solution times of different bounding methods in seconds for the capital budgeting problem (continued)}
    \label{tab:sol-time-CB}%
  \setlength{\tabcolsep}{3mm}
  \scalebox{0.8}{
    \begin{tabular}{rrrrrrrrrrr}
    \toprule
    \multirow{2}[4]{*}{Instance} & \multirow{2}[4]{*}{$(T,I)$} & \multirow{2}[4]{*}{$B$} & \multicolumn{4}{c}{Fractional loan} & \multicolumn{4}{c}{Integral loan} \\
\cmidrule(rl){4-7}  \cmidrule(rl){8-11}         &       &       & \cellcolor[rgb]{ .949,  .949,  .949}$\nu^{\texttt{K}}$ & $\nu^{\Omega}$ & \cellcolor[rgb]{ .949,  .949,  .949}$\bar{\nu}^{\textsc{NA}}_R$ & $\bar{\nu}^{\textsc{DNA}}_R$ & \cellcolor[rgb]{ .949,  .949,  .949}$\nu^{\texttt{K}}$ & $\nu^{\Omega}$ & \cellcolor[rgb]{ .949,  .949,  .949}$\bar{\nu}^{\textsc{NA}}_R$ & $\bar{\nu}^{\textsc{DNA}}_R$ \\
\midrule
    109   & \multirow{9}[2]{*}{(6,5)} & 0     & \cellcolor[rgb]{ .949,  .949,  .949}3016.8 & 7432.8 & \cellcolor[rgb]{ .949,  .949,  .949}5463.7 & 9377.8 & \cellcolor[rgb]{ .949,  .949,  .949}3240.6 & 13927.0 & \cellcolor[rgb]{ .949,  .949,  .949}6806.2 & 13651.4 \\
    110   &       & 50    & \cellcolor[rgb]{ .949,  .949,  .949}5783.4 & 10153.4 & \cellcolor[rgb]{ .949,  .949,  .949}7404.2 & 13493.4 & \cellcolor[rgb]{ .949,  .949,  .949}5940.9 & 11026.0 & \cellcolor[rgb]{ .949,  .949,  .949}9714.2 & 18790.5 \\
    111   &       & 100   & \cellcolor[rgb]{ .949,  .949,  .949}5794.8 & 10968.8 & \cellcolor[rgb]{ .949,  .949,  .949}8336.0 & 15323.8 & \cellcolor[rgb]{ .949,  .949,  .949}6040.9 & 18911.0 & \cellcolor[rgb]{ .949,  .949,  .949}11222.6 & 25385.0 \\
    112   &       & 150   & \cellcolor[rgb]{ .949,  .949,  .949}4342.2 & 9224.2 & \cellcolor[rgb]{ .949,  .949,  .949}6689.4 & 11467.2 & \cellcolor[rgb]{ .949,  .949,  .949}4508.2 & 13912.7 & \cellcolor[rgb]{ .949,  .949,  .949}8315.9 & 16565.2 \\
    113   &       & 200   & \cellcolor[rgb]{ .949,  .949,  .949}4884.6 & 10710.6 & \cellcolor[rgb]{ .949,  .949,  .949}7782.0 & 14405.6 & \cellcolor[rgb]{ .949,  .949,  .949}5035.6 & 18792.8 & \cellcolor[rgb]{ .949,  .949,  .949}10012.1 & 23365.9 \\
    114   &       & 250   & \cellcolor[rgb]{ .949,  .949,  .949}4446.6 & 8944.6 & \cellcolor[rgb]{ .949,  .949,  .949}6575.8 & 11688.6 & \cellcolor[rgb]{ .949,  .949,  .949}4707.9 & 16686.1 & \cellcolor[rgb]{ .949,  .949,  .949}7770.8 & 17262.4 \\
    115   &       & 300   & \cellcolor[rgb]{ .949,  .949,  .949}5827.2 & 11153.2 & \cellcolor[rgb]{ .949,  .949,  .949}8325.5 & 14773.2 & \cellcolor[rgb]{ .949,  .949,  .949}6062.7 & 17358.7 & \cellcolor[rgb]{ .949,  .949,  .949}11622.5 & 25515.4 \\
    116   &       & 350   & \cellcolor[rgb]{ .949,  .949,  .949}4345.8 & 8927.8 & \cellcolor[rgb]{ .949,  .949,  .949}6816.8 & 12029.8 & \cellcolor[rgb]{ .949,  .949,  .949}4605.7 & 16729.0 & \cellcolor[rgb]{ .949,  .949,  .949}9140.3 & 17563.5 \\
    117   &       & 400   & \cellcolor[rgb]{ .949,  .949,  .949}4180.2 & 9836.2 & \cellcolor[rgb]{ .949,  .949,  .949}7160.2 & 12836.2 & \cellcolor[rgb]{ .949,  .949,  .949}4367.5 & 14595.2 & \cellcolor[rgb]{ .949,  .949,  .949}8522.6 & 19430.3 \\
    \midrule
    118   & \multirow{9}[2]{*}{(6,10)} & 0     & \cellcolor[rgb]{ .949,  .949,  .949}10104.6 &$>10h$ & \cellcolor[rgb]{ .949,  .949,  .949}10200.6 & 24872.6 & \cellcolor[rgb]{ .949,  .949,  .949}10546.4 & 27686.4 & \cellcolor[rgb]{ .949,  .949,  .949}13701.4 &$>10h$\\
    119   &       & 50    & \cellcolor[rgb]{ .949,  .949,  .949}8583.0 &$>10h$ & \cellcolor[rgb]{ .949,  .949,  .949}9351.0 & 21351.0 & \cellcolor[rgb]{ .949,  .949,  .949}9209.9 & 22440.2 & \cellcolor[rgb]{ .949,  .949,  .949}12711.7 & 35626.3 \\
    120   &       & 100   & \cellcolor[rgb]{ .949,  .949,  .949}11280.6 &$>10h$ & \cellcolor[rgb]{ .949,  .949,  .949}11376.6 & 22130.6 & \cellcolor[rgb]{ .949,  .949,  .949}11747.3 & 29362.6 & \cellcolor[rgb]{ .949,  .949,  .949}16320.9 & 35307.2 \\
    121   &       & 150   & \cellcolor[rgb]{ .949,  .949,  .949}9476.4 &$>10h$ & \cellcolor[rgb]{ .949,  .949,  .949}9860.4 & 23221.4 & \cellcolor[rgb]{ .949,  .949,  .949}9962.1 & 31165.4 & \cellcolor[rgb]{ .949,  .949,  .949}12165.8 &$>10h$\\
    122   &       & 200   & \cellcolor[rgb]{ .949,  .949,  .949}9439.2 &$>10h$ & \cellcolor[rgb]{ .949,  .949,  .949}9919.2 & 20601.2 & \cellcolor[rgb]{ .949,  .949,  .949}9959.5 & 23400.1 & \cellcolor[rgb]{ .949,  .949,  .949}13121.1 &$>10h$\\
    123   &       & 250   & \cellcolor[rgb]{ .949,  .949,  .949}8892.6 &$>10h$ & \cellcolor[rgb]{ .949,  .949,  .949}9084.6 & 17979.6 & \cellcolor[rgb]{ .949,  .949,  .949}9535.3 & 24697.2 & \cellcolor[rgb]{ .949,  .949,  .949}10903.3 & 26120.8 \\
    124   &       & 300   & \cellcolor[rgb]{ .949,  .949,  .949}9023.4 &$>10h$ & \cellcolor[rgb]{ .949,  .949,  .949}9983.4 & 18604.4 & \cellcolor[rgb]{ .949,  .949,  .949}9466.8 & 30363.2 & \cellcolor[rgb]{ .949,  .949,  .949}11968.1 & 26395.9 \\
    125   &       & 350   & \cellcolor[rgb]{ .949,  .949,  .949}11331.6 &$>10h$ & \cellcolor[rgb]{ .949,  .949,  .949}11427.6 & 27149.6 & \cellcolor[rgb]{ .949,  .949,  .949}11807.6 & 28773.7 & \cellcolor[rgb]{ .949,  .949,  .949}13932.5 &$>10h$\\
    126   &       & 400   & \cellcolor[rgb]{ .949,  .949,  .949}11216.4 &$>10h$ & \cellcolor[rgb]{ .949,  .949,  .949}12080.4 & 28663.4 & \cellcolor[rgb]{ .949,  .949,  .949}11832.9 & 26791.1 & \cellcolor[rgb]{ .949,  .949,  .949}14588.3 &$>10h$\\
    \midrule
    127   & \multirow{9}[2]{*}{(6,15)} & 0     & \cellcolor[rgb]{ .949,  .949,  .949}18997.2 &$>10h$ & \cellcolor[rgb]{ .949,  .949,  .949}$>10h$ &$>10h$ & \cellcolor[rgb]{ .949,  .949,  .949}19922.6 &$>10h$ & \cellcolor[rgb]{ .949,  .949,  .949}$>10h$ &$>10h$\\
    128   &       & 50    & \cellcolor[rgb]{ .949,  .949,  .949}18051.6 &$>10h$ & \cellcolor[rgb]{ .949,  .949,  .949}$>10h$ &$>10h$ & \cellcolor[rgb]{ .949,  .949,  .949}19040.8 &$>10h$ & \cellcolor[rgb]{ .949,  .949,  .949}$>10h$ &$>10h$\\
    129   &       & 100   & \cellcolor[rgb]{ .949,  .949,  .949}16756.2 &$>10h$ & \cellcolor[rgb]{ .949,  .949,  .949}$>10h$ &$>10h$ & \cellcolor[rgb]{ .949,  .949,  .949}17648.3 &$>10h$ & \cellcolor[rgb]{ .949,  .949,  .949}$>10h$ &$>10h$\\
    130   &       & 150   & \cellcolor[rgb]{ .949,  .949,  .949}20196.0 &$>10h$ & \cellcolor[rgb]{ .949,  .949,  .949}$>10h$ &$>10h$ & \cellcolor[rgb]{ .949,  .949,  .949}21130.6 &$>10h$ & \cellcolor[rgb]{ .949,  .949,  .949}$>10h$ &$>10h$\\
    131   &       & 200   & \cellcolor[rgb]{ .949,  .949,  .949}19142.4 &$>10h$ & \cellcolor[rgb]{ .949,  .949,  .949}$>10h$ &$>10h$ & \cellcolor[rgb]{ .949,  .949,  .949}20092.3 &$>10h$ & \cellcolor[rgb]{ .949,  .949,  .949}$>10h$ &$>10h$\\
    132   &       & 250   & \cellcolor[rgb]{ .949,  .949,  .949}19387.2 &$>10h$ & \cellcolor[rgb]{ .949,  .949,  .949}$>10h$ &$>10h$ & \cellcolor[rgb]{ .949,  .949,  .949}20388.0 &$>10h$ & \cellcolor[rgb]{ .949,  .949,  .949}$>10h$ &$>10h$\\
    133   &       & 300   & \cellcolor[rgb]{ .949,  .949,  .949}16451.4 &$>10h$ & \cellcolor[rgb]{ .949,  .949,  .949}$>10h$ &$>10h$ & \cellcolor[rgb]{ .949,  .949,  .949}17483.9 &$>10h$ & \cellcolor[rgb]{ .949,  .949,  .949}$>10h$ &$>10h$\\
    134   &       & 350   & \cellcolor[rgb]{ .949,  .949,  .949}16513.2 &$>10h$ & \cellcolor[rgb]{ .949,  .949,  .949}$>10h$ &$>10h$ & \cellcolor[rgb]{ .949,  .949,  .949}17490.8 &$>10h$ & \cellcolor[rgb]{ .949,  .949,  .949}$>10h$ &$>10h$\\
    135   &       & 400   & \cellcolor[rgb]{ .949,  .949,  .949}15548.4 &$>10h$ & \cellcolor[rgb]{ .949,  .949,  .949}$>10h$ &$>10h$ & \cellcolor[rgb]{ .949,  .949,  .949}16621.3 &$>10h$ & \cellcolor[rgb]{ .949,  .949,  .949}$>10h$ &$>10h$\\
    \midrule
    136   & \multirow{9}[2]{*}{(6,20)} & 0     & \cellcolor[rgb]{ .949,  .949,  .949}26038.2 &$>10h$ & \cellcolor[rgb]{ .949,  .949,  .949}$>10h$ &$>10h$ & \cellcolor[rgb]{ .949,  .949,  .949}28075.0 &$>10h$ & \cellcolor[rgb]{ .949,  .949,  .949}$>10h$ &$>10h$\\
    137   &       & 50    & \cellcolor[rgb]{ .949,  .949,  .949}25179.0 &$>10h$ & \cellcolor[rgb]{ .949,  .949,  .949}$>10h$ &$>10h$ & \cellcolor[rgb]{ .949,  .949,  .949}26638.6 &$>10h$ & \cellcolor[rgb]{ .949,  .949,  .949}$>10h$ &$>10h$\\
    138   &       & 100   & \cellcolor[rgb]{ .949,  .949,  .949}27861.0 &$>10h$ & \cellcolor[rgb]{ .949,  .949,  .949}$>10h$ &$>10h$ & \cellcolor[rgb]{ .949,  .949,  .949}29812.5 &$>10h$ & \cellcolor[rgb]{ .949,  .949,  .949}$>10h$ &$>10h$\\
    139   &       & 150   & \cellcolor[rgb]{ .949,  .949,  .949}27286.2 &$>10h$ & \cellcolor[rgb]{ .949,  .949,  .949}$>10h$ &$>10h$ & \cellcolor[rgb]{ .949,  .949,  .949}29240.4 &$>10h$ & \cellcolor[rgb]{ .949,  .949,  .949}$>10h$ &$>10h$\\
    140   &       & 200   & \cellcolor[rgb]{ .949,  .949,  .949}26205.0 &$>10h$ & \cellcolor[rgb]{ .949,  .949,  .949}$>10h$ &$>10h$ & \cellcolor[rgb]{ .949,  .949,  .949}27575.8 &$>10h$ & \cellcolor[rgb]{ .949,  .949,  .949}$>10h$ &$>10h$\\
    141   &       & 250   & \cellcolor[rgb]{ .949,  .949,  .949}26002.2 &$>10h$ & \cellcolor[rgb]{ .949,  .949,  .949}$>10h$ &$>10h$ & \cellcolor[rgb]{ .949,  .949,  .949}27503.4 &$>10h$ & \cellcolor[rgb]{ .949,  .949,  .949}$>10h$ &$>10h$\\
    142   &       & 300   & \cellcolor[rgb]{ .949,  .949,  .949}28728.0 &$>10h$ & \cellcolor[rgb]{ .949,  .949,  .949}$>10h$ &$>10h$ & \cellcolor[rgb]{ .949,  .949,  .949}30580.6 &$>10h$ & \cellcolor[rgb]{ .949,  .949,  .949}$>10h$ &$>10h$\\
    143   &       & 350   & \cellcolor[rgb]{ .949,  .949,  .949}24949.8 &$>10h$ & \cellcolor[rgb]{ .949,  .949,  .949}$>10h$ &$>10h$ & \cellcolor[rgb]{ .949,  .949,  .949}26952.8 &$>10h$ & \cellcolor[rgb]{ .949,  .949,  .949}$>10h$ &$>10h$\\
    144   &       & 400   & \cellcolor[rgb]{ .949,  .949,  .949}29950.2 &$>10h$ & \cellcolor[rgb]{ .949,  .949,  .949}$>10h$ &$>10h$ & \cellcolor[rgb]{ .949,  .949,  .949}31406.3 &$>10h$ & \cellcolor[rgb]{ .949,  .949,  .949}$>10h$ &$>10h$\\
    \midrule
    145   & \multirow{9}[2]{*}{(6,25)} & 0     & \cellcolor[rgb]{ .949,  .949,  .949}$>10h$ &$>10h$ & \cellcolor[rgb]{ .949,  .949,  .949}$>10h$ &$>10h$ & \cellcolor[rgb]{ .949,  .949,  .949}$>10h$ &$>10h$ & \cellcolor[rgb]{ .949,  .949,  .949}$>10h$ &$>10h$\\
    146   &       & 50    & \cellcolor[rgb]{ .949,  .949,  .949}$>10h$ &$>10h$ & \cellcolor[rgb]{ .949,  .949,  .949}$>10h$ &$>10h$ & \cellcolor[rgb]{ .949,  .949,  .949}$>10h$ &$>10h$ & \cellcolor[rgb]{ .949,  .949,  .949}$>10h$ &$>10h$\\
    147   &       & 100   & \cellcolor[rgb]{ .949,  .949,  .949}$>10h$ &$>10h$ & \cellcolor[rgb]{ .949,  .949,  .949}$>10h$ &$>10h$ & \cellcolor[rgb]{ .949,  .949,  .949}$>10h$ &$>10h$ & \cellcolor[rgb]{ .949,  .949,  .949}$>10h$ &$>10h$\\
    148   &       & 150   & \cellcolor[rgb]{ .949,  .949,  .949}$>10h$ &$>10h$ & \cellcolor[rgb]{ .949,  .949,  .949}$>10h$ &$>10h$ & \cellcolor[rgb]{ .949,  .949,  .949}$>10h$ &$>10h$ & \cellcolor[rgb]{ .949,  .949,  .949}$>10h$ &$>10h$\\
    149   &       & 200   & \cellcolor[rgb]{ .949,  .949,  .949}$>10h$ &$>10h$ & \cellcolor[rgb]{ .949,  .949,  .949}$>10h$ &$>10h$ & \cellcolor[rgb]{ .949,  .949,  .949}$>10h$ &$>10h$ & \cellcolor[rgb]{ .949,  .949,  .949}$>10h$ &$>10h$\\
    150   &       & 250   & \cellcolor[rgb]{ .949,  .949,  .949}$>10h$ &$>10h$ & \cellcolor[rgb]{ .949,  .949,  .949}$>10h$ &$>10h$ & \cellcolor[rgb]{ .949,  .949,  .949}$>10h$ &$>10h$ & \cellcolor[rgb]{ .949,  .949,  .949}$>10h$ &$>10h$\\
    151   &       & 300   & \cellcolor[rgb]{ .949,  .949,  .949}$>10h$ &$>10h$ & \cellcolor[rgb]{ .949,  .949,  .949}$>10h$ &$>10h$ & \cellcolor[rgb]{ .949,  .949,  .949}$>10h$ &$>10h$ & \cellcolor[rgb]{ .949,  .949,  .949}$>10h$ &$>10h$\\
    152   &       & 350   & \cellcolor[rgb]{ .949,  .949,  .949}$>10h$ &$>10h$ & \cellcolor[rgb]{ .949,  .949,  .949}$>10h$ &$>10h$ & \cellcolor[rgb]{ .949,  .949,  .949}$>10h$ &$>10h$ & \cellcolor[rgb]{ .949,  .949,  .949}$>10h$ &$>10h$\\
    153   &       & 400   & \cellcolor[rgb]{ .949,  .949,  .949}$>10h$ &$>10h$ & \cellcolor[rgb]{ .949,  .949,  .949}$>10h$ &$>10h$ & \cellcolor[rgb]{ .949,  .949,  .949}$>10h$ &$>10h$ & \cellcolor[rgb]{ .949,  .949,  .949}$>10h$ &$>10h$\\
    \midrule
    154   & \multirow{9}[2]{*}{(6,30)} & 0     & \cellcolor[rgb]{ .949,  .949,  .949}$>10h$ &$>10h$ & \cellcolor[rgb]{ .949,  .949,  .949}$>10h$ &$>10h$ & \cellcolor[rgb]{ .949,  .949,  .949}$>10h$ &$>10h$ & \cellcolor[rgb]{ .949,  .949,  .949}$>10h$ &$>10h$\\
    155   &       & 50    & \cellcolor[rgb]{ .949,  .949,  .949}$>10h$ &$>10h$ & \cellcolor[rgb]{ .949,  .949,  .949}$>10h$ &$>10h$ & \cellcolor[rgb]{ .949,  .949,  .949}$>10h$ &$>10h$ & \cellcolor[rgb]{ .949,  .949,  .949}$>10h$ &$>10h$\\
    156   &       & 100   & \cellcolor[rgb]{ .949,  .949,  .949}$>10h$ &$>10h$ & \cellcolor[rgb]{ .949,  .949,  .949}$>10h$ &$>10h$ & \cellcolor[rgb]{ .949,  .949,  .949}$>10h$ &$>10h$ & \cellcolor[rgb]{ .949,  .949,  .949}$>10h$ &$>10h$\\
    157   &       & 150   & \cellcolor[rgb]{ .949,  .949,  .949}$>10h$ &$>10h$ & \cellcolor[rgb]{ .949,  .949,  .949}$>10h$ &$>10h$ & \cellcolor[rgb]{ .949,  .949,  .949}$>10h$ &$>10h$ & \cellcolor[rgb]{ .949,  .949,  .949}$>10h$ &$>10h$\\
    158   &       & 200   & \cellcolor[rgb]{ .949,  .949,  .949}$>10h$ &$>10h$ & \cellcolor[rgb]{ .949,  .949,  .949}$>10h$ &$>10h$ & \cellcolor[rgb]{ .949,  .949,  .949}$>10h$ &$>10h$ & \cellcolor[rgb]{ .949,  .949,  .949}$>10h$ &$>10h$\\
    159   &       & 250   & \cellcolor[rgb]{ .949,  .949,  .949}$>10h$ &$>10h$ & \cellcolor[rgb]{ .949,  .949,  .949}$>10h$ &$>10h$ & \cellcolor[rgb]{ .949,  .949,  .949}$>10h$ &$>10h$ & \cellcolor[rgb]{ .949,  .949,  .949}$>10h$ &$>10h$\\
    160   &       & 300   & \cellcolor[rgb]{ .949,  .949,  .949}$>10h$ &$>10h$ & \cellcolor[rgb]{ .949,  .949,  .949}$>10h$ &$>10h$ & \cellcolor[rgb]{ .949,  .949,  .949}$>10h$ &$>10h$ & \cellcolor[rgb]{ .949,  .949,  .949}$>10h$ &$>10h$\\
    161   &       & 350   & \cellcolor[rgb]{ .949,  .949,  .949}$>10h$ &$>10h$ & \cellcolor[rgb]{ .949,  .949,  .949}$>10h$ &$>10h$ & \cellcolor[rgb]{ .949,  .949,  .949}$>10h$ &$>10h$ & \cellcolor[rgb]{ .949,  .949,  .949}$>10h$ &$>10h$\\
    162   &       & 400   & \cellcolor[rgb]{ .949,  .949,  .949}$>10h$ &$>10h$ & \cellcolor[rgb]{ .949,  .949,  .949}$>10h$ &$>10h$ & \cellcolor[rgb]{ .949,  .949,  .949}$>10h$ &$>10h$ & \cellcolor[rgb]{ .949,  .949,  .949}$>10h$ &$>10h$\\
    \bottomrule
    \end{tabular}%
    }
\end{table}%
\end{appendices}

\end{document}